%% file: singl-K3-mth.tex
\documentclass[11pt, amssymb]{article}
\usepackage{amsmath}
\usepackage{amsfonts}
\usepackage{amscd}
\usepackage{amssymb}
\usepackage{graphicx}
\usepackage{color}
\input{colordvi.tex}
\newcommand{\vev}[1]{ \left\langle {#1} \right\rangle }

\newcommand{\nin}{\mbox{\ooalign{\hfil/\hfil\crcr$\in$}}}

\def\SU{\mathop{\rm SU}}

\def\Z{\mathbb{Z}}
\def\Q{\mathbb{Q}}
\def\R{\mathbb{R}}
\def\C{\mathbb{C}}
\def\P{\mathbb{P}}
\def\F{\mathbb{F}}

\newcommand{\nn}{ {\nonumber} }

\newcommand{\sta}[3]{
\begin{array}{ccccc}
             & &  {#1} & &  \\
             & \nearrow & & \searrow & \\
             {#3} & & \leftarrow  & &{#2}
              \end{array}
              }

\begin{document}

%%%%%%%%%%
%%%%%%%%%%      title page
%%%%%%%%%%

\begin{titlepage}
 
\begin{flushright}
KCL-mth-13-11 \\
YITP-13-127 \\
IPMU13-0236
\end{flushright}
 
\vskip 1cm
\begin{center}
 
 {\large \bf 
On the Classification of Elliptic Fibrations modulo Isomorphism \\
on K3 Surfaces with large Picard Number} \\
--- including a review on classifications of elliptic fibrations on K3 surfaces ---
 
 \vskip 1.2cm
 
 Andreas P. Braun$^1$, Yusuke Kimura$^2$ and Taizan Watari$^{3}$

 \vskip 0.4cm
 
 {\it $^1$Department of Mathematics, King's College, London WC2R 2LS, UK
\\[2mm]

 $^2$Yukawa Institute for Theoretical Physics, Kyoto University, 
Kyoto 606-8502, Japan  \\[2mm]
 
 $^3$Kavli Institute for the Physics and Mathematics of the Universe, 
University of Tokyo, Kashiwano-ha 5-1-5, 277-8583, Japan
 }
 \vskip 1.5cm
 
\abstract{Motivated by a problem originating in string theory, we study elliptic fibrations on K3 surfaces 
with large Picard number modulo isomorphism. We give methods to determine upper bounds for the number of 
inequivalent K3 surfaces sharing the same frame lattice. For any given Neron--Severi lattice $S_X$, such a bound
on the `multiplicity' can be derived by investigating the quotient of the isometry group of $S_X$ by the automorphism group. 
The resulting bounds are strongest for large Picard numbers and multiplicities of unity do indeed occur for a number of
K3 surfaces with Picard number 20. Under a few extra conditions, a more refined analysis is also possible by explicitly studying the 
embedding of $S_X$ into the even unimodular lattice ${\rm II}_{1,25}$ and exploiting the detailed structure of the isometry 
groups of $S_X$ and ${\rm II}_{1,25}$. We illustrate these methods in examples and derive bounds for the number of elliptic 
fibrations on Kummer surfaces of Picard numbers 17 and 20. As an intermediate step, we also discuss coarser classification schemes 
and review known results.} 
 
\end{center}
\end{titlepage}

\tableofcontents

%%%%%%%%%%%%%%%%%%%%%%%%%%%%%%%%%%%%%%%%%%%%%%%%%%%%%%%%%%%%%%%%%%
\section{Introduction and Summary}
%%%%%%%%%%%%%%%%%%%%%%%%%%%%%%%%%%%%%%%%%%%%%%%%%%%%%%%%%%%%%%%%%%

Let $X$ be a smooth complex K3 surface, and 
$(\pi_X, \sigma; X, \P^1)$ a set of data defining an elliptic fibration:
$\pi_X: X \longrightarrow \P^1$ is a morphism with the generic fibre 
a smooth curve of genus one and $\sigma: \P^1 \longrightarrow X$ a section 
satisfying $\pi_X \cdot \sigma = {\bf id}_{\P^1}$. One can think of 
classification problems of such data of elliptic fibrations by introducing 
various equivalence relations among them. 
One classification is by the type of elliptic fibration, and closely related 
is another classification by the isometry class of the frame lattice of an elliptic 
fibration. These two are referred to as ${\cal J}^{(\rm type)}(X)$ and 
${\cal J}_2(X)$ classifications, respectively, in this article 
(see section \ref{sec:3-classifications} for more information). 
The main theme of this article is a classification of elliptic fibrations 
by isomorphism. 
\begin{quote}
{\bf Definition}: Two elliptic fibrations 
$(\pi_X, \sigma; X, \P^1)$ and $(\pi'_X, \sigma'; X, \P^1)$ are said to be 
{\it isomorphic} if there is a bijective morphism from $X$ to itself 
({\it automorphism}) $f \in {\rm Aut}(X)$ and a
$g \in {\rm Aut}(\P^1) = PGL(2; \C)$ such that 
$\pi'_X = g \cdot \pi_X \cdot f$ and $\sigma = f \cdot \sigma' \cdot g$.
\end{quote}
This introduces an equivalence relation among the data of elliptic fibration 
that $X$ admits, and individual equivalence classes are called \emph{isomorphism 
classes of elliptic fibrations} on $X$. 
Let ${\cal J}_1(X)$ denote the set of such isomorphism classes.
The equivalence relation for ${\cal J}_1(X)$ is smaller than that of 
${\cal J}_2(X)$, and hence the former classification is finer than 
the latter. 
In this article, mainly in section \ref{sec:J1-classification}, 
we discuss how to derive upper bounds on
the number of isomorphism classes of elliptic fibration that $X$ admits 
for $X$ with large Picard number $\rho_X$. 

The second theme of this article is to provide a review on issues 
associated with elliptic fibrations on K3 surface in general, with 
contemporary string theorists in mind as our target readers.  
This purpose has been served by \cite{Aspinwall} for already more than a 
decade. More can be learned by consulting reviews purely dedicated to mathematics,
such as \cite{Sch}. However, the problems faced by the present authors in \cite{BKW-phys},
require more material to be covered. This article is therefore not only written 
to report new mathematical results on the issue mentioned above, but also to serve (when combined with 
sections 3--4.2 of \cite{BKW-phys}) as a review article 
supplementary to \cite{Aspinwall, Sch}.\footnote{We included more review
material than needed for an original article in mathematics, in order to 
make this article accessible to string theorists. Although jargon
appears without explanation in this summary section, we have tried to provide 
explanations sufficient for string theorists in the following sections.} 

Sections \ref{sec:prelim}--\ref{ssec:Niemeier} should be regarded 
as review material for the most part. 
Section \ref{sec:3-classifications} explains four different classifications 
of elliptic fibrations on $X$---${\cal J}_0(X)$, ${\cal J}_2(X)$, 
${\cal J}^{({\rm type})}(X)$ and ${\cal J}_1(X)$---corresponding to 
different choices of equivalence relations. 
Exploiting the Torelli theorem for K3 surfaces, these classification 
problems are completely translated into the language of lattice theory. 
Although this subject is in principle considered to be well-known, i) we are 
interested in a version of this problem when the existence of a section $\sigma$ is required, 
and ii) a careful attention is paid to the choice of the modular group. 
Lemmas A and B in the logical chain, the proofs of which we have not found in the 
literature, are also provided. 
Sections \ref{ssec:J2-Nishiyama}--\ref{ssec:Niemeier} are devoted
to a review of the Kneser--Nishiyama method \cite{Nish01, Nish02}, which 
makes it possible to work out the ${\cal J}_2(X)$ classification of 
elliptic fibrations on $X$ systematically (for relatively large Picard 
numbers $\rho_X$).

The Kneser--Nishiyama method \cite{Nish01, Nish02} computes a negative definite 
rank $(26-\rho_X)$ lattice $T_0$ from the transcendental lattice $T_X$, 
and then uses $T_0$ to determine the ${\cal J}_2(X)$ classification of $X$. 
Although this method has been used to determine the ${\cal J}_2(X)$ 
classification for some K3 surfaces with large Picard number, all studies 
that we are aware of are for cases where $T_0$ is either a direct sum of
root lattices of $A$--$D$--$E$ type, or an overlattice of a sum of root lattices \cite{Nish01,Nish02, Bertin:2011}. 

We have computed the lattice $T_0$ for twenty-four K3 surfaces with 
$\rho=20$ (the companion article \cite{BKW-phys} explains the motivation to study them), 
the results are found in Table \ref{tab:T0}. The $T_0$ lattice can be chosen as 
a sum of $A$--$D$--$E$ root lattices for 11 K3 surfaces among them, while 
this is not the case for the other 13. One K3 surface, which we name 
$X_{[3~0~2]}$, is chosen from the latter 13 K3 surfaces as an example for which 
$T_0$ is not a root lattice. Tables \ref{tab:J2classify-[abc=302]} 
and \ref{tab:J2classify-[abc=302]_Etype} in section \ref{ssec:new&known} 
show how to carry out the ${\cal J}_2(X)$ classification for this 
K3 surface. These two tables contain only a part of the ${\cal J}_2(X)$ classification of $X = X_{[3~0~2]}$, but 
they already exhausted all cases for which ${\rm II}^*$, ${\rm III}^*$ or ${\rm IV}^*$-type singular fibres appear. 
It turns out that there are 43 such entries in ${\cal J}_2(X)$ for $X = X_{[3~0~2]}$.  

The subject of section \ref{sec:J1-classification} is the modulo--isomorphism classification 
of elliptic fibrations, which we denote by ${\cal J}_1(X)$. 
The ${\cal J}_1(X)$ classification is built on top of the 
${\cal J}_2(X)$ classification \cite{Oguiso}. Since the former is finer 
than the latter, the ${\cal J}_1(X)$ classification can be described by 
specifying how many isomorphism classes are contained in a given 
equivalence class of the ${\cal J}_2(X)$ classification. We call 
this number the {\it multiplicity} or the {\it number of isomorphism classes} in
this article.
Oguiso completely worked out the ${\cal J}_2(X)$ and ${\cal J}_1(X)$ 
classifications for Kummer surfaces $X = {\rm Km}(E \times F)$ associated 
with the product of two elliptic curves $E$ and $F$ with generic 
complex structures \cite{Oguiso}. By generalizing ideas of
\cite{Oguiso}, we first derive an upper bound on the multiplicity 
for a given K3 surface $X$ which holds for any one of the equivalence 
classes in ${\cal J}_2(X)$ (Proposition C and Corollary D). The results
recorded in Table \ref{tab:coset-alt-rho20} can be used with 
Proposition C and Corollary D to derive such upper bounds for 
thirty-four K3 surfaces $X$ with $\rho_X = 20$. In particular, 
we found ten K3 surfaces with $\rho_X = 20$ where 
the multiplicity cannot be larger than 1, which means that 
${\cal J}_1(X) = {\cal J}_2(X)$. 

It is also possible to derive a stronger upper bound on the
multiplicity (Proposition E and Corollary F), by working on 
individual equivalence classes in ${\cal J}_2(X)$. This is done 
by studying orbits of a finite subgroup of the isometry group 
${\rm Isom}(S_X)$ of the Neron--Severi lattice $S_X$ of $X$ 
acting on a divisor defining an elliptic fibration of a given 
equivalence class in ${\cal J}_2(X)$, as in \cite{Oguiso}. 
For practical computations involving ${\rm Isom}(S_X)$, 
we use the theory of \cite{Borcherds2, Kondo-auto-JacobKumm}, which 
is applicable to K3 surfaces where $T_0$ lattice is a direct sum of 
root lattices of $A$--$D$--$E$ type (or an overlattice of a direct sum of root lattices). 
Using Proposition E and Corollary F, we derive upper bounds on the multiplicity 
of elliptic fibrations individually for all the equivalence class of ${\cal J}_2(X)$ for two K3 surfaces: 
$X = {\rm Km}(A)$ (where $\rho_X =17$) and 
$X = {\rm Km}(E_\omega \times E_\omega)$ (where $\rho_X = 20$). The
results are shown in Table \ref{tab:Kummer-Kumar}, Example G, 
Table \ref{tab:KmEoxEo-multiplicity}, and Example J.

Appendices \ref{ssec:Leech-review}--\ref{ssec:review-Kummer} are included 
in order to make this article self-contained even for readers with a physics 
background. In particular, the Appendices \ref{ssec:Leech-review} and 
\ref{ssec:review-Kummer} explain basic material necessary 
for section \ref{sec:J1-classification}, while 
Appendix \ref{ssec:X3-app} will serve as an exercise problem that 
helps to digest the theory of \cite{Borcherds2, Kondo-auto-JacobKumm}.

%%%%%%%%%%%%%%%%%%%%%%%%%%%%%%%%%%%%%%%%%%%%%%%%%%%%%%%%%%%%%%%%%%
\section{Preliminaries}
\label{sec:prelim}
%%%%%%%%%%%%%%%%%%%%%%%%%%%%%%%%%%%%%%%%%%%%%%%%%%%%%%%%%%%%%%%%%%

Lattice theory, the Torelli theorem for K3 surfaces and the structure 
theory of automorphism groups of K3 surfaces lie at the heart of 
all that is discussed in this article. 
For the convenience of readers with a physics background, 
sections \ref{ssec:prelim-lattice} and \ref{ssec:prelim-K3} 
(along with section 3.1 of \cite{BKW-phys}) provide a minimum version 
of explanation concerning this material, and quote important theorems 
that will be used frequently in this article. 
We recommend the readers to refer to original articles such as \cite{Nik01} (for lattice theory)   
and \cite{PSS} (for automorphism group) or textbooks/lecture notes such as \cite{Barthetal, Huybrechts}
for a more thorough treatment. In section \ref{ssec:finite_groups}, we also give a brief summary 
of notations used in expressing substructure of groups.

%%%%%%%%%%%%%%%%%%%%%%%%%%%%%%%%%%%%%%%%%%%%%%%%%%
\subsection{Basic Notions in Lattice Theory}
\label{ssec:prelim-lattice}
%%%%%%%%%%%%%%%%%%%%%%%%%%%%%%%%%%%%%%%%%%%%%%%%%%%

%%%%%%%%%%%%%%%%%%%%%%%%%%%%%%%%%%%%%%%%
\subsubsection{Lattice, Primitive Embedding, Orthogonal Complement}
%%%%%%%%%%%%%%%%%%%%%%%%%%%%%%%%%%%%%%%%
\indent

{\bf Definition} A {\it lattice} $L$ is a free abelian group of finite 
rank, i.e. $L \cong \mathbb{Z}^r$ for some $r \in\mathbb{N}$, 
with {\it non-degenerate symmetric} pairing
$(\cdot,\cdot):L\times L\longrightarrow\mathbb{R}$. The value 
$(x,x) \in \R$ for an element $x \in L$ is called its {\it norm} or 
{\it self-intersection} and we sometimes abbreviate $(x,x)$ by $x^2$. 
If there is an isomorphism 
$\phi: L \longrightarrow L'$ of Abelian groups between a pair of
lattices $L$ and $L'$ and $\phi$ preserves the symmetric pairing, 
the two lattices $L$ and $L'$ are said to be 
{\it isometric} and $\phi$ is said to be an {\it isometry} or 
a {\it lattice isomorphism}.

\vspace{5mm}

{\bf Definition}
A lattice $L$ is said to be {\it
positive-definite} (resp. {\it negative-definite}) if $(x,x)=x^2 > 0$
(resp. $x^2 < 0$) for all {\it non-zero} $x$ in $L$. A lattice that is 
neither positive definite nor negative definite is called indefinite. 
When a non-degenerate lattice $L$ has an intersection form (matrix
representation of the symmetric pairing) with $r_+$ positive
``eigenvalues'' and $r_-$ negative ``eigenvalues'', the pair of integers
$(r_+, r_-)$ are said to be the {\it signature} of $L$. 
A lattice $L$ is said to be {\it integral}, if 
$(x,y)\in\mathbb{Z}$ for any $x,y\in$L. 
An {\it Even} lattice, or equivalently {\it Type II} lattice\footnote{
Lattices that are not even are said to be {\it odd}, or equivalently
{\it Type I}. } 
is an integral lattice $L$ with $x^2$ an even integer for all $x\in L$. 

\vspace{5mm}

{\bf Definition} 
For a lattice $L$, ${\rm discr}(L)$---the {\it discriminant of} $L$---is the
determinant of the intersection form (symmetric pairing) of $L$. 
An integral lattice $L$ is said to be {\it unimodular} 
when ${\rm discr}(L)=\pm 1$. 

\vspace{5mm}

{\bf Definition}
For a non-degenerate integral lattice $L$, its dual lattice $L^*$ is 
defined by 
\begin{equation}
 L^* := \left\{ y \in L \otimes \Q \; | \; (y,x) \in \Z {\rm
	 ~for~}{}^\forall x \in L \right\}
\end{equation}
as an Abelian group. Its symmetric pairing is given by naturally extending 
the symmetric pairing of $L$ to $L \otimes \Q$ first, and then by
restricting it to $L^* \subset L \otimes \Q$. A lattice $L$ is said to
be {\it self-dual}, if $L^*$ is isometric to $L$. From these
definitions one can see that unimodularity and self-duality
are equivalent, and the two words can be used interchangeably.

\vspace{5mm}

There is a strong classification theorem for even unimodular
lattices modulo isometry. For a lattice $L$ with indefinite signature, 
i.e., $r_+ > 0$ {\it and} $r_- > 0$, an even unimodular lattice exists if
$r_+ \equiv r_-$ modulo 8. In this case, it is uniquely determined by the 
signature and its rank. Furthermore, $L$ must be isometric to a direct product of $U$, $E_8$ 
and\footnote{For a lattice $L$, we denote a lattice that is isomorphic to $L$ as an Abelian group, but has a
symmetric pairing which is $n$ times larger (with $n \in \Z$) than that of $L$ by $L[n]$.} $E_8[-1]$. 
See section \ref{sssec:prelim-lattice-miscl} for the definition of 
$U$ and $E_8$ lattices. Lattices characterized uniquely (modulo
isometry) by their signature $(r_+,r_-)$ are denoted by 
${\rm II}_{r_+, \; r_-}$:
\begin{equation}
 {\rm II}_{r_+, \; r_-} \cong U^{\oplus r_+} \oplus E_8^{\oplus m} \quad 
{\rm if~} r_- = r_+ + 8m, \qquad 
 {\rm II}_{r_+, \; r_-} \cong U^{\oplus r_-} \oplus E_8[-1]^{\oplus m} \quad 
{\rm if~} r_+ = r_- + 8m \, . 
\end{equation}

For positive definite (or negative definite) even unimodular lattices, 
however, this uniqueness is lost. Although $E_8$ (root lattice) is the 
only rank-8 even unimodular lattice, there are two isometry classes of 
rank-16 even unimodular lattice and twenty-four rank-24 even 
unimodular lattices that are not mutually isometric \cite{CS-book}. 

\vspace{5mm}

{\bf Definition} A lattice $M$ is said to be a {\it sublattice} of $L$, 
when $M$ can be identified with a free Abelian subgroup of $L$ and 
the symmetric pairing on $M$ is given by restricting the symmetric 
pairing of $L$.

\vspace{5mm}

{\bf Definition} An {\it embedding} $\phi$ of a lattice $M$ into $L$ 
is an injective homomorphism $\phi:M \longrightarrow L$ of Abelian
groups such that the lattice $M$ and the sublattice $\phi(M)$ of $L$
are isometric under $\phi$. A lattice embedding $\phi: M
\longrightarrow L$ is said to be {\it isomorphic} to another embedding 
$\phi': M \longrightarrow L'$ if there is an isometry $f: L
\longrightarrow L'$ so that $\phi' = f \cdot \phi: M \longrightarrow L'$.

\vspace{5mm}

{\bf Definition} A sublattice $M$ in $L$ is {\it primitive} when the
quotient $L/M$ becomes a torsion-free Abelian group. A lattice embedding 
$\phi: M \longrightarrow L$ is said to be {\it primitive}, if $\phi(M)$ 
is a primitive sublattice of $L$.

\vspace{5mm}

Suppose that $M$ is a primitive sublattice of $L$. 
Then the short exact sequence 
$0 \longrightarrow M \longrightarrow L \longrightarrow L/M \longrightarrow 0$ 
of Abelian groups always splits and there exists an isomorphism (of Abelian groups)
\begin{equation}
 L \cong M \oplus (L/M) \, .
\end{equation}
Note, however, that the symmetric pairing (intersection form) does not necessarily respect this direct 
sum decomposition. 

It is useful to note that, for a lattice $L$ and its sublattice $M$, 
the two conditions 
\begin{itemize}
 \item $M$ is primitive in $L$,
 \item $(M \otimes \Q \subset L \otimes \Q) \cap L = M$,
\end{itemize}
are equivalent. 

\vspace{5mm}

{\bf Definition} For a sublattice $M$ of a lattice $L$, 
we can always define another sublattice $M^\perp$ of $L$: 
\begin{equation}
M^{\perp}:=\{x\in L|(x,y)=0 \hskip 2mm {}^{\forall}y\in M\}.
\end{equation}
This new lattice $M^\perp \subset L$ is called the {\it orthogonal complement}
of $M$ in $L$. By definition, a sublattice of a lattice $L$ obtained as the orthogonal complement 
of some sublattice of $L$ is primitive in $L$. The lattice homomorphism 
\begin{equation}
 i: M \oplus M^{\perp} \hookrightarrow L
\end{equation}
is injective, and this embedding is always of finite index (i.e., 
$[L: i(M \oplus M^\perp)] < \infty$). When a sublattice $M$ of a lattice $L$ is unimodular, then 
there even is an isometry (lattice isomorphism)
\begin{equation}
 L \simeq M \oplus M^{\perp} \, .
\end{equation}

{\bf Example} Let $L = \Z\vev{e_1} \oplus \Z \vev{e_2}$, with
symmetric pairing $(e_i, e_j) = 2\delta_{ij}$. For a sublattice 
$M$ generated by $(e_1 + e_2) \in L$, i.e., 
$M = \Z \vev{e_1 + e_2} \subset L$, the orthogonal complement is given
by $M^\perp = \Z \vev{e_1-e_2}$. 
$M \oplus M^\perp$ generates a sublattice of $L$, but 
$L \neq M \oplus M^\perp$.

\vspace{5mm} 

{\bf Definition} A lattice $L$ is said to be an {\it overlattice} of a 
lattice $M$, when $M$ is a sublattice of $L$ with the index $[L:M]$ 
being finite. (So $L$ and $M$ have the {\it same} rank.) 
\vspace{5mm}

\noindent
The dual lattice $L^*$ is always an overlattice of $L$, 
and $L$ is an overlattice of $M \oplus M^\perp$.

%%%%%%%%%%%%%%%%%%%%%%%%%%%%%%%%%%%%%%%%%
\subsubsection{Discriminant Group and Discriminant Form}
%%%%%%%%%%%%%%%%%%%%%%%%%%%%%%%%%%%%%%%%%

Although there is no strong classification theorem for lattices that are 
not unimodular, the discriminant group of a lattice and the 
discriminant form provide powerful tools in studying various 
properties associated with even non-unimodular lattices \cite{Nik01}. 

\vspace{5mm}

{\bf Definition} A {\it discriminant group} $G_L$ of an even 
lattice $L$ is a finite Abelian group given by
\begin{equation}
 G_L := L^* / L,
\end{equation}
which is non-trivial for lattices that are not unimodular.
A {\it discriminant form} $q_L$ on $G_L$ is a quadratic form 
$q_L: G_L \rightarrow \Q / 2\Z$ defined by 
\begin{align}
q_L:G_L & \longrightarrow \mathbb{Q}/2\mathbb{Z}\\
      x & \longmapsto x^2 \qquad {\rm mod~} 2\mathbb{Z} .
\end{align}
From this quadratic form, a {\it bilinear form} 
$b_L: G_L \times G_L \rightarrow \Q / \Z$ can be defined by 
\begin{equation}
b_L(x,y) = \frac{1}{2}\{ q_L(x+y,x+y) - q_L(x,x) - q_L(y,y)\} \in \Q/\Z 
\end{equation}
for $(x,y) \in G_L \times G_L$. Obviously $q_L(x) = b_L(x,x)$.

\vspace{5mm}

Note that a finite Abelian group $G_L$ can be written in the form 
$\prod_i (\Z/p_i^{k_i} \Z)$, where the $p_i$'s are prime numbers, and the $k_i$'s are
positive integers. This is because $\Z/m\Z \cong \prod_j (\Z/p_j^{k_j} \Z)$, when $m \in \Z_{\geq 0}$ 
is decomposed into primes as $m = \prod_j p_j^{k_j}$. 
Let $\sigma_i$ and $\sigma_j$ be generators of the first two factors in 
$(\Z/p_i^{k_i} \Z) \times (\Z/p_j^{k_j} \Z) \times \cdots \subset G_L$,
where $p_i \neq p_j$.
Then $b_L(\sigma_i, \sigma_j) \equiv 0$ mod $\Z$, because 
$b_L(\sigma_i, \sigma_j) \in \Q/\Z$ has to simultaneously be an integral multiple of 
$1/p_i^{k_{i}}$ and $1/p_j^{k_j}$. 
When $G_L$ contains only one factor of the form of $\Z/p^{k_i} \Z$ 
for a given prime number $p$, then the discriminant form on 
the factor $\Z/p^k\Z\vev{\sigma}$ is always in the form of 
$q(\sigma) = b(\sigma,\sigma) = a/p^k$ with $a \in \Z$, 
since $b(p^k \sigma,\sigma) \equiv 0$ mod $\Z$. 
Such a discriminant form on $G_L \cong (\Z/p^k \Z)$ may be denoted 
by $q_a(p^k)$.
When $G_L$ contains more than one factor of subgroups of the form 
$\Z/p^{k_i} \Z$ for a prime number $p$, then the discriminant form 
can be more complicated. 

\vspace{5mm}

Let ${\rm Isom}(G_L,q_L)$ or just simply ${\rm Isom}(q_L)$ be the group 
of isomorphisms of the Abelian group $G_L$ preserving the discriminant
form. There is a natural homomorphism 
\begin{equation}
 p_L: {\rm Isom}(L) \longrightarrow {\rm Isom}(G_L,q_L).
\end{equation}

\vspace{5mm}

{\bf Definition} Any non-zero primitive element\footnote{
An element $x \in L$ is called a {\it primitive element}, if 
$\Z \vev{x}$ is a primitive sublattice of $L$. This is precisely when $x$ is 
one of the elements in $\Q x \cap L$ closest to the origin.} 
$x$ in a lattice $L$ is called a {\it root}, 
if and only if the following reflection isometry 
\begin{equation}
 s_x: L \otimes \Q \longrightarrow L \otimes \Q, \qquad 
     y \longmapsto y - 2\frac{(x,y)}{(x,x)} x,
\end{equation}
maps the subspace $L$ of $(L \otimes \Q)$ to $L$ itself. 
If $L$ is an even lattice, all the elements $x \in L$ satisfying 
$x^2 = -2$ are roots, and hence they are called $(-2)$ roots.
Roots that are not $(-2)$-roots are called $(-n)$-roots 
if $x^2 = -n \neq -2$.
If $L$ is an even unimodular lattice, then all the roots of $L$ 
are $(-2)$-roots. 

\vspace{5mm}

The subgroup of ${\rm Isom}(L)$ generated by reflections is denoted by 
$W(L)$. Its subgroup generated only by reflections associated with
$(-2)$-roots is denoted by $W^{(2)}(L)$. 
Whereas the $W^{(2)}(L)$ subgroup is contained in the kernel of the homomorphism 
$p_L: {\rm Isom}(L) \longrightarrow {\rm Isom}(G_L,q_L)$ above, this is not 
necessarily true for $W(L)$.

Here, we quote some results, mostly in \cite{Nik01}, that are quite 
important and will also be used in the rest of this article.

{\bf Proposition} $\alpha$ (\cite{Nik01}, Prop. 1.4.1): 
For an even lattice $M$, and 
for an isotropic subgroup $H$ of $G_M$ (meaning that $q_M|_H = 0$), 
we can define a lattice $(M;H)$
\begin{equation}
 M; H := \left\{ m' \in M^* \; |
       \; [m' {\rm ~mod~}M] \in H \subset G_M \right\},
\end{equation}
which is still an even lattice. $M$ is a sublattice of $(M;H)$, and 
$(M;H)$ in turn, is a sublattice of $M^*$. The discriminant group 
$G_{(M; H)}$ of the new lattice $(M; H)$ is given by 
\begin{equation}
  G_{(M;H)} = \left\{ y {\rm ~mod~} H \in G_M \; | \; b_L(y, m')=0 
  {\rm ~for~}{}^{\forall} [m' {\rm ~mod~}M] \in H \right\} \subset G_M/H.
\end{equation}
Conversely, for any even sublattice $M'$ of $M^*$, with $M$ being a
sublattice of $M'$, a corresponding isotropic subgroup $H$ of $G_M$ 
can be found such that $M' = M; H$. 

\vspace{5mm}

{\bf Proposition} $\beta$ (\cite{Nik01}, Prop. 1.6.1): 
Let $M$ and $N$ be even lattices for which there is an isomorphism 
$\gamma: G_M \cong G_N$ such that $q_M = - q_N \cdot \gamma$. Then 
\begin{equation}
 \Delta := \left\{ \left( m', \gamma(m') \right) \in G_M \times G_N \; | \;
    {}^{\forall} m' \in G_M \right\} 
\end{equation}
is an isotropic subgroup of $G_M \times G_N$, the discriminant group of 
an even lattice $M \oplus N$. Thus, Proposition $\alpha$ above introduces 
an even lattice $L = (M \oplus N); \Delta$. This lattice $L$ is 
unimodular. $M$ and $N$ are primitive sublattices of $L$, and are 
mutually orthogonal complements in $L$. 
Conversely, for a primitive sublattice $M$ of an even unimodular lattice 
$L$, the orthogonal complement $N := [ M^\perp \subset L]$ and $M$ have
isomorphic discriminant groups, ${}^\exists \gamma: G_M \cong G_N$, 
such that $q_M = - q_N \cdot \gamma$.

\vspace{5mm}

{\bf Proposition} $\gamma$ (\cite{Nik01}, special case of
Cor. 1.5.2): Suppose that $M$ and $N$ are primitive sublattices of 
an even unimodular lattice $L$ and are mutually orthogonal
complements, as in Proposition $\beta$ above. 
% 1.6.1 of \cite{Nik01}. 
A pair of isometries $\varphi \in {\rm Isom}(M)$ and 
$\psi \in {\rm Isom}(N)$ can be lifted to an isometry of $L$ 
that restricts to $(\varphi, \psi)$ on $(M \oplus N) \subset L$, if and 
only if $\gamma \cdot p_M(\varphi) = p_N(\psi) \cdot \gamma$, where 
$\gamma: G_M \cong G_N$ is the isomorphism that appeared already in
Proposition $\beta$.
% Proposition 1.6.1 of \cite{Nik01}.

\vspace{5mm}

{\bf Theorem} $\delta$ (\cite{Nik01}, Thm. 1.12.4): 
This is a sufficient condition for existence of a primitive embedding of 
an even lattice $M$ of signature $(m_+,m_-)$ into an even unimodular lattice $L$ of signature $(l_+, l_-)$:
A primitive embedding $\phi: M \hookrightarrow L$ exists if 
${\rm rank}(M) \leq ({\rm rank}(L))/2$ and $l_+ \equiv l_-$ (mod 8) 
as well as $m_+ \leq l_+$, $m_- \leq l_-$.

\vspace{5mm}

{\bf Theorem} $\epsilon$ (\cite{Nik01}, Thm. 1.14.4; 
\cite{Mor}, Thm. 2.8): This is a sufficient condition 
for the existence and uniqueness (modulo isometry) of a primitive embedding 
of an even lattice $M$ into an even unimodular lattice $L$: 
$m_+ < l_+$, $m_- < l_-$ and $l(G_M) \leq {\rm rank}(L)- {\rm
rank}(M)-2$, where $l(G_M)$ is the smallest number of generators 
of $G_M$.

\vspace{5mm}

{\bf Theorem} $\zeta$ (\cite{Mor}, Cor. 2.10):
There exists a primitive embedding of an even lattice $T$ with 
signature $(2,20-\rho)$ into the even unimodular lattice ${\rm
II}_{3,19}$ if $\rho \geq 12$. Furthermore, such an embedding 
is unique modulo isometry of ${\rm II}_{3,19}$. The uniqueness 
of this embedding guarantees that the orthogonal complement 
$S:=[T^\perp \subset {\rm II}_{3,19}]$ is also determined uniquely 
modulo isometry. 

%%%%%%%%%%%%%%%%%%%%%%%%%%%%%%%%%%%%%%%%%%%%%%%%%%%
\subsubsection{Miscellany}
\label{sssec:prelim-lattice-miscl}
%%%%%%%%%%%%%%%%%%%%%%%%%%%%%%%%%%%%%%%%%%%%%%%%%%%

The rank-2 signature $(1,1)$ even unimodular lattice is called
hyperbolic plane and is denoted by $U$. It has a set of generators 
$U = {\rm Span}_\Z \{ f_1, f_2 \}$ so that the intersection form is
\begin{equation}
 \left[ \begin{array}{cc} 
    (f_1, f_2) & (f_1, f_2) \\ (f_2,f_1) & (f_2,f_2)
 \end{array} \right] = 
 \left[ \begin{array}{cc} 0 & 1 \\ 1 & 0 \end{array} \right].
\label{eq:U-int-form}
\end{equation}
We also denote the quadratic form above with this choice of basis by $U$. 
When we refer to the hyperbolic plane lattice $U$, however, one is free 
to choose a set of generators different from $\{f_1, f_2\}$, and the 
symmetric pairing of the lattice may correspondingly be written as a matrix different 
from the one above. We do not expect confusion to arise from this abuse of notation.
The lattice $U$ has an isometry group ${\rm Isom}(U) \cong 
\Z/2\Z \vev{-{\rm id.}_U} \times \Z /2\Z \vev{\sigma_U}$, where 
\begin{equation}
 (-{\rm id.}_U): f_{1,2} \mapsto -f_{1,2}, \qquad 
 \sigma_U: f_{1,2} \mapsto f_{2,1}.
\end{equation}
The isometry group of $U \oplus U$ is studied and described in detail in \cite{Hosono:2002yb}.
\vspace{3mm}

\noindent
A lattice denoted by $(n)$ is a rank-1 lattice $\cong \Z\vev{e}$ whose 
symmetric pairing is given by $(e,e) = n$.
\vspace{3mm}

\noindent
Root lattices of the $A_n$, $D_n$ and $E_n$ Lie algebras are also denoted by 
$A_n$, $D_n$ and $E_n$, respectively, in this article. The intersection
form (= symmetric pairing) of these lattices are set to be the 
\emph{negative} of their Cartan matrices, which means that 
the diagonal entries are all $(-2)$. 
They are negative definite and even integral lattices. 
The dual lattice of a root lattice of one of the $A$-$D$-$E$ types is 
its weight lattice, and the discriminant group is given by 
\begin{align}
 & 
 G_{A_n} \cong \Z/(n+1)\Z, \qquad 
 G_{D_n} \cong \Z/2\Z \times \Z/2\Z \quad (n {\rm ~even}), \quad 
 \Z/4\Z \quad (n {\rm ~odd}), \nonumber \\
 & 
 G_{E_6} \cong \Z/3\Z, \quad 
 G_{E_7} \cong \Z/2\Z.
\end{align}
For a set of simple roots $\alpha_i$ ($i = 1, \cdots, r$) of a rank-$r$ Lie
algebra of A--D--E type, a set of fundamental weights 
$\omega^i$ ($i=1,\cdots, r$) are the elements of the weight lattice 
satisfying $(\alpha_i, \omega^j) = - \delta_i^{\; j}$. 
As summarized e.g. in [\S 1 of \cite{Nish01}]:

\begin{itemize}
 \item $G_{A_n}$ is generated by the mod-$A_n$ equivalence class of the 
weights of the ``defining representation'' of $\SU(n+1)$, $[\omega^1]$. 
This generator is also denoted by $a_{n+1}$ in this article. For this generator, 
the value of the discriminant form is $q([\omega^1]) = -n/(n+1)$ modulo 2.
 \item $G_{D_n}$ is isomorphic to $\Z/2\Z \times \Z/2\Z$ for an even
       $n$, where one of the two $\Z/2\Z$'s is generated by the mod-$D_n$
       equivalence class of the weights of a spinor representation
       (denoted by ${\bf sp}$ or $(1,0)$), and the other by that of 
       the weights of the other spinor representation (denoted by 
       $\overline{\bf sp}$ or $(0,1)$). ${\bf sp}+\overline{\bf sp}$ 
       in $G_{D_n} \cong \Z/2\Z \times \Z/2\Z$  is denoted by ${\bf v}$. 
       The discriminant bilinear form $b_L$ on this $\{ {\bf sp},
       \overline{\bf sp}\} $ basis is given by 
\begin{equation}
 b_{D_n} = \left( \begin{array}{cc}
	    -n/4 & -(n-2)/4 \\ -(n-2)/4 & -n/4
		  \end{array}\right). \qquad \qquad 
 \left\{ \begin{array}{ll}
  {\rm diag.~entries} & {\rm mod~}2 \\
  {\rm off~diag.~entries} & {\rm mod~}1 
	 \end{array}\right.
\end{equation} 
 \item $G_{D_n} \cong \Z/4\Z$ for an odd $n$ is generated by the mod-$D_n$
       equivalence class of the weights of a spinor representation,
       denoted by ${\bf sp}$ . $q({\bf sp}) = -n/4$ modulo 2.
 \item $G_{E_6} \cong \Z/3\Z$ is generated by the weights of a  
       27-dimensional representation of $E_6$, and hence the generator 
       of $G_{E_6}$ is denoted by ${\bf 27}$. 
       $q({\bf 27}) = - 4/3$ modulo 2.
 \item $G_{E_7} \cong \Z/2\Z$ is generated by the weights of the 
       56-dimensional representation of $E_7$, and hence the generator 
       of $G_{E_7}$ is denoted by ${\bf 56}$. $q({\bf 56}) = -3/2$
       modulo 2.
\end{itemize}

{\bf Definition} A {\it root lattice} $L_{\rm root}$ is the sublattice
of a lattice $L$ generated by all the $(-2)$-roots of $L$. 

%%%%%%%%%%%%%%%%%%%%%%%%%%%%%%%%%%%%%%%%%%%%%%%%%%%%%%
\subsection{Substructure of Groups}
\label{ssec:finite_groups}
%%%%%%%%%%%%%%%%%%%%%%%%%%%%%%%%%%%%%%%%%%%%%%%%%%%%%%

When there is a short exact sequence of groups
\begin{equation}
 1 \longrightarrow N \longrightarrow 
 G \longrightarrow O  \longrightarrow 1 \, ,
\end{equation}
where $N$ is a normal subgroup 
of $G$, we may write $G = N.O$ in a shorthand notation 
(but not $G = O.N$).\footnote{It is often said that $G$ is an extension 
of $O$ by $N$, but there are also authors with the opposite conventions saying 
that $G$ is an extension of $N$ by $O$. We are not using any of these expressions, 
and, in order to avoid confusion, use $G=N.O$ instead.
Reference \cite{ATLAS} suggests to say either that 
$G$ is an upward extension of $N$ by $O$ or that 
$G$ is a downward extension of $O$ by $N$.} 

If the above short exact sequence splits, i.e., there exists 
a homomorphism $\phi: O \longrightarrow G$ in addition to the short 
exact sequence 
\begin{equation}
 1 \rightarrow N \rightarrow 
 G \leftrightarrows_\psi^\phi O  \rightarrow 1 \, , 
\end{equation}
and the condition 
$\psi \cdot \phi = {\bf id}._O: O \longrightarrow O$ 
is satisfied, we say that $G$ is a semi-direct
product of $N$ and $O$ (or $O$ and $N$), and this situation is expressed
in shorthand notations such as  
\begin{equation}
 G = N \rtimes \phi(O), \quad  G = \phi(O) \ltimes N \quad 
 ({\rm but~not~} N \ltimes \phi(O)),  \qquad 
 G = N:O, \quad {\rm and} \quad G = N:{}^{{\rm Ad} \cdot \phi} O.
\end{equation}
When a group $G$ has a structure $G = N.O$, but not $N:O$, one may write 
$G = N\raisebox{6pt}{.}O$ to make this situation explicit.\footnote{An easy
example of this situation is the case $G = \Z/4\Z$ with $N$ the
subgroup generated by 2 mod 4.} 
$G=N.O$ can imply either $N:O$ or $N\raisebox{6pt}{.}O$.

All those notations above are also defined and explained, along with others, 
in section 5.2 of \cite{ATLAS}.

%%%%%%%%%%%%%%%%%%%%%%%%%%%%%%%%%%%%%%%%%%%%%%%%%%%%%%
\subsection{Automorphism Groups of K3 Surfaces}
\label{ssec:prelim-K3}
%%%%%%%%%%%%%%%%%%%%%%%%%%%%%%%%%%%%%%%%%%%%%%%%%%%%%%

%%%%%%%%%%%%%%%%%%%%%%%%%%%%%%%%%%%%%%%%%%%%%%%%%%%%%%
\subsubsection{The Neron--Severi Lattice and the Transcendental Lattice}
%%%%%%%%%%%%%%%%%%%%%%%%%%%%%%%%%%%%%%%%%%%%%%%%%%%%%%

Throughout this article, let $X$ denote a K3 surface: a smooth surface over $\C$ with trivial
canonical bundle and $h^1(X, {\cal O}_X) = 0$. As is well-known, the second cohomology group $H^2(X; \Z)$ 
along with its intersection form is a lattice isometric to 
\begin{equation}
 \Lambda_{\rm K3} = U \oplus U \oplus U \oplus E_8 \oplus E_8 \cong 
   {\rm II}_{3,19}\, .
\end{equation}

Whenever we refer to a K3 surface $X$ in this article, we understand 
that it is equipped with a certain complex structure; if we are referring 
to a family of K3 surface, we will explicitly say so. 
Thus, a period vector
\begin{equation}
 [\Omega_X] \in \P \left[ \left\{ \omega \in H^2(X; \C) \; | \; 
     \omega \wedge \omega = 0, \; \omega \wedge \overline{\omega} > 0
     \right\}\right]
\end{equation} 
is fixed and given. For a period vector, and hence for $X$, the
Neron--Severi lattice\footnote{The Neron--Severi group is defined to be
the set of divisors (algebraic curves) on $X$ modulo algebraic
equivalence. Algebraic equivalence, however, is equivalent to 
numerical equivalence in the case of a K3 surface. Hence the
Neron--Severi group does not have a torsion part, i.e., it is a free
Abelian group. The intersection number between a pair of divisors
introduces a symmetric pairing to this free Abelian group. On algebraic
curves, the integral of the period vector $\Omega_X$ vanishes.}---denoted 
by $S_X$---is given by
\begin{equation}
 S_X = \left\{x \in H^2(X; \Z) \; | \; (x,\Omega_X) = 0 \right\}.
\end{equation}
Because of this characterization, $S_X$ is a primitive 
sublattice of $H^2(X; \Z)$. Its rank,
\begin{equation}
 \rho_X := {\rm rank}(S_X)\, ,
\end{equation}
is called the Picard number, which ranges from 0 to 20, depending on the 
complex structure ($[\Omega_X]$) of $X$.
The Neron--Severi lattice $S_X$ has signature $(1,\rho_X-1)$ and
its discriminant group and discriminant form are denoted by $G_{S_X}$ and 
$q_{S_X}$, respectively.

The transcendental lattice of a K3 surface $X$---denoted by $T_X$---is 
defined to be the orthogonal complement of $S_X$ in $H^2(X; \Z)$. 
It has signature $(2,20-\rho_X)$ and its discriminant group and
discriminant form are denoted by $G_{T_X}$ and $q_{T_X}$, respectively. 
$T_X$ is a primitive sublattice of $H^2(X; \Z)$, because of its
definition. Proposition $\beta$ (quoted in section \ref{ssec:prelim-lattice}) 
guarantees that there always exists an isomorphism of Abelian groups
$\gamma: G_{S_X} \cong G_{T_X}$ with $q_{T_X} \cdot \gamma = -q_{S_X}$ holds.
The lattice $H^2(X; \Z)$ can be regarded as a subset of $S_X^* \oplus T_X^*$ characterized by 
\begin{equation}
 \left\{ (s, \gamma(s)) \in G_{S_X} \times G_{T_X} \; | \; s \in G_{S_X}\right\}
  \subset G_{S_X} \times G_{T_X}.
\end{equation}
The identification between $G_{T_X}$ and $G_{S_X}$ by $\gamma$ also 
establishes a canonical isomorphism 
\begin{equation}
 {\rm Isom}(G_{T_X}, q_{T_X}) \cong {\rm Isom}(G_{S_X},q_{S_X}).
\end{equation}
%

%%%%%%%%%%%%%%%%%%%%%%%%%%%%%%%%%%%%%%%%%%%%%%%%%%%%%%
\subsubsection{Automorphism Groups of K3 Surfaces}
%%%%%%%%%%%%%%%%%%%%%%%%%%%%%%%%%%%%%%%%%%%%%%%%%%%%%%

Automorphisms of a K3 surface $X$ form a group denoted by 
${\rm Aut}(X)$. Any automorphism $f \in {\rm Aut}(X)$ induces 
$f_*: H_2(X; \Z) \longrightarrow H_2(X; \Z)$, which is an isometry 
of the lattice $H^2(X; \Z)$. In particular, the homomorphism 
${\rm Aut}(X) \longrightarrow {\rm Isom}(H^2(X; \Z))$ is injective
(\cite{PSS}, \S 2 Prop. 2).

\vspace{5mm}

{\bf Definition} For a K3 surface $X$, a {\it Hodge isometry} is an 
isometry $\phi: H^2(X; \Z) \longrightarrow H^2(X; \Z)$ that maps 
the complex line represented by $[\Omega_X] \in \P[H^2(X; \C)]$ 
to itself. A {\it Hodge and effective isometry} of an algebraic K3
surface is a Hodge isometry that maps all the classes represented by 
effective curves to the same set of classes. 

\vspace{5mm}

{\bf Proposition} $\eta$ (\cite{PSS}, \S 7 Prop.): This is one of several different versions of the Global Torelli Theorem. 
It states that the image of the injective homomorphism 
${\rm Aut}(X) \longrightarrow {\rm Isom}(H_2(X; \Z))$ for an algebraic 
K3 surface $X$ is the group of Hodge and effective isometries.
Thus, the automorphism group can be identified with this subgroup of 
${\rm Isom}(H_2(X; \Z))$.

\vspace{5mm}

{\bf Definition} The {\it positive cone} of an algebraic surface $X$ is the
subspace of $S_X \otimes \R$ given by 
\begin{equation}
 {\rm Pos}_X := \left\{ x \in S_X \otimes \R \; | \; x^2 > 0 \right\}.
\end{equation}
Because $S_X$ has signature $(1,\rho_X-1)$, it consists of two connected components. 
The connected component containing the classes of effective curves are denoted by ${\rm Pos}_X^+$. \\
From the definition of {\rm ample cone} commonly adopted in algebraic
geometry, it follows that the {\it ample cone} of an algebraic surface
$X$ is a subspace of $S_X \otimes \R$ given by 
\begin{equation}
 {\rm Amp}_X = \left\{ x \in {\rm Pos}_X^+ \; | \; (x,C) > 0 
  {\rm ~for~any~curves~(effective~divisors)~} C \right\}.
\end{equation}

\vspace{5mm}

{\bf Theorem} 1 of \S 6 in \cite{PSS} states that the index 2 subgroup
of ${\rm Isom}(S_X)$ that preserves ${\rm Pos}_X^+$---called group of 
autochronos isometries of $S_X$ and denoted by ${\rm Isom}^+(S_X)$---has
the following structure, for a K3 surface $X$:
\begin{equation}
 {\rm Isom}^+(S_X) \cong W^{(2)}(S_X) \rtimes {\rm Isom}(S_X)^{({\rm
  Amp})}\, . 
\label{eq:IsomSx-W-Amp}
\end{equation}
Here, $W^{(2)}(S_X)$ is the subgroup of ${\rm Isom}^+(S_X)$ generated by 
reflections associated with $(-2)$-curves, and 
${\rm Isom}(S_X)^{({\rm Amp})}$ is the subgroup of ${\rm Isom}^+(S_X)$
preserving the ample cone. Note that the reflection hyperplane of
an irreducible genus $g$ curve passes through ${\rm Pos}_X^+$ only when 
$g = 0$. Thus, the ample cone is the fundamental region of the action of $W^{(2)}(S_X)$ on 
$\overline{{\rm Pos}^+_X}$.

\vspace{5mm}

Proposition $\eta$ implies that the image of the injective
homomorphism from ${\rm Aut}(X)$ to ${\rm Isom}(H^2(X ; \Z))$ is
contained in the following subgroup
\begin{equation}
 {\rm Isom}(T_X)^{({\rm Hodge})} \times {\rm Isom}(S_X)^{({\rm Amp})}
  \subset {\rm Isom}(H^2(X ; \Z)), 
\label{eq:TxS-subgroup-IsomH2XZ}
\end{equation}
where ${\rm Isom}(T_X)^{({\rm Hodge})}$ is the subgroup of 
${\rm Isom}(T_X)$ preserving the period $[\Omega_X] \in \P [T_X \otimes \C]$.

The image of ${\rm Aut}(X)$ under the injective homomorphism into 
${\rm Isom}(T_X)^{({\rm Hodge})} \times {\rm Isom}(S_X)^{({\rm Amp})}$ 
can be specified precisely by using the language of discriminant form
(reviewed in section \ref{ssec:prelim-lattice}). Consider group 
homomorphisms 
\begin{eqnarray}
 p_T: && {\rm Isom}(T_X) \longrightarrow {\rm Isom}(G_{T_X}, q_{T_X}),
\label{eq:def-pT} \\
 p_S: && {\rm Isom}(S_X) \longrightarrow {\rm Isom}(G_{S_X}, q_{S_X}).
\label{eq:def-pS}
\end{eqnarray}
The latter homomorphism factors through ${\rm Isom}(S_X)^{({\rm Amp})}
\times \{ \pm {\rm id.}_{S_X}\}$, because $W^{(2)}(S_X)$ is in the kernel of 
$p_S$. Proposition $\eta$ implies that ${\rm Aut}(X)$ is characterized 
as the fibre product:
\begin{equation}
 \begin{array}{rcl}
  & {\rm Aut}(X) & \\
 \swarrow & & \searrow \\
{\rm Isom}(T_X)^{({\rm Hodge})} & & {\rm Isom}(S_X)^{({\rm Amp})}, \\
 \searrow & & \swarrow \\
 & {\rm Isom}(G,q) &
 \end{array}
\end{equation}
where the isomorphic groups ${\rm Isom}(G_{T_X}, q_{T_X})$ and 
${\rm Isom}(G_{S_X}, q_{S_X})$ are simply denoted by ${\rm Isom}(G,q)$. 
That is, 
\begin{equation}
 {\rm Aut}(X) = {\rm Isom}(T_X)^{({\rm Hodge})} 
\times_{{\rm Isom}(G,q)} {\rm Isom}(S_X)^{({\rm Amp})}. 
\label{eq:AutX-fibr-prod}
\end{equation}
See e.g. section 1.5 of \cite{Vin-2most} for a proof of this characterization. The images of the projection 
\begin{eqnarray}
\pi_T: && {\rm Aut}(X) \longrightarrow {\rm Isom}(T_X)^{({\rm
 Hodge})}, \\
\pi_S: && {\rm Aut}(X) \longrightarrow {\rm Isom}(S_X)^{({\rm Amp})}, 
\label{eq:def-piS}
\end{eqnarray}
associated with the fibre product (\ref{eq:AutX-fibr-prod}) are denoted 
by ${\rm Isom}(T_X)^{({\rm Hodge~Amp})}$ and 
${\rm Isom}(S_X)^{({\rm Amp~Hodge})}$, respectively. 

Finally, let us list up a couple of relations among the various groups 
discussed above. First, it follows from the arguments above that 
the automorphism group ${\rm Aut}(X)$ has the following structures: 
\begin{eqnarray}
 {\rm Aut}(X) & = &
    {\rm Ker} \; .  \; {\rm Isom}(S_X)^{({\rm Amp~Hodge})}, \\
 {\rm Aut}(X) & = & 
    {\rm Aut}_N(X) \; . \; {\rm Isom}(T_X)^{({\rm Hodge~Amp})} \, .
\end{eqnarray}
These relations are obtained from the projections $\pi_S$ and $\pi_T$, 
respectively, accompanied by the injectivity of the homomorphism 
from ${\rm Aut}(X)$ to the group (\ref{eq:TxS-subgroup-IsomH2XZ}). Here, 
${\rm Ker}$ is the kernel of 
$p_T: {\rm Isom}(T_X)^{({\rm Hodge~Amp})}
  \longrightarrow {\rm Isom}(q)$, 
and ${\rm Aut}_N(X)$ the kernel\footnote{This ${\rm Aut}_N(X)$ subgroup 
of ${\rm Isom}(S_X)^{({\rm Amp~Hodge})}$ is also regarded as a subgroup 
of ${\rm Aut}(X)$ given by the fibre product (\ref{eq:AutX-fibr-prod}), 
and is called {\it group of symplectic automorphisms of} $X$. 
Another definition of the ${\rm Aut}_N(X)$ subgroup of ${\rm Aut}(X)$ is 
that of automorphisms acting trivially on the complex line 
$[\Omega_X] \subset H^2(X; \C)$. These two definitions are equivalent 
because of \cite{Nik-auto}, Thm. 3.1.} of 
$p_S: {\rm Isom}(S_X)^{({\rm Amp~Hodge})} \longrightarrow {\rm Isom}(q)$. 
It also follows that 
\begin{equation}
 {\rm Isom}(S_X)^{({\rm Amp~Hodge})} = {\rm Aut}_N(X) \; . \; 
   \left[p_T\left({\rm Isom}(T_X)^{({\rm Hodge~Amp})}\right)\right].
\end{equation}
The group of Hodge isometries of $H_2(X; \Z)$, 
${\rm Isom}\left( H_2(X; \Z)\right)^{({\rm Hodge})}$, is isomorphic to 
\begin{equation}
 {\rm Isom}\left(H_2(X; \Z)\right)^{({\rm Hodge})} \cong 
 \left\{ \pm 1 \right\} \times 
 \left[ W^{(2)}(S_X) \rtimes {\rm Aut}(X) \right].
\end{equation}
%

%%%%%%%%%%%%%%%%%%%%%%%%%%%%%%%%%%%%%%%%%%%%%%%%%%%%%%%%%%%%%%%%%%%%
\section{Classifications of Elliptic Fibrations on a K3 Surface}
\label{sec:3-classifications}
%%%%%%%%%%%%%%%%%%%%%%%%%%%%%%%%%%%%%%%%%%%%%%%%%%%%%%%%%%%%%%%%%%%%

%%%%%%%%%%%%%%%%%%%%%%%%%%%%%%%%%%%%%%%%%%%%%%%%%%%%%%%%%%%%%%%
\subsection{The ${\cal J}_0(X)$ Classification}
%%%%%%%%%%%%%%%%%%%%%%%%%%%%%%%%%%%%%%%%%%%%%%%%%%%%%%%%%%%%%%%

Consider a K3 surface $X$ with a given complex structure, so
that $T_X$ and $S_X$ are already determined in $H_2(X; \Z)$.
Suppose there is an elliptic fibration\footnote{Throughout this
article, we only consider elliptic fibrations that have a section 
$\sigma: \P^1 \longrightarrow X$, even when we just refer to them as
``elliptic fibrations''.} $(\pi_X, \sigma; X, \P^1)$.
Then one can define a sublattice 
$U_* := {\rm Span}_\Z \{[F], [\sigma + F]\}$ of $S_X$, where $[F]$ and 
$[\sigma]$ are the classes represented by the fibre class and the 
section $\sigma$. This rank-2 lattice is isomorphic to the hyperbolic 
plane lattice $U$. Since $U$ is unimodular, an orthogonal decomposition 
can be introduced to the Neron--Severi lattice $S_X$: 
\begin{equation}
 S_X \cong U_* \oplus W_{{\rm frame}*}.
\label{eq:canon-frame-in-NS}
\end{equation}
Here, $W_{{\rm frame}*}:= [U_*^\perp \subset S_X]$ 
is called the {\it canonical frame (sub)lattice of the elliptic fibration} 
$(\pi_X,\sigma; X, \P^1)$. 

For an elliptic fibration $(\pi_S, \sigma; S, C)$ of an algebraic surface, 
it is common to define the frame lattice $W_{\rm frame}$ as follows:
\begin{equation}
 W_{\rm frame} := \left[ [F]^\perp \subset S_S \right] / \vev{ [F] },
\label{eq:frame-def}
\end{equation}
where $S_S$ is the Neron--Severi lattice of $S$ and $[F]$ the fibre class. 
The canonical frame sublattice $W_{{\rm frame}*}$ in $S_X$ for a K3 
surface $X$ is isomorphic to the $W_{\rm frame}$ defined as above.

\vspace{5mm}
Let us call any embedding $\phi_*$ of the hyperbolic plane lattice $U$ 
into $S_X$ of a K3 surface $X$ a {\it canonical embedding} when 
the following conditions are satisfied: for generators $f_1$ and $f_2$ 
of $U$ satisfying the symmetric pairing (\ref{eq:U-int-form}), 
$\phi_*(f_1)$ is in $\overline{{\rm Amp}_X}$, and $\phi_*(f_2-f_1)$ 
is a class of an irreducible $(-2)$ curve. For any elliptic fibration 
$(\pi_X, \sigma; X, \P^1)$, there is a canonical embedding of $U$ given 
by $\phi: f_1 \longmapsto [F]$ and $f_2-f_1 \longmapsto [\sigma]$, so 
that $\phi: U \longmapsto U_* \subset S_X$. The converse is also true:

\vspace{5mm}

{\bf Theorem} $\theta$ (\cite{PSS}, \S 3 Thm. 1; 
\cite{Kondo-triv-Pic}, Lemma 2.1):
Whenever one finds a canonical embedding of $U$ into $S_X$, 
one finds that $h^0(X; {\cal O}_X(\phi(f_1))) = 2$, and 
an elliptic fibration morphism associated with the 
complete linear system of $\phi_*(f_1)$ is provided by 
$\Phi_{|\phi(f_1)|}: X \longrightarrow \P^{(2-1)=1}$. Here,
$\phi(f_1)$ becomes the fibre class and the class $\phi_*(f_2-f_1)$ can be taken as a section of this
fibration.
\vspace{5mm}

For a given K3 surface $X$, let ${\cal J}_{0}(X)$ be the set of all 
elliptic fibrations, $(\pi_X,\sigma; X, \P^1)$. Two elliptic 
fibration morphisms $\pi_X: X\longrightarrow \P^1$ and 
$\pi'_X: X\longrightarrow \P^1$ are regarded equivalent in the
 ${\cal J}_{0}(X)$ classification, if they are different only by 
$\pi'_X = g \cdot \pi_X$ for some ${}^{\exists} g \in PGL(2; \C)$ 
on $\P^1$. Theorem $\theta$ implies that ${\cal J}_{0}(X)$ is 
characterized in terms of canonical embeddings of the hyperbolic plane lattice
as follows:
\begin{equation}
 {\cal J}_{0}(X) \cong \left\{ {\rm canon.~embedding~} \phi_*:
  U \hookrightarrow S_X \right\}. 
%  /
%   \left(\Z/2\Z \vev{\sigma_U}\right);  
\label{eq:J0-1}
\end{equation}

The classification problem of ${\cal J}_0(X)$ can be translated into 
pure lattice theory language, so that the problem no longer involves 
geometric conditions such as $\phi(f_1)$ being in $\overline{{\rm Amp}_X}$ 
or $\phi_*(f_2-f_1)$ being the class of an irreducible curve. Consider 
\begin{equation}
 {\cal J}_{0}(X) \longrightarrow W^{(2)}(S_X) \backslash 
  \left\{ {\rm embedding~} \phi: U \hookrightarrow S_X \right\}
 / \left\{ \pm {\rm id.}_U \right\}. % (\Z/2\Z \vev{-{\rm id.}_U}).
\label{eq:J0-J0}
\end{equation}
This map is surjective. In fact
\begin{itemize}
\item {[\cite{PSS}, \S 6 Thm. 1]} For any embedding $\phi$, one can 
exploit the freedom to choose $\delta \in \{\pm {\rm id.}_U\}$ to make 
sure that $\phi(\delta(f_1))$ is in $\overline{{\rm Pos}^+_X}$, and then 
there exists an element in $w \in W^{(2)}(S_X)$ such that 
$w \cdot \phi \cdot \delta(f_1)$ is in $\overline{{\rm Amp}_X}$. 
\item {[\cite{PSS}, \S 6 Thm. 1; \cite{Kondo-triv-Pic}, 
Lemma 2.1]}: 
When the choice $\delta$ and $w$ are made as above, the class 
$w \cdot \phi \cdot \delta (f_1)$ contains a smooth curve of genus 1, 
and $w \cdot \phi \cdot \delta (f_2-f_1)$ is an effective curve with 
self-intersection $(-2)$. Although this $(-2)$ curve is not necessarily 
irreducible, one can always find an irreducible component $C_0$ 
of $w \cdot \phi \cdot \delta (f_2-f_1)$ such that 
$C_0$ is a class of smooth curve of genus 0, satisfies 
$(C_0, w \cdot \phi \cdot \delta (f_1))=1$, and the multiplicity of 
$C_0$ in $w \cdot \phi \cdot \delta (f_2-f_1)$ is 1.  
\item {[{\bf Lemma A} (see below for a proof)]} 
There exists an element in $w' \in W^{(2)}(S_X)$ 
that keeps $w \cdot \phi \cdot \delta (f_1)$ invariant, while  
it maps $w \cdot \phi \cdot \delta(f_2 - f_1)$ to $C_0$.  
\end{itemize}
All of the above combined implies that for any embedding 
$\phi: U \hookrightarrow S_X$, one can always find an element 
$w' \cdot w \cdot \phi \cdot \delta$ in the orbit 
$W^{(2)}(S_X) \cdot \phi \cdot \{\pm {\rm id.}_U\}$ that satisfies 
all the properties required for a canonical embedding  
$\phi_*: U\hookrightarrow S_X$. 

\begin{quote}
{\bf Proof} of Lemma A: We show that there is an algorithm of finding an 
appropriate $w' \in W^{(2)}(S_X)$. \\
Because of Theorem $\theta$, 
% \S 3 Theorem 1 of \cite{PSS} and 
% Lemma 2.1 of \cite{Kondo-triv-Pic}, 
we can assume that there is an 
elliptic fibration $\pi: X \rightarrow \P^1$, where 
$w \cdot \phi\cdot \delta(f_1) =: [F]$ is the fibre class, and $C_0$ is
a section. Now, note that all the irreducible components of the effective 
divisor $w\cdot \phi \cdot \delta (f_2-f_1)$ except $C_0$ are
either the fibre class or irreducible components in the singular fibres 
of the elliptic fibration. Thus we can write 
\begin{equation}
 w\cdot \phi \cdot \delta (f_2-f_1)  = C_0 +
   \sum_I \left[ m^{(I)} [F] +  \lambda^{(I)} \right],
\label{eq:wphidel-f2-f1-dcmp}
\end{equation} 
where $I$ runs over all the singular fibres of $A$-$D$-$E$ type in the
elliptic fibration, and $\lambda^{(I)}$ is an element of the $A$-$D$-$E$ type
root lattice generated by irreducible components of the $I$-th fibre not 
meeting the section $C_0$. (We do not have to include the irreducible 
components {\it meeting} the section in the decomposition above, 
because they are given by linear combinations of elements in the root 
lattice and $[F]$.) Because $(f_2-f_1)^2 = -2$, 
the coefficient $m^{(I)}$ above should be chosen as 
$(\lambda^{(I)})^2 = -2m^{(I)}$. \\
Let us now consider a subgroup of $W^{(2)}(S_X)$ generated by reflections 
associated with $\alpha^{(I)}$ and $[F]-\alpha^{(I)}$, where $\alpha^{(I)}$ 
are any roots in the $A$-$D$-$E$ root lattice of the $I$-th fibre. Under
this subgroup of $W^{(2)}(S_X)$, the fibre class $[F]$ remains unchanged and 
$w \cdot \phi \cdot \delta (f_2-f_1)$ also remains in the form 
of (\ref{eq:wphidel-f2-f1-dcmp}), except that $\lambda^{(I)}$ (and
correspondingly $m^{(I)}$) may change. In the following, we show that
$\lambda^{(I)}$ can be brought to zero by the subgroup of reflections we
have just introduced, so that the assertion for the existence
of $w' \in W^{(2)}(S_X)$ is proved. This can be done separately for each $I$-th
fibre. \\
 To this end, it is sufficient to note that a reflection by
 $\alpha^{(I)}$ followed by another one associated with $[F]-\alpha^{(I)}$ 
turns (\ref{eq:wphidel-f2-f1-dcmp}) into 
\begin{equation}
 C_0 + \sum_I \left[ (m^{(I)}+(\alpha^{(I)},\lambda^{(I)})+1) [F] +
  (\lambda^{(I)}-\alpha^{(I)}) \right],
\end{equation}
which is to change $\lambda^{(I)}$ into 
$\lambda^{(I)'} := \lambda^{(I)}-\alpha^{(I)}$, and 
correspondingly $2m^{(I)} := - (\lambda^{(I)})^2$ into 
$2(m^{(I)})' := - (\lambda^{(I)'})^2 = 
2m^{(I)}+2(\alpha^{(I)},\lambda^{(I)})+2$.
Since $\lambda^{(I)}$ is given by an integer linear combination of 
root vectors in the root lattice for the $I$-th fibre, 
$\lambda^{(I)'}$ becomes zero eventually. $\blacksquare$
\end{quote}

{\bf Lemma B}: The map (\ref{eq:J0-J0}) is also injective.
\begin{quote}
{\bf proof}: If there are two canonical embeddings 
$\phi_* :U \longrightarrow S_X$ and $\phi_*':U \longrightarrow S_X$ 
that fall within a common orbit of 
$W^{(2)}(S_X) \cdot \phi \cdot \{\pm {\rm id.}_U\}$ for some $\phi$, 
then there is an element $w'' \in W^{(2)}(S_X)$ such that
$w'' \cdot \phi_* = \phi'_*$. In the following, we prove that 
$w'' \cdot \phi_* = \phi'_*$ is the same embedding as $\phi_*$.
\\
Let us first see that $\phi'_*(f_1) = \phi_*(f_1)$, when $\phi_*$ and 
$\phi'_* = w'' \cdot \phi_*$ are both canonical embeddings of $U$ into
$S_X$. This can be proved by contradiction as follows:
assume that $\phi_*(f_1)$ and $\phi'_*(f_1)$ are different and
consider the following subsets\footnote{
Here, the group $W^{(2)}(S_X)$ is generated by a set of reflections 
$\{ r_i \}$ (called {\it simple reflections}) associated with a set of 
$(-2)$-roots, $\{\alpha_i\}$ (called {\it simple roots}), and the
fundamental chamber of this discrete reflection symmetry group is 
bounded by hyperplanes $\{H_i \}$ corresponding to those reflections.
Any element $w \in W^{(2)}(S_X)$ is written as a product of simple
reflections and $l(w)$ is the minimum number of simple reflections 
needed in obtain $w$.
} of $W^{(2)}(S_X)$ for 
$n \in \mathbb{N}_{>0}$:
\begin{equation}
 {\cal W}_n := \left\{ w \in W^{(2)}(S_X) \; | \; 
    w \cdot \phi_*(f_1) = \phi'_*(f_1), \quad 
    l(w) = n \right\}.
\end{equation}
Let $n_0$ be the smallest number of $n$ where ${\cal W}_n$ is not empty, 
and choose any element $w_0 \in {\cal W}_{n_0}$. $w_0$ can be written as a
successive application of $n_0$ simple reflections, 
$r_{i_{n_0}} \cdot \cdots \cdot r_{i_2} \cdot r_{i_1}$. 
Since we have assumed that 
$\phi_*(f_1) \neq \phi'_*(f_1)$, at least one of the $n_0$ reflection
planes does not contain $\phi_*(f_1)$. Let $H_{i_k}$ be the first 
reflection plane of that kind appearing in the sequence of simple 
reflections in $w_0$. 
Then both $\phi_*(f_1)$ and $\phi'_*(f_1) = w_0 \cdot \phi_*(f_1)$ are 
on the same side of $H_{i_k}$ (because both $\phi_*$ and $\phi'_*$ are 
a canonical embedding of $U$), yet the path connecting them passes through 
the other side of $H_{i_k}$ during the way. This means that there is 
a short-cut path; there is another $w \in W^{(2)}(S_X)$ with 
$0 < l(w) < n_0$. This is a contradiction\footnote{This proof is an easy modification of the proof of simple transitivity of 
reflection groups. See e.g. the proof for Theorem 11.6 
in \cite{mirror-reflection} or any other textbook on Coxeter groups.}. 
\\
The element $w'' \in W^{(2)}(S_X)$ therefore belongs to the stabilizer 
subgroup of $\phi_*(f_1) = \phi'_*(f_1)$. 
This stabilizer subgroup is
% see section ``note on Coxeter group'' } 
once again a reflection group. As a set of generators, we can take 
the reflections $s_\alpha$ associated with $(-2)$ curves $\alpha$
satisfying $(\alpha, \phi_*(f_1)) = 0$.
In the present context, this means that we can take reflections associated with $\alpha^{(I)}$'s 
and $[F]-\alpha^{(I)}$'s (in the elliptic fibration corresponding to the canonical embedding $\phi_*(U)$) 
as the generators of this stabilizer group. It is now obvious that 
this group cannot map a section $\phi_*(f_2-f_1)$ of this elliptic 
fibration to another section. Hence 
$w'' \cdot \phi_*(f_2-f_1)=\phi_*(f_2-f_1)$.  
Now, we have seen that $\phi_*$ and $\phi'_* = w'' \cdot \phi_*$ are the same 
embedding of $U$ into $S_X$. $\blacksquare$
\end{quote}

Therefore, the classification ${\cal J}_{0}(X)$ has been translated 
into pure lattice theory language:
\begin{equation}
  {\cal J}_{0}(X) \cong W^{(2)}(S_X) \backslash 
  \left\{ {\rm embedding~} \phi: U \hookrightarrow S_X \right\}
 / \left\{ \pm {\rm id.}_U \right\}. % (\Z/2\Z \vev{-{\rm id.}_U}).
\label{eq:J0-2}
\end{equation}
%

%%%%%%%%%%%%%%%%%%%%%%%%%%%%%%%%%%%%%%%%%%%%%%%%%%%%%%%%%%%%%%%%
\subsection{The ${\cal J}_2(X)$ and ${\cal J}^{({\rm type})}(X)$ 
Classifications}
%%%%%%%%%%%%%%%%%%%%%%%%%%%%%%%%%%%%%%%%%%%%%%%%%%%%%%%%%%%%%%%%

For any elliptic fibration $(\pi_X, \sigma; X,\P^1)$ of a K3 surface
$X$, there is an associated frame lattice $W_{{\rm frame}*} \cong W_{\rm frame}$, which
can equivalently be defined by \eqref{eq:canon-frame-in-NS} or \eqref{eq:frame-def}. Any frame lattice
$W_{\rm frame}$ for a K3 surface $X$ has the following three properties:
i) it is an even lattice, ii) it has signature $(0,\rho_X-2)$, and
finally iii) there exists an Abelian group isomorphism 
$G_{W_{\rm frame}} \cong G_{S_X}$ such that $q_W = q_{S_X}$. 
This motivates to consider a set of lattices with these
three properties modulo isometry for a K3 surface $X$:
\begin{equation}
 {\cal J}'_{2}(X) := {\rm Isom} \backslash 
 \left\{ W \; | \; 
  {\rm even~}
  {\rm sgn}(W)=(0,\rho_X-2), \quad
  (G_{W},q_W) {}^{\exists} \cong (G_{S_X},q_{S_X}) \right\}. 
\label{eq:J2-prime}
\end{equation}
A map ${\cal J}_{0}(X) \longrightarrow {\cal J}'_{2}(X)$ between the two
classifications is given by 
\begin{eqnarray}
 {\cal J}_{0}(X) \ni 
 \left( \phi_*: U \hookrightarrow S_X \cong \phi_*(U) \oplus W_{{\rm
  frame}*} \right) & \longmapsto & [W_{{\rm frame}*}] \in {\cal J}'_2(X), \\
 {\cal J}_{0}(X) \ni \left[ \phi: U \hookrightarrow S_X \right] 
  & \longmapsto & 
  \left[ (\phi(U))^\perp \subset S_X \right] \in {\cal J}'_{2}(X),
\end{eqnarray}
where the first line used the language of (\ref{eq:J0-1}) and the second
line that of (\ref{eq:J0-2}).

The map ${\cal J}_{0}(X) \longrightarrow {\cal J}'_2(X)$ is not
necessarily injective, but always factors through 
\begin{equation}
 {\cal J}_2(X) := {\rm Isom}^+(S_X) \backslash
   \left\{ {\rm embedding~} \phi: U \hookrightarrow S_X \right\} / 
  \left\{ \pm {\rm id.}_U \right\} \cong 
 {\rm Isom}(S_X)^{({\rm Amp})} \backslash {\cal J}_{0}(X)\, ,
\label{eq:J2}
\end{equation}
where we used (\ref{eq:IsomSx-W-Amp}) in order to obtain the last
expression. That is, 
\begin{equation}
 {\cal J}_{0}(X) \twoheadrightarrow 
 \left[ {\rm Isom}(S_X)^{({\rm Amp})} \backslash {\cal J}_{0}(X)
  = {\cal J}_{2}(X) \right] \hookrightarrow {\cal J}'_{2}(X)\, .
\end{equation}
It is obvious that two elements of ${\cal J}_{0}(X)$ identified in 
${\cal J}_{2}(X)$ provide frame lattices that are isometric to each other. 
Conversely, for two embeddings $\phi: U \hookrightarrow S_X$ and 
$\phi': U \hookrightarrow S_X$, if there is an isometry $\psi: W
\rightarrow W'$, where $W := [\phi(U)^\perp \subset S_X]$ and 
$W' := [\phi'(U)^\perp \subset S_X]$, we can construct an isometry of 
$S_X$ by $(\phi' \cdot \phi^{-1},\psi): S_X \cong \phi(U)\oplus W
\longrightarrow \phi'(U) \oplus W' \cong S_X$.

Furthermore, any lattice $W$ characterized as in (\ref{eq:J2-prime})
can be isomorphic to the frame lattice of some elliptic fibration in
$X$, when ${\rm rank}(T_X) \leq {\rm rank}(S_X)- 2$, or equivalently
$\rho_X \geq 12$ \cite{Nish01, Sch}. Hence 
${\cal J}_{2}(X) \hookrightarrow {\cal J}'_2(X)$ is also surjective 
for a K3 surface with $\rho_X \geq 12$, and 
\begin{equation}
 {\cal J}_2(X) \cong {\cal J}'_2(X).
\end{equation}
This is because, for a given $W$ in (\ref{eq:J2-prime}), 
$T_X \oplus (U \oplus W)$ can be embedded into the even 
unimodular lattice ${\rm II}_{3,19} \cong \Lambda_{\rm K3}$ 
(Proposition $\beta$),
% (Prop. 1.6.1 of \cite{Nik01}), 
from which a primitive embedding $T_X \hookrightarrow \Lambda_{\rm K3}$ 
is given. Theorem $\zeta$
% Cor 2.10 of \cite{Mor} 
guarantees that such a primitive embedding is unique modulo 
${\rm Isom}(\Lambda_{\rm K3})$ if ${\rm rank}(T_X)\leq {\rm rank}(S_X)-2$. Thus there must be 
an isometry $[T_X^\perp \subset \Lambda_{\rm K3}] =: 
S_X \cong (U \oplus W)$, and this $W$ is characterized as the orthogonal
complement of some embedding $U \hookrightarrow S_X$.

Note that both classifications ${\cal J}_{2}(X)$ in (\ref{eq:J2}) and 
${\cal J}'_{2}(X)$ in (\ref{eq:J2-prime}) are characterized purely in the language of 
lattice theory without reference to geometry.
Note also
that when two elliptic fibrations $\pi_{X}$ and $\pi'_X$ 
fall into the same element in the 
${\cal J}_2(X) \hookrightarrow {\cal J}'_2(X)$ classification, 
the lattice isometry between the canonical frame lattices of $\pi_X$
and $\pi'_X$ is always lifted to an ample-cone preserving isometry of 
$S_X$ such that $U_*$ in $U_* \oplus W_* \cong S_X$ (resp. $W_*$ in
$S_X$) for $\pi_X$ is mapped to $U'_*$ in $U'_* \oplus W'_* \cong S_X$ 
(resp. $W'_*$ in $S_X$).

Closely related to the ${\cal J}_2(X)$ classification of elliptic
fibration is the ${\cal J}^{({\rm type})}(X)$ 
classification. For an elliptic fibration $(\pi_S, \sigma; S, C)$, 
its {\it type} 
is the data of how many singular fibres of a given type in the
Kodaira classification are present in $\pi_S: S \longrightarrow C$. This can be
expressed in the form of $n_1 I_1 + n_2 I_2 + \cdots + m_1 I_1^* + \cdots$.
Two elliptic fibrations $(\pi_S, \sigma)$ and $(\pi'_S, \sigma')$ for 
a given $(S, C)$ are regarded to {\it be of the same type}, when the type
data of these two elliptic fibrations are the same and
${\cal J}^{({\rm type})}$ is the set of all possible types of elliptic 
fibration of a K3 surface $X$. 

The ${\cal J}_2(X)$ and ${\cal J}_2'(X)$ classifications are
very close to the ${\cal J}^{({\rm type})}(X)$ classification, in that 
type data can almost be read out from the root lattice part of the frame 
lattice, $W_{\rm root}$.
The root lattice $W_{\rm root}$ is in the form of $\oplus_a R_a$, where 
$R_a$ is one of $A$--$D$--$E$ lattice. An $I_{n+1}$ type singular fibre
gives rise to an $A_n$ component in $W_{\rm root}$, $I_{n-4}^*$ type singular fibre 
to a $D_n$ component, ${\rm IV}^*$ to $E_6$ ,${\rm III}^*$ to $E_7$, ${\rm II}^*$ to 
$E_8$, ${\rm III}$ to $A_1$ and ${\rm IV}$ to $A_2$. On the one hand, the frame lattice
hence misses certain information, because both $I_2$ and ${\rm III}$ 
type singular fibres appear as an $A_1$ component in $W_{\rm root}$, and both $I_3$ and ${\rm IV}$ type singular 
fibres as an $A_2$ component. It is not possible either, at least immediately, 
to read out the numbers of $I_1$ and ${\rm II}$ type singular fibres
from the frame lattice. The frame lattice, on the other hand, contains 
more data---such as the Mordell--Weil group---than those we can read out 
from its sublattice $W_{\rm root} \subset W_{\rm frame}$. 
In section \ref{ssec:new&known}, we will come across examples 
of mutually non-isometric frame lattices $W_{\rm frame}$ of a given K3
surface $X$ sharing the same $W_{\rm root}$. 

Apart from these subtle differences, the two classifications of 
elliptic fibrations ${\cal J}_2(X)$ and ${\cal J}^{({\rm type})}(X)$ 
are still quite close. For this reason, there are some authors
where individual isometry classes of frame lattices in ${\cal J}_2(X)$
and ${\cal J}_2'(X)$ are referred to as types.

%%%%%%%%%%%%%%%%%%%%%%%%%%%%%%%%%%%%%%%%%%%%%%%%%%%%%%%%%%%%%%%%%%
\subsection{The ${\cal J}_1(X)$ Classification}
%%%%%%%%%%%%%%%%%%%%%%%%%%%%%%%%%%%%%%%%%%%%%%%%%%%%%%%%%%%%%%%%%%

It is definitely a question of mathematical interest to consider the
classification of elliptic fibrations $\left(\pi_X, \sigma; X, \P^1\right)$
modulo automorphisms of $X$. An equivalence relation 
is introduced to the set of all possible elliptic fibration data 
$(\pi_X, \sigma; X, \P^1)$ for a K3 surface $X$, as defined 
at the beginning of the Introduction. The quotient space with respect to 
this equivalence relation is denoted by ${\cal J}_1(X)$.
This modulo-isomorphism classification of elliptic fibrations is not 
just interesting as a problem in mathematics, but also the relevant 
classification scheme in F-theory compactification on a K3 surface $X$ 
in the context of string theory \cite{BKW-phys}. 

Since the image of 
$\pi_S: {\rm Aut}(X) \longrightarrow {\rm Isom}(S_X)^{({\rm Amp})}$
is defined to be ${\rm Isom}(S_X)^{({\rm Amp~Hodge})}$, the
classification ${\cal J}_{1}(X)$ is equivalent to 
\begin{equation}
 {\cal J}_{1}(X) = {\rm Isom}(S_X)^{({\rm Amp~Hodge})} \backslash
  {\cal J}_{0}(X),  
\label{eq:J1-1}
\end{equation}
which can be phrased purely in lattice theory language:\footnote{ 
Such conditions as ``ample-cone preserving'' or ``Hodge'' on the 
subgroup ${\rm Isom}(S_X)^{({\rm Amp~Hodge})}$ are already well-defined 
without invoking geometric intuition; the ample cone is one of the fundamental
region of the reflection group $W^{(2)}(S_X)$ acting on the lattice
$S_X$, and the ${\rm Isom}(T_X)^{({\rm Hodge})} \subset {\rm Isom}(T_X)$ 
subgroup is characterized as the stabilizer of a complex plane 
$[\Omega_X] \in \P[T_X^* \otimes \C]$.}
\begin{equation}
 {\cal J}_{1}(X) = \left[ 
   W^{(2)}(S_X) \rtimes {\rm Isom}(S_X)^{({\rm Amp~Hodge})} \right]
 \backslash
  \left\{ {\rm embed.~}\phi: U \hookrightarrow S_X \right\} /
 \{ \pm {\rm id.}_U \}.
\label{eq:J1-2}
\end{equation}

The quotient maps ${\cal J}_{0}(X) \longrightarrow {\cal J}_{2}(X)$ and 
${\cal J}_{0}(X) \longrightarrow {\cal J}_{1}(X)$ already determine a map 
\begin{equation}
 {\cal J}_{1}(X) \longrightarrow {\cal J}_{2}(X)
\end{equation} 
automatically. The fibre of this map is the decomposition of 
a single orbit of ${\rm Isom}(S_X)^{({\rm Amp})}$ into orbits of 
the subgroup ${\rm Isom}(S_X)^{({\rm Amp~Hodge})}$.
Since ${\rm Isom}(S_X)^{({\rm Amp~Hodge})}$ is not
necessarily a normal subgroup of ${\rm Isom}(S_X)^{({\rm Amp})}$, 
however, ${\cal J}_{2}(X)$ is not necessarily regarded as a quotient 
of ${\cal J}_{1}(X)$.

%%%%%%%%%%%%%%%%%%%%%%%%%%%%%%%%%%%%%%%%%%%%%%%%%%%%%%%%%%%%%%%%%%%%%
\section{The ${\cal J}_2(X)$ Classification}
\label{sec:J2-classification}
%%%%%%%%%%%%%%%%%%%%%%%%%%%%%%%%%%%%%%%%%%%%%%%%%%%%%%%%%%%%%%%%%%%

%%%%%%%%%%%%%%%%%%%%%%%%%%%%%%%%%%%%%%%%%%%%%%%%%%%%
\subsection{The ${\cal J}_{2}(X)$ Classification using Niemeier Lattices \\
(the Kneser--Nishiyama method)}
\label{ssec:J2-Nishiyama}
%%%%%%%%%%%%%%%%%%%%%%%%%%%%%%%%%%%%%%%%%%%%%%%%%%%%

In \cite{Nish01, Nish02} a systematic way to study the ${\cal J}'_{2}(X)$ classification 
which is quite convenient in dealing with K3 surfaces with large Picard number was introduced.
Nishiyama's approach starts with a following observation:

\vspace{5mm}

For a K3 surface $X$ with transcendental lattice $T_X$, 
there exists a lattice $T_0$ with the three following properties: 
i) it is an even lattice, ii) its signature is $(0,26-\rho)$, and
iii) $(G_{T_0}, q_{T_0}) {}^{\exists}\cong (G_{T_X}, q_{T_X})$. 
This is because Theorem $\delta$ 
% Theorem 1.12.4 of \cite{Nik01} 
guarantees that there always exists a primitive embedding 
\begin{equation}
 T_X[-1]^{(20-\rho,2)} \hookrightarrow 
  \left(E_8 \oplus U^{\oplus (20-\rho)}\right)^{(20-\rho,28-\rho)},
\label{eq:embed-TX[-]-E8U}
\end{equation}
and the orthogonal complement of the $T_X[-1]$ sublattice satisfies
all the three properties of $T_0$ we mentioned above \cite{Nish02}.

\vspace{5mm}

Let us now pick one such a $T_0$ and consider the map 
\begin{equation}
 \amalg_{I=\alpha}^{\omega} \left[{\rm Isom}(L^{(I)}) \backslash 
  \left\{ {\rm prim.~embed.~}\phi_{T_0}: T_0 \hookrightarrow L^{(I)} \right\}
  / {\rm Isom}(T_0) \right] \longrightarrow {\cal J}'_{2}(X)\, ,
\label{eq:Niemeier-classify-4-J2'}
\end{equation}
where $L^{(I)}$ labelled by Greek letters 
$I = \alpha, \beta, \cdots, \omega$ are even unimodular negative
definite lattices of rank 24. As reviewed briefly 
in section \ref{ssec:Niemeier}, there are twenty-four mutually
non-isometric lattices of that kind, and they are called \emph{Niemeier
lattices}. For any primitive embedding of $T_0$ into any one of the
Niemeier lattices $L^{(I)}$, the orthogonal complement 
$[(\phi_{T_0}(T_0))^\perp \subset L^{(I)}]$ satisfies the properties required 
for an element in ${\cal J}'_2(X)$, see (\ref{eq:J2-prime}), and this 
is how the above map is determined. 
This map is well-defined and furthermore surjective (use Proposition $\beta$). 
This allows us to study the left-hand side 
of (\ref{eq:Niemeier-classify-4-J2'}) instead in order to find 
all the elements in ${\cal J}'_2(X)$ (see \cite{Nish01}, Lemma 6.3).

One can also show that the map (\ref{eq:Niemeier-classify-4-J2'}) is
injective if $p_{T_0}: {\rm Isom}(T_0) \longrightarrow {\rm Isom}(G_{T_0},q_{T_0})$
is surjective. This condition is satisfied in many examples of singular 
K3 surfaces ($\rho=20$) that we will be interested 
in (see Table \ref{tab:T0}). 
Thus, it is not too inefficient to study the left-hand side 
of (\ref{eq:Niemeier-classify-4-J2'}) instead of ${\cal J}'_2(X)$ itself. 

Overall, the problem of finding the following embeddings and working out 
orthogonal complements, 
\begin{equation}
 T_X^{(2,20-\rho)} \hookrightarrow H_2(X;\Z)^{(3,19)}, \qquad 
\phi_U: U^{(1,1)} \hookrightarrow S_X^{(1,\rho - 1)}
   := [(T_X)^\perp \subset H_2(X; \Z)] 
\label{eq:PSS-embedding}
\end{equation}
modulo isometry is now translated into the similar problem:
\begin{equation}
 T_X[-1]^{(20-\rho,2)} \hookrightarrow 
 \left(E_8\oplus U^{\oplus (20-\rho)}\right)^{(20-\rho,28-\rho)}, \quad 
 \phi_{T_0}:[T_0^{(0,26-\rho)} := (T_X[-1])^\perp] \subset (L^{(I)})^{(0,24)}.
\label{eq:Nishiyama-embedding}
\end{equation}
The (isometry class of) frame lattice $W_{\rm frame}$ is obtained either as 
$[(\phi_U(U))^\perp \subset S_X]$ or as 
$[(\phi_{T_0}(T_0))^\perp \subset L^{(I)}]$.  
The latter problem is suitable for systematic calculations for various K3
surfaces with large Picard number. The target of the embedding 
$\phi_{T_0}$, $L^{(I)}$, does not depend on individual choices of $T_X$ and
is also technically easier because of the (negative) definite 
symmetric pairing of $L^{(I)}$.

%%%%%%%%%%%%%%%%%%%%%%%%%%%%%%%%%%%%%%%%%%%%%%%%%%%%%%%%%%%%%%
\subsection{A Brief Review on Niemeier Lattices}
\label{ssec:Niemeier}
%%%%%%%%%%%%%%%%%%%%%%%%%%%%%%%%%%%%%%%%%%%%%%%%%%%%%%%%%%%%%%

Even unimodular lattices with signature $(r_+, r_-)$ are known to be 
unique modulo isometry, if both $r_+$ and $r_-$ are non-zero, as
reviewed in section \ref{ssec:prelim-lattice}. 
However, there can be more than one isometry classes of negative
definite (i.e., $r_+=0$) even unimodular lattices of given rank 
$r_- = 16, 24, \cdots$ (see \cite{CS-book} for more information).
 
For example, there are two mutually non-isometric negative definite 
even unimodular lattices of rank $r_- = 16$. One is $L = E_8 \oplus E_8$. 
The other one, $L'$, is characterized as an overlattice of $D_{16}$.
The discriminant group $G_{D_{16}} = D_{16}^* / D_{16}$ is isomorphic to
$\Z/2\Z \vev{{\bf sp}} \times \Z/2\Z\vev{\overline{\bf sp}}$. 
Both $\Z/2\Z$ subgroups generated by the weights of the spinor/cospinor
representations are isotropic subgroups of $(G_{D_{16}},q_{D_{16}})$  
($D_n$ with $n$ divisible by 8).
The lattice $L'=D_{16}; \Delta$ with $\Delta = \Z/2\Z\vev{{\bf sp}}$ or 
$\Delta=\Z/2\Z\vev{\overline{\bf sp}}$ (see Proposition $\alpha$) 
becomes the other even unimodular negative definite lattice of rank-16. 
$D_{16}; \Delta$ with the two different choices of $\Delta$ are  
isometric. That is, $L'$ is obtained by adding the weights of either the
spinor or cospinor representation to the root lattice of $D_{16}$.

Similarly, there are twenty-four isometry classes of negative definite 
even unimodular lattices of rank 24 that are called \emph{Niemeier lattices}.
Following \cite{CS-book} we denote them by $L^{(I)}$ with 
$I = \alpha, \beta, \cdots, \psi, \omega$. Out of those twenty-four 
rank-24 lattices, twenty-three ($I \neq \omega$) allow a
description like the one just given for $L'$. They can be specified by
% Those twenty-three lattices, $L^{(I)}$ with a label $I = \alpha, \beta, 
% \cdots, \psi$, are all characterized as an overlattice of some root
choosing of a rank-24 lattice given as the direct sum of
$A$--$D$--$E$ lattices $\oplus_a R_a = L^{(I)}_{\rm root}$, and an appropriate 
isotropic subgroup $\Delta$ of $G_{\oplus_a R_a}$. For a fixed Niemeier lattice, 
the dual Coxeter numbers of all components $R_a$ are equal and denoted by $h^{(I)}$. 
This information is displayed in Table \ref{tab:Niemeier-root}.
There is one more rank-24 even unimodular negative definite lattice, which is called the Leech lattice and
is denoted by $L^{(I=\omega)} = \Lambda_{24}$, that does not 
have any norm $(-2)$ points and correspondingly cannot be described in a similar fashion to
the other Niemeier lattices. More information on Leech lattice is 
provided in our review in the appendix \ref{ssec:Leech-review}.
%%%%%%%%%%%%%%%%%%%%%%%%%%%%%%%%%%%%%%%%%%%%%%%%%%%%%%%%%%%%%%%
\begin{table}[tb]
\begin{center}
%%%%%%%%%%%%%%%%%%%%%%%%%%%%%%%%%%%%%%%%%%%%%%%%%%%%%%%%%%
  \begin{tabular}{|c|r|c|c||c|r|c|c|}
 \hline 
$I$ & $L^{(I)}_{\rm root}$ & $h^{(I)}$ & $\Delta$ &
$I$ & $L^{(I)}_{\rm root}$ & $h^{(I)}$ & $\Delta \subset G_{L_{\rm root}}$ \\
 \hline
 \hline
 $\alpha$ & $D_{24}$ & 46 & $({\bf sp})_2$ &
 $\nu$ & $A_9 A_9 D_6$ & 10 & $\Z_5(\Z_2)^2 \subset (\Z_{10})^2(\Z_2)^2$ \\
 $\beta$ & $D_{16}E_8$ & 30 & $({\bf sp})_2$ & 
 $\xi$ & $D_6^{\oplus 4}$ & 10 & $(\Z_2)^4\subset (\Z_2)^8$ \\
 $\gamma$ & $E_8 E_8 E_8$ & 30 & &  
 $o$ & $A_8 A_8 A_8$ & 9 & $\Z_9\Z_3 \subset (\Z_9)^3$ \\
 $\delta$ & $A_{24}$ & 25 & $(5)_5$ & 
 $\pi$ & $A_7^{\oplus 2}D_5^{\oplus 2}$ & 8 &
     $(\Z_8\Z_4) \subset (\Z_8)^2(\Z_4)^2$ \\
 $\epsilon$ & $D_{12}D_{12}$ & 22 &
   $({\bf sp},{\bf v})_2, \; ({\bf v},{\bf sp})_2$ & 
 $\rho$ & $A_6^{\oplus 4}$ & 7 & $(\Z_7)^2 \subset (\Z_7)^4$ \\
 $\zeta$ & $A_{17} E_7$ & 18 & $(3,1)_6$ &
 $\sigma$ & $A_5^{\oplus 4} D_4$ & 6 &
     $(\Z_2)^3(\Z_3)^2 \subset (\Z_6)^4(\Z_2)^2$ \\
 $\eta$ & $D_{10} E_7 E_7$ & 18 &
      $({\bf sp},1,0)_2, \; (\overline{\bf sp},0,1)_2$ &
 $\tau$ & $D_4^{\oplus 6}$ & 6 & $(\Z_2)^6 \subset (\Z_2)^{12}$ \\
 $\theta$ & $A_{15} D_9$ & 16 & $(2,{\bf sp})_8$ &
 $\upsilon$ & $A_4^{\oplus 6}$ & 5 & $(\Z_5)^3 \subset (\Z_5)^6$ \\
 $\iota$ & $D_{8} D_8 D_8$ & 14 &
	       $({\bf sp},{\bf v},{\bf v}),({\bf v},{\bf v},{\bf
	       sp}),({\bf v},{\bf sp},{\bf v})$ & 
 $\phi$ & $A_3^{\oplus 8}$ & 4 &  $(\Z_4)^4 \subset (\Z_4)^8$ \\
 $\kappa$ & $A_{12}A_{12}$ & 13 & $(1,5)_{13}$ &
 $\chi$ & $A_2^{\oplus 12}$ & 3 &  $(\Z_3)^{6} \subset (\Z_3)^{12}$ \\
 $\lambda$ & $A_{11} D_7 E_6$ & 12 & $(1,{\bf sp},1)_{12}$ & 
 $\psi$ & $A_1^{\oplus 24}$ & 2 & $(\Z_2)^{12} \subset (\Z_2)^{24}$ \\
 $\mu$ & $E_6^{\oplus 4}$ & 12 & $(1,1,1,0)_3, \; (0,2,1,1)_3$ &
 $\omega$ & ------ & --- & ($\Lambda_{24}$)  % & --------- 
\\
 \hline
 \hline
  \end{tabular}
  
  \caption{\label{tab:Niemeier-root}
The Niemeier lattices: out of the 24 Niemeier lattices $L^{(I)}$ with
$I=\alpha, \beta, \cdots, \psi, \omega$, twenty-three 
have a rank-24 root lattice, $L^{(I)}_{\rm root} = \oplus_a R_a$, while 
the root lattice for $L^{(I=\omega)}$ is empty. The root lattices are given 
in the second column.
The quotient $\Delta = L^{(I)}/L^{(I)}_{\rm root}$ must be an isotropic 
subgroup of the discriminant group $G_{\oplus_a R_a}$ for all the cases
except $I=\omega$.
Generators of $\Delta \subset G_{\oplus_a R_a}$ are specified in detail 
for the first twelve Niemeier lattices (in the left-hand side) 
in this table; subscripts attached to the generators, $m$ of $(g)_m$, 
carry the information that they generate a $\Z/m\Z\vev{g}$ factor in 
$\Delta$. 
For the next eleven Niemeier lattices, the limited space in this table is
not enough to include detailed information on the generators,
see Table 16.1 of \cite{CS-book} for more information. 
In this table, we have introduced abbreviated notations 
$R_1R_2 \cdots$ for $R_1 \oplus R_2 \oplus \cdots$, 
$\Z_m$ for $\Z/m\Z$, and $G_1 G_2 G_3 \cdots $ for 
$G_1 \times G_2 \times G_3 \times \cdots$, to save space.}
 \end{center}
\end{table}
%%%%%%%%%%%%%%%%%%%%%%%%%%%%%%%%%%%%%%%%%%%%%%%%%%%%%%%%%%%%%%%%

%%%%%%%%%%%%%%%%%%%%%%%%%%%%%%%%%%%%%%%%%%%%%%%%%%%%%%%%%
\subsection{Computing $T_0$}
%%%%%%%%%%%%%%%%%%%%%%%%%%%%%%%%%%%%%%%%%%%%%%%%%%%%%%%%%

The first step of the program reviewed 
in section \ref{ssec:J2-Nishiyama} (the Kneser--Nishiyama method) is to compute the lattice $T_0$ 
for a given K3 surface $X$ with the transcendental lattice $T_X$.
This is done by first embedding $T_X[-1]$ primitively into 
$E_8 \oplus U^{\oplus (20-\rho)}$ and then using $T_0=T_X[-1]^\perp$. 
Theorem $\epsilon$
% Theorem 1.14.4 of \cite{Nik01} and Theorem 2.8 of \cite{Mor} 
does not guarantee uniqueness (modulo isometry of 
$E_8 \oplus U^{\oplus (20-\rho)}$) of primitive embeddings of
$T_X[-1]$, but this is not a problem. The crucial point is that
\emph{any} $T_0$ satisfying the properties i)--iii) does the job of 
(\ref{eq:Niemeier-classify-4-J2'}) for a given K3 surface $X$ 
\cite{Nish01, Nish02}.  Hence we only need to find one primitive embedding 
of $T_X[-1]$ and can then use the resulting $T_0$ for a given $X$.

\newpage\noindent
For some K3 surfaces, possible choices of $T_0$ can be found in the literature:\footnote{The appendix \ref{ssec:review-Kummer} provides a quick
summary of various definitions, explanations and facts on Kummer
surfaces that we need in this article.} 
% [Nishiyama '96 $X_{4,3}, {\rm Km}(Eli \times Eli)$, 
% Kondo '98 ${\rm Km}(A)$, Museum-Kondo '01 ]
%
\begin{itemize}
 \item Reference \cite{Nish01} contains examples of 
       $T_0$ for the four singular K3 surfaces $(\rho_X = 20)$, 
       $X = X_4$, $X_3$, ${\rm Km}(E_i \times E_i)$ and 
       ${\rm Km}(E_\omega \times E_\omega)$, where 
       $T_0=D_6$, $E_6$, $A_3\oplus A_3$ and $D_4 \oplus A_2$, respectively.
       The same reference also contains the result $T_0 = D_4 \oplus A_3$ for the 
       slightly more general class of Kummer surfaces 
       $X = {\rm Km}(E \times E)$ (isogenous case), where $\rho_X = 19$ and
       $T_0=D_4 \oplus D_4$ for the even more general class of Kummer surfaces 
       $X = {\rm Km}(E \times F)$ (product type), where $\rho_X = 18$.
 \item The most general class of Kummer surfaces, $X = {\rm Km}(A)$,   
       are associated with generic Abelian surfaces $A$, and 
       $\rho_X = 17$. In this case, we can take $T_0$ as 
       $[A_3 \oplus A_1^{\oplus 6}]; H$ with the isotropic subgroup  
       given by \cite{Kondo-auto-JacobKumm}  
\begin{equation}
G_{A_3 \oplus A_1^{\oplus 6}} \cong (\Z/4\Z) \times (\Z/2\Z)^6 \supset 
  H \cong \Z_2.  
\end{equation}
\end{itemize}

In this article, we are interested in K3 surfaces with large Picard
number, particularly in singular K3 surfaces $(\rho_X = 20)$.
The present authors came up with a list of thirty-four singular K3
surfaces that can be used for a particular class of string 
compactifications \cite{BKW-phys}. Those thirty-four K3 surfaces 
are listed in the first column of Table \ref{tab:coset-alt-rho20}. 
To explain the notation used in the table and elsewhere in this article, 
note that when the transcendental lattice 
$T_X$ is given an oriented basis\footnote{We call $\{f_2, f_1\}$ an {\it oriented basis} of $T_X = {\rm Span}_\Z\{f_1,f_2 \}$ if 
${\rm Im} \left[ \langle \Omega_X, f_2 \rangle / \langle \Omega_X, f_1
\rangle \right] > 0$ for the complex structure of $X$, 
$[\Omega_X] \in \P[T_X^* \otimes \C]$. } $\{f_2, f_1\}$, and the symmetric 
pairing of $T_X$ is given by\footnote{We have adopted a convention of writing this matrix in the order of $f_1$ and $f_2$
(rather than $f_2$ and $f_1$) and parametrizing the matrix 
with $a$, $b$ and $c$ as in the main text, by following a literature 
in string theory \cite{AK}; thus, $[a~b~c]$ in \cite{AK} and also 
in this article should be read as $[c~b~a]$ in \cite{SI}. } 
\begin{equation}
 \left[ \begin{array}{cc}
    (f_1,f_1) & (f_1, f_2) \\ (f_2,f_1) & (f_2,f_2)
 \end{array} \right] =
 \left[ \begin{array}{cc}
  2a & b \\ b & 2c
	\end{array}\right], 
\end{equation}
this determines a unique singular K3 surface $X$ \cite{SI}. 
This K3 surface $X$ is denoted by $X_{[a~b~c]}$ in this article, and 
the transcendental lattice of $X_{[a~b~c]}$ may sometimes be abbreviated 
by a symbol $[a~b~c]$. 

We have computed a possible choice for $T_0$ for twenty-four singular K3 surfaces 
among the thirty-four listed in Table \ref{tab:coset-alt-rho20} 
and the result is presented in Table \ref{tab:T0}. Four among them are 
the same already obtained in \cite{Nish01}.\footnote{
$X_{[1~0~1]} = X_4$, $X_{[1~1~1]} = X_3$, $X_{[2~0~2]} = {\rm Km}(E_i
\times E_i)$ and $X_{[2~2~2]} = {\rm Km}(E_\omega \times E_\omega)$.}
There are eleven singular K3 surfaces for which $T_0$ can be taken to 
be a direct sum of $A$--$D$--$E$ root lattice. For the other thirteen 
cases, the result of the computation is presented by describing 
a sublattice given by $(T_0)_{\rm root}$ and its orthogonal complement 
in $T_0$ and supplying information of the quotient by 
this sublattice:
\begin{equation}
  \left( (T_0)_{\rm root} \oplus [(T_0)_{\rm root}^\perp \subset T_0] 
  \right); H, \qquad H := 
  T_0/  \left( (T_0)_{\rm root} \oplus [(T_0)_{\rm root}^\perp \subset T_0] 
  \right).
\end{equation}
We find this way of describing lattices convenient (partially with
application in physics in mind), and will use this also in 
Tables \ref{tab:J2classify-[abc=302]} and 
\ref{tab:J2classify-[abc=302]_Etype}.
%%%%%%%%%%%%%%%%%%%%%%%%%%%%%%%%%%%%%%%%%%%%%%%%%%%%%%%%%%%%%%%%%%
\begin{table}[tbp]
 \begin{center}
  \begin{tabular}{|c|c|c|c|c|}
\hline 
 [a b c] & $T_0$ & ${\rm Isom}(T_0)/W^{(2)}(T_0)$ & 
Isom($q_{T_0}$) & $p_{T_0}$surj.? \\
\hline
\hline
 [1 0 1] & $D_6$ & $\Z_2$ & $\Z_2$ & yes \\
{}[1 1 1] & $E_6$ & $\Z_2$ & $\Z_2$ & yes \\
{}[2 0 1] & $D_5A_1$  & $\Z_2$ & $\Z_2$ & yes \\
{}[2 1 1] & $A_6$ & $\Z_2$ & $\Z_2$ & yes \\
{}[3 0 1] & $A_5A_1$  & $\Z_2$ & $\Z_2$ & yes \\
{}[3 1 1] & $[D_5(-44)]; \Z_4$  &&&\\
{}[4 0 1] & $[A_5(-24)]; \Z_3$  &&&\\
{}[4 1 1] & $A_4A_2$  & $\Z_2 \times  \Z_2$ & $\Z_2 \times \Z_2$ & yes \\
{}[5 0 1] & $[D_4A_1(-10)]; \Z_2$  &&&\\
{}[5 1 1] & $[A_5(-114)]; \Z_6$  &&&\\
{}[6 0 1] & $A_3 A_2 A_1$  & $\Z_2 \times \Z_2$ & $\Z_2 \times \Z_2$
		   & yes \\
{[6 1 1]} & $[A_4A_1(-230)]; \Z_{10}$ &&& \\
{[2 0 2]} & $A_3 A_3$  & $D_4$ & $D_4$ & yes \\
{[2 1 2]} & $[A_5(-10)]; \Z_2$  &&&\\
{}[2 2 2] & $D_4A_2$  & $S_3 \times \Z_2$ & $\Z_2 \times S_3$ & yes \\
{}[3 0 2] & $A_5(-4)$  &&&\\
{}[3 1 2] & $[A_4A_1(-230)]; \Z_{10}$  &&&\\
{}[3 2 2] & $[D_5(-20)]; \Z_2$ &&& \\
{}[3 2 2] & $A_4 A_1^{\oplus 2}$ & $\Z_2 \times \Z_2$  & $\Z_2 \times \Z_2$ & yes\\
{}[4 0 2] & $[A_3A_1A_1(-8)]; \Z_2$ &&& \\
{}[4 2 2] & $[A_3A_2(-84)]; \Z_6$  &&&\\
{}[6 0 2] & $A_3A_2(-4)$  &&&\\
{}[6 2 2] & $[A_2A_1A_1A_1(-66)]; \Z_6$ &&& \\
{}[3 0 3] & $(A_2A_1)^{\oplus 2}$  & $D_4 \times \Z_2$ & $\Z_2
		   \times D_4$ & yes \\
{}[6 0 6] & $(A_2(-4))^{\oplus 2}$  &&&\\
\hline
\hline  
\end{tabular}
\caption{\label{tab:T0} The $T_0$ lattice (the 2nd column) for 
twenty-four singular K3 surfaces (the 1st column), with some additional 
information (the 3rd--5th columns) to be used 
in section \ref{sec:J1-classification}.
We use the same abbreviated notation as in Table \ref{tab:Niemeier-root}.
Out of the thirty-four singular K3 surfaces appearing in
\cite{BKW-phys} (and also in Table \ref{tab:coset-alt-rho20}), we have 
included all cases for which $T_0$ can be taken to be a root lattice. }
\end{center}
\end{table}
%%%%%%%%%%%%%%%%%%%%%%%%%%%%%%%%%%%%%%%%%%%%%%%%%%%%%%%%%%%%%%%%%

Let us briefly sketch how to compute $T_0$ in practice. The first step 
is to find a primitive embedding $\phi: T_X[-1] \hookrightarrow E_8$, which is not difficult for 
a singular K3 surface $X$ because the transcendental lattice $T_X$ is just 
a rank-2 lattice in this case. The $E_8$ lattice can be described 
as\footnote{The negative definite symmetric pairing in this $E_8$
lattice is given by $(L_i, L_j) = -\delta_{ij}$.} 
\begin{equation}
 E_8 \cong \left\{ \left. \sum_{i=1}^8 n_i L_i \; \right| \; 
   n_i \in \Z, \; \sum_i n_i \equiv 0 {\rm ~mod~}2 \right\} \cup
  \left\{ \left. \sum_{i=1}^8 \left(\frac{1}{2}+m_i\right) L_i \; \right| \; 
   m_i \in \Z, \; \sum_i m_i \equiv 0 {\rm ~mod~}2  \right\}.
\label{eq:E8-description}
\end{equation}
In case $2a$ (or $2c$) is $2$, the images of $f_1$ (or $f_2$) should be chosen from the 
the $(-2)$-roots of $E_8$. If $f_1$ or $f_2$ are norm 4 vectors in $T_X$ (i.e., $2a$ or $2c$ is 4), 
then the image should be chosen from one of the vectors of the form 
\begin{equation}
 (\pm 2,0^7), \quad ((\pm 1)^4,0^4), \quad (\pm 3,(\pm 1)^7)/2.
\end{equation}
Norm $(-6)$ elements in $E_8$ are of the form 
\begin{equation}
 ( \pm 2,(\pm 1)^2, 0^5), \quad ((\pm 1)^6,0^2), \quad 
 ( (\pm 3)^2, (\pm 1)^6)/2,
\end{equation}
which can be used as the images of $f_{1,2}$ if either $2a$ or $2c$ 
is 6. To take the case [a~b~c]=[5~1~1] as an example, we can take 
\begin{eqnarray}
\phi(f_2) & = &  (L_1+\cdots + L_8))/2 = (1,\cdots,1)/2, \nonumber \\
\phi(f_1) & = & 2L_1 + (L_2+L_3+L_4)-(L_5+L_6+L_7) = (2,1^3,(-1)^3,0), 
\end{eqnarray}
so that $(\phi(f_2),\phi(f_1)) = -1$.
This embedding $\phi: T_X[-1] \hookrightarrow E_8$ is primitive, as 
$\phi(f_2)$ is chosen from the second half of (\ref{eq:E8-description}), 
while $\phi(f_1)$ from the first half. Note that 
$\phi(f_2)$ has coefficient $\pm 1/2$ for some $L_i$'s ($L_8$ in this
case), while $\phi(f_1)$ has vanishing coefficient.
Similar choices of $\phi(f_2)$ and $\phi(f_1)$ for other singular K3
surfaces easily produce primitive embeddings of $T_X[-1]$ into $E_8$.

The second step is to work out the orthogonal complement $T_0$, which 
is also almost straightforward. 
In the case of [a b c]=[5 1 1] we first see that the five roots 
\begin{equation}
(0,1,-1,0,0^3,0), \; 
(0,0,1,-1,0^3,0), \;
(0,(-1)^2,1,-1,1^2,-1)/2, \; 
(0,0^3,1,-1,0,0), \;
(0,0^3,0,1,-1)
\end{equation}
of $E_8$ are orthogonal to $\phi(T_X[-1])$ and generate an $A_5$
lattice. A primitive vector 
\begin{equation}
 \bar{e}_{114} := \frac{1}{2}(9,(-7)^3,(-1)^3,15)
\end{equation}
of $E_8$ is also in $[\phi(T_X[-1])^\perp \subset E_8]$, yet it is also 
orthogonal to $A_5$. The lattice $T_0 = [\phi(T_X[-1])^\perp \subset E_8]$ 
contains a rank-6 lattice $A_5 \oplus \Z\vev{\bar{e}_{114}} = 
A_5 \oplus (-114)$, but there are some elements left over. 
The lattice $T_0$ as the orthogonal complement of $\phi(T_X[-1])$ 
contains an element 
\begin{equation}
 E_8 \ni \frac{1}{2}(-1,1^3,-1,1^2,-3) = 
\frac{\alpha_1+2\alpha_2+3\alpha_3-2\alpha_4-\alpha_5}{6}
  - \frac{1}{6}\bar{e}_{114}, 
\end{equation}
which is not within the sublattice $A_5 \oplus (-114)$; 
here, $\alpha_{1,\cdots,5}$ are the five simple roots of $A_5$.
Modulo $A_5 \oplus \Z\vev{\bar{e}_{114}}$, this element can be regarded as 
an isotropic element\footnote{$q((-1,-19)) = q_{A_5}(-1)+q_{114}(-19) =
 (-5/6)(-1)^2 + (-1/114)(-19)^2 = -5/6-19/6 \equiv 0$ mod 2.} 
$(-1,-19) \in G_{A_5+(-114)} \cong 
 \Z/6\Z\vev{[\omega^1]} \times \Z_{114}\vev{[\bar{e}_{114}/114]}$, 
and generates a subgroup $\Z/6\Z \subset \Z/6\Z \times \Z/114 \Z$.
This result is simply described as $[A_5(-114)];\Z_6$ 
in Table \ref{tab:T0}.
It serves as a sanity check to compute\footnote{
In the case of [a b c]=[5 1 1], ${\rm discr}(T_0) = [{\rm discr}(A_5)]
\times [{\rm discr}(-114)]/ |\Z/6\Z|^2 = (-6)\times(-114)/6^2=19$,
which agrees with ${\rm discr}(T_X[-1]) = (-10)(-2)-(-1)^2 = 19$.} 
${\rm discr}(T_0)$ to see that it agrees with ${\rm discr}(T_X[-1])$.

It turns out from the computation of $T_0$, that the $T_0$ lattice for 
the singular K3 surface $X_{[6~1~1]}$ and that of $X_{[3~1~2]}$ are mutually 
isometric\footnote{There can be more than one primitive embeddings of 
this $T_0 = [A_4A_1(-230)]; \Z_{10}$ into $E_8$, since Theorem
$\epsilon$ for uniqueness of primitive embeddings (modulo ${\rm
Isom}(E_8)$) cannot be applied to this case.}.
This is enough to conclude that 
\begin{equation}
 {\cal J}_2(X_{[6~1~1]}) = {\cal J}_2(X_{[3~1~2]}).
\end{equation}

One will also note from Table~\ref{tab:T0} that there exist two inequivalent embeddings of $[3~2~2]$ into
the root lattice of $E_8$. The first one was found with the technique explained above, whereas the second one,
$T_0= A_4\oplus A_1^{\oplus 2}$, can be found in the literature \cite{Nish01}. In the second case, $T_0$
is a sum of $A$--$D$--$E$ root lattices. The existence such a choice of $T_0$
is irrelevant for the ${\cal J}_2$ classification (except for maybe some practical simplifications), 
where having one $T_0$ satisfying the conditions spelled out at the beginning of
Section \ref{ssec:J2-Nishiyama} is sufficient. However, this question becomes relevant
in the light of the results presented in Section \ref{ssec:Borcherds-Kondo}, which are
only applicable to cases in which $T_0$ is a sum of root lattices\footnote{We expect that this 
should be extendable to cases where $T_0$ is an overlattice of a sum of $A$--$D$--$E$ root lattices, but we 
do not discuss this here.}. This motivates to ask for all possibilities where $T_0$ can be 
chosen to be a sum of root lattices. We can answer this question by embedding by compiling a 
list of all rank-6 lattices of this type and studying their primitive embeddings into $E_8$.
The resulting orthogonal complements then tell us about the corresponding $T_X$. 

There are sixteen rank-6 lattices given by a direct sum of $A$--$D$--$E$ 
root lattices, out of which 11 already occur in Table~\ref{tab:T0}. The remaining five
are $D_4A_1^{\oplus 2}$, $A_3A_1^{\oplus 3}$, $A_2^{\oplus 3}$, $A_2A_1^{\oplus 4}$ and 
$A_1^{\oplus 6}$. It follows from Table 2.3 in \cite{Nish01} that none of these allow a
primitive embedding into $E_8$. Note that this also means that they can never have the
same discriminant form as one of the $T_X=[a~b~c]$, as this would imply the existence of 
such an embedding. Furthermore, using results of \cite{Nish01} one can show that all primitive
embeddings of the eleven (sums of) $A$--$D$--$E$ roots lattices appearing in Table~\ref{tab:T0} are unique,
so that we have found all singular K3 surfaces for which $T_0$ can be a sums of  $A$--$D$--$E$ 
root lattices.

%%%%%%%%%%%%%%%%%%%%%%%%%%%%%%%%%%%%%%%%%%%%%%%%%%%%%%%%%
\subsection{The ${\cal J}_2(X)$ Classification for Singular K3 Surfaces:
  New and Known Results}
\label{ssec:new&known}
%%%%%%%%%%%%%%%%%%%%%%%%%%%%%%%%%%%%%%%%%%%%%%%%%%%%%%%%%%%%%

The ${\cal J}_{2}(X)$ classification of elliptic fibrations has already 
been worked out for some K3 surfaces $X$ with large Picard number 
$\rho_X$. To name a few, 
\begin{itemize}
 \item For $X = {\rm Km}(E \times F)$, where $\rho_X = 18$, 
       there are 11 entries in ${\cal J}_2(X)$. 
       See Ref.~\cite{Oguiso} or Table \ref{tab:Oguiso-result} 
       in the next section (which carries partial information of the
       results obtained in \cite{Oguiso}).
 \item For four singular K3 surfaces $X = X_4$, $X_3$, 
       ${\rm Km}(E_i \times E_i)$ and ${\rm Km}(E_\omega \times
       E_\omega)$, and for ${\rm Km}(E \times E)$, where $\rho_X=19$, 
       Ref.~\cite{Nish01} worked out the ${\cal J}_2(X)$ classification 
       through the procedure reviewed in section
       \ref{ssec:J2-Nishiyama}. ${\cal J}_2(X)$ consists of 13, 6, 63
       and 30 entries for these four singular K3 surfaces above,
       respectively, and $|{\cal J}_2(X)| = 34$ for $X = {\rm Km}(E \times E)$.
 \item For $X = {\rm Km}(A)$, where $\rho_X = 17$, there are 
       25 entries in ${\cal J}_2(X)$. 
       See Ref.~\cite{Kumar}, or Table \ref{tab:Kummer-Kumar} in the
       next section of this article (which carries partial information of 
       the results in \cite{Kumar}).
 \item For $X=X_{[2~0~1]}$, one can choose $T_0 \cong D_5 \oplus A_1$,
       and $|{\cal J}_2(X)|=30$ as worked out in \cite{Bertin:2011}.
\end{itemize}
%

%%%%%%%%%%%%%%%%%%%%%%%%%%%%%%%%%%%%%%%%%%%%%%%%%%%%%%%%%%%%%%%%%%%%%%%%%%%%
\begin{table}[tb]
\begin{center}
 \begin{tabular}{|c||c|c|c|c|}
\hline 
  $T_X$ & $L^{(I)}$ & $W_{\rm root}$ & $MW(X)$ & S-Z \\
\hline
\hline
 [1 1 1] & % $E_6$ & $\{ 1 \}$ &
   $\gamma$ % $E_8^{\oplus 3}$
   & $A_2 \oplus E_8^{\oplus 2}$ & 0 & 297 \\
 & % & & 
   $\beta$ % $E_8 \oplus D_{16}$
   & $A_2 \oplus D_{16}$ & $\Z/2\Z$ & 216 \\
 & % & & 
   $\eta$ % $E_7^{\oplus 2} \oplus D_{10}$ 
   & $D_{10} \oplus E_7$ & $\Z \times \Z/2\Z$ & N/A \\
 & % & & 
   $\zeta$ % $E_7 \oplus A_{17}$ 
   & $A_{17}$ & $\Z \times \Z/3\Z$ & N/A \\
 & % & & 
   $\mu$ % $E_6^{\oplus 4}$ 
   & $E_6^{\oplus 3}$ & $\Z/3\Z$ & 219 \\
 & % & & 
   $\lambda$ % $E_6 \oplus D_7 \oplus A_{11}$ 
   & $A_{11} \oplus D_7$ & $\Z/4\Z$ & 166. \\
\hline 
\hline 
 \end{tabular}
\caption{\label{tab:Nishiyama-Table11}
The ${\cal J}_2(X)$ classification for $X_3$, quoted from \cite{Nish01}. Each entry is obtained by 
embedding $T_0 \cong E_6$ primitively into one of the Niemeier lattices 
$L^{(I)}$ specified in the second column; Greek letters correspond to 
those in Table \ref{tab:Niemeier-root}.
The last column indicates the type of `extremal elliptic K3' 
classified in \cite{SZ}. An {\it extremal elliptic K3 surface} is a K3
surface with elliptic fibration satisfying $\rho=20$ and ${\rm rank}(MW) = 0$.}
\end{center} 
\end{table}
%%%%%%%%%%%%%%%%%%%%%%%%%%%%%%%%%%%%%%%%%%%%%%%%%%%%%%%%%%%%%%%%%%%%%%%%%%%%

Each entry of ${\cal J}_2(X)$ for a given K3 surface $X$ is
characterized by an isometry class of frame lattice $W_{\rm frame}$. 
Once the frame lattice $W_{\rm frame}$ is given, one can extract information about
various objects.  
Its sublattice $W_{\rm root}$ generated by all the norm $(-2)$
elements of $W_{\rm frame}$ corresponds to the collection of singular fibres. 
The Mordell--Weil group of can be computed as the Abelian quotient
group (see \cite{Shioda}, Thm. 1.3)
\begin{equation}
 {\rm MW} := W_{\rm frame} / W_{\rm root}. 
\end{equation}
Table 1.1 of \cite{Nish01}---the result of $W_{\rm root}$ and the
Mordell--Weil group for the six distinct isometry classes of $W_{\rm frame}$ 
in ${\cal J}_2(X)$ for a singular K3 surface $X=X_3=X_{[1~1~1]}$---is reproduced here
as Table \ref{tab:Nishiyama-Table11} for the convenience of the reader.
With string theory applications in mind, however, more interesting objects are 
\begin{equation}
 W_{\rm gauge} := W_{\rm root} \oplus L(X), \qquad 
 L(X) := \left[ (W_{\rm root})^\perp \subset W_{\rm frame} \right] 
\end{equation}
and the subgroup $(W_{\rm frame}/W_{\rm gauge}) \subset G_{W_{\rm gauge}}$
\cite{BKW-phys}. For $X=X_3$, this information is summarized 
in section 4.1 of \cite{BKW-phys}. Here, $L(X)$ is called the 
{\it essential lattice} of an elliptic fibration on a K3 surface $X$. It is
also isometric to $MW(X)^0[-1]$, where $MW(X)^0$ is the {\it narrow Mordell--Weil lattice}
(see \cite{Shioda}, Thm. 8.9).

In all the cases referred to above, the lattice $T_0$ can be chosen 
to be either a direct sum of root lattices of A--D--E type, or an
overlattice of a direct sum of root lattices of A--D--E type.
Not all the $T_0$ lattices in Table \ref{tab:T0} are of that form.
However, even when $T_0$ does not contain a rank-6 root lattice, 
there is nothing preventing us from carrying out the procedure described 
in section \ref{ssec:J2-Nishiyama}.

%%%%%%%%%%%%%%%%%%%%%%%%%%%%%%%%%%%%%%%%%%%%%%%%%%%%%%%%%%%%%%%%%%%
\begin{table}[tbp]
\begin{center}
\begin{tabular}{||c|c|c|c|c|c|}
\hline 
 $L^{(I)}$ & $A_5$ & (-4) & $W_{\rm frame}=[T_0^\perp \subset L^{(I)}]$
 & $W_{\rm root}$ & MW \\
\hline
\hline 
$\alpha$ & $D_{24}$ & $\wedge^4 {\bf v} \in D_{18} \subset A_5^\perp$
   & $[A_3D_{14}(-6)];\Z_2$ & $A_3D_{14}$ & $\Z$ \\
\cline{3-6}
 & & ${\rm sym}^2 {\bf v} \in D_{18} \subset A_5^\perp$
   & $D_{17}(-6)$ & $D_{17}$ & $\Z$ \\
\hline
\hline 
$\beta$ & $E_8$ & $\wedge^4 {\bf v} \in D_{16}$ & $A_2A_1+[A_3D_{12}];\Z_2$ 
    & $A_1A_2A_3D_{12}$ & $\Z_2$ \\
\cline{3-6}
 & & ${\rm sym}^2 {\bf v} \in D_{16}$ & $A_2A_1D_{15}$ &
		 $A_1A_2D_{15}$ & $\{1 \}$ \\
\cline{3-6}
 & & ${\bf sp} \in D_{16}; \Z_2\vev{{\bf sp}}$ & $A_2A_1+A_{15};\Z_2$ &
		 $A_1A_2A_{15}$ & $\Z_2$ \\
\cline{3-6}
 & & $A_2+A_1=A_5^\perp$ & $(-6)(-4)+D_{16};\Z_2$ & $D_{16}$ &
		     $\Z\times \Z \times \Z_2$ \\
\cline{3-6}
 & & $[A_1\subset A_5^\perp] + D_{16}$ & $(-4)A_2+[D_{14}A_1];\Z_2$
   & $A_1A_2D_{14}$ & $\Z\times \Z_2$ \\
\cline{3-6}
 & & $[A_2 \subset A_5^\perp]+D_{16}$ &
   $A_1+[D_{14}A_1(-6)(-4)];(\Z_2)^2$ & $A_1^{\oplus 2}D_{14}$  
   & $\Z \times \Z \times \Z_2$ \\
\cline{2-6}
 & $D_{16}$ & $(1^4,0^4) \in E_8$ & $[A_3D_4];\Z_2+[D_{10}(-6)];\Z_2$ 
   & $A_3D_4D_{10}$ & $\Z\times \Z_2$ \\
\cline{3-6}
 & & $(2,0^7) \in E_8$ & $D_7+[D_{10}(-6)];\Z_2$ & $D_7D_{10}$ & $\Z$ \\
\cline{3-6}
 & & $(-3,1^7)/2 \in E_8$ & $[A_6(-28)]; \Z_7 + [D_{10}(-6)];\Z_2$  
   & $A_6D_{10}$ & $\Z \times \Z$ \\
\cline{3-6}
 & & $\wedge^4 {\bf v} \in D_{10} \subset A_5^\perp$
   & $E_8+[A_3D_6(-6)];\Z_2$ & $A_3D_6E_8$ & $\Z$ \\
\cline{3-6}
 & & ${\rm sym}^2{\bf v} \in D_{10} \subset A_5^\perp$ 
   & $D_9E_8(-6)$ & $D_9E_8$ & $\Z$ \\
\cline{3-6}
 & & ${\bf sp} \in D_{16}; \Z_2\vev{{\bf sp}}$ & $E_8+[A_9(-60)]; \Z_5$ 
   & $A_9E_8$ & $\Z$ \\
\cline{3-6}
 & & $E_8+A_5^\perp$ & $[E_7(-4)(-6)D_8A_1];(\Z_2)^2$
   & $A_1D_8E_7$ & $\Z \times \Z$ \\
\hline
\hline
$\gamma$ 
 & $E_8^{(1)}$ & $(1^4,0^4) \in E_8^{(2)}$
   & $E_8^{(3)}A_2A_1+[A_3D_4];\Z_2$ & $A_1A_2A_3D_4E_8$ & $\Z_2$ \\
\cline{3-6}
 & & $(2,0^7) \in E_8^{(2)}$ & $E_8^{(3)}A_2A_1D_7$
   & $A_1A_2D_7E_8$ & $\{ 1 \}$ \\
\cline{3-6}
 & & $(-3,1^7)/2 \in E_8^{(2)}$ & $E_8^{(3)}A_2A_1+[A_6(-28)];\Z_7$
   & $A_1A_2A_6E_8$ & $\Z$ \\
\cline{3-6} 
 & & $E_8^{(2)}+E_8^{(3)}$ & $A_2A_1+[E_7^{(2)}E_7^{(3)}(-4)];\Z_2$ 
   & $A_1A_2E_7E_7$ & $\Z$ \\
\cline{3-6}
 & & $A_5^\perp$ & $E_8^{(2)}E_8^{(3)}(-6)(-4)$
   & $E_8E_8$ & $\Z \times \Z$ \\
\cline{3-6} 
 & & $E_8^{(2)}+[A_1 \subset A_5^\perp]$ & $E_8^{(3)}E_7A_2(-4)$ 
   & $A_2E_7E_8$ & $\Z$ \\
\cline{3-6}
 & & $E_8^{(2)}+[A_2 \subset A_5^\perp]$
   & $A_1E_8^{(3)}E_7+[(-6)(-4)];\Z_2$ & $A_1E_7E_8$ & $\Z \times \Z$ \\   
\hline
\hline 
\end{tabular}
\caption{\label{tab:J2classify-[abc=302]}
${\cal J}_{2}(X)$ classification for a K3 surface $X=X_{[3~0~2]}$. 
% whose transcendental lattice symmetric pairing is given by [a b c]=[3 0 2].
This table only shows the entries of ${\cal J}_2(X)$ obtained from 
primitive embeddings of $T_0 = A_5 \oplus (-4)$ into Niemeier lattices 
$L^{(I)}$ with $I=\alpha, \beta, \gamma$.
The first column indicates which one of twenty-three Niemeier lattices 
% $L^{(I)}$ $(I=\alpha,\beta,\cdots)$ the rank-6 lattice 
% $T_0 \cong A_5+(-4)$
$T_0$ is embedded into and the second and third columns specify the
particular irreducible root lattices $R_a$ in $(L^{(I)})_{\rm root} \cong \oplus_a
R_a$ into which $A_5$ and  $(-4)$ are embedded. 
% $e_4$ of the rank-1 lattice $(-4)$ lies. 
}
\end{center}
\end{table}
%%%%%%%%%%%%%%%%%%%%%%%%%%%%%%%%%%%%%%%%%%%%%%%%%%%%%%%%%%%%%%%%%%%%
As another example, let us consider the singular K3 surface $X_{[3~0~2]}$.
As in Table \ref{tab:T0}, we can take $T_0 = A_5 \oplus (-4)$ for 
this singular K3 surface. The $A_5$ component has to be embedded
primitively into the sublattice $(L^{(I)})_{\rm root} \subset L^{(I)}$
of the Niemeier lattices, and the results of \S 4.1 of \cite{Nish01} 
can be used for this part of the calculation. One still has to work out 
all possible embeddings of the generator of the $(-4)$ part of $T_0$, 
and make sure that the embedding into $L^{(I)}$ is primitive.

Instead of carrying out the ${\cal J}_{2}(X)$ classification for 
$X_{[3~0~2]}$ completely, we first carried out the part obtained 
from embedding into the Niemeier lattices $L^{(I)}$ with 
$I=\alpha, \beta$ and $\gamma$; the results are found in 
Table \ref{tab:J2classify-[abc=302]}. 
As a second part, we have studied all embeddings
into the remaining Niemeier lattices for which the frame lattice contains 
a root lattice of type $E_6$, $E_7$ or $E_8$, the results are found
in Table \ref{tab:J2classify-[abc=302]_Etype}. Together, the two tables 
hence contain all the isometry classes of the frame lattice available in 
$X_{[3~0~2]}$ % cases of $T_X=$ [3 0 2] 
for which singular fibres of type ${\rm IV}^*$, ${\rm III}^*$ or 
${\rm II}^*$ occur. Overall, there are 43 distinct isometry classes in the 
${\cal J}_2(X_{[3~0~2]})$ classification containing either $E_6$, $E_7$ or $E_8$.

The resulting frame lattice $W_{\rm frame}$ is expressed 
in Table \ref{tab:J2classify-[abc=302]} in the form 
of $W_{\rm gauge}; (W_{\rm frame}/W_{\rm gauge})$, in order to make 
it easier to read out the information necessary in string theory 
applications.

Note that this list already contains some interesting subtleties. 
Let us focus on cases for which $W_{\rm root}= A_1D_8E_7$;
the 15th entry of Table \ref{tab:J2classify-[abc=302]} and 
the 12th and 18th entries of Table \ref{tab:J2classify-[abc=302]_Etype} 
share this property. Although the Mordell--Weil group is of rank 2 for all 
the three cases, the first one is without a torsion part, while the last two 
have $\Z_2$ torsion. 
Furthermore, the last two cases share the same $W_{\rm root}$ and the 
Mordell--Weil group, yet their frame lattices are not mutually isometric.
The same relation holds between the 4th and 5th entries of Table 
\ref{tab:J2classify-[abc=302]_Etype}. 

%%%%%%%%%%%%%%%%%%%%%%%%%%%%%%%%%%%%%%%%%%%%%%%%%%%%%%%%%%%%%%%%%%%
\begin{table}[tbp]
\begin{center}
\begin{tabular}{||c|c|c|c|c|c|}
\hline 
 $L^{(I)}$ & $A_5$ & (-4) & $W_{\rm frame}=[T_0^\perp \subset L^{(I)}]$
 & $W_{\rm root}$ & MW \\
\hline
\hline 
$\zeta$ & $A_{17}$ & $(0^{14},1^2,-1^2) \in A_{17}$ & $ E_7 A_7A_{1}A_1 (-36)(-24);\Z_{12}\times \Z_{4}$ & $A_{1}A_1A_7E_7  $ 
& $\Z\times \Z$ \\
\cline{3-6}
 & & $g_{\zeta}$ & $E_6A_8A_2(-90)(-240); \Z_{90}\times \Z_{3}$ & $A_2A_8E_6$ & $\Z \times \Z$  \\
 \cline{3-6}
 & & $2g_{\zeta}$ mod $A_{17}\oplus E_7$ & $E_7A_5A_5(-12); \Z_{6}$ & $A_5^{\oplus 2}E_7$ & $\Z$  \\
\hline
\hline 
$\eta$ & $D_{10}$ & $(0^6,1^4) \in D_{10}$ & $E_7+[A_3(-6)E_7];\Z_2$  & $A_3E_7^2$ & $\Z$ \\
\cline{3-6}
 & & $(0^6,2,0^3) \in D_{10}$ & $A_3(-6)E_7^2;\Z_2 $ & $A_3E_7^2$ & $\Z$ \\
\cline{3-6}
 & & $g_{\eta 1}$ & $A_3 E_6 E_7 (-12)(-48);\Z_4 \times \Z_6$ &
		 $A_3 E_6 E_7 $ & $\Z\times\Z$ \\
\cline{3-6}
 & & $D_{10}+E_7$ & $(-4)(-6)A_1^3D_6E_7;\Z_2^3$
   & $A_1^3D_6E_7$ & $\Z\times \Z$ \\
\cline{3-6}
 & & $(1^2,-1^2,0^4)\in E_7$ &
   $D_4D_5A_1E_7(-6);(\Z_2)^2$ & $A_1D_4D_5E_7$  
   & $\Z$ \\
\cline{3-6}
 & & $g_{\eta 1}-g_{\eta 2}$ &
   $A_3(-6)(-12)(-48)E_6E_6;\Z_2^3\times\Z_{3}^2$ & $A_3E_6E_6$  
   & $\Z\times\Z\times\Z$ \\
\cline{2-6}
 & $E_7$,$1$ & $(1^2,-1^2,0^6) \in D_{10}$ & $D_6A_3A_1(-6)E_7;\Z_2\times\Z_2$ 
   & $A_1A_3D_6E_7$ & $\Z\times \Z_2$ \\
\cline{3-6}
 & & $(2,0^9) \in D_{10}$ & $D_9+[A_{1}(-6)E_7];\Z_2$ & $A_{1}D_9E_7$ & $\Z$ \\
\cline{3-6}
 & & $D_{10} + E_7$ & $[A_1D_8E_7;\Z_2+(-6)];\Z_2+(-4)$  
   & $A_1D_8E_7$ & $\Z \times \Z \times \Z_2$ \\
\cline{3-6}
 & &  $g_{\eta 1}$
   & $A_9A_1E_7(-240);\Z_2\times \Z_{10}$ & $A_1A_9E_7$ & $\Z\times\Z_2$ \\
\cline{3-6}
 & & $g_{\eta 2}$
   & $A_9A_1(-6)E_6(-240);\Z_6\times\Z_{10}$ & $A_1A_9E_6$ & $\Z\times\Z$ \\
\cline{3-6}
 & & $g_{\eta 1} - g_{\eta 2}$
   & $D_9A_1E_6(-12)(-48);\Z_{24}$ & $A_1D_9E_6$ & $\Z\times\Z$ \\
   \cline{2-6}
 & $E_7$,$2$ & $(1^2,-1^2,0^6) \in D_{10}$ & $A_3D_6E_7;\Z_2+A_2$ 
   & $A_2A_3D_6E_7$ & $\Z_2$ \\
\cline{3-6}
 & &$(2,0^9) \in D_{10}$ & $D_9A_2E_7$ 
   & $A_2D_9E_7$ & $\{1\}$ \\
\cline{3-6}
 & &$D_{10}+ ( A_2 = A_5^\perp \subset E_7)$ & $[A_1D_8E_7;\Z_2+(-4)(-6)];\Z_2$ 
   & $A_1D_8E_7$ & $\Z\times\Z\times\Z_2$ \\
   \cline{3-6}
 & & $g_{\eta 2}$ & $A_2+[A_9E_6(-240)];\Z_{20}$ 
   & $A_2A_9E_6$ & $\Z$ \\
\hline
\hline
$\lambda$ 
 & $A_{11}$ & $(1^2,-1^2,0^8) \in A_{11}$
   & $A_1^3(-12)(-12)D_7E_6;\Z_2\times \Z_{12}$ & $A_1^3D_7E_6$ & $\Z\times\Z$ \\
\cline{3-6}
 & & $(1,1,-1,-1,0^3) \in D_7$ & $A_5(-12)A_3^{\oplus 2}E_6;\Z_{12}$
   & $A_3^{\oplus 2}A_5E_6$ & $\Z$ \\
\cline{3-6}
 & & $(2,0^6) \in D_7$ & $A_5(-12)E_6;\Z_6+D_6$
   & $A_5D_6E_6$ & $\Z$ \\
\cline{3-6} 
 & & $A_{11}+D_7$ & $A_{3}D_5A_1E_6(-12)(-12)(-4);\Z_{12}\times \Z_4$ 
   & $A_1A_{3}D_5E_6$ & $\Z\times\Z\times\Z$ \\
\cline{3-6} 
 & & $3g_{\lambda}$ mod $E_6\oplus D_7$ & $[A_2A_2A_6(-18)(-4032)+E_6;\Z_3];\Z_9\Z_4\Z_7$
   & $A_2^{\oplus 2}A_6E_6$ & $\Z\times\Z\times\Z_3$ \\   
\cline{3-6} 
 & & $6g_{\lambda}$ mod $A_{11}\oplus E_6\oplus D_7$ & $[A_5D_6E_6(-48)];\Z_{12}$
   & $A_5D_6E_6$ & $\Z$ \\   
\cline{2-6}
 & $D_7$ & $A_{11}$ & $[A_1^2A_7;\Z_2+(-24)(-6)(-4)E_6];\Z_{24}$ 
   & $A_1^2A_7E_6$ & $\Z^3\times\Z_2$ \\
\cline{3-6}
 & & $D_7$ & $(-6)+A_{11}E_6;\Z_3$
   & $A_{11}E_6$ & $\Z \times \Z_3$ \\
\cline{3-6} 
 & & $3g_{\lambda}$ mod $(E_6\oplus D_7)$ & $[A_8A_2\left[\begin{array}{cc}
                                           -90 & 36 \\
                                           36 & -360
                                          \end{array}\right]+E_6;\Z_3];\Z_{108}$
   & $A_2^{\oplus 2}A_6E_6$ & $\Z\times\Z\times\Z_3$ \\  
\cline{3-6} 
 & & $6g_{\lambda}$ mod $A_{11}\oplus E_6\oplus D_7$ & $A_5A_5E_6(-6)(-48);\Z_{12}\times\Z_3$
   & $A_5^{\oplus 2}E_6$ & $\Z\times\Z$ \\   
\hline
\hline 
$\mu$ 
 & $E_6$ & $(0,1^2,-1^2,0^3) \in E_{6}^{(2)}$
   & $A_1+(-12)D_4E_6^2;\Z_6$ & $A_1 D_4E_6^2$ & $\Z$ \\
\cline{3-6}
 & & $E_6^{(2)}+E_6^{(3)}$ & $[(-4)A_1+E_6A_5^2;\Z_3]\Z_2$
   & $A_1A_5^2E_6$ & $\Z\times\Z_3$ \\
\cline{3-6} 
 & & $A_{1}+E_6^{(2)}$ & $(-4)+A_5E_6^2;\Z_{3}$ 
   & $A_5E_6^2$ & $\Z\times\Z_3$ \\
\hline 
\end{tabular}
\caption{\label{tab:J2classify-[abc=302]_Etype}
${\cal J}_{2}(X)$ classification for a K3 surface $X=X_{[3~0~2]}$. 
This table only shows the entries of ${\cal J}_2(X)$ obtained from 
primitive embeddings of $T_0 = A_5 \oplus (-4)$ into Niemeier lattices (apart from those already
considered in Table \ref{tab:J2classify-[abc=302]}) such that $T_0^{\perp}$ 
contains one of the lattices $E_6$, $E_7$ or $E_8$. 
The columns are as in Table \ref{tab:J2classify-[abc=302]}.
The $g_I$'s appearing in the 3rd column are glue vectors of the indicated Niemeier lattices, i.e. 
generators of $\Delta \subset G_{L^{(I)}_{\rm root}}$ shown in Table \ref{tab:Niemeier-root}.}
\end{center}
\end{table}
%%%%%%%%%%%%%%%%%%%%%%%%%%%%%%%%%%%%%%%%%%%%%%%%%%%%%%%%%%%%%%%%%%%%

%%%%%%%%%%%%%%%%%%%%%%%%%%%%%%%%%%%%%%%%%%%%%%%%%%%%%%%%%%%%%%%%%%%
\section{The ${\cal J}_1(X)$ Classification}
\label{sec:J1-classification}
%%%%%%%%%%%%%%%%%%%%%%%%%%%%%%%%%%%%%%%%%%%%%%%%%%%%%%%%%%%%%%%%%%%

%%%%%%%%%%%%%%%%%%%%%%%%%%%%%%%%%%%%%%%%%%%%%%%%%%%%%%%%%%%%%
\subsection{Uniform Upper Bounds on Multiplicity}
%%%%%%%%%%%%%%%%%%%%%%%%%%%%%%%%%%%%%%%%%%%%%%%%%%%%%%%%%%%%%

The classification ${\cal J}_1(X)$ of elliptic fibrations $(\pi_X,
\sigma; X, \P^1)$ contains finer information than the
${\cal J}_2(X)$ classification. 
For a given isometry class of frame lattice $[W] \in {\cal J}_2(X)$, 
there may be more than one isomorphism classes of elliptic fibrations 
$(\pi_X, \sigma; X, \P^1)$ in ${\cal J}_1(X)$. In this section, 
we make an attempt at deriving upper bounds 
on this number, {\it the number of isomorphism classes for a given} 
$[W] \in {\cal J}_2(X)$, which we call its {\it multiplicity}.
The ${\cal J}_1(X)$ classification is not only an interesting mathematical 
question, but also relevant to string theory applications \cite{BKW-phys}.

The ${\cal J}_1(X)$ classification has been worked out completely 
for $X={\rm Km}(E \times F)$ in \cite{Oguiso}, the result is summarized 
in Table \ref{tab:Oguiso-result}.
%%%%%%%%%%%%%%%%%%%%%%%%%%%%%%%%%%%%%%%%%%%%%%%%%%%%%%%%%%%%%%%%%%%%%%%%%%%%
\begin{table}[tbp]
  \begin{tabular}{|c||c|c|c|c|c|c|c|c|c|c|c|}
\hline
  ${\cal J}^{({\rm type})}(X)$ & $2I_8+8I_1$ & $I_4+I_{12}+8I_1$ & 
      $2{\rm IV}^*+(8-2b)I_1+b{\rm II}$ & $4I_0^*$ & $I_6^* + 6I_2$ \\  
\hline 
  nmb.isom.clss. & 9 & 6 & 1 & 2 & 1 \\
\hline 
\hline
\end{tabular}
\begin{tabular}{|c||c|c|c|}
\hline
 ${\cal J}^{({\rm type})}(X)$ &
     $2I_2^*+4I_2$ & $I_4^*+2I_0^*+2I_1$ &
      ${\rm III}^* + I_2^* +(3-b)I_2 + (1-b)I_1+b{\rm III}$ \\  
\hline
  nmb.isom.clss. & 9 & 9 & 6 \\
\hline
\hline
\end{tabular}
\begin{tabular}{|c||c|c|c|}
\hline
  ${\cal J}^{({\rm type})}(X)$ & 
      ${\rm II}^*+2I_0^*+(2-2b)I_1+b{\rm II}$ & 
      $I_8^*+I_0^*+(4-2b)I_1+b{\rm II}$ & 
      $2I_4^*+(4-2b)I_1 + b {\rm II}$ \\
\hline
  nmb.isom.clss. & 1 & 6 & 9 \\
\hline
\hline 
  \end{tabular}
\caption{\label{tab:Oguiso-result} Part of the results obtained 
in \cite{Oguiso}. There are 11 elements in the 
${\cal J}^{({\rm type})}(X) \approx {\cal J}_2(X)$ classification 
of elliptic fibrations for $X = {\rm Km}(E \times F)$ with a pair of 
mutually non-isogenous elliptic curves $E$ and $F$. The lower rows present 
the number of isomorphism classes in the individual isometry classes ${\cal J}_2(X)$, 
i.e. the multiplicities. }
\end{table}
%%%%%%%%%%%%%%%%%%%%%%%%%%%%%%%%%%%%%%%%%%%%%%%%%%%%%%%%%%%%%%%%%%%%%%%%%
The study of \cite{Oguiso} relies, to some extent, on things that are 
specific to the particular family of K3 surfaces  
$X = {\rm Km}(E \times F)$, but there are also ideas and structures in 
\cite{Oguiso} that can also be applied to other K3 surfaces with large Picard
number. In these cases we can apply them for deriving 
upper bounds on the multiplicities rather than for precisely determining them. 

The first thing to notice is 
\begin{quote}
{\bf Proposition $C$}: for a K3 surface $X$, the number of isomorphism
classes of elliptic fibration (multiplicity) is bounded uniformly 
from above by the number of elements of the coset space
\begin{equation}
  {\rm Isom}(S_X)^{({\rm Amp~Hodge})} \backslash {\rm Isom}(S_X)^{({\rm Amp})} 
 \cong 
    \left[ W^{(2)}(S_X) \rtimes {\rm Isom}(S_X)^{({\rm
     Amp~Hodge})}\right] \backslash  {\rm Isom}^+(S_X) 
\label{eq:the-coset}
\end{equation}
for {\it any} isometry class $[W] \in {\cal J}_2(X)$,   
\end{quote}
because the difference between the two classifications 
${\cal J}_2(X)$ and ${\cal J}_1(X)$ is only in the choice of the
quotient group. $\blacksquare$

In practice, though, it is not easy to compute the coset 
space (\ref{eq:the-coset}) for many different K3 surfaces 
without a great deal of knowledge about the geometry and 
the Neron--Severi lattice. Easier to use for K3 surfaces with 
large Picard number is 
\begin{quote}
{\bf Proposition $C'$}: for a K3 surface $X$, the number of isomorphism
 classes of elliptic fibration (multiplicity) is bounded uniformly 
from above by the number of elements of the coset space 
\begin{equation}
 p_T\left( {\rm Isom}(T_X)^{({\rm Hodge})}\right) \backslash 
   {\rm Isom}(q) 
\label{eq:coset-alt1}
\end{equation}
for {\it any} isometry class $[W] \in {\cal J}_2(X)$.
This upper bound is generally weaker than that of Proposition C, 
but the two are the same if the homomorphism 
$p_S: {\rm Isom}^+(S_X) \longrightarrow {\rm Isom}(q)$ is surjective. 
\end{quote}
{\bf proof}: Let 
$G^{\rm tot} = {\rm Isom}(q)$, 
$G^{(s)} = p_S[{\rm Isom}(S_X)^{({\rm Amp})}]$, 
$G^{(t)} = p_T[{\rm Isom}(T_X)^{({\rm Hodge})}]$, 
$H = G^{(t)} \cap G^{(s)}$, and 
$G^{\rm relev}$ be the subgroup of $G^{\rm tot}$ generated by all the
elements in $G^{(s)}$ and $G^{(t)}$.
First, note that the homomorphism $p_S$ maps the coset space 
(\ref{eq:the-coset}) one-to-one to another coset $H \backslash G^{(s)}$ 
defined in $G^{\rm tot}$. Thus, the upper bound in Proposition C is given by 
$|H \backslash G^{(s)}|$. We claim, then, that 
\begin{equation}
 | H \backslash G^{(s)}| \leq | G^{(t)} \backslash G^{\rm relev} |
  \leq |G^{(t)} \backslash G^{\rm tot} |;  
\end{equation}
since the last inequality is obvious, only the first inequality needs
to be verified. Now, let $\{s_i\}_{i \in I}$ be a set of
representatives of the coset space $H \backslash G^{(s)}$. 
If two representatives $s_i$ and $s_j$ ($i \neq j$) were to be in the 
same orbit of $G^{(t)}$ in $G^{\rm relev}$, i.e., 
${}^\exists t \in G^{(t)}$ such that $t \cdot s_i = s_j$, then 
$t = s_j \cdot s_i^{-1} \in G^{(s)}$, and hence $t \in H$; this is a 
contradiction. $\blacksquare$

In the case of $X = {\rm Km}(E \times F)$, 
\begin{equation}
  {\rm Isom}(q) \cong (S_3 \times S_3) \rtimes (\Z/2\Z), \qquad 
 p_T\left( {\rm Isom}(T_X)^{({\rm Hodge})}\right) = \{ 1 \}, 
\end{equation}
and $p_S: {\rm Isom}^+(S_X) \longrightarrow {\rm Isom}(q)$ is 
surjective \cite{Oguiso}. Since the coset space (\ref{eq:coset-alt1}) is 
a group of order 72, Proposition C' implies that the number of
isomorphism classes must be no more than 72 for any one of the eleven 
isometry classes of frame lattice ${\cal J}_2(X)$ for 
$X = {\rm Km}(E \times F)$. 
From Table \ref{tab:Oguiso-result} this is indeed seen to be the case. 

We computed the group ${\rm Isom}(q)$ and the coset space 
$p_T( {\rm Isom}(T_X)^{({\rm Hodge})} ) \backslash {\rm Isom}(q)$ for 
the thirty-four singular K3 surfaces that showed up in the study of 
\cite{BKW-phys}. The result is summarized in Table \ref{tab:coset-alt-rho20}.
%%%%%%%%%%%%%%%%%%%%%%%%%%%%%%%%%%%%%%%%%%%%%%%%%%%%
\begin{table}[tbp]
\begin{center}
  \begin{tabular}{|l|c|c|c|r|r|} \hline
    $T_X$ & $\Z/m\Z$ & $G_{T_X}$  
    & ${\rm Isom}(G_{T_X}, q_{T_X})$ & $p_T(\theta_m)$ & Coset \\ \hline \hline
 {}[1 0 1] & $\Z_4$ & $\Z_2 \times \Z_2$ & $\Z_2\vev{\rm exch}$ && $\{ 1 \}$ \\
 {}[1 1 1] & $\Z_6$ &  $\Z_3$ & $\Z_2\vev{-_3}$ && $\{ 1 \}$ \\
 {}[2 0 1] & $\Z_2$ & $\Z_4 \times \Z_2$ & $\Z_2\vev{-_4}$
		   && $\{ 1 \}$ \\
 {}[2 1 1] & $\Z_2$ & $\Z_7$ & $\Z_2\vev{-_7}$ && $\{ 1 \}$ \\
 {}[3 0 1] & $\Z_2$ & $\Z_6 \times \Z_2$ & $\Z_2\vev{-_3}$ &&
		       $\{ 1 \}$\\
 {}[3 1 1] & $\Z_2$ & $\Z_{11}$ & $\Z_2\vev{-_{11}}$ && $\{ 1 \}$ \\
 {}[4 0 1] & $\Z_2$ & $\Z_8 \times \Z_2$ & $\Z_2\vev{-_8}$ && $\{ 1 \}$ \\
 {}[4 1 1] & $\Z_2$ & $\Z_{15}$ & $\Z_2\vev{-_3} \times \Z_2\vev{-_5}$ 
                    & $(-_3,-_5)$ & $\Z_2$ \\
 {}[5 0 1] & $\Z_2$ & $\Z_{10} \times \Z_2$ & $\Z_2\vev{-_5} \times
	       \Z_2\vev{{\rm exch}}$ & $(-_5,1)$  & $\Z_2$ \\
 {}[5 1 1] & $\Z_2$ & $\Z_{19}$ & $\Z_2\vev{-_{19}}$ && $\{ 1 \}$ \\
 {}[6 0 1] & $\Z_2$ & $\Z_{12} \times \Z_2$ &
     $\Z_2\vev{-_3}\times \Z_2\vev{-_4}$ & $(-_3,-_4)$  & $\Z_2$ \\
 {}[6 1 1] & $\Z_2$ & $\Z_{23}$ & $\Z_2\vev{-_{23}}$ && $\{ 1\}$ \\
\hline
 {}[2 0 2] & $\Z_4$ & $\Z_4 \times \Z_4$ & $\Z_4\vev{\tau} \rtimes 
          \Z_2\vev{\sigma} \cong D_4$ & $\tau$ & $\Z_2\vev{\sigma}$ \\
 {}[2 1 2] & $\Z_2$ & $\Z_{15}$ & $\Z_2 \vev{-_3} \times \Z_2\vev{-_5}$ 
                    & $(-_3,-_5)$ & $\Z_2$ \\
 {}[2 2 2] & $\Z_6$ & $(\Z_6 \times \Z_2)$ &
     $\Z_2\vev{-_3} \times S_3\vev{\sigma_2,\tau_2}$ 
    & $(-_3,\tau_2)$ & $\Z_2 \vev{\sigma_2}$ \\
 {}[3 0 2] & $\Z_2$ & $\Z_6 \times \Z_4$ &
     $\Z_2 \vev{-_3} \times \Z_2 \vev{-_4}$ & $(-_3,-_4)$ & $\Z_2$ \\
 {}[3 1 2] & $\Z_2$ & $\Z_{23}$ & $\Z_2\vev{-_{23}}$ & & $\{ 1 \}$ \\
 {}[3 2 2] & $\Z_2$ & $(\Z_{10} \times \Z_2)$ & $\Z_2\vev{-_5}
	       \times \Z_2\vev{{\rm mix}_{22}}$ & $(-_5,1)$ &
	       $\Z_2\vev{{\rm mix}_{22}}$ \\
 {}[4 0 2] & $\Z_2$ & $\Z_8 \times \Z_4$ &
     $\Z_2\vev{-_8} \times \Z_2\vev{-_4} 
        \times \Z_2\vev{{\rm mix}_{84}}$ & $(-_8,-_4,1)$
     & $\Z_2 \times \Z_2\vev{{\rm mix}_{84}} $ \\
 {}[4 2 2] & $\Z_2$ & $(\Z_{14} \times \Z_2)$ & $\Z_2\vev{-_7} 
       \times \Z_2\vev{{\rm mix}_{22}}$ & $(-_7, 1)$ &
       $\Z_2\vev{{\rm mix}_{22}}$ \\
 {}[6 0 2] & $\Z_2$ & $\Z_{12} \oplus \Z_4$ &
     $\Z_2\vev{-_3} \times \Z_2\vev{-_{4;12}} \times
	       \Z_2\vev{-_{4;4}}$ & $(-,-,-)$ 
    & $\Z_2 \times \Z_2$ \\
 {}[6 2 2] & $\Z_2$ & $(\Z_{22} \times \Z_2)$ &
      $\Z_2\vev{-_{11}} \times S_3\vev{\tau_2,\sigma_2}$ & 
      $(-_{11},1)$ & $S_3\vev{\tau_2,\sigma_2}$ \\
\hline 
 {}[3 0 3] & $\Z_4$ & $\Z_6 \times \Z_6$
	       & $\Z_2\vev{{\rm exch}_2} \times D_4\vev{\tau_3,\sigma_3}$ &
   $({\rm exch}_2, \tau_3)$ & $\Z_2 \times \Z_2\vev{\sigma_3}$ \\
 {}[3 3 3] & $\Z_{6}$ & $\left(\Z_{9} \times \Z_3\right)$ &
          $\Z_2\vev{-} \times S_3\vev{\tau,\sigma}$ & $(-,\tau)$
          & $\Z_2\vev{\sigma}$  \\
 {}[6 0 3] & $\Z_2$ & $\Z_{12} \times \Z_6$ &
    $\Z_2\vev{-_4} \times \Z_2\vev{-_{3;12}} \times \Z_2\vev{-_{3;6}}$ &
    $(-,-,-)$ & $\Z_2 \times \Z_2$ \\
 {}[6 3 3] & $\Z_2$ & $(\Z_{21} \times \Z_3)$ &
      $\Z_2\vev{-_7} \times D_4\vev{\tau_3,\sigma_3}$
	       & $(-_7,\tau^2_3)$ & 
       $D_4\vev{\tau,\sigma_3}$  \\
\hline
 {}[4 0 4] & $\Z_4$ & $\Z_8 \times \Z_8$ &
     $D_4\vev{\tau_8,\sigma_8} \times \Z_2\vev{{\rm exch}'_8}$ & 
     $(\tau_8,{\rm exch}'_8)$ & $\Z_2 \times \Z_2$ \\
 {}[4 4 4] & $\Z_6 $ & $(\Z_{12} \times \Z_4)$ &
     $\Z_2\vev{-_3}\times \Z_2\vev{-_4}\times S_3\vev{\tau_4,\sigma_4}$
     & $(-,-,\tau)$ & $\Z_2 \times \Z_2$ \\
 {}[5 0 5] & $\Z_4$ & $\Z_{10} \times \Z_{10}$ & $\Z_2\vev{{\rm exch}_2}
	       \times D_4^{(5)}$ &
    $({\rm exch}_2,\tau_5)$ & $\Z_2 \times \Z_2$ \\
 {}[5 5 5] & $\Z_6$ & $(\Z_{15} \times \Z_5)$ & 
     $\Z_2\vev{-_3}\times \Z_2\vev{-_5}\times S_3\vev{\tau_5,\sigma_5}$
     & $(-,-,\tau)$ & $\Z_2 \times \Z_2$ \\
 {}[6 0 6] & $\Z_4$ & $\Z_{12} \times \Z_{12}$ & $D_4^{(3)} \times
	       D_4^{(4)}$ & $(\tau_3, \tau_4)$
		   & $(D_4\times D_4)/\Z_4$ \\
 {}[6 6 6] & $\Z_6$ & $(\Z_{18} \times \Z_6)$ &
     $S_3\vev{\tau_2,\sigma_2} \times \Z_2 \vev{-_3} 
      \times S_3\vev{\tau_3,\sigma_3}$ & $(\tau,-,\tau)$ &
      $(S_3 \times S_3)/\Z_3$ \\
 {}[7 7 7] & $\Z_6$ & $(\Z_{21} \times \Z_7)$ &
     $\Z_2\vev{-_3}\times \Z_2\vev{-_7}\times S_3\vev{\tau_7,\sigma_7}$
     & $(-,-,\tau)$ & $\Z_2 \times \Z_2$ \\
 {}[8 8 8] & $\Z_6$ & $(\Z_{24} \times \Z_8)$ & 
     $\Z_2\vev{-_3}\times \Z_2\vev{(-)_8}\times \Z_2\vev{(3)_8} \times 
      S_3\vev{\tau_8,\sigma_8}$ & $(-,-,1,\tau)$ &
      $\Z_2 \times \Z_2 \times \Z_2$ \\
\hline    
  \end{tabular}
\caption{\label{tab:coset-alt-rho20}
We study here thirty-four singular K3 surfaces $X_{[a~b~c]}$ specified by 
$[a~b~c]$ in the first column. 
Their ${\rm Isom}(T_X)^{({\rm Hodge})} \cong \Z/m\Z$ and 
the discriminant group $G_{T_X} \cong T_X^*/T_X$ are shown in the 
2nd and 3rd columns, respectively. In the 3rd column, expressions such as  
$G_1 \times G_2$ and $(G_1 \times G_2)$ mean that the discriminant
bilinear form is block diagonal in the former case, and it is not in the latter. 
$\theta_m$ is the generator of 
${\rm Isom}(T_X)^{({\rm Hodge})} \cong \Z/m\Z$ 
(angle $2\pi/m$ rotation in $T_X \otimes \R$). The coset space 
(\ref{eq:coset-alt1}) in the last column is determined by the
information in the 4th and 5th columns. The 5th column is left empty, when $p_T$ is surjective.
}
\end{center}
\end{table}
%%%%%%%%%%%%%%%%%%%%%%%%%%%%%%%%%%%%%%%%%%%%%%%%%%%%
In the table, 
\begin{equation}
 D_4\vev{\tau,\sigma} := 
 \Z_4\vev{r_1\cdot s=\tau} \rtimes \Z_2\vev{r_2=\sigma} \cong 
 (\Z_2\vev{r_1} \times \Z_2\vev{r_2})\rtimes \Z_2\vev{s} 
\end{equation}
is the dihedral group of order 8. We have, 
\begin{quote}
{\bf Corollary D}: 
For any one of the thirty-four singular K3 surfaces $X_{[a~b~c]}$ studied in 
Table \ref{tab:coset-alt-rho20}, 
the number of isomorphism classes of elliptic fibration (multiplicity)
is bounded from above uniformly for any $[W] \in {\cal J}_2(X)$ 
by the number of elements of the coset space in the last column of the table. 
In particular, there are ten singular K3 surfaces in the table, 
\begin{eqnarray}
&&
 X_{[1~0~1]}, \quad X_{[1~1~1]}, \quad X_{[2~0~1]}, \quad 
 X_{[2~1~1]}, \quad X_{[3~0~1]}, \nonumber \\
&&
 X_{[3~1~1]}, \quad  X_{[4~0~1]}, \quad X_{[5~1~1]}, \quad
 X_{[6~1~1]}, \quad  X_{[3~1~2]}, \nonumber 
\end{eqnarray}
where ${\cal J}_1(X) = {\cal J}_2(X)$. The multiplicity is at most 2 
for any $[W] \in {\cal J}_2(X)$ for ten other singular K3 surfaces $X$ 
in the table. 
Even for the fourteen other singular K3 surfaces in the table, the 
multiplicity may only maximally be as large as $|(D_4 \times D_4)/\Z_4| = 16$
for some $[W] \in {\cal J}_2(X)$, which happens in the case of $X_{[6~0~6]}$. $\blacksquare$
\end{quote}

We also remark here that the coset space (\ref{eq:coset-alt1}) is not
necessarily a group, because $p_T({\rm Isom}(T_X)^{({\rm Hodge})})$ 
is not always a normal subgroup of ${\rm Isom}(q)$. In fact, 
for singular K3 surfaces $X=X_{[6~0~6]}$ and $X_{[6~6~6]}$ in 
Table \ref{tab:coset-alt-rho20}, the coset space (\ref{eq:coset-alt1}) 
is not a group.  

%%%%%%%%%%%%%%%%%%%%%%%%%%%%%%%%%%%%%%%%%%%%%%%%%%%%%%%%%%%%%
\subsection{Type-dependent Upper Bounds on Multiplicity Using
  the ${\rm II}_{1,25}$ Lattice}
\label{ssec:Borcherds-Kondo}
%%%%%%%%%%%%%%%%%%%%%%%%%%%%%%%%%%%%%%%%%%%%%%%%%%%%%%%%%%%%%

Propositions C and C' provide an upper bound on the number of 
isomorphism classes of elliptic fibrations (multiplicity) for a given 
K3 surface $X$ which is applicable for any isometry class 
$[W] \in {\cal J}_2(X)$. 
The number of isomorphism classes, however, can be different for 
different isometry classes $[W] \in {\cal J}_2(X)$ for a given $X$, and 
there must be room for deriving an upper bound   
on the number of isomorphism classes for individual isometry classes 
$[W] \in {\cal J}_2(X)$. Furthermore, the universal bound derived in the last section
is expected to be a relatively weak bound. 
In this section, we derive Proposition E 
and Corollary F by exploiting the theory on the structure of 
${\rm Isom}(S_X)$ developed by \cite{Borcherds2} and \cite{Kondo-auto-JacobKumm}.

The basic idea is this. First of all, since we consider the classification of 
elliptic fibrations $(\pi_X, \sigma; X, \P^1)$ modulo automorphism 
of $X$, an explicit choice of the section $\sigma$ is no longer important.
When $\sigma'$ is another section of $\pi_X: X \longrightarrow \P^1$, 
$(\pi_X, \sigma; X, \P^1)$ and $(\pi_X, \sigma'; X, \P^1)$ are 
mutually isomorphic because the translation automorphism (the group law
sum in the Mordell--Weil group extended also to singular fibres) maps 
one to the other \cite{Oguiso}. Thus, when an elliptic fibration is
given by a primitive embedding of $U \cong {\rm Span}_\Z\{f_1, f_2\}$ 
into $S_X$, $f_2$ does not play any other role than providing a divisor
intersecting $f_1$ once. We call divisors satisfying $f_1^2=0$ which 
also span a lattice $U$ with any other divisor in $S_X$ {\it elliptic divisors}.
\footnote{One might further impose one more condition---$f_1 \in \overline{\rm Amp}_X$---
for a divisor $f_1$ to be called an {\it elliptic divisor}.}
When we study the quotient space (\ref{eq:J1-2}), we only have to 
focus on elliptic divisors under the action of the groups 
${\rm Isom}(S_X)^{({\rm Amp~Hodge})}$ and 
${\rm Isom}(S_X)^{({\rm Amp})}$. 

We can always map an elliptic divisor for some isometry class 
$[W] \in {\cal J}_2(X)$
into $\overline{\rm Amp}_X$ by the $W^{(2)}(S_X)  \times \{ \pm 1_{S_X} \}$ 
subgroup of ${\rm Isom}(S_X)$ \cite{PSS}. For a given isometry class  
$[W] \in {\cal J}_2(X)$, such elliptic divisors within $\overline{\rm Amp}_X$
form a single orbit of ${\rm Isom}(S_X)^{({\rm Amp})}$. This orbit 
% under the group ${\rm Isom}(S_X)^{({\rm Amp})}$ 
is decomposed into a set of orbits of the subgroup 
${\rm Isom}(S_X)^{({\rm Amp~Hodge})}$.
The number of isomorphism classes for $[W] \in {\cal J}_2(X)$ is the number of 
${\rm Isom}(S_X)^{({\rm Amp~Hodge})}$ orbits within the single orbit 
of ${\rm Isom}(S_X)^{({\rm Amp})}$. The trouble in formulating the problem 
in this way, however, is that the groups ${\rm Isom}(S_X)^{({\rm Amp~Hodge})}$
and ${\rm Isom}(S_X)^{({\rm Amp})}$ may not be finite groups, and have 
complicated structures. 

The theory of \cite{Borcherds, Borcherds2, Kondo-auto-JacobKumm} 
is a powerful tool 
in studying the structure of these two groups. The crucial part
is an algorithm for setting a fundamental region under the action 
of these groups that is much smaller than ${\rm Amp}_X$---the fundamental 
region of the reflection group $W^{(2)}(S_X)$ generated only by $(-2)$ roots. 
In the following, we provide a brief review of the theory 
of \cite{Borcherds2, Kondo-auto-JacobKumm}, with a focus primarily 
on the aspect of setting a smaller fundamental region, while explaining 
notations and clarifying a set of sufficient conditions for the theory 
of \cite{Borcherds, Borcherds2, Kondo-auto-JacobKumm} to work.\footnote{ 
The appendix \ref{ssec:X3-app} will serve as a side reader on this subject. 
It deals with this theory applied to a particular (and the simplest)
example $X = X_3 = X_{[1~0~1]}$. Example 5.3 of
\cite{Borcherds2}---a single-page dense description---is extended into a 7
page-long pedagogical presentation there, so that even those without a firm
background in mathematics (including the present authors) can understand.}

\begin{center}
...........................................................................
\end{center}

Let us first assume that the lattice 
$T_0 \subset E_8 \oplus U^{\oplus (20-\rho)}$ (introduced 
in \cite{Nish01, Nish02}
and reviewed in section \ref{ssec:J2-Nishiyama}) is a direct sum of 
root lattices of $A$--$D$--$E$ type---({\bf as-1}).\footnote{
In the case discussed in \cite{Kondo-auto-JacobKumm}, $T_0$ itself is 
not a direct sum of root lattices of $A$--$D$--$E$ type, but contains 
such a lattice as an index 2 sublattice. There is room for a 
choice of $T_0$, more general than just being a direct sum of root lattices of 
$A$--$D$--$E$ type for the theory of \cite{Borcherds2, Kondo-auto-JacobKumm} 
to be applicable, but we do not try to present it in its most general form possible.} The isomorphism 
$\gamma: G_{T_0} \cong G_{T_X[-1]} \cong G_{S_X}$ consistent with the
discriminant form determines an embedding 
$\phi_{(T_0,S_X)*}: T_0 \oplus S_X \hookrightarrow {\rm II}_{1,25}$.  
Upon restriction to $T_0$, this is regarded as a primitive embedding 
of $T_0$: $\phi_{T_0*}: T_0 \hookrightarrow {\rm II}_{1,25}$. 
The Neron--Severi lattice $S_X$ is now regarded as the orthogonal
complement $[\phi_{T_0*}(T_0)^\perp \subset {\rm II}_{1,25}]$.
Theorem $\epsilon$ 
% Theorem 1.14.4 of \cite{Nik01} and Theorem 2.8 of \cite{Mor}
also guarantees, as long as $\rho_X \geq 12$, that for any primitive 
embedding $\phi_{T_0}: T_0 \hookrightarrow {\rm II}_{1,25}$, there
exists an isometry $f \in {\rm Isom}({\rm II}_{1,25})$ such that 
$f \cdot \phi_{T_0} = \phi_{T_0*}$. Thus, we can choose any primitive 
embedding $\phi_{T_0}: T_0 \hookrightarrow {\rm II}_{1,25}$, and 
regard the orthogonal complement as $S_X$, 
as long as $\rho_X \geq 12$---({\bf as-2}).
 
Let $J$ denote the simple roots of $T_0$ and $W_J$ the Weyl group of $T_0$.
If we find a map $\phi_{T_0}: J \longrightarrow \Pi$ preserving the 
Coxeter matrix (the intersection form),\footnote{$\Pi$ is the set of 
Leech roots, a set of simple roots of ${\rm II}_{1,25}$. 
See the appendix \ref{ssec:Leech-review} for more.}  then an embedding 
$\phi_{T_0}:T_0 \hookrightarrow {\rm II}_{1,25}$ is obtained 
by extending the map $\phi_{T_0}: J \longrightarrow \Pi$ linearly. In an
abuse of notation we use $\phi_{T_0}$ for the embeddings of both $J$ and $T_0$. 
We focus our attention to primitive embeddings $\phi_{T_0}$
satisfying this property. 

Suppose that the homomorphism $p_{T_0}: {\rm Isom}(T_0) \longrightarrow
{\rm Isom}(G_{T_0},q_{T_0}) \cong {\rm Isom}(G_{S_X},q_{S_X})$ is 
surjective. Since ${\rm Isom}(T_0) \cong W_J \rtimes {\rm Aut}(J)$, where
${\rm Aut}(J)$ is the group of automorphism of the Dynkin diagram of
$T_0$, this assumption---({\bf as-3}) is equivalent to the surjectiveness of 
\begin{equation}
 p_{T_0}: {\rm Aut}(J) \longrightarrow {\rm Isom}(q) \, ,
\end{equation}
where we use $p_{T_0}$ also for this homomorphism.
Hence all isomorphisms of $S_X$ (resp. ${\rm Isom}^+(S_X)$) are 
obtained by restricting certain isomorphism of ${\rm II}_{1,25}$ 
(resp. ${\rm Isom}^+({\rm II}_{1,25})$).
This assumption is satisfied by all the ten cases in Table \ref{tab:T0}
for which $T_0$ is a direct sum of root lattices of $A$--$D$--$E$ type.
The group of autochronous isomorphisms of ${\rm II}_{1,25}$ preserving 
the two subspaces $\phi_{T_0}(T_0)$ and $S_X$ has a structure 
(see \cite{Borcherds2}, Lemma 2.1)
\begin{equation}
 {\rm Isom}^+({\rm II}_{1,25})^{(\phi(T_0),S_X)} =
    W_{\phi_{T_0}(J)} \rtimes  W'_{\phi_{T_0}(J)}\, .
\end{equation}
Here, the Coxeter system $(W_{\phi_{T_0}(J)}, \phi_{T_0}(J))$ acts 
on\footnote{The representation of $W_\Pi$ on ${\rm II}_{1,25} \otimes
\R$ is regarded as the reduced representation of the canonical linear 
representation of $W_\Pi$ on $V := \{ \sum_{\lambda \in \Pi} 
x_\lambda e_\lambda \; | \; x_\lambda \in \R \}$ in the sense 
of \cite{Vin-lin-discr}.} 
${\rm II}_{1,25} \otimes \R$, and $W'_{\phi_{T_0}(J)}$ is the subgroup of 
${\rm Isom}^+({\rm II}_{1,25})$ mapping the fundamental chamber of 
the Coxeter group $W_{\phi_{T_0}(J)}$ to itself. 
Because the homomorphism ${\rm Aut}(J) \longrightarrow {\rm Isom}(q)$ 
is injective (which follows from the injectivity of 
${\rm Aut}(J) \longrightarrow {\rm Isom}(q)$ for any one of the root
lattices of $A$--$D$--$E$ type), restriction of this subgroup 
on $S_X$ induces an identification 
\begin{equation}
 {\rm Isom}^+({\rm II}_{1,25})^{(T_0,S_X)} \longrightarrow 
  W'_{\phi_{T_0}(J)} \cong {\rm Isom}^+(S_X).
\end{equation}

The $W^{(2)}(S_X)$ subgroup of ${\rm Isom}^+(S_X)$ can also be regarded 
as a subgroup of $W'_{\phi_{T_0}(J)} \subset {\rm Isom}^+({\rm II}_{1,25})$.
The fundamental chamber of the reflection group $W^{(2)}(S_X)$ is the 
ample cone of $X$, ${\rm Amp}_X \subset S_X \otimes \R$, but the approach of 
\cite{Borcherds2} and \cite{Kondo-auto-JacobKumm} is to 
exploit much a larger subgroup of ${\rm Isom}^+(S_X)$ in order to obtain 
a smaller fundamental region. 

Let us now briefly explain a few concepts and introduce some notation 
in order to spell out the statements obtained from the works of 
\cite{Borcherds2, Kondo-auto-JacobKumm}.
First, the Coxeter group $(W_\Pi, \Pi)$ in 
\begin{equation}
 W_\Pi \subset {\rm Isom}^+({\rm II}_{1,25}) \cong 
   W_\Pi \rtimes Co_\infty
\end{equation}
acts on ${\rm II}_{1,25} \otimes \R$ as a reflection group and 
the fundamental chamber is denoted by $C_\Pi$:  
\begin{equation}
C_\Pi = 
 \left\{ x \in {\rm II}_{1,25} \otimes \R \; | \; (x,\lambda) > 0 \; 
   {\rm for~} {}^\forall \lambda \in \Pi \right\}. 
\end{equation}
The half of $\{ x \in {\rm II}_{1,25} \otimes \R \; | \; x^2 > 0 \}$
containing $C_\Pi$ is called the {\it positive cone of} ${\rm II}_{1,25}$.
The interior of 
\begin{equation}
  (S_X \otimes \R) \cap \overline{C_\Pi}
\end{equation}
in $S_X \otimes \R$ is denoted by $D'$. 
We choose an isometry between 
$[\phi_{T_0}(J)^\perp \subset {\rm II}_{1,25}]$ and $S_X$ so that 
this $D'$ is contained in the ample cone ${\rm Amp}_X$ of $S_X\otimes \R$.

Secondly, in the Coxeter diagram\footnote{Nodes of the Coxeter diagram 
correspond to individual Leech roots in $\Pi$. A pair of nodes corresponding to $\lambda_i$ and $\lambda_j$ in $\Pi$ 
with $(\lambda_i, \lambda_j)_{{\rm II}_{1,25}} = \delta m_{ij}$ are 
joined by a single line if $\delta m_{ij} = +1$, 
by a thick line if $\delta m_{ij}=+2$, and by 
a dotted line if $\delta m_{ij} = +3, 4, \cdots $, following the conventions 
of \cite{Vin-Lobacevskii}. } 
of the Coxeter system $(W_\Pi, \Pi)$, 
any subdiagram is said to be {\it spherical} if it corresponds to one of 
the Dynkin diagrams of an $A$--$D$--$E$ root 
system. For any one of the $A$--$D$--$E$ root lattices, there is a unique 
element in the Weyl group that maps the Weyl chamber $D$ to $-D$ (note the 
simple transitive action on the chambers). This element 
is called {\it opposition involution} of a root lattice/system  
$R=A_n, D_n, E_{6,7,8}$, and is denoted by $\sigma_R$ or $\sigma(R)$.

Finally, the following group 
\begin{equation}
 {\rm Aut}(D') := \left\{ g \in Co_\infty \; | \; 
    g \left(\phi_{T_0}(J)\right) = \phi_{T_0}(J) \right\}
  = W'_{\phi_{T_0}(J)} \cap Co_\infty
\end{equation}
plays an important role. Now, we are ready to spell out the following 

\vspace{5mm}

\begin{quote}
{\bf Theorem} $\iota$   
(\cite{Borcherds, Borcherds2};  
\cite{Kondo-auto-JacobKumm}, Lemma 7.3): 
One can find a subgroup $N$ of 
${\rm Isom}(S_X)^{({\rm Amp~Hodge})}$ % of a K3 surface $X$ 
acting on $\overline{\rm Amp}_X$ so that the images of $\overline{D'}$ 
under $N$ cover the entire $\overline{\rm Amp}_X$ if, in addition to the assumptions (as-1), (as-2) 
and (as-3) that are stated already, the two following 
assumptions are satisfied:
\begin{itemize}
\item ({\bf as-4}): for any spherical subdiagram of $\Pi$ in the form of 
$R' := R \cup r$, $R:= \phi_{T_0}(J)$ for some $r \in \Pi$, either 
$\phi_{T_0}(J)$ is mapped to itself by $\sigma_{R'} \cdot \sigma_R$, or 
there exists an element $g' \in Co_\infty$ such that it is mapped 
to itself by $g' \cdot \sigma_{R'} \cdot \sigma_{R}$, and 
\item ({\bf as-5}): the homomorphism 
${\rm Aut}(D') \longrightarrow {\rm Aut}(J)$ is surjective.
\end{itemize}

When all of these assumptions are satisfied, it also follows that 
any element $g \in {\rm Isom}(S_X)^{({\rm Amp})}$ 
(resp. $g \in {\rm Isom}(S_X)^{({\rm Amp~Hodge})}$) can be written 
as $g = g_d \cdot n$ for some $n \in N$ and $g_d \in {\rm Aut}(D')$ 
(resp. $g_d \in {\rm Aut}(D')^{({\rm Hodge})}$). 
${\rm Aut}(D')^{({\rm Hodge})}$ is the inverse image of 
$p_T({\rm Isom}(T_X)^{({\rm Hodge})}) \subset {\rm Isom}(q)$.
\end{quote}

\vspace{5mm}

In the process of examining whether the assumptions (as-1)--(as-5) 
are all satisfied, one also has to carry out the following tasks:
\begin{itemize}
 \item [I] find an embedding $\phi_{T_0}: J \longrightarrow \Pi$ 
so that $\phi_{T_0}: T_0 \hookrightarrow {\rm II}_{1,25}$ is primitive,
 \item [II] list up Leech roots $r \in \Pi$ where $\phi_{T_0}(J) \cup
       r$ forms a spherical subdiagram of the Coxeter diagram of the
       Leech roots $\Pi$, 
 \item [III] compute 
       $\sigma_{\phi_{T_0}(J) \cup r} \cdot \sigma_{\phi_{T_0}(J)}$ for the Leech 
   roots $r$ listed up in II, and 
 \item [IV] compute the group ${\rm Aut}(D')$.
\end{itemize}

\begin{center}
.........................................................................
\end{center}

Let us now return to the problem set at the beginning of this 
section. If all of the assumptions 
(as-1)--(as-5) are satisfied, then the theorem above implies 
for any isomorphism class of elliptic fibrations 
$[(\pi_X, \sigma; X, \P^1)] \in {\cal J}_1(X)$, 
that its elliptic divisor (modulo ${\rm Aut}(X)$) can be chosen 
not just within $\overline{\rm Amp}_X$, but even within $\overline{D'}$.

Let us focus on an isometry class $[W] \in {\cal J}_2(X)$ and ask for 
its multiplicity in ${\cal J}_1(X)$. The ${\rm Aut}(D')$ group acts on the 
set of all the elliptic divisors $f$ in $\overline{D'}$ for the class $[W]$. 
Let 
\begin{equation}
 \amalg_a F^{[W]}_a := \amalg_a \left\{ {\rm Aut}(D') \cdot f_a \right\}
\end{equation}
denote the orbit decomposition of such elliptic divisors.
The action of ${\rm Aut}(D')$ is not necessarily transitive (explicit 
examples of the non-transitive action are found in \cite{Kumar}), and 
the index $a$ above is meant to label such different ${\rm Aut}(D')$ orbits.

\begin{quote}
{\bf Proposition E}: Suppose that the assumptions (as-1)--(as-5) are 
satisfied for a K3 surface $X$. Then for an isometry class 
$[W] \in {\cal J}_2(X)$, each isomorphism class of elliptic fibrations 
in ${\cal J}_1(X)$ falling into $[W]$ has its elliptic divisor in any
of the ${\rm Aut}(D')$ orbits, $\{ F^{[W]}_a \}$.
Thus, for any ${\rm Aut}(D')$ orbit $F^{[W]}_a$, its decomposition into 
the orbits of the action of the ${\rm Aut}(D')^{({\rm Hodge})}$ subgroup is 
regarded as a complete set of representatives of the ${\cal J}_1(X)$ 
classification for $[W] \in {\cal J}_2(X)$.
\end{quote}

{\bf proof}: If the assumptions (as-1)--(as-5) are satisfied, 
we can find an elliptic divisor $f$ in $\overline{D'}$ for any isomorphism 
class in ${\cal J}_1(X)$, as we have already seen above. This divisor 
$f$ must be in one of the ${\rm Aut}(D')$-orbits, say, $F^{[W]}_{a*}$. 
For any other ${\rm Aut}(D')$-orbit, say $F^{[W]}_{a'}$ (if there is any),
there must be a transformation $\phi \in {\rm Isom}(S_X)^{({\rm Amp})}$ 
mapping $f$ to $F^{[W]}_{a'}$, because the elliptic divisors in $F^{[W]}_{a'}$ 
belong to the same element $[W] \in {\cal J}_2(X)$ as $f \in F^{[W]}_{a*}$. 
Thanks to the assumption (as-5), there must be a transformation 
$g \in {\rm Aut}(D')$ such that 
$g \cdot \phi \in {\rm Isom}(S_X)^{({\rm Amp~Hodge})}$, which means that 
both the elliptic divisors $f$ and $g \cdot \phi ( f )$ define elliptic 
fibrations in the same isomorphism class, i.e. the same element 
in ${\cal J}_1(X)$). $\blacksquare$

\begin{quote}
{\bf Corollary F}: For an isometry class $[W] \in {\cal J}_2(X)$, 
the number of isomorphism classes of elliptic fibrations (= multiplicity = 
the number of elements in the inverse image of $[W]$ in ${\cal J}_1(X)$) 
is bounded from above by the number of ${\rm Aut}(D')^{({\rm Hodge})}$-orbits 
within any one of $F^{[W]}_a$, and in particular, by the smallest one 
among them. $\blacksquare$
\end{quote}

We remark that the complete sets of representatives in Proposition E are  
not necessarily minimal complete sets of representatives of isomorphism 
classes. In order to find a minimal complete set of representatives 
one would have to exploit more information of $X$ as in the original work 
of \cite{Oguiso} for $X = {\rm Km}(E \times F)$.

%%%%%%%%%%%%%%%%%%%%%%%%%%%%%%%%%%%%%%%%%%%%%%%%%%%%%%%%%%%%%%%%%%
\subsection{An Example: $X = {\rm Km}(A)$}
%%%%%%%%%%%%%%%%%%%%%%%%%%%%%%%%%%%%%%%%%%%%%%%%%%%%%%%%%%%%%%%%%%%

As the first application of Corollary F, we choose a $\rho = 17$ family of 
K3 surfaces $X = {\rm Km}(A)$, for which the tasks I--IV in 
section \ref{ssec:Borcherds-Kondo} have been carried out and all 
the assumptions (as-1)--(as-5) verified already \cite{Kondo-auto-JacobKumm}. 
Furthermore, one representative elliptic divisor $f^{[W]}_a$ 
has been determined for any of the ${\rm Aut}(D')$ orbits $F^{[W]}_a$ 
of any $[W] \in {\cal J}_2(X)$ in \cite{Kumar} for this 
family of K3 surfaces. 

There are 25 isometry classes in the ${\cal J}_2(X)$ classification, and the
representative elliptic divisors obtained by \cite{Kumar} are shown 
in the 2nd column of Table \ref{tab:Kummer-Kumar}. 
We now only have to work out the action of \cite{Kondo-auto-JacobKumm}
\begin{equation}
   {\rm Aut}(D') = (\Z/2\Z)^5 \rtimes S_6, \qquad 
   {\rm Aut}(D')^{({\rm Hodge})} = (\Z/2\Z)^5   
\end{equation}
on the elliptic divisors $f^{[W]}_a$ to see how many ${\rm
Aut}(D')^{({\rm Hodge})}$-orbits individual ${\rm Aut}(D')$-orbits 
$F^{[W]}_a$ are decomposed into. Necessary technical details 
on these two groups as well as $S_X$ of this family of K3 surfaces are 
summarized in the appendix \ref{ssec:review-Kummer}, more interested readers may prefer to consult 
to \cite{Barthetal, Nik-Kummer, Keum-auto-JacobKumm, 
Kondo-auto-JacobKumm}.

%%%%%%%%%%%%%%%%%%%%%%%%%%%%%%%%%%%%%%%%%%%%%%%%%%%%%%%%%%%%%%%%%%%%%%%%%%%%
\begin{table}[tbp]
\begin{center}
\begin{tabular}{|c||cc|c|}
\hline
$[W_{\rm root}]$ & & Elliptic Divisor & \# orbits \\
\hline
\hline 
$D_4^2A_1^6$ & $f^{[1]}$ &
    $(-1, -1, 0, 0, 0, 0, 0, 0, 0, 0, 0, 0, 0, 0, 0, 0, 1)$ &  15 \\
\hline
$D_6D_4A_1^4$ & $f^{[2]}$ & 
    $(-2,-1,0,0,0,-1,0,0,0,0,0,-1,-1,0,0,0,2)$ &  360 \\
\hline
$D_6A_3A_1^6$ & $f^{[3]}$ &
    $(-2,-1,0,0,-1,0,0,0,0,0,-1,0,-1,0,0,0,2)$ & 90 \\
\hline
$D_4^3A_3$ & $f^{[4]}$ &
    $(-2,-1,0,0,0,0,0,0,0,0,-1,-1,0,-1,0,0,2)$ & 60 \\
\hline
$A_7A_3^2$ & $f^{[5]}$ &
    $\frac{1}{2} (-3,-3,-1,-1,0,0,-1,-1,0,0,-1,-2,0,-2,0,-1,4)$ & 90 \\
\hline
$A_7A_3A_1^2$ & $f^{[6]}$ & 
    $\frac{1}{2} (-4,-2,-1,0,0,-1,-1,0,0,-1,-1,-1,0,-2,-1,-1,4)$ & 45  \\
\hline
$A_3^4$ & $f^{[7]}$ & 
    $\frac{1}{2} \{-1,-1,0,0,0,0,0,0,0,0,-1,-1,-1,-1,-1,-1,2)$ & 15 \\
\hline
$D_6^2A_2$ & $f^{[8]}_0$ &
    $(-3,-2,-1,0,0,0,0,0,-1,0,-1,-1,0,-1,0,0,3)$ & 360 \\
% \cline{2-4}
 & $f^{[8]}_A$ & 
    $(-2,-2,0,0,0,-1,0,0,0,-1,0,-2,-2,0,0,0,3)$ & 180  \\
\hline
$D_6D_5A_1^4$ & $f^{[9]}$ & 
    $(-3,-2,-1,0,0,0,0,0,0,-1,-1,-1,0,-1,0,0,3)$ & 360 \\
\hline
$A_5^2A_1^2$ & $f^{[10]}$ &
    $\frac{1}{2} (-2,-2,-1,0,0,-1,0,-1,-1,0,-2,0,-1,-1,0,0,3)$ & 60\\
\hline
$D_8A_1^6$ & $f^{[11]}_0$ & 
    $(-3,-2,0,0,-1,0,0,0,0,-1,-1,-1,-1,0,0,0,3)$ & 360 \\
% \cline{2-4}
  & $f^{[11]}_A$ & $(-2,-2,0,0,0,-2,0,-1,0,0,0,-1,-2,0,0,0,3)$ & 180 \\
\hline
$E_6D_5$ & $f^{[12]}$ &
    $\frac{1}{2} (-5,-4,0,0,0,-3,-1,-1,-1,0,-1,-3,-2,-1,0,-2,6)$ & 360 \\
\hline 
$E_6D_4$ & $f^{[13]}$ & 
    $\frac{1}{2} (-6,-3,0,0,-3,0,-1,-1,0,-1,-2,-1,-2,-1,-2,-1,6)$ & 60 \\
\hline
$A_5^2$ & $f^{[14]}$ & 
    $\frac{1}{2} (-3,-1,0,0,-1,-1,-1,-1,0,0,0,-1,-1,-1,-1,0,3)$ & 10 \\
\hline
$D_8D_4A_3$ & $f^{[15]}_0$ & 
    $(-4,-2,0,0,0,-2,0,0,-1,-1,0,-2,-1,-1,0,0,4)$ & 720 \\
% \cline{2-4}
  & $f^{[15]}_A$ & $(-3,-3,0,0,0,-2,0,-1,0,-1,0,-2,-2,0,0,0,4)$ & 720 \\
\hline
$E_7D_4A_1^3$ & $f^{[16]}_0$ & 
    $(-4,-2,0,0,0,-2,-1,0,-1,0,0,-2,-1,-1,0,0,4)$ & 720 \\
  & $f^{[16]}_A$ & $(-3,-3,0,0,0,-2,-1,0,0,0,0,-2,-2,0,0,-1,4)$ & 720  \\
\hline
$D_7D_4^2$ & $f^{[17]}$ & 
    $(-4,-2,0,0,-1,-1,0,-1,0,0,0,-2,-2,0,0,-1,4)$ & 360 \\
\hline
$E_7A_3A_1^5$ & $f^{[18]}_0$ & 
    $(-3,-3,0,0,0,-2,-1,-1,0,0,0,-2,-2,0,0,0,4)$ & 360 \\
 & $f^{[18]}_A$ & $(-4,-2,0,0,0,-2,0,0,-1,0,-1,-2,-1,-1,0,0,4)$ & 360 \\
\hline
$D_5^2$ & $f^{[19]}$ & 
    $\frac{1}{2} (-3,-3,-2,0,0,0,0,0,0,-2,-1,-1,-1,-1,-1,-1,4)$ & 45 \\
\hline
$D_{10}A_1^4$ & $f^{[20]}_0$ & 
    $(-5,-2,0,-1,0,-2,0,0,-1,-1,0,-3,-2,-1,0,0,5)$ & 720 \\
 & $f^{[20]}_A$ & $(-4,-4,0,0,-2,0,0,-1,0,-1,-2,-2,-2,0,0,0,5)$ & 360 \\
 & $f^{[20]}_B$ & $(-4,-3,0,0,-1,-2,0,-1,0,-1,0,-3,-3,0,0,0,5)$ & 720 \\
\hline
$A_9A_1^3$ & $f^{[21]}$ & 
    $\frac{1}{2} (-4,-3,0,0,-1,-2,0,-2,-1,0,-1,-2,-3,0,-1,0,5)$ & 180 \\
\hline 
$A_9A_1$ & $f^{[22]}$ & 
    $\frac{1}{2} (-5,-2,0,-1,0,-2,0,-1,-2,0,-1,-2,-2,-1,-1,0,5)$ & 72 \\
\hline 
$D_9A_1^6$ & $f^{[23]}$ & 
    $(-5,-2,0,0,-2,-1,0,-1,-1,0,-3,0,-2,-1,0,0,5)$ & 360 \\
\hline
$E_8A_1^6$ & $f^{[24]}_0$ & $(-6,-5,0,0,-3,0,0,-1,0,-2,-3,-3,-2,-1,0,0,7)$ & 720 \\
 & $f^{[24]}_A$ & $(-4,-5,0,0,0,-5,0,-3,-1,0,0,-2,-4,-1,-1,0,7)$ & 720 \\
 & $f^{[24]}_B$ & $(-6,-4,0,0,-2,-2,0,-2,-1,-1,0,-4,-4,0,0,0,7)$ & 360 \\
 & $f^{[24]}_C$ & $(-3,-2,-4,0,-3,-1,0,-1,0,-1,-7,0,0,0,-2,-2,7)$ & 720  \\
\hline
$D_9$ & $f^{[25]}$ &  
    $\frac{1}{2} (-7,-5,0,-2,-2,0,0,0,0,-4,-3,-3,-3,-1,-1,-1,8)$ & 360 \\
\hline
\hline
\end{tabular}
\caption{\label{tab:Kummer-Kumar} The first two columns quote the results of 
${\cal J}_2(X)$ classification for $X = {\rm Km}(A)$ from \cite{Kumar}.
A representative elliptic divisor in each ${\rm Aut}(D')$-orbit is 
described as linear combinations of divisors 
$N_{00}, N_{01}, \cdots, N_{05}, N_{12}, \cdots, N_{15},
 N_{23}, \cdots, N_{25}, N_{34}, N_{35}, N_{45}$ and $H$ in $S_X$, and their 17
 coefficients are shown in the 2nd column.
The last column shows our computation on the number of 
${\rm Aut}(D')^{({\rm Hodge})}$-orbits within the individual 
${\rm Aut}(D')$-orbits. }
\end{center}
\end{table}
%%%%%%%%%%%%%%%%%%%%%%%%%%%%%%%%%%%%%%%%%%%%%%%%%%%%%%%%%%%%%%%%%%%%%%%%%%

We have studied the decomposition of each of the 
${\rm Aut}(D')$-orbits (each row in Table \ref{tab:Kummer-Kumar}) into 
the orbits under the ${\rm Aut}(D')^{({\rm Hodge})}$ subgroup.
Let us take the first isometry class in ${\cal J}_2(X)$, where 
$W_{\rm root} = D_4^{\oplus 2} \oplus A_1^{\oplus 6}$, as an example.
It turns out that the ${\rm Aut}(D')$-orbit, ${\rm Aut}(D') f^{[1]}$ 
consist of 240 elliptic divisors in $\overline{D'}$, which are grouped into 
15 distinct ${\rm Aut}(D')^{({\rm Hodge})}$-orbits, each one of 
which consists of 16 elliptic divisors. 
Similar computations have been carried out for all other elliptic
divisors, and the result is presented in the last column of 
Table \ref{tab:Kummer-Kumar}.
Note that, in the case of $X = {\rm Km}(A)$, all the 
${\rm Aut}(D')^{({\rm Hodge})}$-orbits within a given 
${\rm Aut}(D')$-orbit consist of the same number of elliptic
divisors because ${\rm Aut}(D')^{({\rm Hodge})}$ is a normal subgroup 
of ${\rm Aut}(D')$. The labels $0,A,B,\cdots$ refer to different elliptic divisors 
equivalent under ${\rm Aut}(X)$, as given in \cite{Kumar}. Interestingly, different
bounds are established starting from different elliptic divisors.
\begin{quote}
 {\bf Example G}: The multiplicity is bounded from 
above by 15 (no.1), 360 (no.2), 90 (no.3), 60 (no.4), 90 (no.5),
45 (no.6), 15 (no.7), 180 (no.8), 360 (no.9), 60 (no.10),  
180 (no.11), 360 (no.12), 60 (no.13), 10 (no.14), 720 (no.15), 
720 (no.16), 360 (no.17), 360 (no.18), 45 (no.19), 360 (no.20), 
180 (no.21), 72 (no.22), 360 (no.23), 360 (no.24), and 360 (no.25), 
respectively, for the twenty-five isometry classes in ${\cal J}_2(X)$
for $X = {\rm Km}(A)$. (We follow the numbering (no.1--25) on the 
isometry classes as used in \cite{Kumar}.) $\blacksquare$
\end{quote}

Note that these upper bounds on the multiplicity for individual 
isometry classes $[W] \in {\cal J}_2(X)$ is stronger than the uniform 
upper bound obtained from Proposition C': $|S_6| = 720$.
Note also that the actual multiplicity for a given $[W] \in {\cal J}_2(X)$ 
may be smaller than the upper bound obtained here. 

%%%%%%%%%%%%%%%%%%%%%%%%%%%%%%%%%%%%%%%%%%%%%%%%%%%%%%%%%
\subsection{A Systematic Way to Find an Elliptic Divisor in $\overline{D'}$}
\label{ssec:systematic-23}
%%%%%%%%%%%%%%%%%%%%%%%%%%%%%%%%%%%%%%%%%%%%%%%%%%%%%%%%%%

In order to apply Corollary F to other K3 surfaces, 
we need to be able to find at least one elliptic divisor in
$\overline{D'}$ for each isometry class $[W] \in {\cal J}_2(X)$. 
In this subsection, we present a systematic procedure to find 
such an elliptic divisor,\footnote{An alternative will be to implement 
this problem on a computer and pick up candidates of elliptic divisors from 
the edges of $\overline{D}'$, as in \cite{Kumar}.} 
where one can combine both a) the Kneser--Nishiyama method for the 
${\cal J}_2(X)$ classification using Niemeier lattices (reviewed in 
section \ref{ssec:J2-Nishiyama}) and 
b) Borcherds and Kondo's theory of handling the isometry group of $S_X$ 
using the ${\rm II}_{1,25}$ lattice 
(reviewed in section \ref{ssec:Borcherds-Kondo}).

Because of the uniqueness of the even unimodular lattice of signature 
(1,25), there must be isometries among all the twenty-four lattices 
of the form $U \oplus L^{(I)}$, where $I=\alpha, \beta, \cdots, \omega$ and 
$L^{(I)}$'s are the twenty-four Niemeier lattices. The isometry between 
the ones with $I=\alpha, \beta, \cdots, \psi$ and the Leech lattice $L^{(\omega)}$
is denoted by  
\begin{equation}
  \phi^{(I)}:  L^{(I)} \oplus U^{(I)} \cong \Lambda_{24} \oplus U^{(\omega)} \cong 
  {\rm II}_{1,25}.
 \label{eq:23-vs-omega}
\end{equation}
It is known (\cite{CS-book}, Chap. 26 Thm. 5) that, 
for a vector
\begin{equation}
  u^{(I)} :=  (\vec{c}^{(I)}, 1, - (\vec{c}^{(I)})^2/2 ) \in 
       \left(\Lambda_{24} \oplus U^{(\omega)} \right) \otimes \Q,
\end{equation}
the twenty-three other Niemeier lattices $L^{(I)}$ are obtained by 
$[(u^{(I)} )^\perp \subset (\Lambda_{24} \oplus U^{(\omega)})]/
(u^{(I)})$. Here, $\vec{c}^{(I)}$ is a centre of a deep hole of type 
$I=\alpha, \beta, \cdots, \psi$, see appendix \ref{ssec:Leech-review} 
for a brief review, or \cite{CS-book} for extensive exposition.
We start off by refining this theorem a little more to the level of 
constructing the isometry $\phi^{(I)}$ in (\ref{eq:23-vs-omega}) explicitly.

Let $(L^{(I)})_{\rm root} \cong \oplus_{a \in A} R_a$, where 
each one of the $R_a$'s is a root lattice of $A$--$D$--$E$ type 
(see Table \ref{tab:Niemeier-root}). 
Let $\vec{c}^{(I)} \in \Lambda_{24} \otimes \Q$ be a centre of a deep hole 
of type $I$, and $\vec{v}^{(I)a}_i$ ($i=0,\cdots, r_a = {\rm rank}(R_a)$) 
the integral points of $\Lambda_{24}$ surrounding it. 
With this data, an embedding 
$\phi^{(I)}: (L^{(I)})_{\rm root} \hookrightarrow \Lambda_{24} \oplus U^{(\omega)}$ 
is specified by mapping the simple roots of $R_a$  
\begin{equation}
 \phi^{(I)} : \alpha^{(I)a}_i \longmapsto \lambda^{(I)a}_i = 
    \left(\vec{v}^{(I)a}_i, 1, -1-\frac{(\vec{v}^{(I)a}_i)^2}{2} \right)
\label{eq:def-phi-I}
\end{equation}
for $i = 1, \cdots, r_a$ ($i=0$ is not included here). 

We claim that this embedding of 
$(L^{(I)})_{\rm root} \hookrightarrow \Lambda_{24} \oplus U^{(\omega)}$ can be 
lifted to an isometry (\ref{eq:23-vs-omega}). To see this, we only need to 
prove that the orthogonal complement of this embedding, 
$[\phi^{(I)}( \oplus_{a\in A} R_a )^\perp \subset (\Lambda_{24} \oplus U^{(\omega)})]$, is
isometric to the hyperbolic plane lattice $U$---(**). 
In order to prove (**), let us first see that the two vectors 
\begin{eqnarray}
 u_1^{(I)} & = & 
   \left(h^{(I)} \vec{c}^{(I)}, h^{(I)}, - h^{(I)} (\vec{c}^{(I)})^2/2  \right), 
      \qquad   \label{eq:f1-as-elliptic-div} \\
 u_2^{(I)} & = & 
   \left(\vec{c}^{(I)}+ \frac{1}{h^{(I)}} \vec{\rho}^{(I)},1,
         - \frac{(\vec{c}^{(I)})^2}{2}
         + \frac{1}{h^{(I)}} \left( 1-(\vec{\rho}^{(I)},\vec{c}^{(I)}) \right)
    \right), 
\end{eqnarray}
are integral elements of $\Lambda_{24} \oplus U^{(\omega)}$. Here, 
\begin{equation}
 \vec{\rho}^{(I)} := \sum_{a \in A} \vec{\rho}^{(I)a}, \qquad 
  \left( \vec{\rho}^{(I)a} , \vec{v}_i^{(I)b}-\vec{c}^{(I)}
  \right)_{\Lambda_{24} \otimes \Q} = - \delta^{ab} \quad (i=1,\cdots, r_b);
\end{equation}
that is, $\vec{\rho}^{(I)a}$ is the Weyl vector of the (negative definite) 
root lattice $R_a$ contained in $(L^{(I)})_{\rm root} \cong \oplus_{a \in A} R_a$.
One can show after a little manipulation (using an arbitrary $a \in A$) that 
\begin{eqnarray}
 - \frac{h^{(I)}}{2} (\vec{c}^{(I)})^2 & = &
     h^{(I)} \left(1 + \frac{(\vec{v}^{(I)a}_0)^2}{2} \right) 
   - (\vec{v}^{(I)a}_0, h^{(I)} \vec{c}^{(I)} ) \in \Z, \\
 - \frac{(\vec{c}^{(I)})^2}{2}
 + \frac{1}{h^{(I)}}(1-(\vec{\rho}^{(I)}, \vec{c}^{(I)}))  & = & 
   2 + \frac{(\vec{v}^{(I)a}_0)^2}{2}
   - \left(\vec{v}^{(I)a}_0, \vec{c}+ \frac{1}{h^{(I)}}\vec{\rho}^{(I)} \right).
\end{eqnarray}
Chap. 24, section 2 of \cite{CS-book} guarantees that 
$\vec{c}^{(I)} + \frac{1}{h^{(I)}} \vec{\rho}^{(I)} \in \Lambda_{24}$, and
hence both $u_1$ and $u_2$ are indeed integral elements in 
$\Lambda_{24} \oplus U^{(\omega)}$. % rather than ${\rm II}_{1,25} \otimes \Q$.
Secondly, it is easy to see that these two vectors are orthogonal to 
$\phi^{(I)}(\oplus_a R_a)$, and hence 
${\rm Span}_\Z \left\{ u_1, u_2 \right\} \subset 
[\phi^{(I)}(\oplus_a R_a)^\perp \subset (\Lambda_{24} \oplus U^{(\omega)})]$. 
Finally, the symmetric pairing on ${\rm Span}_\Z \{ u_1, u_2\}$ turns out 
to be $\left[ \begin{array}{cc} 0 & 1 \\ 1 & -2 \end{array}\right]$, where 
the following relation \cite{Borcherds-thesis}
\begin{equation}
 (\vec{\rho}^{(I)})^2 = - 2 h^{(I)} (h^{(I)} + 1)
\end{equation}
is used to compute $(u_2^{(I)})^2$.
This means that ${\rm Span}_\Z \{ u_1, u_2 \}$ is isometric to $U$, and 
forms a primitive sublattice of $\Lambda_{24} \oplus U^{(\omega)}$. 
It thus follows that 
\begin{equation}
 \left[ \phi^{(I)}(\oplus_a R_a)^\perp \subset 
        \left(\Lambda_{24} \oplus U^{(\omega)} \right)  \right] = 
  {\rm Span}_\Z \left\{ u_1^{(I)}, u_2^{(I)} \right\} \cong U.
\end{equation}
The claim (**) is now proven, and we have 

\begin{quote}
{\bf Lemma H}: An isometry $\phi^{(I)}$ in (\ref{eq:23-vs-omega}) is obtained 
by mapping $L^{(I)}$ by (\ref{eq:def-phi-I}) and embedding $U^{(I)}$ into 
the hyperbolic plane ${\rm Span}_\Z\{u_1^{(I)}, u_2^{(I)} + u_1^{(I)}\}$. 
Furthermore, the $u_1^{(I)}$ can be written as a positive
coefficient sum of Leech roots under this isometry: 
\begin{equation}
  u_1^{(I)} = \sum_{i=0}^{r_a} n_i^{(I)a} \lambda_i^{(I)a} \qquad 
   {\rm for~}{}^\forall  a \in A, 
\label{eq:lmmH-2}
\end{equation}
where $\lambda_0^{(I)a} := (\vec{v}^{(I)a}_0, 1,
 -1-(\vec{v}^{(I)a}_0)^2/2)$.
\end{quote}

Suppose that the lattice $T_0$ for a K3 surface $X$ is a direct sum of 
root lattices of $A$--$D$--$E$ type, and a set of its simple roots is denoted 
by $J$, as in section \ref{ssec:Borcherds-Kondo}. 
For an isometry class $[W] \in {\cal J}_2(X)$ associated with 
a primitive embedding $\phi'_{T_0}: T_0 \hookrightarrow L^{(I)}$, suppose that 
the primitive embedding $\phi_{T_0}: T_0 \hookrightarrow \Lambda_{24} \oplus 
U^{(\omega)} \cong {\rm II}_{1,25}$ is given by a combination of 
$\phi'_{T_0}: J \longrightarrow ({\rm simple~roots~of})~L^{(I)}_{\rm root}$ 
and $\phi^{(I)}$ in Lemma H. Then the Neron--Severi lattice 
$S_X = [\phi_{T_0}(T_0)^\perp \subset {\rm II}_{1,25}]$ is isomorphic to 
\begin{equation}
{\rm Span}_\Z\{u_1^{(I)}, u_2^{(I)} + u_1^{(I)} \} \oplus \phi^{(I)} \left( 
    W = [ \phi'_{T_0}(T_0)^\perp \subset L^{(I)} ] \right).
\end{equation}
Therefore, $u_1^{(I)}$ and $u_2^{(I)}+u_1^{(I)}$ in 
$\Lambda_{24} \oplus U^{(\omega)}\cong{\rm II}_{1,25}$ can be regarded as 
elliptic divisors of this isometry class $[W] \in {\cal J}_2(X)$.

\begin{quote}
{\bf Lemma I}: The elliptic divisor $u_1^{(I)}$ is in $\overline{D'}$.
\end{quote}

{\bf proof}: It is in $S_X$, and also in $\overline{C_\Pi}$ 
because $(u_1^{(I)}, \lambda_i^{(I)a}) = 0$ for all of $i=0,\cdots,r_a$, 
$a \in A$, and $(u_1^{(I)},\lambda)$ is positive for all other Leech roots 
$\lambda \in \Pi$. $\blacksquare$

%%%%%%%%%%%%%%%%%%%%%%%%%%%%%%%%%%%%%%%%%%%%%%%%%%%%%%%%%%%%%
\subsection{Another Example: $X={\rm Km}(E_\omega \times E_\omega)$}
\label{ssec:KmEoxEo-multiplicity}
%%%%%%%%%%%%%%%%%%%%%%%%%%%%%%%%%%%%%%%%%%%%%%%%%%%%%%%%%%%%%%

In this section, we take 
a singular K3 surface $X = {\rm Km}(E_\omega \times E_\omega) = X_{[2~2~2]}$ 
as an example, and apply Corollary F. 
The tasks I--IV in section \ref{ssec:Borcherds-Kondo} have been carried
out and all the assumptions (as-1)--(as-5) are verified also for this 
K3 surface \cite{KK-auto-prod}. The ${\cal J}_2(X)$ classification has 
also been worked out in \cite{Nish01} for this K3 surface, and there are 
30 different isometry classes. However, we are not aware of an identification
of the elliptic divisors, in $\overline{D'}$ in particular, for these isometry classes
in the literature.  
Therefore, we combine Lemma I and Corollary F with all that is  
known from \cite{Nish01} and \cite{KK-auto-prod} to derive an upper 
bound on the number of isomorphism classes of elliptic fibrations 
for each isometry class $[W] \in {\cal J}_2(X)$ of 
$X = {\rm Km}(E_\omega \times E_\omega)$.

The lattice $T_0$ for this K3 surface is $D_4 \oplus A_2$ \cite{Nish01}, 
and for all the thirty types in ${\cal J}_2(X)$, we stick to the
following embedding of the simple roots of $D_4 \oplus A_2$ to the 
Leech roots. First, we define the following 6 vectors in $\Lambda_{24}$:
\begin{eqnarray}
  \vec{v}_1 & = & \nu_\Omega + 4 \nu_{\infty}=: X, 
\label{eq:D4-A2-embed-II-1-25-1}\\
  \vec{v}_2 & = & \vec{0} =: Z, \\
  \vec{v}_3 & = & \nu_{\Omega} + 4 \nu_{0} =:Y, \\
  \vec{v}_4 & = & \nu_\Omega - 2\nu_{K_U} + 4 \nu_{\{\infty, 0,1\}} =: U, \\
  \vec{v}_p & = & 4 \nu_{\{0,\infty\}} =: P, \\
  \vec{v}_q & = & \nu_\Omega-4\nu_2 =: Q_2,
\label{eq:D4-A2-embed-II-1-25-q}
\end{eqnarray}
where $K_U$ is a codeword in ${\cal C}_{24}(8)$ containing 
$\{ \infty,0,1,2\}$ as a subset. Leech roots that correspond to 
$\vec{v}_{1,2,3,4}$ and $\vev{v}_{p,q}$ in $\Lambda_{24}$ 
through the relation (\ref{eq:L24-Pi-relation}) are denoted by 
$\lambda_{1,2,3,4}$ and $\lambda_{p,q}$, respectively. The embedding 
$\phi_{T_0}: J \hookrightarrow \Pi$ of simple roots is given by 
assigning $\lambda_{1,2,3,4}$ to the simple roots $\alpha_{1,2,3,4}$ of 
$D_4 \subset T_0$, and $\lambda_{p,q}$ to those of $A_2 \subset T_0$ 
\cite{KK-auto-prod}. For $K_U$ we use the codeword given 
in (\ref{eq:Ku}).\footnote{This choice (\ref{eq:Ku}) for $K_U$ is not
the same as $R_3$ in \cite{KK-auto-prod}. We use $K_U$ (\ref{eq:Ku})
for detailed computations in this section (and in 
the appendix \ref{ssec:review-Kummer}) only because some codewords in the MOG
representation looks nice with (\ref{eq:Ku}), no other reasons.
Our choice $K_U$, (\ref{eq:Ku}), is equivalent to $R_3$ in \cite{KK-auto-prod} 
modulo $M_{24}$. In particular, the permutation in $M_{24}$ mapping $K_U$ to 
$R_3$ is such that it leaves all other vectors except $\vec{v}_4$ 
in (\ref{eq:D4-A2-embed-II-1-25-1}--\ref{eq:D4-A2-embed-II-1-25-q}) 
invariant. Hence our embedding is equivalent to the one used 
in \cite{KK-auto-prod}, and it is not necessary to repeat the Tasks
I--IV and verify (as-4) independently. 
}
Once this embedding is fixed, then $S_X$ for 
$X = {\rm Km}(E_\omega \times E_\omega)$ is obtained as the orthogonal 
complement of $\phi_{T_0}(T_0)$ in $\Lambda_{24} \oplus U^{(\omega)}
\cong {\rm II}_{1,25}$. ${\rm Aut}(D')$ and ${\rm Aut}(D')^{({\rm
Hodge})}$ were determined for this set-up in \cite{KK-auto-prod}, whose 
result is quoted in appendix \ref{ssec:review-Kummer}.

Since each isometry class $[W] \in {\cal J}_2(X)$ is obtained in the 
form of $W := [ \phi'_{T_0}(T_0)^\perp \subset L^{(I)}]$ for some 
primitive embedding into one of the Niemeier lattices, Lemma I can be used to 
determine an elliptic divisor for this $[W] \in {\cal J}_2(X)$, as 
we have stated prior to Lemma I. In practice, however, we need a centre 
$\vec{c}^{(I)} \in \Lambda_{24} \otimes \Q$ of a deep hole of type $I$ 
(and vectors $\vec{v}_i^{(I)a}$'s around the centre) in order to
determine the isometry $\phi^{(I)}$. It looks as if we can use the
centre $\vec{c}^{(I)}$ and vectors $\vec{v}^{(I)a}_i$  given explicitly 
in Chapt. 23 of \cite{CS-book}, but the lattice $T_0$ is now embedded into 
${\rm II}_{1,25}$ through $\phi^{(I)} \cdot \phi'_{T_0}$, which is 
often different from the embedding determined 
by (\ref{eq:D4-A2-embed-II-1-25-1}--\ref{eq:D4-A2-embed-II-1-25-q}).
It is one possibility to work out how the groups ${\rm Aut}(D')$ 
and ${\rm Aut}(D')^{({\rm Hodge})}$ act on the orthogonal complement of 
$\phi^{(I)} \cdot \phi'_{T_0} (T_0)$ for each isometry class 
$[W] \in {\cal J}_2(X)$. An alternative, however, is to find a centre 
of deep hole of type $I$ (and the vectors $\vec{v}_i^{(I)a}$ around it) 
for each 
$[W] = [\phi'_{T_0}(T_0)^\perp \subset L^{(I)}] \in {\cal J}_2(X)$ so
that the vectors 
(\ref{eq:D4-A2-embed-II-1-25-1}--\ref{eq:D4-A2-embed-II-1-25-q}) are
always included. We will take the latter approach in the following.

We have determined an elliptic divisor for each one of thirty isometry 
classes in ${\cal J}_2(X)$ for $X={\rm Km}(E_\omega \times E_\omega)$. 
See Table \ref{tab:KmEoxEo-multiplicity} for the results.
Nishiyama assigned numbers (type identification number) from 1 to 30 
to the thirty isometry classes $[W] \in {\cal J}_2(X)$ in Table 1.3 of 
\cite{Nish01}. We use the same type identification number in Table 
\ref{tab:KmEoxEo-multiplicity}. 
Sections \ref{sssec:D4toDn-A2toOther}--\ref{sssec:D4A2toE8} provide 
detailed information of the determination process of elliptic divisors.
Upper bounds on the multiplicities for individual isometry classes 
are discussed in section \ref{sssec:KmEwxEw-multiplicity}.

%%%%%%%%%%%%%%%%%%%%%%%%%%%%%%%%%%%%%%%%%%%%%%%%%%%%%%%%%%%%%%%%%%%
\subsubsection{Determination of Elliptic Divisors: 
$D_4\hookrightarrow D_n$ ,\, \,$A_2\hookrightarrow {\rm other}$.} 
\label{sssec:D4toDn-A2toOther}
%%%%%%%%%%%%%%%%%%%%%%%%%%%%%%%%%%%%%%%%%%%%%%%%%%%%%%%%%%%%%%%%%%%

Let us begin with the first 13 types\footnote{In the rest of this
article, we refer to each {\it isometry class} $[W] \in {\cal J}_2(X)$
as each {\it type (of elliptic fibration)} for brevity, given the close
relation between ${\cal J}_2(X)$ and ${\cal J}^{({\rm type})}(X)$ 
classification.} in ${\cal J}_2(X)$ of 
$X = {\rm Km}(E_\omega \times E_\omega)$ listed up in Table 
\ref{tab:KmEoxEo-multiplicity}. For these 13 types
the frame lattice $[W]$ is obtained as the orthogonal complement of 
a primitive embedding of $T_0 \cong D_4 \oplus A_2$ into $L^{(I)}$, as 
in section \ref{sec:J2-classification}. In this subsection, we study
cases where $D_4$ is embedded
into a $D_n$ component of the irreducible decomposition of 
$(L^{(I)})_{\rm root} \cong \oplus_a R_a$, and $A_2 \subset T_0$ into 
another irreducible component. Some definitions of octads (codewords)
used in this section can be found in the appendices \ref{section:octads} 
and \ref{sssec:NS-for-KmEwEw}.

The types [no. 29] and [no. 30] use an embedding into Niemeier lattices $L^{(I)}$ 
whose root lattices contain $D_4$ as an irreducible component.
Thus, we need to find the centre of a deep hole of $\Lambda_{24} \otimes \R$ 
that is surrounded by $\vec{v}_{1,2,3,4}$ and one more vector $\vec{v}'_4$. 
From the conditions that $(\vec{v}'_4-\vec{v}_{1,3,4})^2 = -4$, 
$(\vec{v}'_4-\vec{v}_2)^2 = -6$ and $(\vec{v}'_4-\vec{v}_{p,q})^2 = -4$, 
we find that this vector should be of either one of the forms 
\begin{eqnarray}
 \vec{v}'_4 & = & \nu_{\Omega} - 2 \nu_{K_U} + 4 \nu_{\{ \infty, 0,i\} }, \qquad 
    i \in K_U \backslash \{ \infty,0,1,2\}, 
    \label{eq:vec-Leech-KmEoxEo-1st-D4inD4-a} \\
 \vec{v}'_4 & = & \nu_\Omega - 2 \nu_{K} + 4 \nu_{\{\infty, 0,1\}}, \qquad 
     K = K_{\epsilon\heartsuit, \epsilon \clubsuit,
            \epsilon \spadesuit, \epsilon \diamondsuit }, 
   \label{eq:vec-Leech-KmEoxEo-1st-D4inD4-b} \\
 \vec{v}'_4 & = & \nu_\Omega -2\nu_K + 4 \nu_{\{ \infty, 0, i\}}, \qquad 
     \{\infty, 0,2\} \subset K \cap K_U, |K \cap K_U|=4, 1 \nin K, i \nin K_U.
   \label{eq:vec-Leech-KmEoxEo-1st-D4inD4-c} 
\end{eqnarray}

Let us take $\vec{v}'_4$ of the form (\ref{eq:vec-Leech-KmEoxEo-1st-D4inD4-a})
with $i = 11$ for now. Then there are 25 vectors $\vec{u} \in \Lambda_{24}$ 
satisfying $(\vec{u}-\vec{v}_{1,2,3,4})^2 = -4$ and 
$(\vec{u}-\vec{v}'_4)^2 = -4$. It turns out that the diagram of these 
25 vectors\footnote{
The 25 vectors $\vec{u}\in \Lambda_{24}$ form five groups of 5 vectors, 
each one of which gives rise to an extended Dynkin diagram of $D_4$.
The five groups of vectors consist of $P$, $Q_{k=2,18,8,12}$; 
$2\nu_K$ with $K \in 
\{K_{\alpha \heartsuit}, K_p, K_{q1}, K_{q\omega}, K_{q\bar{\omega}} \}$; 
$2\nu_K$ with $K \in 
\{K_{\alpha \clubsuit}, K_{r1}, K_{\bar{r}1}, K_{s1}, K_{\bar{s}1} \}$;
$2\nu_{K}$ with $K \in 
\{K_{\alpha \spadesuit}, K_{r\omega}, K_{\bar{r}\omega}, K_{s\omega},
K_{\bar{s}\omega} \}$; and finally, 
$2\nu_K$ with $K \in \{K_{\alpha \diamondsuit}, 
K_{r\bar{\omega}}, K_{\bar{r}\bar{\omega}}, K_{s\bar{\omega}}, 
K_{\bar{s}\bar{\omega}} \}$.}
drawn under the rule explained in the appendix \ref{ssec:Leech-review} 
is the collection of five extended Dynkin diagrams of $D_4$. 
Hence a deep hole of type $\tau$ is surrounded by the vectors 
(see Table \ref{tab:Niemeier-root}), and the choice of 
(\ref{eq:D4-A2-embed-II-1-25-1}--\ref{eq:D4-A2-embed-II-1-25-q}) 
and $\vec{v}^{(\tau)a}_i$'s corresponds to an embedding of 
$D_4 \subset D^{(1)}_4$ and $A_2 \subset D^{(2)}_4$ into 
$L^{(\tau)}_{\rm root} = \oplus_{a=1}^6 D_4^{(a)}$. This is the type
[no. 29] in Table 1.3 of \cite{Nish01}. 
% (Table 1.3 of \cite{Nish01}). 
The centre of this deep hole is 
\begin{equation}
 \vec{c}^{(I=\tau)} = \frac{1}{h^{(\tau)}} \left( 4\nu_{K_{\alpha \heartsuit}}+
  2\nu_{K_p}+2\nu_{K_{q1}}+2\nu_{K_{q\omega}}+2\nu_{K_{q\bar{\omega}}} \right) = 
\frac{1}{h^{(\tau)}}  \begin{tabular}{|cc|cc|cc|}
 \hline
  12 & 12 & 4 & 4 & 0 & 4 \\
   4 &  4 & 4 & 4 & 4 & 0 \\
   4 &  4 & 4 & 4 & 4 & 0 \\
   4 &  4 & 4 & 4 & 4 & 0 \\
  \hline    
   \end{tabular}.
\end{equation}
The elliptic divisor is given by $u^{(\tau)}_1 \in {\rm II}_{1,25}$ 
in (\ref{eq:f1-as-elliptic-div}), with $h^{(I=\tau)}=6$. 
Using the 4 groups of $\vec{v}^{(\tau)a}_{i=0,1,2,3,4}$ forming an extended
Dynkin diagram of $D_4$ (excluding those containing 
(\ref{eq:D4-A2-embed-II-1-25-1}--\ref{eq:D4-A2-embed-II-1-25-q})), 
we obtain 
\begin{eqnarray}
 u^{(\tau)}_1 & = & 2G_{11} + E_1 + F_1 + C_1 + D_1, \\
              & = & 2G_{24} + E_4 + F_2 + C_3 + D_3, \\
              & = & 2G_{43} + E_3 + F_4 + C_2 + D_2, \\
              & = & 2G_{32} + E_2 + F_3 + C_4 + D_4; 
\end{eqnarray} 
they describe 4 distinct singular fibres of $I_0^*$ type, and are 
algebraically equivalent. One of these is used in Table \ref{tab:KmEoxEo-multiplicity}. 
Choosing $\vec{v}'_4$ of the form \eqref{eq:vec-Leech-KmEoxEo-1st-D4inD4-b} also
implies an embedding of $A_4 \oplus A_2$ into $L^{(\tau)}$, so we do not discuss this
choice further.

Let us now take $\vec{v}'_4$ of the 
form (\ref{eq:vec-Leech-KmEoxEo-1st-D4inD4-c}) with 
\begin{equation}
  K = 
\begin{array}{|cc|cc|cc|}
\hline
   *&*&{\color{white} *} & {\color{white} *} &*&* \\
    & & & & {\color{white} *}&  \\
    & & & & {\color{white} *}&  \\
   *&*& & &*&* \\
\hline 
\end{array} \; 
\end{equation}
instead. Then there are 24 vectors $\vec{u} \in \Lambda_{24}$ satisfying 
$(\vec{u}-\vec{v}_{1,2,3,4})^2 = -4$ and $(\vec{u}-\vec{v}'_4)^2 = -4$.
The diagram of these 24 vectors turns out to be the collection of 4 extended 
Dynkin diagrams of $A_5$. The vectors $P$ and $Q_2$, associated with the simple roots 
of $A_2 \subset T_0$, are part of one of the 4 $A_5$'s. This is a deep hole of type $\sigma$ (see Table
\ref{tab:Niemeier-root}), and correspond to the 
type [no.30: $D_4 \subset D_4$, $A_2 \subset A_5$; $A_5^{\oplus 3}$] in
Nishiyama's ${\cal J}_2(X)$ classification. The centre of this deep hole 
is located at
\begin{equation}
 h^{(I=\sigma)} \vec{c}^{(\sigma)} = 
   \left(2\nu_{K_{\alpha \heartsuit}}+2\nu_{K_{q1}}
 + 2\nu_{K_{\beta \heartsuit}} + 2\nu_{K_{\bar{r}\bar{\omega}}}
 + 2\nu_{K_{\gamma \heartsuit}} + 2\nu_{K_{q\omega}} \right), 
\end{equation}
and the elliptic divisor in $\overline{D'} \subset S_X \subset {\rm II}_{1,25}$ 
can be obtained by using Lemma I:
\begin{eqnarray}
 u_1^{(\sigma)}
 & = & G_{1 1} + C_{1} + G_{3 4} + F_{3} + G_{3 3} + D_{1},  \qquad (A_5, I_6) \\
 & = & G_{24} + D_3 + G_{13} + C_4 + G_{21} + F_2, \qquad (A_5, I_6) \\
 & = & G_{43} + F_4 + G_{41} + D_4 + G_{14} + C_2. \qquad (A_5, I_6),  
\end{eqnarray}
which describe three singular fibres of $I_6$ type.

Let us now turn to the types [no. 27] and [no. 28] in ${\cal J}_2(X)$, 
where the Niemeier lattice $L^{(\pi)}$ is used, and $D_4 \subset T_0$ is 
embedded into an irreducible component $D_5$ in the root lattice of 
$L^{(\pi)}$. Thus, we look for $\vec{v}_5$, and choose 
\begin{equation}
 \vec{v}_5 = \nu_\Omega - 2\nu_{K_L} +4\nu_0 =: X_L, \qquad 
   \{0,2\} \subset K_L, \; \infty, 1 \nin K_L,  \;  |K_L \cap K_U|=4\, .
\end{equation}
For $K_L \in {\cal C}_{24}(8)$, we use the one in (\ref{eq:Kl-Kr}).
The centre of deep hole of $I=\pi$ type, 
$\vec{c}^{(\pi)} \in \Lambda_{24} \otimes \R$, should be surrounded by 
$\vec{v}_{1,\cdots,5}$ and one more vector $\vec{v}'_5 \in \Lambda_{24}$.
The vector $\vec{v}'_5$ should be of the same form as $\vec{v}_5$ 
with $K_L \in {\cal C}_{24}(8)$ for $\vec{v}_5$ replaced by another 
$K \in {\cal C}_{24}(8)$ satisfying the same conditions as $K_L$ plus 
one more condition, $|K \cap K_L| = 4$. 

For two different choices of this $K \in {\cal C}_{24}(8)$ for $\vec{v}'_5$,
say, 
\begin{equation}
 K = 
 \begin{array}{|cc|cc|cc|}
   \hline
     *& & &*&*&  \\
      & & & &*&* \\
      & & & & &  \\
      &*&*& &*&  \\ 
   \hline
 \end{array}\; ,  \qquad 
 \begin{array}{|cc|cc|cc|}
   \hline
     *& &{\color{white} *} &{\color{white} *} &*&  \\
      &*& & & &* \\
      &*& & & &* \\
     *& & & &*&  \\ 
   \hline
 \end{array} \; ,
\end{equation}
the vectors $\vec{u} \in \Lambda_{24}$ satisfying 
$(\vec{u}-\vec{v}_{1,\cdots,5})^2 = -4$ and 
$(\vec{u}-\vec{v}'_5)^2 = -4$ form the extended Dynkin diagram of 
$A_7+A_7+D_5$. 
The two vectors associated with the simple roots of $A_2 \subset T_0$, 
however, are contained in $D_5$ when the first $K \in {\cal C}_{24}(8)$ is 
used for $\vec{v}'_5$, and they are in one of $A_7$ when the second 
$K \in {\cal C}_{24}(8)$ is used. Therefore, these two cases correspond to 
the type [no.28] and [no.27] in Nishiyama's classification,
respectively. 
In type [no.28: $D_4 \subset D_5$, $A_2 \subset D_5$; $A_7^{\oplus 2}$], 
the centre of the deep hole is given by 
\begin{equation}
 h^{(\pi)} \vec{c}^{(\pi)} = 2\nu_{K_{\gamma \heartsuit}} + 2\nu_{K_{q\omega}} +
    2 \nu_{K_{\delta \clubsuit}} + 2\nu_{K_{\bar{s}\bar{\omega}}} + 
    2 \nu_{K_{\gamma \diamondsuit}} + 2\nu_{K_p} + 
    2 \nu_{K_{\delta \spadesuit}} + 2\nu_{K_{\bar{r}\bar{\omega}}},
\end{equation}
and elliptic divisors are given by 
\begin{eqnarray}
u_1^{(\pi)} 
 & = & G_{21} + E_1 + G_{31} + F_3 + G_{33} + D_1 + G_{44} + C_4 , 
     \qquad (A_7, I_8), \\
 & = & G_{14} + F_1 + G_{12} + E_2 + G_{42} + C_1 + G_{23} + D_4. 
     \qquad (A_7, I_8). 
\end{eqnarray}
For type [no.27: $D_4 \subset D_5$, $A_2 \subset A_7$; $A_7 \oplus D_5$], 
\begin{eqnarray}
  h^{(\pi)} \vec{c}^{(\pi)} & = &
   2\nu_{K_{\alpha \clubsuit}} + 2\nu_{K_{\bar{s}1}}
 + 2\nu_{K_{\delta \spadesuit}} + 2\nu_{K_{\bar{r}\bar{\omega}}} 
 + 2\nu_{K_{\gamma \heartsuit}} + 2\nu_{K_{q\omega}} 
 + 2\nu_{K_{\delta \clubsuit}} + 2\nu_{K_{r1}}. \\
  u_1^{(\pi)} 
 & = & G_{24} + D_3 + G_{31} + F_3 + G_{33} + D_1 + G_{44} + E_4, 
   \qquad (A_7, I_8) \\
 & = & G_{43} + G_{42} + 2(D_2 + G_{12}) + F_1 + E_2. 
   \qquad (D_5, I_1^*)
\end{eqnarray}

\begin{figure}[tbp]
\begin{center}
   \scalebox{.5}{ 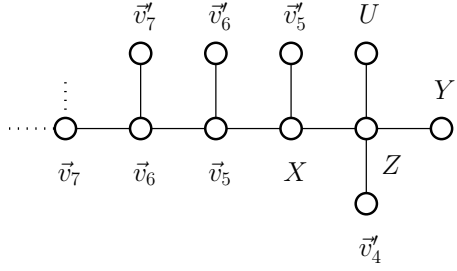 }
\caption{A possible embedding of $D_4 \subset T_0$ into $D_n, n\geq 6$. 
The vectors $\vec{v}_i$ are those given in \eqref{basisD4Dn}. 
The extended Dynkin diagram of $D_n$ is formed by the skeleton made of
$U$, $Y$, $Z$, $X$ and $\vec{v}_{5, \cdots, n}$ and one more 
$\vec{v}'_n$ or $\vec{v}''_n$. \label{dnchain}}
\end{center} 
\end{figure}

There are 9 more types in ${\cal J}_2(X)$ for 
$X = {\rm Km}(E_\omega \times E_\omega)$ by Nishiyama where 
the embedding $\phi'_{T_0}: T_0 \hookrightarrow L^{(I)}_{\rm root} \cong 
\oplus_{a \in A} R_a$ is given by sending $D_4 \subset T_0$ into 
$R_a = D_n$ and $A_2$ into another irreducible component $R_{a'}$.
The first 9 entries of Table \ref{tab:KmEoxEo-multiplicity} are 
in this category. In order to determine the elliptic divisors 
for those types, the required task is to find a deep hole of 
type $I$ so that all the vectors 
(\ref{eq:D4-A2-embed-II-1-25-1}--\ref{eq:D4-A2-embed-II-1-25-q}) are 
among the surrounding vectors. This task can be carried out systematically, 
just like such a systematic approach was possible for a series of 
$L^{(I)}$ containing $D_n$ sublattice in Chapt. 23 of \cite{CS-book}.
Defining vectors in $\Lambda_{24}$ as follows,  
\begin{eqnarray}
 \vec{v}_6 & = & 2\nu_{K_{\bar{r}1}}, \qquad 
     \vec{v}'_6 = 2\nu_{K_{\alpha \heartsuit}} \; [{\rm no.}\,26], \quad 
     \vec{v}''_6 = 2\nu_{K_{s1}} \; [{\rm no.}\,25],\nn \\
 \vec{v}_7 & = & 2\nu_{K_{\delta\heartsuit}}, \qquad 
     \vec{v}'_7 = 2\nu_{K_{\gamma \heartsuit}} \; [{\rm no.}\,22], \quad 
     \vec{v}''_7 = 2\nu_{K_{\alpha\clubsuit}} \; [{\rm no.}\,20],\nn \\
 \vec{v}_8 & = & 2\nu_{K_{r\bar{\omega}}}, \qquad 
     \vec{v}'_8 = 2\nu_{K_{q\bar{\omega}}} \; [{\rm no.}\,16], \nn\\
 \vec{v}_9 & = & 2\nu_{K_{\gamma \clubsuit}}, \qquad 
     \vec{v}'_9 = 2\nu_{K_{\beta \diamondsuit}} \; [{\rm no.}\,17], \nn\\
 \vec{v}_{10} & = & 2\nu_{K_{\bar{s}1}}, \qquad 
     \vec{v}'_{10} = 2\nu_{K_{q1}} \; [{\rm no.}\,8],\nn \\
 \vec{v}_{11} & = & 2\nu_{K_{\delta\spadesuit}},\nn \\
 \vec{v}_{12} & = & 2\nu_{K_{p}}, \qquad 
     \vec{v}'_{12} = 2\nu_{K_{\bar{r}\bar{\omega}}} \; [{\rm no.}\,14],  
   % ({\rm altn.~notation~} \vec{v}_0^{(\epsilon)D_{12}}?)
 \nn \\
 \vec{v}_{13} & = & 2\nu_{K_{\gamma\diamondsuit}},\nn \\
 \vec{v}_{14} & = & 2\nu_{K_{\bar{s}\bar{\omega}}},\nn \\
 \vec{v}_{15} & = & 2\nu_{K_{\delta \clubsuit}},\nn \\
 \vec{v}_{16} & = & 2\nu_{K_{q\omega}}, \qquad 
 \vec{v}'_{16} = 2\nu_{K_{r1}} \; [{\rm no.}4],  
 \label{basisD4Dn}
\end{eqnarray}
we see that the vectors $\vec{v}_{1, \cdots, n}$ along with one more vector 
$\vec{v}'_n$ (or $\vec{v}''_n$) form an extended Dynkin diagram of $D_n$
for $n=6,7,8,9,10,12,16$, where those vectors are arranged in the diagram 
as in Figure \ref{dnchain}. 
%%%%%%%%%%%%%%%%%%%%%%%%%%%%%%%%%%%%%%%%%%%%%%%%%%%%%%%%%%%%%%%%%%%
%%%%%%%%%%%%%%%%%%%%%%%%%%%%%%%%%%%%%%%%%%%%%%%%%%%%%%%%%%%%%%%%%%%%
By working out other irreducible components of the diagram, we see that 
they indeed correspond to an irreducible component of vectors surrounding 
a centre of deep hole of some type. The corresponding type identification 
numbers in the ${\cal J}_2(X)$ classification are already given above, 
and the type $I$ of Niemeier lattices/deep holes are found in 
Table \ref{tab:KmEoxEo-multiplicity}. Elliptic divisors are determined 
for all those 9 types, as in the previous 4 types, and the results 
are shown in Table \ref{tab:KmEoxEo-multiplicity}.

For types [no.4], [no.14] and [no.17], where $T_0$ is embedded into 
$L^{(I)}_{\rm root} \cong D_{16}\oplus E_8$, $D_{12}\oplus D_{12}$ and 
$D_9\oplus A_{15}$, respectively, there are only two irreducible components 
of the root lattice $L^{(I)}_{\rm root}$. The elliptic divisors obtained in these
cases inevitably have expressions containing $\lambda_p$, $\lambda_{q2}$. Of course
$S_X$ is orthogonal to $T_0$ by construction, so this does not imply that these
are fibre components. 
If we choose a basis of $T_0$ and $S_X$ without redundancy, the elliptic divisors can 
be solely written in terms of the basis of $S_X$.

%%%%%%%%%%%%%%%%%%%%%%%%%%%%%%%%%%%%%%%%%%%%%%%%%%%%%%%%%%%%%%%%%%%%%
\subsubsection{Determination of Elliptic Divisors: $D_4 \hookrightarrow E_n\, ,\, A_2\hookrightarrow {\rm other}$} 
% (type 2,5,9--11, 19, 23, 24)
\label{sssec:D4toEn-A2toOther}
%%%%%%%%%%%%%%%%%%%%%%%%%%%%%%%%%%%%%%%%%%%%%%%%%%%%%%%%%%%%%%%%%%%%%%%

%%%%%%%%%%%%%%%%%%%%%%%%%%%%%%%%%%%%%%%%%%%%%%%%%%%%%%%%%%%%%%%%%%%%
\begin{figure}[!h]
\begin{center}
    \scalebox{.5}{ 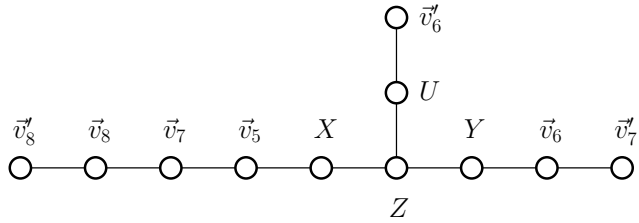 }
\caption{A possible embedding of $D_4 \subset T_0$ into 
$E_n$ \label{enchain}. The vectors in the diagram are defined 
in \eqref{v5v6d4en} and below. The extended Dynkin diagram of 
a) $E_6$, b) $E_7$ and c) $E_8$ are formed by adding a) $\vec{v}'_6$, 
b) $\{\vec{v}_7, \vec{v}'_7 \}$  and c) $\{ \vec{v}_{7,8}, \vec{v}'_8
 \}$, respectively, on top of a common subdiagram that consists of the 
six vectors $\{U, X, Y, Z\}$ and $\vec{v}_{5,6}$ in (\ref{v5v6d4en}).
}
\end{center} 
\end{figure}
%%%%%%%%%%%%%%%%%%%%%%%%%%%%%%%%%%%%%%%%%%%%%%%%%%%%%%%%%%%%%%%%%%%%%%

Let us now move on to the second group in Table \ref{tab:KmEoxEo-multiplicity}
consisting of 8 types (no.$2, 5, \cdots, 24$) in ${\cal J}_2(X)$ 
for $X = {\rm Km}(E_\omega \times E_\omega)$. The frame lattice of the  
types in this group is obtained by embedding $D_4 \subset T_0$ into an 
irreducible component $E_n \subset (L^{(I)})_{\rm root}$, and 
$A_2 \subset T_0$ into another irreducible component of 
$(L^{(I)})_{\rm root}$. For all of these types, we use 
\begin{eqnarray}
 \vec{v}_5 & = & \nu_\Omega - 2\nu_{K_L} + 4 \nu_0 =: X_L,\nn \\
 \vec{v}_6 & = & \nu_\Omega - 2\nu_{K_R} + 4\nu_\infty =: Y_R, \label{v5v6d4en}
\end{eqnarray}
with $K_L$ and $K_R$ given in (\ref{eq:Kl-Kr}).

The types [no. 19], [no. 23] and [no. 24] in Table \ref{tab:KmEoxEo-multiplicity} 
involve an embedding of $D_4 \subset T_0$ into an irreducible component 
$E_6 \subset (L^{(I)})_{\rm root}$. In order to find elliptic divisors 
for these types, we look for one more remaining vector 
$\vec{v}'_6 \in \Lambda_{24}$
which surrounds the centre of a deep hole of the corresponding type along
with 6 other vectors $\{X_L, X, Z, Y, Y_R, U\}$. This remaining 
vector should be of the form $\vec{v}'_6 = 2\nu_K$ for 
$K \in {\cal C}_{24}(8)$. If we take
\begin{equation}
 K= 
\begin{array}{|cc|cc|cc|}
 \hline
  *&*& & & &* \\
   & &*&*&*&  \\
   & & & & &* \\
   & & & & &* \\
 \hline
\end{array} \; , \quad 
 \begin{array}{|cc|cc|cc|}
  \hline
   *&*& &*& &  \\
   *& & & &*&  \\
   *& & & & &* \\
    & &*& & &  \\ 
  \hline
 \end{array} \; , \quad 
 \begin{array}{|cc|cc|cc|}
  \hline
    *&*&{\color{white} *} &  {\color{white} *}& &  \\
    *&*& & & &  \\
     & & & &*&* \\
     & & & &*&* \\
  \hline
 \end{array} \; , 
\end{equation}
it turns out that $\vec{v}_{1,\cdots,6}$ and $\vec{v}'_6$ indeed 
form an extended Dynkin diagram of $E_6$, as shown in Figure \ref{enchain}. 
Examining the rest of the diagram of vectors, we see that the three
different choices of the codeword $K \in {\cal C}_{24}(8)$ lead to 
three different embedding, [no.19: $D_4 \subset E_6^{(1)}$, 
$A_2 \subset E_6^{(2)}$; $E_6^{\oplus 2}$], 
[no.23: $D_4 \subset E_6$, $A_2 \subset A_{11}$; $D_7$] and 
[no.24: $D_4 \subset E_6$, $A_2 \subset D_7$; $A_{11}$], respectively.
Lemma I is used to determine elliptic divisors of those three types, and
the results are recorded in Table \ref{tab:KmEoxEo-multiplicity}.

For the types [no.9], [no.10] and [no.11], where an embedding 
$D_4 \subset E_7$ is used, we take 
$\vec{v}_7 = 2\nu_{K_{\alpha \diamondsuit}}$. 
The centres of deep holes with an appropriate embedding of 
the vectors $\{U, X, Y, Z, \vec{v}_{5,6,7}\}$ and $\vec{v}_{p,q2}$
for these three types are found when we take one more vector
$\vec{v}'_7$ to complete the extended Dynkin diagram of $E_7$ as   
\begin{equation}
  \vec{v}'_7 = 2\nu_{K_{\alpha\spadesuit}}, \qquad 
  2\nu_{K_{\beta\diamondsuit}}, \qquad 
  2\nu_{K_{\bar{s}1}},
\end{equation}
respectively.

The types [no. 2] and [no. 5] are associated with the embedding of 
$D_4$ into $E_8$. Along with $\{ U, X, Y, Z, \vec{v}_{5,6,7}\}$ 
we can take $\vec{v}_8 = 2\nu_{K_{r\bar{\omega}}}$ for both types, and
complete the extended Dynkin diagram by adding
\begin{equation}
 \vec{v}'_8 = 2\nu_{K_{\gamma\clubsuit}}, \qquad 
   2\nu_{K_{\delta\heartsuit}}\, ,
\end{equation}
respectively. Working out the other irreducible components of the
diagram of the vectors around the deep hole, we see that these two 
types correspond to the embedding 
[no.2: $D_4 \subset E_8^{(1)}$, $A_2 \subset E_8^{(2)}$; $E_8^{(3)}$] and 
[no.5: $D_4 \subset E_8$, $A_2 \subset D_{16}$], as 
desired. 

%%%%%%%%%%%%%%%%%%%%%%%%%%%%%%%%%%%%%%%%%%%%%%%%%%%%%%%%%%%%%%%%%%%%
\subsubsection{Determination of Elliptic Divisors: 
$D_4 \oplus A_2 \hookrightarrow D_n\, , n\geq 7$}
\label{sssec:D4A2toDn}
%%%%%%%%%%%%%%%%%%%%%%%%%%%%%%%%%%%%%%%%%%%%%%%%%%%%%%%%%%%%%%%%%%%%%

\begin{figure}[!h]
\begin{center}
   \scalebox{.5}{ 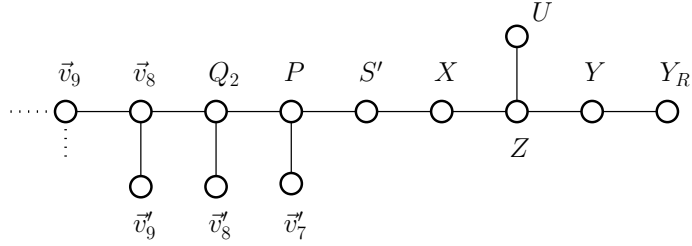 }
\caption{A possible embedding of $D_4\oplus A_2 \subset T_0$ into 
$D_n,n\geq 7$ or $E_8$ \label{dnchaindn7e8}. The vectors appearing
in the figure are defined in \eqref{d4a2dnvecs} and \eqref{v5v6d4en}. 
The extended Dynkin diagram is composed of $\{U, Y, Z, X, S', P, Q_2, 
\vec{v}_{8, \cdots, n}\}$ and one more vector $\vec{v}'_n \in
 \Lambda_{24}$. The extended Dynkin diagram of $E_8$, on the other hand,
 is made of $\{U, Y, Z, X, S', P, Q_2, Y_R\}$ and $\vec{v}'_8$.}
\end{center} 
\end{figure}

We now discuss the third group in Table \ref{tab:KmEoxEo-multiplicity}, 
type no.$12, 3, \cdots, 21$ of ${\cal J}_2(X)$, for which both $D_4$ and
$A_2$ are embedded in a $D_n$ component of $(L^{(I)})_{\rm root}$. 
We start with the octad
\begin{equation}
 K_S = \begin{array}{|cc|cc|cc|}
  \hline
    *& &*&{\color{white} *} &{\color{white} *} &  \\
    *& &*& & &  \\
    *& &*& & &  \\
    *& &*& & &  \\
  \hline
\end{array} \; 
\label{eq:Ks}
\end{equation}
which allows us to define $\vec{v} = 2\nu_{K_S} =: S'$, sitting 
in between $X$ and $P$, see Figure \ref{dnchaindn7e8}.
%%%%%%%%%%%%%%%%%%%%%%%%%%%%%%%%%%%%%%%%%%%%%%%%%%%%%%
% \vspace{5mm}
%%%%%%%%%%%%%%%%%%%%%%%%%%%%%%%%%%%%%%%%%%%%%%%%%%%%%%
The 4 vectors $\{U, Y, Z, X\}$ for $D_4 \subset T_0$, 2 vectors 
$\{P, Q_2\}$ for $A_2 \subset T_0$ and $S'$ already form 
the Dynkin diagram of $D_7$. These 7 vectors form the stem for the 7 types 
of ${\cal J}_2(X)$ studied in section \ref{sssec:D4A2toDn} as well as 
for the 2 types to be covered in section \ref{sssec:D4A2toE8}.

\begin{figure}[!h]
\begin{center}
   \scalebox{.5}{ 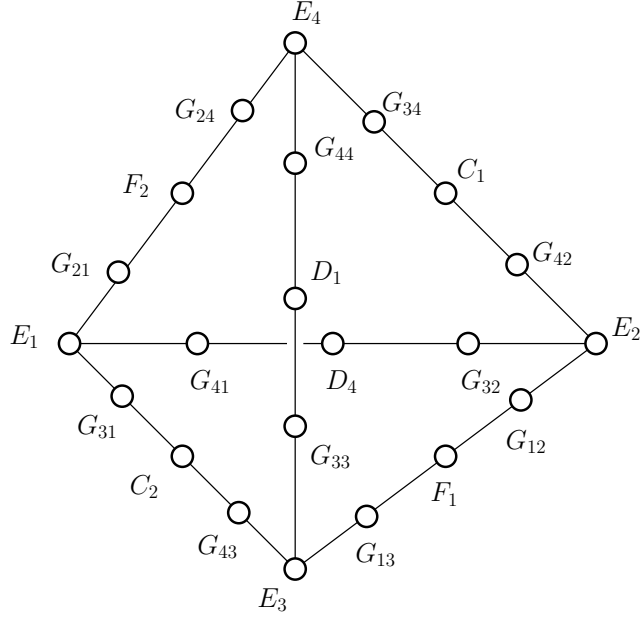 }
\caption{\label{fig:Ew-Coxeter} Part of the Coxeter diagram of $S_X$ for 
$X={\rm Km}(E_\omega\times E_\omega)$. 
The 22 vertices are labelled using the identifications \eqref{gijvsK} and \eqref{efcdvsK}. The displayed subset 
is tailored towards the study of embeddings $D_4\oplus A_2 \hookrightarrow D_n$ or $E_n$. As explained in the main
text, this diagram allows to read off the vectors of the remaining components around deep holes for these cases. }
\end{center} 
\end{figure}

The type [no.21] should be associated with the embedding 
[$(D_4 \oplus A_2) \subset D_7$; $E_6 \oplus A_{11}$]. 
For the other vector $\vec{v}'_7$ forming an extended Dynkin diagram 
of $D_7$, we choose $\vec{v}'_7 = \nu_\Omega - 4\nu_k \equiv Q_k$ for 
$k \in \{ 11,18,8,12\}$ (e.g., 11), see Figure \ref{dnchaindn7e8}. 
The other irreducible components of $L^{(I=\lambda)}_{\rm root}$, namely 
$E_6 \oplus A_{11}$ should be made of vectors $\vec{u} \in \Lambda_{24}$
that are at norm-$(-4)$ distance from all the seven vectors 
$U, Y, Z, X, S', P, Q_2$---(***), as well as from $\vec{v}'_7$; here, 
the condition $(\vec{u}-U)^2 = (\vec{u}- Y)^2 =  \cdots = -4$ is
implied. There are 22 vectors satisfying (***), they are shown in 
Figure~\ref{fig:Ew-Coxeter}. Note that the labels on the vertices 
are not $\vec{u} \in \Lambda_{24}$ but their corresponding Leech roots 
$\lambda = (\vec{u}, 1, -1 - (\vec{u})^2/2) \in S_X \subset \Lambda_{24}
\oplus U^{(\omega)}$. Among those 22 vectors, three vectors
corresponding to $G_{43}$, $G_{24}$ and $G_{32}$ are not at the
norm-$(-4)$ distance from $\vec{v}'_7$. Eliminating those three
vertices from the diagram in Figure \ref{fig:Ew-Coxeter}, we see that 
the extended Dynkin diagrams of $E_6$ and $A_{11}$ are left indeed. 
This is how we see that this choice of $\vec{v}^{(I)a}_i$'s is for 
a deep hole of type $I=\lambda$, and for the embedding [no.21: 
$(D_4 \oplus A_2) \subset D_7$; $E_7 \oplus A_{11}$].
Elliptic divisors can be computed by using Lemma I, as before.
%%%%%%%%%%%%%%%%%%%%%%%%%%%%%%%%%%%%%%%%%%%%%%%%%%%%%%%%%%%%%%%%
%%%%%%%%%%%%%%%%%%%%%%%%%%%%%%%%%%%%%%%%%%%%%%%%%%%%%%%%%%%%%%%%%%%

Similarly to the study in sections \ref{sssec:D4toDn-A2toOther} and
\ref{sssec:D4toEn-A2toOther}, and to Chapt. 23 of \cite{CS-book}, 
we can proceed systematically after this. 
The following vectors are used to complete the extended Dynkin diagram 
of $D_n$ to which $T_0 = D_4 \oplus A_2$ is embedded:
\begin{eqnarray}
 \vec{v}_8 & = & 2\nu_{K_{\epsilon \clubsuit}}, \qquad 
    \vec{v}'_8 = 2\nu_{K_{\epsilon \diamondsuit}} \; [{\rm ~no.}\,15],\nn \\ 
 \vec{v}_9 & = & 2\nu_{K_{q\bar{\omega}}}, \qquad 
    \vec{v}'_9 = 2\nu_{K_{s1}} \; [{\rm ~no.}18], \nn \\ 
 \vec{v}_{10} & = & 2\nu_{K_{\beta \clubsuit}}, \qquad 
    \vec{v}'_{10} = 2\nu_{K_{\delta \heartsuit}} \; [{\rm ~no.}\,7], \nn \\
 \vec{v}_{11} & = & 2\nu_{K_{r\omega}}, \nn \\
 \vec{v}_{12} & = & 2\nu_{K_{\gamma \heartsuit}}, \qquad 
    \vec{v}'_{12} = 2\nu_{K_{\alpha \spadesuit}} \; [{\rm ~no.}\,13], \nn \\
 \vec{v}_{13} & = & 2\nu_{K_{q\omega}}, \nn \\
 \vec{v}_{14} & = & 2\nu_{K_{\delta \clubsuit}},\nn  \\
 \vec{v}_{15} & = & 2\nu_{K_{r1}}, \nn \\
 \vec{v}_{16} & = & 2\nu_{K_{\beta \heartsuit}}, \qquad 
    \vec{v}'_{16} = 2\nu_{K_{\alpha \clubsuit}} \; [{\rm ~no.}\,3],\nn  \\
 \vec{v}_{17} & = & 2\nu_{K_{q1}}, \qquad
 \vec{v}_{18} = 2\nu_{K_{\gamma \clubsuit}}, \qquad 
 \vec{v}_{19} = 2\nu_{K_{r\bar{\omega}}}, \qquad 
 \vec{v}_{20} = 2\nu_{K_{\alpha \diamondsuit}}, \nn \\
&&  \qquad \qquad \quad \vec{v}_{21} = 2\nu_{K_{s\bar{\omega}}}, \qquad 
 \vec{v}_{22} = 2\nu_{K_{\beta \spadesuit}}, \qquad 
 \vec{v}_{23} = 2\nu_{K_{p}}, \nn \\
 \vec{v}_{24} & = & 2\nu_{K_{\gamma \diamondsuit}}, \qquad 
    \vec{v}'_{24} = 2\nu_{K_{\delta \spadesuit}} \; [{\rm ~no.}\,12], \label{d4a2dnvecs}
\end{eqnarray} 
With this data and Lemma I, one can compute the elliptic divisors for
all of these types. An explicit identification with the Leech roots 
in $S_X \subset \Lambda_{24} \oplus U^{(\omega)}$ and curves 
$\{G_{ij}, E_i, F_i, C_i, D_i \}$ (see appendix
\ref{ssec:review-Kummer}) allows us to write those divisors 
in the form presented in Table \ref{tab:KmEoxEo-multiplicity}.

Before moving on to section \ref{sssec:D4A2toE8}, let us explain how 
Figure \ref{fig:Ew-Coxeter} is used to read off the type of deep holes
and embedding of $T_0$ into Niemeier lattices, as well as to determine 
the vectors $\vec{v}_n$'s systematically for large $n$. 
In the following, we frequently use the identification among 
the vectors in $\Lambda_{24}$ listed above, 
their corresponding Leech roots, and elements in $S_X$.

As for type [no.15: $T_0 \subset D_8$; $D_8 \oplus D_8$], vectors
corresponding to the two vertices $F_1$ and $F_2$ in 
Figure \ref{fig:Ew-Coxeter} are not at norm-$(-4)$ distance from 
$\vec{v}_8$, and the vectors corresponding to $D_1$ and $D_4$ in the
figure are not from $\vec{v}'_8$. Removing these four vertices from Figure
\ref{fig:Ew-Coxeter}, the extended Dynkin diagram of $D_8$ appears twice.

For the type [no.18: $T_0 \subset D_9$; $A_{15}$], 
we keep the vertices $D_1$ and $D_4$ (because $\vec{v}'_8$ is not in
$D_9 \subset L^{(I)}_{\rm root}$), but still excise $F_2$ and $F_1$. 
Furthermore, since $\vec{v}_9=F_1$ and $\vec{v}_9'=F_2$, we also have 
to delete the nearest neighbours (not at norm-$(-4)$ distance) of $F_2$
and $F_1$, namely the four vertices $G_{21}$, $G_{24}$, $G_{12}$ and $G_{13}$. 
The rest of the diagram in Figure \ref{fig:Ew-Coxeter} is in the shape
of the extended Dynkin diagram of $A_{15}\subset (L^{(\theta)})_{\rm root}$.
For the type [no.7: $T_0 \subset D_{10}$; $E_7 \oplus E_7$], 
we have $\vec{v}_{10}=G_{13}$ and $\vec{v}_{10}'=G_{12}$ completing the
extended Dynkin diagram of the vectors $D_{10}$.
% $T_0$ to $D_{10}\subset(L^{(\eta)})_{\rm root}$. 
Hence we have to remove $F_2$, $F_1$, $G_{12}$, $G_{13}$ and two more
vertices $E_2$ and $E_3$ (because they are not at norm-$(-4)$ distance 
from $\vec{v}_{10}$ and $\vec{v}'_{10}$, respectively) from 
Figure \ref{fig:Ew-Coxeter}. 
As expected from the list of Niemeier lattices 
(the $D_{10}$ irreducible component is found only in $L^{(\eta)}$; 
see table \ref{tab:Niemeier-root}), 
this reproduces twice the extended Dynkin diagram of $E_7$. 
Using \eqref{d4a2dnvecs}, the remaining cases are treated in the same fashion.

Diagrams similar to Figure~\ref{fig:Ew-Coxeter} can be used to determine
the vectors (\ref{basisD4Dn}) and types of embeddings systematically 
in sections \ref{sssec:D4toDn-A2toOther} and \ref{sssec:D4toEn-A2toOther}.

%%%%%%%%%%%%%%%%%%%%%%%%%%%%%%%%%%%%%%%%%%%%%%%%%%%%%%%%%%%%%%%%%%%%%%%%%%%%%
\subsubsection{Determination of Elliptic Divisors: 
$D_4 \oplus A_2 \hookrightarrow E_8$}
\label{sssec:D4A2toE8}
%%%%%%%%%%%%%%%%%%%%%%%%%%%%%%%%%%%%%%%%%%%%%%%%%%%%%%%%%%%%%%%%%%%%%%%%%%%%

Finally, we discuss the cases where $T_0 = (D_4\oplus A_2)$ are embedded 
into $E_8$ (no other embeddings into lattices of type $E_n$ 
are possible). We use $(\vec{v}_6=Y_R), Y, Z, X, S', P, Q_2$ and $U$ for 
$E_8$ and enlarge it to the extended diagram by adding a vector
$\vec{v}'_8$, as shown in Figure \ref{dnchaindn7e8}. 
The only two candidates for $\vec{v}'_8$ are $2\nu_K$ with 
$K_{\epsilon\diamondsuit}$, and $K_{\epsilon \spadesuit}$. 
The two possible choices correspond to [no. 1], and [no. 6], respectively.
This can be seen from Figure \ref{fig:Ew-Coxeter} as follows:
first we note---through computations using codewords---that $\vec{v}_6$ 
are not at norm-$(-4)$ distance from $D_4$, $G_{43}$ and $G_{34}$, 
so that these vertices must be removed from Figure \ref{fig:Ew-Coxeter} 
for both cases. In case [no. 1], 
where $\vec{v}'_8=2\nu_{K_{\epsilon\diamondsuit}}$, we furthermore 
have to remove $D_1$, so that twice the extended Dynkin diagram of $E_8$ 
appears. Hence this choice of $\vec{v}'_8$ corresponds to the embedding 
$T_0\hookrightarrow E_8 \hookrightarrow L^{(\gamma)}_{\rm root}$, type 
[no.1] of \cite{Nish01}. For case [no. 6], where 
$\vec{v}'_8=2\nu_{K_{\epsilon \spadesuit}}$, we need to remove $C_1$ 
instead of $D_1$ and find the extended Dynkin diagram of 
$D_{16}\subset L^{(\beta)}_{\rm root}$. Elliptic divisors are computed 
by using Lemma I, and the results are presented in 
Table \ref{tab:KmEoxEo-multiplicity}.

%%%%%%%%%%%%%%%%%%%%%%%%%%%%%%%%%%%%%%%%%%%%%%%%%%%%%%%%%%%%%%%
\subsubsection{Multiplicities in $X = {\rm Km}(E_\omega \times E_\omega)$}
\label{sssec:KmEwxEw-multiplicity}
%%%%%%%%%%%%%%%%%%%%%%%%%%%%%%%%%%%%%%%%%%%%%%%%%%%%%%%%%%%%%%%

%%%%%%%%%%%%%%%%%%%%%%%%%%%%%%%%%%%%%%%%%%%%%%%%%%%%%%%%%%%%%
\begin{table}[tbp]
\begin{center}
%%%%%%%%%%%%%%%%%%%%%%%%%%%%%%%%%%%%%%%%%%%%%%%%%%%%%%%%%%
\begin{tabular}{|c|c|c|c|c|c|}
\hline
no. & $I$ & $D_4$ & $A_2$ & elliptic divisor & n \\
\hline
4 & $\beta$ & $D_{16}$ & $E_8$ & $D_4 + 2 G_{23} + 3 E_{3} + 4 G_{4 3} + 5 \lambda_{q_{11}} + 6 \lambda_p + 3 \lambda_{q_{18}} + 
  4 \lambda_{q_2} + 2 \lambda_{K_{\epsilon}\clubsuit}$ & 2\\
14 & $\epsilon$ & $D_{12}^{(1)}$ & $D_{12}^{(2)}$ & $D_{1} + C_{4} + \lambda_{q_{18}} + \lambda_{q_2} + 
 2 (\lambda_p + \lambda_{q_{11}} + G_{4 3} + E_{3} + G_{2 3} + D_{4} + G_{1 4} + 
    E_{4} + G_{4 4})$ & 2\\
8 & $\eta$ & $D_{10}$ & $E_7$ & $E_{1} + D_{4} + 2 (G_{2 1} + G_{1 4} + D_{1}) + 
 3 (C_{4} + E_{4}) + 4 G_{4 4}$ & 2\\
17 & $\theta$ & $D_9$ & $A_{15}$ & $\lambda_p + \lambda_{q_2} + \lambda_{K_{\epsilon}\clubsuit} + F_{3} + G_{3 1} + E_{1} + G_{2 1} + C_{4} + 
 G_{4 4} + E_{4} $ & \\
 &&&&$+ G_{1 4} + D_{4} + G_{2 3} + E_{3} + 
 G_{4 3} + \lambda_{q_{11}}$ & 2\\
16 & $\iota$ & $D_8^{(1)}$ & $D_8^{(2)}$& $E_{4} + D_{1} + D_{3} + F_{3} + 
 2 (G_{4 4} + C_{4} + G_{2 1} + E_{1} + G_{3 1})$&  2\\
20 & $\lambda$ & $D_7$ & $A_{11}$ &$ G_{1 4} + G_{4 2} + G_{4 3} + 2 (D_{4} + C_{1} + E_{3}) + 
 3 G_{2 3}$ & 2\\ 
22 & $\lambda$ & $D_7$ & $E_6$ & $G_{1 4} + D_{4} + G_{2 3} + C_{1} + G_{4 2} + D_{3} + 
 E_{1} + G_{3 1} + C_{4} + G_{2 1} + E_{4} + G_{4 4} $& 1\\
25 & $\xi$ & $D_6^{(1)}$ & $D_6^{(2)}$ & $E_{1} + F_{3} + 2 (G_{3 1} + D_{3} + G_{4 2}) + C_{1} + E_{2}$
& 1\\
26 & $\nu$ & $D_6$ & $A_9$ & $G_{4 3} + D_{2} + G_{2 1} + C_{4} + G_{4 4} + E_{4} + 
 G_{1 4} + D_{4} + G_{2 3} + E_{3}$ & 2\\
27 & $\pi$ & $D_5$ & $A_7$ & $G_{2 4} + D_{3} + G_{3 1} + F_{3} + G_{3 3} + D_{1} + 
 G_{4 4} + E_{4}$ & 1\\
28 & $\pi$ & $D_5^{(1)}$ & $D_5^{(2)}$ & $G_{1 4} + F_{1} + G_{1 2} + E_{2} + G_{4 2} + C_{1} + 
 G_{2 3} + D_{4} $&1 \\
29 & $\tau$ & $D_4^{(1)}$ & $D_4^{(2)}$ & $2 G_{1 1} + E_{1} + F_{1} + C_{1} + D_{1}$ &1\\
30 & $\sigma$ & $D_4$ & $A_5$ & $G_{1 1} + C_{1} + G_{3 4} + F_{3} + G_{3 3} + D_{1}$& 1\\
\hline
2 & $\gamma$ & $E_8^{(1)}$ & $E_8^{(2)}$ &$ G_{2 3} + 2 (E_{3} + F_{1}) + 3 (G_{3 3} + G_{2 4}) + 
 4 (D_{1} + G_{1 4}) + 5 G_{4 4} + 6 E_{4}$ & 2\\
5 & $\beta$ & $E_8$ & $D_{16}$ & $D_{1} + E_{3} + F_{1} + D_{2} + 
 2 (G_{3 3} + F_{3} + G_{3 1} $& \\
 &&&&$ + E_{1} + G_{4 1} + \lambda_{q_{18}} + \lambda_p + \lambda_{q_2} +
     \lambda_{K_{\epsilon}\spadesuit} + C_{1} + G_{4 2} + E_{2} + G_{1 2})$ & 2\\
9 & $\eta$ & $E_7$ & $D_{10}$ &$ 2 G_{2 4} + G_{3 3} + 2 D_{1} + 3 G_{4 4} + 4 E_{4} + 
 3 G_{1 4} + 2 F_{1} + G_{1 2} $& 1\\
10 & $\eta$ & $E_7^{(1)}$ & $E_7^{(2)}$ & $\lambda_p + G_{3 3} + 2 (\lambda_{q_2} + E_{3} + G_{4 2}) + 3 (\lambda_{K_{\epsilon}\spadesuit} + G_{2 3}) + 
 4 C_{1}$ & 1\\
11 & $\zeta$ & $E_7$ & $A_{17}$ &$ C_{1} + \lambda_{K_{\epsilon}\spadesuit} + \lambda_{q_2} + \lambda_p + \lambda_{q_{18}} + G_{4 1} + E_{1} + G_{2 1} + D_{2} + 
 G_{1 2} $ &2 \\ 
 &&&& $+F_{1} + G_{1 4} + E_{4} + G_{4 4} + D_{1} + 
 G_{3 3} + E_{3} + G_{2 3}$ &\\
19 & $\mu$ & $E_6^{(1)}$ & $E_6^{(2)}$ & $G_{3 1} + G_{2 3} + G_{4 4} + 2 (F_{3} + E_{3} + D_{1}) + 
 3 G_{3 3}$ & 1\\
23 & $\lambda$ & $E_6$ & $A_{11}$ & $G_{2 4} + G_{1 4} + F_{3} + E_{3} + 
 2 (E_{4} + G_{4 4} + D_{1} + G_{3 3})$ &2 \\
24 & $\lambda$ & $E_6$ & $D_7$ &$ \lambda_{q_{18}} + \lambda_{q_2} + G_{4 4} + G_{1 4} + 2 (\lambda_p + \lambda_{q_{11}} + G_{2 4} + E_{4})$&1 \\
\hline
\hline
12 & $\alpha$ & $D_{24}$ & $D_{24}$ & $\lambda_{1}+\lambda_{4} + G_{2 1} + G_{3 1} + 
 2 (\lambda_{2} + \lambda_{3} + 2\nu_{K_S} + \lambda_p + \lambda_{q_2} + \lambda_{K_{\epsilon}\clubsuit} + F_{1} + G_{1 3} + E_{3}$ & \\  
  &&&&  $+G_{3 3} + D_{1} + G_{4 4} + E_{4} + G_{3 4} + C_{1} + 
    G_{4 2} + E_{2} + G_{3 2} + D_{4} + G_{4 1} + E_{1}) $& 2\\
3 & $\beta$ & $D_{16}$ & $D_{16}$ & $G_{4 2} + 2 (C_{2} + E_{2}) + 3 (G_{3 2} + G_{2 1}) + 
 4 (G_{3 1} + D_{4}) + 5 G_{4 1} + 6 E_{1} $& 2\\
13 & $\epsilon$ & $D_{12}^{(1)}$ & $D_{12}^{(1)}$ & $G_{2 4} + G_{4 4} + G_{2 1} + G_{3 1} + 
 2 (E_{4} + G_{3 4} + C_{1}$ & \\
 &&&&$+ G_{4 2} + E_{2} + G_{3 2} + 
    D_{4} + G_{4 1} + E_{1})$ &1 \\ 
7 & $\eta$ & $D_{10}$ & $D_{10}$ & $G_{3 3} + G_{4 2} + 2 (D_{1} + C_{1} + G_{2 4}) + 
 3 (G_{4 4} + G_{3 4}) + 4 E_{4}$& 1\\
18 & $\theta$ & $D_9$ & $D_9$ &$ E_{1} + G_{4 1} + D_{4} + G_{3 2} + E_{2} + G_{4 2} + 
  C_{1} + G_{3 4} $& \\
  &&&&$+ E_{4} + G_{4 4} + D_{1} + G_{3 3} + 
  E_{3} + G_{4 3} + C_{2} + G_{3 1} $&1\\
15 & $\iota$ & $D_8^{(1)}$ & $D_8^{(1)}$ & $G_{2 4} + G_{4 4} + G_{3 2} + G_{1 2} + 
 2 (E_{4} + G_{3 4} + C_{1} + G_{4 2} + E_{2})$ & 1\\
21 & $\lambda$ & $D_7$ & $D_7$ & $D_{4} + F_{2} + C_{2} + 2 (G_{4 1} + G_{2 1} + G_{3 1}) + 
 3 E_{1} $& 1\\
\hline
1 & $\gamma$ & $E_8^{(1)}$ & $E_8^{(1)}$ & $G_{3 3} + 2 E_{3} + 3 G_{1 3} + 4 F_{1} + 5 G_{1 2} + 
 6 E_{2} + 4 G_{4 2} + 2 C_{1} + 3 G_{3 2}$ &2\\
6 & $\beta$ & $E_8$ & $E_8$ &$ G_{3 1} + G_{4 1} + G_{4 2} + G_{3 2} + 
 2 (E_{1} + G_{2 1} + F_{2} + G_{2 4}$  & \\ 
&&&&   $+ E_{4} + G_{4 4} + 
    D_{1} + G_{3 3} + E_{3} + G_{1 3} + F_{1} + G_{1 2} + 
    E_{2})$& 2\\
\hline 
  \end{tabular}
%%%%%%%%%%%%%%%%%%%%%%%%%%%%%%%%%%%%%%%%%%%%%%%%%%%%%%%%%%
\caption{\label{tab:KmEoxEo-multiplicity}
Type-dependent upper bounds on the multiplicity for the thirty types for
${\cal J}_2(X)$ of $X = {\rm Km}(E_\omega \times E_\omega)$ 
(shown in the last column). 
The type id. number in the first column is that of Table 1.3 
in \cite{Nish01}. In the second column, the Niemeier lattice 
which $T_0$ is embedded into is specified. The third and fourth columns 
give more details of this embedding. In the fifth column, we give 
a possible choice of an elliptic divisor as an element of 
${\rm II}_{1,25}$, followed by the upper bound on the multiplicity in the last column.}
 \end{center}
\end{table}
%%%%%%%%%%%%%%%%%%%%%%%%%%%%%%%%%%%%%%%%%%%%%%%%%%%%%%%%%%%%%%%%%
The ${\rm Aut}(D')$ and ${\rm Aut}(D')^{({\rm Hodge})}$ group 
action on $S_X \subset {\rm II}_{1,25}$ is known in the literature 
for $X = {\rm Km}(E_\omega \times E_\omega) =
X_{[2~0~2]}$, as summarized in the appendix \ref{ssec:review-Kummer}.  
Now that the study in 
sections \ref{sssec:D4toDn-A2toOther}--\ref{sssec:D4A2toE8}
specified at least one elliptic divisor within $\overline{D'}$ for 
individual isometry classes (type) in ${\cal J}_2(X)$, 
we are ready to use Corollary F to derive type-dependent upper bounds 
on the multiplicity. The results are found in the last
column of Table \ref{tab:KmEoxEo-multiplicity}.

\begin{quote}
 {\bf Example J}: Among the 30 different types of elliptic fibrations 
in ${\cal J}_2(X)$ for $X = {\rm Km}(E_\omega \times E_\omega)$, there 
are at least 15 types where there is a unique isomorphism class of
elliptic fibration. They are \\
no.7 [$W_{\rm root} = A_3E_7^{\oplus 2}$], no.9 [$A_1^{\oplus 3}D_7E_7$], 
no.10 [$A_1^{\oplus 3}A_5D_{10}$], no.13 [$D_5D_{12}$], 
no.15 [$D_8^{\oplus 2}$], no.18 [$A_1^{\oplus 2}A_{15}$], 
no.19 [$A_2^{\oplus 2}E_6^{\oplus 2}$], 
no.21 [$A_{11}E_{6}$], no.22 [$A_2^{\oplus 2}A_3A_{11}$], 
no.24 [$A_{11}D_4$], no.25 [$A_1^{\oplus 2}A_3D_6^{\oplus 2}$], 
no.27 [$A_4A_7D_5$], no.28 [$A_1^{\oplus 2}A_7^{\oplus 2}$], 
no.29 [$D_4^{\oplus 4}$], no.30 [$A_2A_5^{\oplus 3}$]. \\
The remaining 15 types have at most two isomorphism classes of 
elliptic fibrations.
\end{quote}

%%%%%%%%%%%%%%%%%%%%%%%%%%%%%%%%%%%%%%%%%%%%%%%%%%%%%%%%%%%%
\section*{Acknowledgements}
%%%%%%%%%%%%%%%%%%%%%%%%%%%%%%%%%%%%%%%%%%%%%%%%%%%%%%%%%%%%

It is our great pleasure to express gratitude to S. Kondo, S. Mukai, 
I. Shimada and T. Shioda for useful comments and communications.
A.~P.~B. likes to thank the IPMU Tokyo for kind hospitality during 
his visit, which made this collaboration possible.
The work of A.~P.~B. was supported by a JSPS postdoctoral fellowship 
under grant PE 12530 and the STFC under grant ST/J002798/1,
and that of T.~W. by WPI Initiative, MEXT, Japan and a Grant-in-Aid for 
Scientific Research on Innovative Areas 2303. 

\appendix

%%%%%%%%%%%%%%%%%%%%%%%%%%%%%%%%%%%%%%%%%%%%%%%%%%%%%%%%%%%%%%
\section{Appendix}
%%%%%%%%%%%%%%%%%%%%%%%%%%%%%%%%%%%%%%%%%%%%%%%%%%%%%%%%%%%%%%%%%%%

% \renewcommand{\thesubsection}{\thesection.A}
%%%%%%%%%%%%%%%%%%%%%%%%%%%%%%%%%%%%%%%%%%%%%%%%%%%%%%%%%%%%%%%%%%
\subsection{The Leech Lattice and the Leech Roots}
\label{ssec:Leech-review}
%%%%%%%%%%%%%%%%%%%%%%%%%%%%%%%%%%%%%%%%%%%%%%%%%%%%%%%%%%%%%%%%%%%

Here, we provide a review on Leech lattice $\Lambda_{24}$, which 
we use extensively in section \ref{sec:J1-classification} in this article. 
Conway--Sloan's textbook \cite{CS-book} is the ideal reference 
for this subject, but this review will explain the minimum prerequisite 
in reading this article, and at the same time, at least serve for 
the purpose of setting notation to be used in this article.

%%%%%%%%%%%%%%%%%%%%%%%%%%%%%%%%%%%%%%%%%%%%%%%%%%%%%%%%%%%%
\subsubsection{The Leech Lattice}
%%%%%%%%%%%%%%%%%%%%%%%%%%%%%%%%%%%%%%%%%%%%%%%%%%%%%%%%%%%%%

The Leech lattice $\Lambda_{24}$ is uniquely characterized as the even unimodular 
negative definite lattice of rank 24 which does not have any element of norm $(-2)$. 

In order to construct Leech lattice $\Lambda_{24}$, we begin by
describing its free Abelian group of rank-24 as a subset of 
$\R^{24} = \{ \sum_{i \in \Omega} x_i \nu_i \; | \; x_i \in \R\}$, 
where $\nu_i$'s labelled by\footnote{The ordering used here anticipates the
introduction of the miracle octad generator, see \eqref{eq:MOG-layout}.} 
\begin{equation}
 \Omega := \{0, 19, 15, 5, \; \; \infty, 3, 6, 9,\; \; 1,20,14,21,\;\; 
11,4,16,13,\; \; 2,10,17,7,\; \; 22,18,8,12\}
\label{eq:order-in-Omega}
\end{equation}
are vectors in $\R^{24}$, mutually linearly independent over $\R$. 
It is convenient to introduce a notation 
$\nu_S := \sum_{i \in S} \nu_i \in \R^{24}$ for a subset 
$S \subset \Omega$. 
With this notation, $\Lambda_{24}^{\rm even}$ is given by 
\begin{equation}
 \Lambda_{24}^{\rm even} = \left\{ \left. 
   2 \nu_C +  \sum_{i \in \Omega} 4 n_i \nu_i \; \right| \; n_i \in \Z, \; 
  \sum_i n_i \equiv 0 \; {\rm mod~}2 \right\}
\end{equation}
where $C$ can be any one of codewords of the extended binary 
Golay code ${\cal C}(24)$ (see below). 
$\Lambda_{24}^{\rm odd}$ consists of 
\begin{equation}
 \left\{ \left. \nu_\Omega + 2 \nu_C + 
    \sum_{i \in \Omega} 4m_i  % \left(1 + 4 m_i \right)
           \nu_i \; \right| \; 
    m_i \in \Z, \; \sum_i m_i \equiv 1 \; {\rm mod~}2 \right\}.
\end{equation}
$\Lambda_{24} = \Lambda_{24}^{\rm even} \cup \Lambda_{24}^{\rm even}$ 
forms a free Abelian group of rank 24 under the ordinary vector sum 
in $\R^{24}$.

A little bit of explanation on the Golay code is in order here. 
Let $\P(\Omega)$ denote a set that consists of 
any subset of $\Omega$. $\P(\Omega)$ consists of $2^{24}$ elements, as 
it can be identified with a vector space $(\F_2)^{24}$ over the field $\F_2$.
The {\it extended binary Golay code of length 24}, ${\cal C}_{24}$, is a specific 
choice of a subset of $\P(\Omega)$. Instead of explaining its construction, we restrict
ourselves to record some properties relevant to this article, see \cite{CS-book} for a detailed treatment.
\begin{itemize}
 \item ${\cal C}_{24}$ is a $\F_2$-linear $12$-dimensional subspace of 
      $\P(\Omega) \cong (\F_2)^{24}$. Thus, $|{\cal C}_{24}| = 2^{12}$.
 \item $\varnothing \in {\cal C}_{24}$, $\Omega \in {\cal C}_{24}$, 
 \item ${\cal C}_{24}$ is decomposed into five subsets 
   $\{ \phi \} \cup \{ \Omega \} \cup 
    {\cal C}_{24}(8) \cup {\cal C}_{24}(12) \cup {\cal
       C}_{24}(16)$, and all the elements of ${\cal C}_{24}(8)$
    (resp. ${\cal C}_{24}(12)$, ${\cal C}_{24}(16)$) are 
    8-element (resp. 12-element, 16-element) subsets of $\Omega$. 
\end{itemize}
Any codeword in ${\cal C}_{24}(8)$ is called an {\it octad}.
\begin{itemize}
 \item $|{\cal C}_{24}(8)| = |{\cal C}_{24}(16)| = 759$ and 
       $|{\cal C}_{24}(12)| = 2576$.
 \item If one chooses four arbitrary elements from $\Omega$, i.e., $P
       \subset \Omega$, $|P| = 4$, then 
       there are five (and not more than five) codewords 
       $K_{1,2,3,4,5} \in {\cal C}_{24}(8)$ 
       that contain the specified 4 elements $P$. $K_i \cap K_j$ consists 
       precisely of the 4 elements $P$, when $i \neq j$. Thus, for a given 
       4-element subset $P \subset \Omega$, there is a unique way to decompose
       $\Omega$ into six 4-element subsets; $\Omega = 
       \amalg_{i=1}^5 (K_i \backslash P) \amalg P$. Such a 
       decomposition is called {\it sextet decomposition}, 
       see below for explicit examples. 
       For an explicit algorithm of finding out this 
       decomposition, \cite{CS-book} Chapter 11 will be the best 
       reference to look at.
\begin{itemize}
 \item A sextet decomposition may be described as 
       $\Omega = \amalg_{a=0,\cdots,5} \Xi_a$. It is a sextet
       decomposition of 4-element subsets $\Xi_a \subset \Omega$ with 
       $a=1,\cdots,5$ as much as for the 4-element subset $\Xi_0 \subset
       \Omega$.
 \item In a given sextet decomposition, $\amalg_{a} \Xi_a$, 
       any 8-element subset of the form $K_{ab} = \Xi_a \amalg \Xi_b$
       ($a \neq b$) is a codeword in ${\cal C}_{24}(8)$.
 \item If one chooses arbitrary 5 elements out of $\Omega$, then there 
       is a unique codeword in ${\cal C}_{24}(8)$ that contains the
       five elements. 
\end{itemize}
\end{itemize}

The free Abelian group $\Lambda_{24}$ becomes th Leech lattice with the 
symmetric pairing $(\nu_i, \nu_j) = - \delta_{ij} /8$. 
It is known that the Leech lattice defined in this way is an even
unimodular negative definite lattice of rank 24.

There is no norm $(-2)$ elements in this lattice. 
Norm $(-4)$ elements of $\Lambda_{24}$ are of the form 
\begin{equation}
 (2^8,0^{16})^T, \qquad (3,1^{23})^T, \qquad (4^2,0^{22})^T
\label{eq:norm-4-in-Leech}
\end{equation}
in the component description $(x_0, x_{19}, x_{15}, \cdots, x_8,
x_{12})^T$ modulo signs of each entries and ordering. 
Norm $(-6)$ elements of $\Lambda_{24}$ are of the form 
\begin{equation}
 (2^{12},0^{12})^T, \quad 
 (3^3,1^{21})^T, \quad 
 (4,2^8,0^{15})^T, \quad 
 (5,1^{23})^T .
\label{eq:norm-6-in-Leech}
\end{equation}

The isometry group of this lattice, ${\rm Isom}(\Lambda_{24})$, 
is often denoted by $(\cdot 0)$, or $Co_0$. Its Affine transformation
group is denoted by $(\cdot \infty)$ or $Co_\infty$:
\begin{equation}
 1 \longrightarrow \Z^{24} \longrightarrow Co_\infty 
   \longrightarrow Co_0 \longrightarrow 1.
\end{equation}
%

%%%%%%%%%%%%%%%%%%%%%%%%%%%%%%%%%%%%%%%%%%%%%%%%%%%%%%%%%%%%
\subsubsection{More on ${\cal C}_{24}(8)$}
\label{section:octads}
%%%%%%%%%%%%%%%%%%%%%%%%%%%%%%%%%%%%%%%%%%%%%%%%%%%%%%%%%%%%

It proves convenient in dealing with codewords of ${\cal C}_{24}$ 
to describe subsets of $\Omega$ (i.e., elements of $\P(\Omega)$) 
as follows. Any subset, let's say $K$, of $\Omega$ is specified by 
whether individual elements are contained in it. So, when the 24
elements in $\Omega$ are laid out as in 
\begin{equation}
 \Omega =
 \begin{array}{|cc|cc|cc|}
   \hline
    0  & \infty & 1  & 11 & 2  & 22 \\
    19 & 3      & 20 & 4  & 10 & 18 \\
    15 & 6      & 14 & 16 & 17 & 8  \\
    5  & 9      & 21 & 13 & 7  & 12 \\
    \hline
  \end{array},
\label{eq:MOG-layout}
\end{equation}
a subset $K_U = \{ 0,\infty,1,11,2,18,8,12\}$, for example, can be denoted by 
\begin{equation}
 \begin{array}{|cc|cc|cc|}
 \hline
 *&*&*&*&*& \\
  & & & & &*\\
  & & & & &*\\
  & & & & &*\\
 \hline
 \end{array}.
\label{eq:Ku}
\end{equation}
We will also use the two following codewords in ${\cal C}_{24}(8)$ 
in the main text of this article:
\begin{equation}
 K_L = 
  \begin{array}{|cc|cc|cc|}
   \hline
    *&{\color{white} *} &{\color{white} *} &*&*&* \\
     & & & & &* \\
     & & & &*&  \\
    *& & &*& &  \\
   \hline
  \end{array}, \qquad 
 K_R = 
  \begin{array}{|cc|cc|cc|}
   \hline
    {\color{white} *}&*&{\color{white} *} &*&*&* \\
    & & & & &* \\
    &*& &*& &  \\
    & & & &*&  \\
   \hline
  \end{array}.
\label{eq:Kl-Kr} 
\end{equation}
The layout of the 24 elements of $\Omega$ in (\ref{eq:MOG-layout})
follows so called ``standard MOG labelling'' (p.309 Fig. 11.7 of \cite{CS-book}).

Because we will use them in the main text of this article, we leave a list of some
of the sextet decompositions we referred to above. The sextet decomposition for 
the following five 4-element subsets of $K_U$ (and hence subsets of $\Omega$), 
$\Xi^{(\epsilon)}_* = \{\infty, 0,1,2\}$, 
$\Xi^{(\alpha)}_* = \{ \infty, 0,1,11\}$, 
$\Xi^{(\beta)}_* = \{ \infty,0,1,18\}$, 
$\Xi^{(\gamma)}_* = \{ \infty,0,1,8\}$, 
$\Xi^{(\delta)}_* = \{ \infty, 0,1,12\}$ are given as follows:\footnote{
We followed the explanations in Chapter 11 of \cite{CS-book} in order 
to determine these sextet decompositions in ${\cal C}_{24}$.}
\begin{align}
&
\amalg_{a \in A} \Xi^{(\epsilon)}_a = 
 \begin{array}{|cc|cc|cc|}
\hline
  * &* &*& \circ &* & \heartsuit \\
  \clubsuit & \diamondsuit & \spadesuit & \heartsuit & \spadesuit &
   \circ \\
  \spadesuit & \clubsuit & \diamondsuit & \heartsuit & \diamondsuit &
   \circ \\
  \diamondsuit & \spadesuit & \clubsuit & \heartsuit & \clubsuit & \circ \\
\hline
 \end{array} \, , \nonumber \\
& 
\amalg_{a \in A} \Xi^{(\alpha)}_a = 
  \begin{array}{|cc|cc|cc|}
 \hline
  * & * & * & * & \circ & \heartsuit \\
  \clubsuit & \clubsuit & \clubsuit & \clubsuit & \heartsuit & \circ \\
  \spadesuit & \spadesuit & \spadesuit & \spadesuit & \heartsuit &
   \circ \\
  \diamondsuit & \diamondsuit & \diamondsuit & \diamondsuit & \heartsuit
   & \circ \\
  \hline
 \end{array} \, ,   \qquad 
\amalg_{a \in A} \Xi^{(\beta)}_a = 
\begin{array}{|cc|cc|cc|}
 \hline 
 * & * & * & \circ & \circ & \heartsuit \\
 \clubsuit & \diamondsuit & \heartsuit & \spadesuit & \spadesuit & * \\
 \clubsuit & \heartsuit & \spadesuit & \diamondsuit & \clubsuit &
  \circ \\
 \heartsuit & \diamondsuit & \spadesuit & \clubsuit & \diamondsuit &
  \circ \\ 
 \hline
\end{array} \, ,  \label{eq:five-sextets} \\
&
\amalg_{a \in A} \Xi^{(\gamma)}_a = 
\begin{array}{|cc|cc|cc|}
\hline
 * & * & * & \circ & \circ & \heartsuit \\
 \heartsuit & \spadesuit & \diamondsuit & \clubsuit & \spadesuit &
  \circ \\
 \clubsuit & \spadesuit & \heartsuit & \diamondsuit & \diamondsuit &
  * \\
 \clubsuit & \heartsuit & \diamondsuit & \spadesuit & \clubsuit & \circ \\
\hline
\end{array} \, , \qquad 
\amalg_{a \in A} \Xi^{(\delta)}_a = 
\begin{array}{|cc|cc|cc|}
\hline
 * & * & * & \circ & \circ & \heartsuit \\
 \clubsuit & \heartsuit & \spadesuit & \diamondsuit & \clubsuit &
  \circ \\
 \heartsuit & \diamondsuit & \spadesuit & \clubsuit & \diamondsuit &
  \circ \\
  \clubsuit & \diamondsuit & \heartsuit & \spadesuit & \spadesuit & * \\
\hline
\end{array} \, , \nonumber 
\end{align}
where 
$a \in A = \{ *, \circ, \heartsuit, \clubsuit, \spadesuit, \diamondsuit\}$.
From these five sextet decompositions, we can find out all the 
codewords $K$ in ${\cal C}_{24}(8)$ satisfying $\{0,\infty,1\} \subset K$
and $|K \cap K_U| = 4$: 
\begin{eqnarray}
 K_{\epsilon a} & = & \Xi^{(\epsilon)}_* \amalg \Xi^{(\epsilon)}_a, 
     \qquad \qquad \qquad \qquad \qquad \qquad 
      a \in \{ \heartsuit, \clubsuit, \spadesuit, \diamondsuit \} \\
 K_{\alpha a} & = & \Xi^{(\alpha)}_* \amalg \Xi^{(\alpha)}_a, \qquad 
 K_{\beta a} = \Xi^{(\beta)}_* \amalg \Xi^{(\beta)}_a, \\ % \quad 
 K_{\gamma a} & = & \Xi^{(\gamma)}_* \amalg \Xi^{(\gamma)}_a, \qquad 
 K_{\delta a} = \Xi^{(\delta)}_* \amalg \Xi^{(\delta)}_a. 
\end{eqnarray}
For example, 
$K_{\epsilon\heartsuit} = \{ 0,\infty,1,4,16,13,2,22\}\subset \Omega$. 
Overall, there are twenty codewords $K \in {\cal C}_{24}(8)$ satisfying 
$\{0,\infty, 1\} \subset K$ and $|K \cap K_U| = 4$.

%%%%%%%%%%%%%%%%%%%%%%%%%%%%%%%%%%%%%%%%%%%%%%%%%%%%%%%%%%%%%%%%%%
\subsubsection{The Centres of Deep Holes in $\Lambda_{24}$}
\label{sssec:deep-hole}
%%%%%%%%%%%%%%%%%%%%%%%%%%%%%%%%%%%%%%%%%%%%%%%%%%%%%%%%%%%%%%%%%%

All that we state in this appendix can also be 
found in Chapt. 23 in \cite{CS-book}. We only quote the results 
that we need in the main text of this article. For more systematic 
exposition on this subject, we recommend to consult \cite{CS-book}.

For any $\vec{u} \in (\Lambda_{24} \otimes \R) \backslash \Lambda_{24}$, 
it is known that 
\begin{equation}
 {\rm Min}\left[ -(\vec{u}-\vec{v})^2 \; | \;
                 \vec{v} \in \Lambda_{24} \right] \leq 2. 
\end{equation}
Points in $\Lambda_{24} \otimes \R$ saturating this inequality---points 
in $\Lambda_{24} \otimes \R$ that can be as far away as possible from
integral points $\Lambda_{24}$---are called {\it centres of deep holes in
} $\Lambda_{24}$. The $Co_\infty$ symmetry group of $\Lambda_{24}$, with its 
action naturally extended linearly to $\Lambda_{24} \otimes \R$, acts 
on these centre of deep holes. It is known that the centres of deep holes
form twenty-three distinct orbits of $Co_\infty$, and are labelled by 
$I = \alpha, \beta, \cdots, \phi,\chi,\psi$, 
as in Table \ref{tab:Niemeier-root}. Centre of deep holes that belong to 
the label-$I$ orbit are said to be {\it type} $I$ in this article.

It is due to the following reason that the twenty-three different labels 
$I = \alpha, \beta, \cdots$ for Niemeier lattices can also be used to 
distinguish $Co_\infty$-inequivalent deep holes.
For the centre of a deep hole $\vec{c}^{(I)} \in \Lambda_{24} \otimes \R$ 
of type $I$, consider all of $\vec{v} \in \Lambda_{24}$ 
satisfying $(\vec{v}-\vec{c}^{(I)})^2 = -2$. A graph is determined 
for this set of points $\{\vec{v} \}$, by assigning one node for each 
$\vec{v}$ in this set, and by drawing, between a pair of nodes for $\vec{v}_i$ 
and $\vec{v}_j$, two lines if $(\vec{v}_i-\vec{c}, \vec{v}_j-\vec{c}) = +2$, 
one line if $(\vec{v}_i-\vec{c}, \vec{v}_j-\vec{c})= +1$, and 
no line if $(\vec{v}_i-\vec{c},\vec{v}_j-\vec{c}) = 0$. 
It is known that the graph becomes the collection of extended Dynkin 
diagrams of A--D--E root system in the combination specified 
in one of the twenty-three ($I=\alpha, \beta, \cdots, \psi$) entries 
in Table \ref{tab:Niemeier-root} (see \cite{CS-book} Chapt. 23). 
The type label $I$ for the $Co_\infty$-orbits of the centres of deep
holes is assigned through the correspondence with the classification of 
Niemeier lattices. 
The points $\vec{v} \in \Lambda_{24}$ satisfying 
$(\vec{v}-\vec{c}^{(I)})^2 = -2$ for a given deep hole of type $I$ are 
denoted by $\vec{v}^{(I)a}_i$, where $a = 1,2,\cdots$ label irreducible 
components of the A--D--E root systems in 
$(L^{(I)})_{\rm root} = \oplus_a R_a$, and $i$ runs from 0 
to the rank $r_a$ of the $a$-th component $R_a$. 
For example, for $I=\beta$, $R_{a=1}$ and $R_{a=2}$ correspond to $D_{16}$ 
and $E_8$, respectively, and $i=0,1,\cdots, r_1 = 16$ for $a=1$ and 
$i=1,\cdots,r_2 = 8$ for $a=2$.
It is further known that 
\begin{equation}
 h^{(I)} \; \vec{c}^{(I)} = \sum_{i=0}^{r_1} n^{(I)a=1}_i \vec{v}^{(I)a=1}_i
\end{equation}
holds \emph{for any} $a$.  Here, $n^{(I)a}_i$ are integers (sometimes referred to as {\it Kac
label} or {\it Dynkin label}) assigned to individual nodes of the
A--D--E extended Dynkin diagram ($n_0 = 1$ for any $A$--$D$--$E$ root 
system), so that the maximal root is given by the linear combination of 
simple roots $\alpha^a_i$ as $\sum_{i=1}^{r_a} n^a_i \alpha_i$. The number
$h^{(I)}$ on the left-hand side is the dual Coxeter number of the A--D--E 
root system, which is common for all irreducible components $R_a$ 
in the individual entries of Table \ref{tab:Niemeier-root}, i.e. there
is one $h^{(I)}$ for any Niemeier lattice.

It is useful for practical computations to note that, 
\begin{equation}
 (\vec{v}_i-\vec{v}_j)^2 = \left\{ \begin{array}{c} 
  -8 \\
  -6 \\
  -4 \end{array} \right. \longleftrightarrow
  \left.
  \begin{array}{c} 2 \\ 1 \\ 0 \end{array} \right\}
  = (\vec{v}_i-\vec{c}, \vec{v}_j-\vec{c}),
\end{equation}
respectively, under the condition that $(\vec{v}-\vec{c})^2 = -2$. 
Thus, in particular, the condition that $\vec{v}_i-\vec{c}$ and 
$\vec{v}_j-\vec{c}$ are orthogonal is equivalent to the norm-$(-4)$ 
distance between the two vectors $\vec{v}_i$ and $\vec{v}_j$ around 
a deep hole.

%%%%%%%%%%%%%%%%%%%%%%%%%%%%%%%%%%%%%%%%%%%%%%%%%%%%%%%%%%%%
\subsubsection{Leech Roots and the ${\rm II}_{1,25}$ Lattice}
%%%%%%%%%%%%%%%%%%%%%%%%%%%%%%%%%%%%%%%%%%%%%%%%%%%%%%%%%%%%

The rank-26 lattice $\Lambda_{24} \oplus U$ is an even (i.e., type II) 
unimodular lattice of signature $(1,25)$. From the classification theorem
for indefinite even unimodular lattices, mentioned in Section \ref{ssec:prelim-lattice},
it follows that it is unique modulo isometry. This lattice---denoted by ${\rm II}_{1,25}$---can be 
described as an Abelian group as 
\begin{equation}
 {\rm II}_{1,25} \cong \Lambda_{24} \oplus U =
   \left\{ (\vec{v}, m, n) \; |
  \vec{v} \in \Lambda_{24}, \; m,n \in \Z \right\}.
\end{equation}

The reflection symmetry group of this even unimodular lattice, 
$W({\rm II}_{1,25}) = W^{(2)}({\rm II}_{1,25})$, can be generated by 
reflections associated with simple roots, $\Pi$, which are called 
Leech roots in this case. It is possible to take a fundamental chamber 
of this reflection symmetry group, $C_\Pi$, so that 
the simple roots  (i.e., Leech roots) are of the form 
\begin{equation}
 \Pi = \left\{ \left.
    \lambda = \left(\vec{v},1,-1 - \frac{(\vec{v})^2}{2} \right) \in
	{\rm II}_{1,25} \; \right| \; \vec{v} \in \Lambda_{24} \right\},
 \label{eq:L24-Pi-relation}
\end{equation}
where $(\vec{v})^2$ is the norm of $\vec{v}$ in $\Lambda_{24}$.
For this reason, there is a one-to-one correspondence between the Leech
roots $\Pi$ and $\Lambda_{24}$.
Note that all these Leech roots are of norm $(-2)$ in ${\rm II}_{1,25}$, 
and that all of them satisfy\footnote{Such a vector is sometimes called a \emph{Weyl vector}} 
$(w,\lambda) = 1$ for $w = (\vec{0},0,1) \in {\rm II}_{1,25}$.
It is also useful to note that 
\begin{equation}
 (\lambda,\lambda')_{{\rm II}_{1,25}} = (\vec{v},\vec{v}') -2 
 - \frac{(\vec{v})^2 + (\vec{v}')^2}{2} =
  -2 - \frac{(\vec{v}-\vec{v}')^2_{\Lambda_{24}}}{2}, 
\label{eq:L24-Pi-relation-ii}
\end{equation}
which takes the values $-2, 0, +1, +2, \cdots$.

The isometry group of ${\rm II}_{1,25}$---${\rm Isom}({\rm II}_{1,25}) \cong 
{\rm Isom}^+({\rm II}_{1,25}) \times \{ \pm 1\}$---has the 
following structure: ${\rm Isom}^+({\rm II}_{1,25}) = W_\Pi \rtimes 
{\rm Isom}({\rm II}_{1,25})^{(C_\Pi)}$, where 
${\rm Isom}({\rm II}_{1,25})^{(C_\Pi)}$ is the group of isometries of 
${\rm II}_{1,25}$ that map $C_\Pi$ to itself. 
Since the Leech roots correspond to the boundary walls (reflection 
hyperplanes) of the chamber $C_\Pi$, 
isometries in ${\rm Isom}({\rm II}_{1,25})^{(C_\Pi)}$ map Leech roots 
$\Pi$ to themselves, and upon identification between $\Pi$ and
$\Lambda_{24}$, the group ${\rm Isom}({\rm II}_{1,25})^{(C_\Pi)}$ 
can be regarded as $Co_\infty$ on $\Lambda_{24}$.

% \renewcommand{\thesubsection}{\thesection.B}

%%%%%%%%%%%%%%%%%%%%%%%%%%%%%%%%%%%%%%%%%%%%%%%%%%%%%%%%%%%%%%%%%%
\subsection{The Neron--Severi Lattice of $X_3$}
\label{ssec:X3-app}
%%%%%%%%%%%%%%%%%%%%%%%%%%%%%%%%%%%%%%%%%%%%%%%%%%%%%%%%%%%%%%%%%

%%%%%%%%%%%%%%%%%%%%%%%%%%%%%%%%%%%%%%%%%%%%%%%%%%%%%%%%%%%%
\subsubsection{Tasks I--IV for $X_3$}
\label{sssec:X3-task}
%%%%%%%%%%%%%%%%%%%%%%%%%%%%%%%%%%%%%%%%%%%%%%%%%%%%%%%%%%%%

In this section, we see how the tasks I--IV introduced in 
section \ref{ssec:Borcherds-Kondo} are carried out and how the
assumptions (as-1)--(as-5) verified in practice by working on 
a specific example in detail. We choose a singular K3 surface 
$X=X_3 = X_{[1~1~1]}$ where the symmetric pairing of $T_X$ is given by 
$\left[ \begin{array}{cc} 2 & 1 \\ 1 & 2 \end{array}\right]$ and
the structure of the Neron--Severi lattice and its isometry group have 
been studied very well (e.g., \cite{Vin-2most, Borcherds2}). The
following discussion in appendices \ref{sssec:X3-task} and \ref{sssec:X3-detl} 
is a kind of calculation note that fills the gaps in these references
which are not obvious for non-experts.\footnote{
The example $X = X_3$ is chosen because it is the easiest example in 
carrying out the tasks I--IV and verify assumptions 1--5 in practice.
Apart from this exercise purpose, though, there is not much meaning 
in calculating ${\rm Aut}(D')$ for this particular example $X=X_3$, 
because we know already that 
$p_T:{\rm Isom}(T_X)^{({\rm Hodge})} \longrightarrow {\rm Isom}(q)$ 
is surjective (Table \ref{tab:coset-alt-rho20} 2nd row). We can conclude 
from Proposition C without any detailed knowledge about the structure of 
${\rm Aut}(D')$ that there is only one isomorphism class of elliptic fibrations  
in ${\cal J}_1(X)$ for $X = X_3$ for any isometry class $[W] \in {\cal J}_2(X)$. }
  
\begin{figure}[tbp]
\begin{center}
   \scalebox{.5}{ 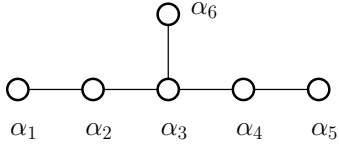 }
\caption{\label{fig:E6-Dynkin} The Dynkin diagram of $E_6$.}
\end{center} 
\end{figure}

In this example, $T_0$ is the $E_6$ root lattice, so that the assumptions (as-1) 
and (as-2) are satisfied. 
Let $J = \left\{ \alpha_1, \cdots, \alpha_6 \right\}$ be the simple roots 
of $T_0 = E_6$, see Figure~\ref{fig:E6-Dynkin}. 
%%%%%%%%%%%%%%%%%%%%%%%%%%%%%%%%%%%%%%%%%%%%
%%%%%%%%%%%%%%%%%%%%%%%%%%%%%%%%%%%%%%%%%%%%
The ${\rm Aut}(J) \cong \Z/2\Z$ group is generated by the left-right
flip of the Dynkin diagram (Figure \ref{fig:E6-Dynkin}), which exchanges 
the weights of the ${\bf 27}$ and $\overline{\bf 27}$ representations. 
Thus ${\rm Aut}(J) \cong \Z/2\Z$ is mapped to ${\rm Isom}(q) \cong
\Z/2\Z$ of the discriminant group $G_{T_0} \cong G_{T_X} \cong \Z/3\Z$
(Table \ref{tab:T0}) and the assumption (as-3) is also satisfied. 

[task I] 
The embedding $\phi_{T_0}: J \longrightarrow \Pi$ is given as follows. 
First, we define 6 points in $\Lambda_{24}$ as follows,
\begin{eqnarray}
 \vec{v}_1& := & 4\nu_{\{0,1\}} =: P',    \label{eq:E6-embed-vect-1} \\
 \vec{v}_2& := & 2\nu_{K_{L'}},           \label{eq:E6-embed-vect-2} \\
 \vec{v}_3& := & \nu_{\Omega} + 4 \nu_1 =: X', \label{eq:E6-embed-vect-3} \\
 \vec{v}_4& := & \vec{0} =: Z,            \label{eq:E6-embed-vect-4} \\
 \vec{v}_5& := & \nu_{\Omega} + 4 \nu_0 =: Y,  \label{eq:E6-embed-vect-5} \\
 \vec{v}_6& := & \nu_{\Omega} - 2\nu_{K_{U'}} + 4\nu_0.
                                          \label{eq:E6-embed-vect-6} 
\end{eqnarray}
$K_{L'}$ and $K_{U'}$ are codewords in ${\cal C}_{24}(8)$ satisfying the
following conditions:
\begin{equation}
 0 \in K_{L'}, K_{U'}, \qquad 
 1 \nin K_{L'}, K_{U'}, \qquad 
 |K_{L'} \cap K_{U'} | = 2.
\end{equation}
As an example, we can take 
\begin{equation}
 K_{L'} =
  \begin{array}{|cc|cc|cc|}
   \hline
   *&*&{\color{white} *} &{\color{white} *} &{\color{white} *} & {\color{white} *} \\
   *&*& & & &  \\
   *&*& & & &  \\
   *&*& & & &  \\
   \hline
  \end{array}, \qquad 
 K_{U'} =
  \begin{array}{|cc|cc|cc|}
   \hline
   *&*& &*&*&* \\
    & &*& & &  \\
    & &*& & &  \\
    & &*& & &  \\
   \hline
  \end{array},
\end{equation}
where we follow the conventions explained 
in appendix \ref{ssec:Leech-review}.
Noting that $(\vec{v}_i-\vec{v}_{i+1})^2 = -6$ for all $i=1,\cdots, 4$ 
and $(\vec{v}_3 - \vec{v}_6)^2 = -6$, but $(\vec{v}_i-\vec{v}_j)^2 = -4$
for $i \neq j$ otherwise, we see that the six simple roots 
$\alpha_{1,\cdots,6}$ of $T_0$ can be embedded into 
the six Leech roots
related to the six points $\vec{v}_i \in \Lambda_{24}$ 
through (\ref{eq:L24-Pi-relation}, \ref{eq:L24-Pi-relation-ii}).  
This embedding $J \longrightarrow \Pi$ as well as its linear extension 
$T_0 \hookrightarrow {\rm II}_{1,25}$ is denoted by $\phi_{T_0}$.

In order to see that $\phi_{T_0}: T_0 \hookrightarrow {\rm II}_{1,25}$ is 
a primitive embedding, it is sufficient to make sure that only 
$T_0 \cong E_6$ within $T_0^* = E_6^*$ are mapped to integral points in 
${\rm II}_{1,25}$. The discriminant group $G_{E_6} \cong \Z/3\Z$ is 
generated by a weight of the 27-dimensional representation of $E_6$, 
$(4\alpha_1+5\alpha_2+6\alpha_3+4\alpha_4+2\alpha_5+3\alpha_6)/3$. Since 
$\phi_{T_0}$ maps this weight to ${\rm II}_{1,25} \otimes \Q$, but not to 
${\rm II}_{1,25}$, this embedding 
$\phi_{T_0}: T_0 \hookrightarrow {\rm II}_{1,25}$ is indeed primitive and
the [task I] is completed.

[task II] Let $R := \phi_{T_0}(J)$. Since $R=E_6$ in this example, the 
only possible spherical subdiagrams $R' = R \cup \lambda$ of the 
Coxeter diagram of $\Pi$ are of the form of either 
type 0) $R'=E_6+A_1$ or type 1) $R' = E_7$.

Leech roots $\lambda$ of type 0), namely 
$\{ \lambda \in \Pi | (\lambda, \lambda_i) = 0 \; {\rm
for~}i=1,\cdots,6\}$,  correspond to 
\begin{equation}
\{\vec{v} \in \Lambda_{24} \; | \;  
(\vec{v}-\vec{v}_i)^2 = -4 {\rm ~for~} i=1,\cdots,6\}. 
\end{equation}
With the embedding given by (\ref{eq:E6-embed-vect-1}--\ref{eq:E6-embed-vect-6}),
the condition for $i=3$ implies that $\vec{v}$ must be 
a norm $(-4)$ element of $\Lambda_{24}$, which is always in the form 
of (\ref{eq:norm-4-in-Leech}). Further imposing the conditions for 
$i=1,2,4,5,6$, it turns out that $\vec{v} \in \Lambda_{24}$ should be 
of the form of $2\nu_K$ for $K \in {\cal C}_{24}(8)$ satisfying 
\begin{equation}
 \{ 0,1 \} \subset K, \qquad |K \cap K_{L'} | = 4, \qquad 
|K \cap K_{U'} | = 2\, .
\label{eq:cond-E6-embed-type0-root-K}
\end{equation}
The corresponding Leech roots of type 0) are of the form $(2\nu_K,1,1)$.
There are 24 codewords $K \in {\cal C}_{24}(8)$ satisfying these conditions, 
and hence there are 24 Leech roots of type 0).
%%%%%%%%%%%%%%%%%%%%%%%%%%%%%%%%%%%%%%%%%%%%%%%%%%%%%%%%%%%%%%
\begin{figure}[!h]
\begin{center}
   \scalebox{.5}{ 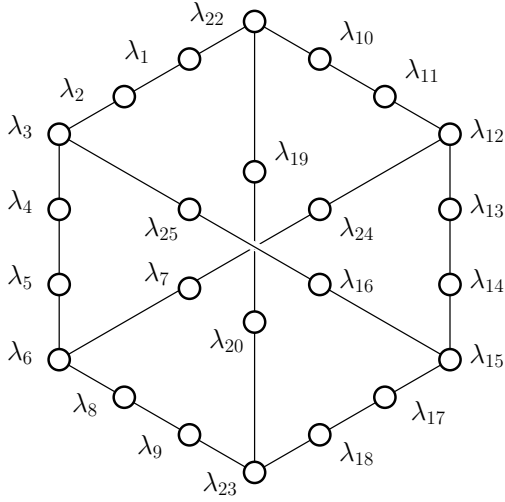 }
\caption{\label{fig:X3-Coxeter} Coxeter diagram of $S_X$ for $X=X_3$. Here,
only nodes associated with $(-2)$ roots are shown. 
The twenty-four vertices $\lambda_1, \cdots, \lambda_{20}$, 
$\lambda_{22} \cdots \lambda_{25}$ correspond to 
$e_{1, \cdots, 20}$, $e_{22, \cdots, 25}$ in \cite{Vin-2most}.
$\lambda_{2,\cdots,9}$ and $\lambda_{11, \cdots,18}$ are
the simple roots of the two $E_8$'s in 
$S_X \cong U \oplus A_2 \oplus E_8^{\oplus 2}$.}
\end{center} 
\end{figure}
%%%%%%%%%%%%%%%%%%%%%%%%%%%%%%%%%%%%%%%%%%%%%%%%%%%%%%%%%%%%%%
Among them, six involve the codewords $K \in {\cal C}_{24}(8)$ 
listed explicitly in the appendix \ref{ssec:Leech-review}:
\begin{equation}
  \lambda_{22} = (2\nu_{K_{\beta \clubsuit}}, 1,1), \quad 
  \lambda_{6} = (2\nu_{K_{\delta \clubsuit}},1,1), \quad 
  \lambda_{15} = (2\nu_{K_{\gamma \clubsuit}},1,1),
\label{eq:Leech-root-X3-triv-up}
\end{equation}
\begin{equation}
 \lambda_{23} = (2\nu_{ K_{\gamma \spadesuit} }, 1,1), \quad 
 \lambda_{3} = (2\nu_{ K_{\delta \diamondsuit} }, 1,1), \quad 
 \lambda_{12} = (2\nu_{ K_{\beta\diamondsuit} }, 1,1).
\label{eq:Leech-root-X3-triv-down}
\end{equation}
The remaining 18 other codewords $K$ are also given 
in Appendix \ref{sssec:X3-detl}.
Because all the Leech roots of type 0) are also norm $(-2)$ roots of 
$S_X$, the Coxeter diagram of the twenty-four Leech roots of type 0) 
becomes a subgraph of the Coxeter diagram of $W(S_X)$.
This graph for $X_3$ is given in Figure \ref{fig:X3-Coxeter}, which 
reproduces the result of \cite{Vin-2most}.

Leech roots of type 1) correspond to vectors $\vec{u} \in \Lambda_{24}$
satisfying either 
\begin{itemize}
 \item $(\vec{u}-\vec{v}_i)^2 = -4$ for $i =1,2,3,4,6$, but 
       $(\vec{u}-\vec{v}_5)^2 = -6$, or 
 \item  $(\vec{u}-\vec{v}_i)^2 = -4$ for $i =2,3,4,5,6$, but 
       $(\vec{u}-\vec{v}_1)^2 = -6$.
\end{itemize}
The vectors of the first form have to be  
\begin{equation}
 \vec{u}_{k'} = \nu_\Omega -2\nu_{K_{k'}} + 4\nu_1  
\end{equation}
for a codeword $K_{k'} \in {\cal C}_{24}(8)$ satisfying 
\begin{equation}
 |K_{k'} \cap K_{U'}|=4, \qquad 1 \in K_{k'}, \quad {\rm and} \quad 
 K_{k'} \cap K_{L'} = \phi\, .
\label{eq:cond-E6-embed-type1-1-root-K}
\end{equation}
There are six codewords satisfying this set of conditions, and they 
are labelled by $k' \in K_{U'} \backslash \{0,\infty\} = \{ 11,2,22,20,14,21\}$.
The explicit form of those codewords is given in Appendix \ref{sssec:X3-detl}.
The vectors $\vec{u}_{k''} \in \Lambda_{24}$ of the second form have to be 
\begin{equation}
 \vec{u}_{k''} = \nu_\Omega -4\nu_{k''}
\end{equation}
for $k'' \in K_{U'} \backslash \{ 0, \infty\}$, so that there are also 
six choices of $k'' \in \{20,14,21,11,2,22\}$. Overall, there are twelve 
Leech roots of type 1), reproducing the result of \cite{Vin-2most}.
The [task II] is now completed.

Section 2 of \cite{Borcherds2} contains all the information 
necessary in carrying out the [task III].
For type 0) roots, $\sigma_{R'} \cdot \sigma_R$ preserves individual simple 
roots or $R=E_6$, and for type 1) roots, $\sigma_{R'} \cdot \sigma_R$ maps 
the simple roots of $R = \phi_{T_0}(J)$ to itself as a whole, but with 
the left-right flip. This also means that the assumption (as-4) is 
automatically satisfied in the case $X=X_3$.

The [task IV] is to compute the ${\rm Aut}(D')$ group, which 
is carried out by following the logic of the latter half of the proof 
of Lemma 4.5 in \cite{Kondo-auto-JacobKumm}.
For the argument in Lemma 4.5 of \cite{Kondo-auto-JacobKumm} to work, 
it is sufficient that the following two conditions are satisfied:
\begin{itemize}
 \item [({\bf as-6})] Leech roots of type 0) generate $S_X$ with 
       $\Z$ coefficients,
 \item [({\bf as-7})] $w'$---the orthogonal projection of the Weyl vector 
$w \in {\rm II}_{1,25}$ onto $S_X$---can be written in the form of 
a linear combination of the type 0) Leech roots that is manifestly 
invariant under the symmetry of the Coxeter graph of the type 0) Leech roots.
\end{itemize}
In the case of $X=X_3$, we will see explicitly in the 
appendix \ref{sssec:X3-detl} that these two conditions are satisfied. They also hold
in the case of $X={\rm Km}(A)$ studied in \cite{Kondo-auto-JacobKumm} and 
$X={\rm Km}(E\times F)$, ${\rm Km}(E \times E)$, 
${\rm Km}(E_\omega \times E_\omega)$ and ${\rm Km}(E_i \times E_i)$ 
in \cite{KK-auto-prod}.
Under the condition (as-6), the symmetry group ${\rm Aut}(\Gamma)$ of the 
Coxeter graph of the type 0) Leech roots $\Gamma$ generates 
isometries of $S_X$ and is identified with a subgroup of ${\rm Isom}(S_X)$.
The condition (as-7) further implies that it is a subgroup of 
${\rm Aut}(D')$. Conversely, any $g \in {\rm Aut}(D')$ maps 
type 0) Leech roots to themselves, inducing a symmetry transformation 
of the graph $\Gamma$. Therefore, 
\begin{equation}
{\rm Aut}(\Gamma) \cong {\rm Aut}(D').
\label{eq:AutGamma-AutD'}
\end{equation}

In the case of $X=X_3$, the symmetry of the Coxeter graph of the type 0) 
Leech roots, i.e., Figure \ref{fig:X3-Coxeter}, is 
$(S_3^{(\vartriangle)} \times S_3^{(\triangledown)}) \rtimes \Z_2\vev{\sigma_*}$, 
where 
$S_3^{(\vartriangle)}\vev{\tau_{\vartriangle}, \sigma_{\vartriangle}}$
is the permutation of the three trivalent vertices 
$\{\lambda_{22}, \lambda_6, \lambda_{15} \}$, and 
$S_3^{(\triangledown)}\vev{\tau_{\triangledown}, \sigma_{\triangledown}}$
that of the other three trivalent vertices 
$\{\lambda_{23}, \lambda_{12}, \lambda_{3} \}$. 
\begin{eqnarray}
 \tau_{\vartriangle}: & & 
 \lambda_{22} \rightarrow \lambda_{6} \rightarrow \lambda_{15} 
\rightarrow \lambda_{22}, \\
 \sigma_{\vartriangle}: & & \lambda_6 \leftrightarrow \lambda_{15}, 
    \qquad \lambda_{22} \quad {\rm remains~invariant}, \\
 \tau_{\triangledown}: & & 
\lambda_{23} \rightarrow \lambda_{12} \rightarrow \lambda_{3}
\rightarrow \lambda_{23}, \\
 \sigma_{\triangledown}: & & 
  \lambda_3 \leftrightarrow \lambda_{12}, 
    \qquad \lambda_{23} \quad {\rm remains~invariant} .
\end{eqnarray}
The generator $\sigma_*$ is the anti-podal transformation of the Coxeter 
diagram Figure \ref{fig:X3-Coxeter}; 
$\sigma_*: \lambda_{22} \leftrightarrow \lambda_{23}$, 
$\lambda_3 \leftrightarrow \lambda_{15}$ and 
$\lambda_6 \leftrightarrow \lambda_{12}$ \cite{Vin-2most}.
With (\ref{eq:AutGamma-AutD'}), now the [task IV] is completed. 

In order to confirm that the assumption (as-5) is satisfied, we
need to know the homomorphism $p_S: {\rm Aut}(D') \longrightarrow 
{\rm Isom}(q)$.In particular, we need to know the kernel of this 
map, which is denoted by ${\rm Aut}(D')_0$. 
Since ${\rm Aut}(D')_0$ is the subgroup of $Co_\infty$ preserving 
all $\vec{v}_{1,\cdots,6}$ 
in (\ref{eq:E6-embed-vect-1}--\ref{eq:E6-embed-vect-6}), 
any elements of ${\rm Aut}(D')_0$ should preserve $Z = \vec{0}$. 
This means that ${\rm Aut}(D')_0 \subset Co_0$.
Furthermore, because any element has to preserve $P'+Y-X=8\nu_0$, 
${\rm Aut}(D')_0 \subset 2^{12}.M_{24}$ (Ref. \cite{CS-book}, 
Chapt. 10 Thm.26). With a little more thought, 
one can see that ${\rm Aut}(D')_0 \subset M_{22}$. 
This great deal of simplification follows from choosing the vectors 
(\ref{eq:E6-embed-vect-1}--\ref{eq:E6-embed-vect-6}) so that 
$X', Y, Z$ and $P'$ are included. This method, introduced in 
\cite{Kondo-auto-JacobKumm} and exploited also in \cite{KK-auto-prod},
can be used for any $T_0$ containing $A_3 \oplus A_1$ as a sublattice. 

In the case of $X=X_3$, ${\rm Aut}(D')_0$ has to be a subgroup 
of $(S_6 \times S_6 \times S_9) \cap M_{22}$, where 
the first two $S_6$'s act as permutation on $K_{L'} \backslash \{0, \infty\}$
and $K_{U'} \backslash \{0,\infty\}$, and the last $S_9$ acts on 
$\Omega \backslash (K_{L'} \cup K_{U'} \cup \{1\})$. The ${\rm Aut}(D')_0$ 
group is also a subgroup of ${\rm Aut}(D') \cong {\rm Aut}(\Gamma)$. 
This is enough to determine ${\rm Aut}(D')_0$. It turns out that 
\begin{equation}
 \left( (\Z/3\Z) \vev{\tau_{\vartriangle}} \times 
        (\Z/3\Z) \vev{\tau_{\triangledown}} \right) \rtimes 
 (\Z/4\Z) \vev{\sigma_{\vartriangle} \cdot \sigma_*},  
\end{equation}
which is an index 2 subgroup of 
${\rm Aut}(D') \cong (S_3 \times S_3) \rtimes (\Z/2\Z)$.
It thus follows that 
\begin{equation}
p_S:  {\rm Aut}(D') \longrightarrow {\rm Aut}(D')/{\rm Aut}(D')_0 \cong 
\Z/2\Z
\end{equation}
is a surjective homomorphism to 
${\rm Isom}(q) \cong {\rm Aut}(J) \cong \Z/2\Z$. 
Now the assumption (as-5) is verified for $X=X_3$.

As a side remark, ${\rm Aut}(D')^{({\rm Hodge})}$ is also the same as 
${\rm Aut}(D')$, because $p_T:{\rm Isom}(T_X)^{({\rm Hodge})} \longrightarrow 
{\rm Isom}(q)$ is surjective in this example (see 
Table \ref{tab:coset-alt-rho20}). 

It is also possible to work out explicitly the subgroup of ${\rm Aut}(D')$
that preserves $R' = \phi_{T_0}(J) \cup \lambda$ for a given Leech root 
$\lambda$ of type 1). 
Let us take $\vec{u}_{11''}$ as an example. Such a subgroup should 
preserve $\{\lambda_{24}, \lambda_{25}, \lambda_{20} \}$ separately 
from all the other Leech roots of type 0), because these three 
type 0) Leech roots are not orthogonal to the type 1) Leech root 
associated with $\vec{u}_{11''}$.
The symmetry group of the three type 0) Leech roots above, and of 
the graph of the remaining 21 Leech roots is a diagonal subgroup of 
$S_3^{(\vartriangle)}$ and $S_3^{(\triangledown)}$. The coset space 
${\rm Aut}(\Gamma)/S_3$ contains 12 elements, which is the number of 
type 1 Leech roots. This calculation roughly reproduces the discussion
in \cite{Borcherds2}.

%%%%%%%%%%%%%%%%%%%%%%%%%%%%%%%%%%
\subsubsection{Supplementary Details }
\label{sssec:X3-detl}
%%%%%%%%%%%%%%%%%%%%%%%%%%%%%%%%%%

Here we present some details that have been omitted in 
the appendix \ref{sssec:X3-task}.

When the Neron--Severi lattice of $S_X$ for $X=X_3$ is realized as the 
orthogonal complement $[\phi_{T_0}(J)^\perp \subset {\rm II}_{1,25}]$ 
for the embedding using (\ref{eq:E6-embed-vect-1}--\ref{eq:E6-embed-vect-6}, 
\ref{eq:L24-Pi-relation}), all the Leech roots of type 0) are of the form 
$(2\nu_K, 1,1)$ with some codewords $K \in {\cal C}_{24}(8)$.  
There are 24 codewords satisfying the 
conditions (\ref{eq:cond-E6-embed-type0-root-K}), six of which 
(those in (\ref{eq:Leech-root-X3-triv-up},
\ref{eq:Leech-root-X3-triv-down}) corresponding to the trivalent
vertices in the Coxeter graph) have already been specified 
in Appendix \ref{sssec:X3-task}.
The remaining 18 codewords correspond to the remaining 18 type 0) Leech roots, 
which are located in between the trivalent vertices of the graph 
in Figure \ref{fig:X3-Coxeter}.

Let $K_{(n)}$ denote the codeword corresponding to the Leech root 
$\lambda_k = (2\nu_{K_{(k)}},1,1)$. The Leech roots in between 
the trivalent vertex $\lambda_{22}$ and $\lambda_{3}$ [resp. 
$\lambda_{12}$ or $\lambda_{23}$] are $\{\lambda_1, \lambda_2\}$ 
[resp. $\{\lambda_{10}, \lambda_{11} \}$ or 
$\{\lambda_{19}, \lambda_{20}\}$], for which the codewords are  
\begin{eqnarray}
 K_{(1)} = 
\begin{array}{|cc|cc|cc|}
\hline
 *& &*&{\color{white} *} & &{\color{white} *} \\
  & & & &*& \\
  &*&*& & & \\
 *&*& & &*& \\
\hline
	   \end{array}, & 
 K_{(10)} = 
\begin{array}{|cc|cc|cc|}
\hline
 *& &*&{\color{white} *} &{\color{white} *} & \\
  &*&*& & & \\
  & & & & &*\\
 *&*& & & &*\\
\hline
	   \end{array}, &
 K_{(19)} = 
\begin{array}{|cc|cc|cc|}
\hline
 *& &*& &{\color{white} *} &{\color{white} *} \\
  &*& &*& & \\
  &*& &*& & \\
 *& &*& & & \\
\hline
	   \end{array}, \\
K_{(2)} = 
\begin{array}{|cc|cc|cc|}
\hline
 *& &*& &*& \\
 *&*& & & & \\
 *& & &*& &*\\
  & & & & & \\
\hline
	   \end{array}, & 
K_{(11)} = 
\begin{array}{|cc|cc|cc|}
\hline
 *& &*& & &*\\
 *& & &*&*& \\
 *&*& & & & \\
  & & & & & \\
\hline
	   \end{array}, & 
K_{(20)} = 
\begin{array}{|cc|cc|cc|}
\hline
 *& &*&*& & \\
 *& & & & & \\
 *& & & & & \\
  &*& & &*&*\\
\hline
	   \end{array}. 
\end{eqnarray}
The Leech roots $\{ \lambda_5, \lambda_4 \}$, $\{ \lambda_8, \lambda_9 \}$ 
and $\{ \lambda_7, \lambda_{24} \}$ lie in between the trivalent vertex 
$\lambda_6$ and three other trivalent vertices $\lambda_3$, $\lambda_{23}$ 
and $\lambda_{12}$ (Figure \ref{fig:X3-Coxeter}), and the corresponding 
codewords are as follows:
\begin{eqnarray}
 K_{(5)} = 
\begin{array}{|cc|cc|cc|}
\hline 
*& &*&{\color{white} *} &{\color{white} *} & \\
 & & & & &*\\
*&*& & & &*\\
 &*&*& & & \\
\hline
\end{array}, & 
 K_{(8)} = 
\begin{array}{|cc|cc|cc|}
\hline 
*& &*&{\color{white} *} & &{\color{white} *} \\
 &*&*& & & \\
*&*& & &*& \\
 & & & &*& \\
\hline
\end{array}, & 
 K_{(7)} = 
\begin{array}{|cc|cc|cc|}
\hline 
*& &*& &{\color{white} *} &{\color{white} *} \\
 &*& &*& & \\
*& &*& & & \\
 &*& &*& & \\
\hline
\end{array}, \\
 K_{(4)} = 
\begin{array}{|cc|cc|cc|}
\hline 
*& &*& & &*\\
*&*& & & & \\
 & & & & & \\
*& & &*&*& \\
\hline
\end{array}, & 
 K_{(9)} = 
\begin{array}{|cc|cc|cc|}
\hline 
*& &*& &*& \\
*& & &*& &*\\
 & & & & & \\
*&*& & & & \\
\hline
\end{array}, & 
 K_{(24)} = 
\begin{array}{|cc|cc|cc|}
\hline 
*& &*&*& & \\
*& & & & & \\
 &*& & &*&*\\
*& & & & & \\
\hline
\end{array}.
\end{eqnarray}
Finally, for the type 0) Leech roots located in between the trivalent vertex 
$\lambda_{15}$ and three other trivalent vertices $\lambda_{12}$, 
$\lambda_{23}$ and $\lambda_3$, 
\begin{eqnarray}
 K_{(14)} = 
\begin{array}{|cc|cc|cc|}
\hline 
*& &*&{\color{white} *} & &{\color{white} *} \\
*&*& & &*& \\
 & & & &*& \\
 &*&*& & & \\
\hline
\end{array}, 
&
 K_{(17)} = 
\begin{array}{|cc|cc|cc|}
\hline 
*& &*&{\color{white} *} &{\color{white} *} & \\
*&*& & & &*\\
 &*&*& & & \\
 & & & & &*\\
\hline
\end{array}, 
&
 K_{(16)} =  
\begin{array}{|cc|cc|cc|}
\hline 
*& &*& &{\color{white} *} &{\color{white} *} \\
*& &*& & & \\
 &*& &*& & \\
 &*& &*& & \\
\hline
\end{array}, 
\\
 K_{(13)} = 
\begin{array}{|cc|cc|cc|}
\hline 
*& &*& &*& \\
 & & & & & \\
*&*& & & & \\
*& & &*& &*\\
\hline
\end{array}, 
& 
 K_{(18)} = 
\begin{array}{|cc|cc|cc|}
\hline 
*& &*& & &*\\
 & & & & & \\
*& & &*&*& \\
*&*& & & & \\
\hline
\end{array}, 
& 
 K_{(25)} = 
\begin{array}{|cc|cc|cc|}
\hline 
*& &*&*& & \\
 &*& & &*&*\\
*& & & & & \\
*& & & & & \\
\hline
\end{array}. 
\end{eqnarray}

Among the twelve Leech roots of type 1), those of the first type 
are given by the codewords $K_{k'}$ satisfying the 
conditions (\ref{eq:cond-E6-embed-type1-1-root-K}).
There are six of them: three are the codewords $K_{k'}$ for 
$k' \in \{11,2,22\}$, which are the collection of the two columns of 1 
and $k'$ within $\Omega$ in the MOG presentation (see below).
The remaining three codewords are $K_{k'}$ for $k=\{20,14,21\}$, which 
consist of the two rows of 1 and $k'$ in the centre and right blocks 
of the MOG presentation. $K_{k'=11'}$ and $K_{k'=20'}$, for example, are 
\begin{equation}
 K_{11'} = 
\begin{array}{|cc|cc|cc|}
  \hline
   {\color{white} *} &{\color{white} *} &*&*&{\color{white} *} &{\color{white} *} \\
    & &*&*& & \\
    & &*&*& & \\
    & &*&*& & \\
  \hline
\end{array}, \qquad 
 K_{20'} = 
\begin{array}{|cc|cc|cc|}
   \hline
    {\color{white} *} &{\color{white} *} &*&*&*&* \\
     & &*&*&*&* \\
     & & & & &  \\
     & & & & &  \\
   \hline
\end{array}\, .
\end{equation}
The $\sigma_{R'} \cdot \sigma_R$ transformation for a type 1) 
Leech root $\lambda$ of $X_3$ is contained in $W'_{\phi_{T_0}(J)}
\subset W_\Pi$ and is hence regarded as an isometry of $S_X$.
This isometry maps $\lambda'$ (the orthogonal projection onto 
$S_X \otimes \Q$) to $-\lambda'$. 
Thus it acts on $S_X$ as a reflection and can hence be regarded as 
a part of the generators of the reflection symmetry group $W(S_X)$. 
With the information above, it is straightforward to compute 
the Coxeter diagram of the reflection symmetry group generated by the 
type 0) and type 1) Leech roots. Such a computation reproduces the result 
of \cite{Vin-2most}. 

\begin{figure}[tbp]
\begin{center}
     \scalebox{.5}{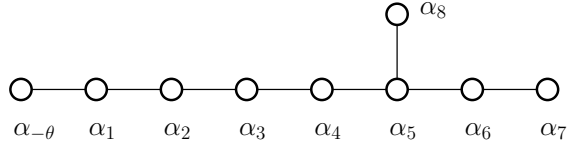} 
\caption{The affine Dynkin diagram of $E_8$.\label{fig:E8-Dynkin}} 
\end{center}
\end{figure}

Let us now verify that the conditions (as-6) and (as-7) in 
Appendix \ref{sssec:X3-task} are indeed satisfied in the example 
$X = X_3$. As for the condition (as-6), note first that 
$S_X \cong U^{(\gamma)} \oplus E_8^{(1)} \oplus E_8^{(2)} \oplus A_2$. 
Let $U^{(\gamma)} = {\rm Span}_\Z \{ u_1^{(\gamma)},
\bar{u}_1^{(\gamma)} \}$ with the symmetric pairing\footnote{ 
$\bar{u}_1^{(\gamma)}$ here corresponds to
$u_2^{(\gamma)}+u_1^{(\gamma)}$ in section \ref{ssec:systematic-23}.} 
$(u_1, \bar{u}_1) = 1$, $(u_1, u_1) = (\bar{u}_1, \bar{u}_1) = 0$. 
The simple roots of $E_8^{(a=1,2)}$ are denoted by
$\alpha_{1,\cdots,8}^{(a)}$ (see Figure~\ref{fig:E8-Dynkin}),
%%%%%%%%%%%%%%%%%%%%%%%%%%%%%%%%%%%%%%%%%%%%%%%%%%%%%%%%%%%%%%%%%%%%%%%%%
%%%%%%%%%%%%%%%%%%%%%%%%%%%%%%%%%%%%%%%%%%%%%%%%%%%%%%%%%%%%%%%%%%%%%%%%%
and the negative of the maximal root by $\alpha_{-\theta}^{(a)}$. 
The simple roots of $A_2$ are denoted by $\alpha_{-\theta}^{(3)}$ and 
$\alpha_1^{(3)}$, because this $A_2$ and an $E_6$ subalgebra generated 
by $\alpha_{i=3,\cdots,8}^{(3)}$ form a maximal subalgebra of $E_8^{(3)}$.
The twenty-four type 0) Leech roots ($(-2)$-roots of $S_X$)
$\lambda_{1,\cdots,\check{21},\cdots,25}$ are written in this basis 
of $S_X$ as \cite{Vin-2most}\footnote{
For the record, we also write down the three remaining Leech roots of
type 0).
\begin{eqnarray} 
&&
 \lambda_{23} \longleftrightarrow 2(u_1^{(\gamma)} + \bar{u}_1^{(\gamma)})
    + \alpha^{(3)}_{-\theta} 
  -(2\alpha_1+4\alpha_2+6\alpha_3+8\alpha_4+10\alpha_5+7\alpha_6+4\alpha_7+5\alpha_8)^{(1)}
  \nonumber \\
&& \qquad \qquad \qquad \qquad \qquad \qquad 
   -(2\alpha_1+4\alpha_2+6\alpha_3+8\alpha_4+10\alpha_5+7\alpha_6+4\alpha_7+5\alpha_8)^{(2)}, \\
&& \lambda_{24} \longleftrightarrow
 3(u_1^{(\gamma)}+\bar{u}_1^{(\gamma)})
   + 2\alpha^{(3)}_{-\theta} + \alpha^{(3)}_1
  -
  (3\alpha_1+6\alpha_2+9\alpha_3+12\alpha_4+15\alpha_5+10\alpha_6+5\alpha_7+8\alpha_8)^{(1)}
  \nonumber \\
&&  \qquad \qquad \qquad \qquad \qquad \qquad \quad   -
  (3\alpha_1+6\alpha_2+8\alpha_3+10\alpha_4+12\alpha_5+8\alpha_6+4\alpha_7+6\alpha_8)^{(2)}, \\
&& \lambda_{25} \longleftrightarrow 
 3(u_1^{(\gamma)}+\bar{u}_1^{(\gamma)})
   + 2\alpha^{(3)}_{-\theta} + \alpha^{(3)}_1
   -
  (3\alpha_1+6\alpha_2+8\alpha_3+10\alpha_4+12\alpha_5+8\alpha_6+4\alpha_7+6\alpha_8)^{(1)},
  \nonumber \\
&&  \qquad \qquad \qquad \qquad \qquad \qquad \quad   -
  (3\alpha_1+6\alpha_2+9\alpha_3+12\alpha_4+15\alpha_5+10\alpha_6+5\alpha_7+8\alpha_8)^{(2)}.
\end{eqnarray}
}
\begin{eqnarray}
&&
 \lambda_1 \longleftrightarrow u_1^{(\gamma)} + \alpha^{(1)}_{-\theta},
      \quad  
 \lambda_{2,\cdots,6} \longleftrightarrow \alpha^{(1)}_{1,\cdots,5}, 
      \quad 
 \lambda_{8,9} \longleftrightarrow \alpha^{(1)}_{6,7}, \quad 
 \lambda_{7} \longleftrightarrow \alpha^{(1)}_8, \\
&&
 \lambda_{10} \longleftrightarrow u_1^{(\gamma)} + \alpha^{(2)}_{-\theta},
      \quad 
 \lambda_{11,\cdots,15} \longleftrightarrow \alpha^{(2)}_{1,\cdots,5},
      \quad 
 \lambda_{17,18} \longleftrightarrow \alpha^{(2)}_{6,7}, \quad 
 \lambda_{16} \longleftrightarrow \alpha^{(2)}_{8}, \\
&&
 \lambda_{19} \longleftrightarrow u_1^{(\gamma)} + \alpha^{(3)}_{-\theta},
      \quad 
 \lambda_{20} \longleftrightarrow \alpha^{(3)}_1, \qquad \qquad 
 \lambda_{22} \longleftrightarrow \bar{u}_1^{(\gamma)} - u_1^{(\gamma)}. 
\end{eqnarray}
Thus, it is easy to see that all the generators 
$\alpha_{i=1,\cdots,8}^{(a=1,2)}$, $\alpha^{(3)}_{-\theta,1}$, 
$u_1^{(\gamma)}$ and $\bar{u}_1^{(\gamma)}$ of $S_X$ are written 
as linear combinations of the type 0) Leech roots with integer coefficients.
Now (as-6) is verified.

As for the condition (as-7), note that the Weyl vector 
$w :=(\vec{0},0,1) \in \Lambda_{24} \oplus U^{(\omega)} \cong 
{\rm II}_{1,25}$ is expressed in the language of 
\begin{equation}
 (\phi^{(\gamma)})^{-1}: 
 {\rm II}_{1,25} \cong U^{(\gamma)} \oplus L^{(\gamma)} = 
  U^{(\gamma)} \oplus E_8^{(1)} \oplus E_8^{(2)} \oplus E_8^{(3)}
\label{eq:II-1-25-gamma}
\end{equation}
as \cite{Borcherds-Leech}
\begin{equation}
 w = -\rho^{(\gamma)1}-\rho^{(\gamma)2}-\rho^{(\gamma)3}
  + 30 \bar{u}_1^{(\gamma)} + 31 u_1^{(\gamma)},
\end{equation}
where $\rho^{(\gamma)a}$ for $a=1,2,3$ are the Weyl vectors of 
the root lattice $E_8^{(a)}$. This means that 
$(\rho^{(\gamma)a}, \alpha_i^{(b)}) = - \delta_{ab}$ 
for any $i=1,\cdots,8$. To be more explicit,  
$\rho^{(\gamma)a} = (29\alpha_1 +57\alpha_2 +84\alpha_3 +110\alpha_4 
+135\alpha_5 +91\alpha_6 +46\alpha_7 +68\alpha_8)^{(a)}$. 
When $T_0 \cong E_6$ is embedded into ${\rm II}_{1,25}$ by mapping 
$\alpha_{1,\cdots,6} \in J$ to $\alpha_{3,\cdots,8}^{(3)}$, $w'$ (the
projection of the Weyl vector $w$ to $S_X \otimes \Q$) is given by 
\begin{equation}
 w' = w + 
  (8\alpha_3+15\alpha_4+21\alpha_5+15\alpha_6+8\alpha_7+11\alpha_8)^{(3)}.
\end{equation}
It is now a straightforward computation to see that this $w'$ is
equal to the sum of $4 \times (\lambda_{22}+\lambda_6+\lambda_{15}+
\lambda_{23}+\lambda_{12}+\lambda_3)$---a sum over all the type 0) Leech
roots at the trivalent vertices in the Coxeter diagram of
Figure~\ref{fig:X3-Coxeter}---and 
$3 \times (\lambda_1+\lambda_2+\lambda_4+\cdots )$---a sum over all the
other 18 Leech roots of type 0). This linear combination is manifestly 
invariant under the action of ${\rm Aut}(\Gamma) \cong (S_3 \times S_3)
\rtimes (\Z/2\Z)$, the symmetry of the graph in Figure \ref{fig:X3-Coxeter}.
Now (as-7) is also verified.

%%%%%%%%%%%%%%%%%%%%%%%%%%%%%%%%%%%%%%%%%%%%%%%%%%%%%%
\subsubsection{Elliptic Divisors for the Six Different Types on $X_3$}
%%%%%%%%%%%%%%%%%%%%%%%%%%%%%%%%%%%%%%%%%%%%%%%%%%%%%

The singular K3 surface $X_3$ admits 6 different types of elliptic
fibrations \cite{Nish01}. Due to Corollary D, they also form the modulo-isomorphism
classification of elliptic fibration on $X_3$.
We have learned from Lemma I and sections 
\ref{sssec:D4toDn-A2toOther}--\ref{sssec:D4A2toE8} that elliptic
divisors can be determined by eliminating some of the vertices from 
the Coxeter diagram of $W^{(2)}(S_X)$. Particularly in the case of 
$X = X_3$, it is obvious from Figure \ref{fig:X3-Coxeter} which 
vertex to eliminate. The result is recorded in the following. 

Type no.1 is associated with the embedding [$E_6 \subset E_8$; $E_8
\oplus E_8$], and $W_{\rm root} = A_2E_8E_8$. The elliptic divisors
is given by 
\begin{eqnarray}
u_1^{(\gamma)}&= & \lambda_1+ 2\lambda_2+ 3\lambda_3+ 4\lambda_4
 + 5\lambda_5+ 6\lambda_6+ 3\lambda_7+ 4\lambda_8+ 2\lambda_9,
   \qquad \qquad (E_8, {\rm II}^*) \\
&=& \lambda_{10}+2\lambda_{11}+3\lambda_{12}+4\lambda_{13}+5\lambda_{14}+
6\lambda_{15}+3\lambda_{16}+4\lambda_{17}+2\lambda_{18}.
   \quad (E_8, {\rm II}^*)
\end{eqnarray}
They correspond to two distinct singular fibres of type ${\rm II}^*$ and
are hence algebraically equivalent. 

The type no.2 is obtained by the embedding 
[$E_6 \subset E_8$; $D_{16}$], and hence $W_{\rm root}=A_2D_{16}$.
\begin{equation}
u_1^{(\beta)} = \lambda_2 + \lambda_{25} + 2(\lambda_{3,4,5,6,8,9} + 
  \lambda_{23} + \lambda_{18,17,15,14,13,12} )
  +\lambda_{24} + \lambda_{11}.  \quad (D_{16}, I_{12}^*)
\end{equation}

The type no.3 is for [$E_6 \subset E_7$; $D_7 \oplus E_7$], where 
$W_{\rm root}=D_{10}E_7$. The elliptic divisor is 
\begin{eqnarray}
 u_1^{(\eta)} & = & \lambda_7+\lambda_5 +
  2(\lambda_{6,8,9}+\lambda_{23}+\lambda_{18,17,15}) 
 + \lambda_{14} + \lambda_{16}, \qquad \quad (D_{10}, I_{6}^*) \\
  & = & \lambda_3+2\lambda_2+3\lambda_1 + 4\lambda_{22} + 
       3\lambda_{10}+2\lambda_{11}+\lambda_{12} + 2\lambda_{19}. \qquad 
    (E_7, {\rm III}^*)
\end{eqnarray}

The type no.4 is obtained by [$E_6 \subset E_7$; $A_{17}$], where 
$W_{\rm root}=A_{17}$. 
\begin{equation}
u_1^{(\zeta)} = \lambda_{22} + \lambda_{1,\cdots,6,8,9} + \lambda_{23} + 
  \lambda_{18,17,15,\cdots,10}.
\end{equation}

For type no.5 based on [$E_6 \subset E_6$; $E_6^{\oplus 3}$] with 
$W_{\rm root}=E_6^{\oplus 3}$,
\begin{eqnarray}
 u_1^{(\theta)}
   & = & 3 \lambda_{22}+ 2\lambda_{1,10,19} + \lambda_{2,11,20},
       \qquad \qquad \qquad \qquad \qquad \qquad (E_6, {\rm IV}^*)\\
   & = & 3 \lambda_{6} + 2\lambda_{5,8,7} + \lambda_{4,9,24}, 
       \qquad \qquad \qquad \qquad \qquad \qquad \qquad (E_6, {\rm IV}^*)\\
   & = & 3\lambda_{15} + 2\lambda_{14,17,16} + \lambda_{13,18,25}.
       \qquad \qquad \qquad \qquad \qquad \qquad (E_6, {\rm IV}^*)
\end{eqnarray}

Finally, for type no.6 from [$E_6 \subset E_6$; $A_{11}\oplus D_7$], 
the elliptic divisor is given by 
\begin{eqnarray}
u_1^{(\lambda)} & = & 
 \lambda_{1, \cdots, 6,8,9} + \lambda_{23} + \lambda_{20,19} +
 \lambda_{22},  \qquad \qquad \qquad \qquad \qquad (A_{11}, I_{12}) \\ 
 & = & \lambda_{24} + \lambda_{11} + 
  2 \lambda_{12,\cdots,15} + \lambda_{17} + \lambda_{16}.
  \qquad \qquad \qquad \qquad \quad  (D_7, I_3^*)
\end{eqnarray}
%

%
% \renewcommand{\thesubsection}{\thesection.C}
%
%%%%%%%%%%%%%%%%%%%%%%%%%%%%%%%%%%%%%%%%%%%%%%%%%%%%%%%%%%%%%%%%%%
\subsection{Kummer surfaces}
\label{ssec:review-Kummer}
%%%%%%%%%%%%%%%%%%%%%%%%%%%%%%%%%%%%%%%%%%%%%%%%%%%%%%%%%%%%%%%%%%

{\bf Definition} A {\it Kummer surface} is the minimal resolution of the sixteen 
$\C^2/(\Z/2\Z)$ singularity points in $A /[-{\bf id.}]$, where 
$A$ is an Abelian surface. Such a Kummer surface $X$ is denoted by 
${\rm Km}(A)$. 

%%%%%%%%%%%%%%%%%%%%%%%%%%%%%%%%%%%%%%%%%%%
\subsubsection{Abelian varieties}
%%%%%%%%%%%%%%%%%%%%%%%%%%%%%%%%%%%%%%%%%%%

{\bf Definition} An {\it Abelian variety} of dimension-$d$ is 
a complex $d$-dimensional torus that has a projective embedding. 
Abelian varieties of $d = 1$ are called elliptic curves, and 
those of $d=2$ Abelian surfaces.

A general treatment of Abelian varieties can be found in the classic textbook 
by Griffith and Harris \cite{Griffiths:1978}. Here, we only highlight those
aspects which are relevant to the further discussion.

\vspace{5mm}
A complex torus of dimension $d$ is given by $T =
\C^d/\Lambda$, where the lattice $\Lambda \subset \C^d$ is defined in terms of 
$\Lambda = {\rm Span}_Z \left\{ \ell_1, \ell_2, \cdots , \ell_{2d} \right\}$.
The generators of the lattice are specified by their coordinates in 
$\C^d$, $\ell_i = (z^1_i, z^2_i,\cdots, z^d_i)$.
These $\ell_i$ form a basis of $H_1(T; \Z)$.
Their dual basis in $H^1(T; \Z)$ in the form of $\R$-valued 1-forms 
are denoted by $\left\{ \lambda^i \right\}_{i=1,2,\cdots,2d}$ and satisfy
$\vev{\lambda^i, \ell_j} = \delta^i_{\; j}$.

Complex valued 1-forms $dz^a$ can be written as linear 
combinations of these 1-forms, $d z^a = \sum_i C^a_{\; i} \lambda^i$.
By taking appropriate linear combinations of the holomorphic coordinates 
and by properly changing the basis $\left\{ \lambda^i \right\}$ 
(and $\left\{\ell_i \right\}$ accordingly), on can always write
\begin{eqnarray}
 dz^1 & = & \lambda^1 + \tau^{11} \lambda^{d+1} + \cdots + \tau^{1d}\lambda^{2d}, \nn\\
 dz^2 & = & \lambda^2 + \tau^{21} \lambda^{d+1} + \cdots + \tau^{2d}\lambda^{2d}, \nn\\
& & \cdots \nn\\
 dz^d & = & \lambda^d + \tau^{d1} \lambda^{d+1} + \cdots + \tau^{dd}\lambda^{2d}. 
\end{eqnarray}

{\bf Theorem} The complex torus $T$ specified by the $\C$-valued 
$d\times d$ matrix $(\tau^{ab})$ is an Abelian variety (i.e. has 
a projective embedding) if and only if $\tau^{ab}$ is symmetric and 
the associated $\R$-valued matrix $H^{ab} = {\rm Im}(\tau^{ab})$ 
is positive definite. 

\vspace{5mm}

{\bf Example}:
In the one-dimensional case $\tau^{ab}$ is a $1\times 1$ matrix
and the condition to be an abelian variety is ${\rm Im}(\tau^{11}) > 0$.
This is the well-known case of an elliptic curve. 
Here, elliptic functions provide a projective embedding 
$T \ni u \mapsto [\wp(u): \wp'(u): 1] \in \P^2$. 

\vspace{5mm}

{\bf Example}:
Another class of examples is given by the Jacobian varieties. 
Consider a curve $\Sigma_g$ of genus $g$.
By integrating the $g$ independent harmonic $(1,0)$ forms over an element of $H_1(\Sigma_g,\Z)$ 
we obtain a vector in $\C^g$. The set of all such vectors defines a lattice $\Lambda\subset\C$.
${\rm Jac}(\Sigma_g)=\C/\Lambda$ is an example of an Abelian variety of dimension 
$d = g$. Once the complex structure of $\Sigma_g$ is given, the corresponding Abelian variety 
${\rm Jac}(\Sigma_g) = {\rm Pic}^0(\Sigma_g)$ is specified. While the moduli space of genus $g$ curves is 
of dimension $(3g-3)$ for $g \geq 2$, which grows only linearly in $g$, the moduli space of $d$-dimensional Abelian varieties has dimension 
$d(d+1)/2$, which grows quadratically in $d$. Thus ${\rm Jac}(\Sigma_g)$ form a very special class of Abelian varieties for large $g$.

\vspace{5mm}

In the case of $d=g=2$, however, any Abelian surface can be obtained 
as the Jacobian variety of a curve $\Sigma_2$ of genus two. In this case,
the moduli spaces of both the genus two curve and Abelian surfaces are 
three-dimensional. When we let $\tau^{12} = \tau^{21} = 0$
(while having ${\rm Im}(\tau^{11})>0$ and ${\rm Im}(\tau^{22})>0$), 
the Abelian variety becomes a product of elliptic curves $E\times F$. 
Such an Abelian variety corresponds to the Jacobian variety of 
$\Sigma_2$ that degenerates into the two elliptic curves $E$ and $F$. 

As Abelian varieties are obtained as the quotient of $\C^n$ by a
lattice, they inherit the additive group structure of $\C^n$. Hence, as
their name suggests, they are manifolds which are equipped with the 
structure of an abelian group. For one-dimensional abelian varieties,
this abelian group is nothing but the classic group law of elliptic curves.

%%%%%%%%%%%%%%%%%%%%%%%%%%%%%%%%%%%%%%%%%%%%%%%%%%%%%
\subsubsection{Algebraic Cycles in Abelian Surfaces}
%%%%%%%%%%%%%%%%%%%%%%%%%%%%%%%%%%%%%%%%%%%%%%%%%%%%%

The algebraic cycles of an Abelian surface form its Neron--Severi
lattice, $S_A = H^{1,1}(A; \R) \cap H^2(A; \Z)$. 
As for K3 surfaces, its rank is denoted by $\rho_A$. The minimal
resolution of $A /[-{\bf id.}]$ introduces $16$ independent classes, so
that the rank of $S_X$ for a Kummer surface $X = {\rm Km}(A)$ 
is $\rho_X = 16 + \rho_A$. 
By definition, Abelian varieties allow for a projective embedding, so
that there is at least one integral (1,1)-form on $A$ and we
generically have $\rho_X = 17$. For special choices of $A$, $\rho_A$
and hence the rank of $S_X$ can be enhanced, which is what we discuss
in the following.

In an Abelian variety $A$ of dimension $d=2$, we can take 
\begin{eqnarray}
  \left( e_1, \; e^1, \; e_2, \; e^2, \; e_3, \; e^3 \right)^T :=
   \left( \lambda^3 \wedge \lambda^2, \; \lambda^1 \wedge \lambda^4, \;
         \lambda^1 \wedge \lambda^2, \; - \lambda^3 \wedge \lambda^4, \;
         \lambda^1 \wedge \lambda^3, \; \lambda^2 \wedge \lambda^4 
  \right)^T
\label{eq:2form-basis-on-A}
\end{eqnarray}
as a basis of $H^2(A; \Z)$. Let us denote the dual basis of $H_2(A; \Z)$ by
\begin{equation}
  \left( C_{32}, C_{14}, C_{12}, C_{43}, C_{13}, C_{24} \right) \, .
\label{eq:2cycle-basis-on-A}
\end{equation}
Let us take $\lambda^1 \wedge \lambda^3 \wedge \lambda^2
\wedge \lambda^4$ as the orientation of the 4-cycle $[A]$, i.e. the volume form
on $A$. As $(e_i, e^j) = \delta_{i}^{\; j}$, the intersection form on $H_2(A; \Z)$ then becomes 
$U \oplus U \oplus U$ in this basis. 

These cycles $C_{ij}$ in $A$ are mapped into cycles $\bar{C}_{ij}$ of ${\rm Km}(A)$. 
Since $( \bar{C}_{ij}, \bar{C}_{kl}) = 2 (C_{ij}, C_{kl})$, their intersection form on ${\rm Km}(A)$ 
is $U[2] \oplus U[2] \oplus U[2]$. 

The period vector $\pi_A \in (H_2(A; \Z))^* \otimes \C $ is obtained by 
expanding $dz^1 \wedge dz^2$ in the basis (\ref{eq:2form-basis-on-A}):
\begin{eqnarray}
 dz^1 \wedge dz^2 & = & 
   \tau^{11} \lambda^3 \wedge \lambda^2
 + \tau^{22} \lambda^1 \wedge \lambda^4
 + \lambda^1 \wedge \lambda^2
 + {\rm det}(\tau^{ab}) \lambda^3 \wedge \lambda^4
 + \tau^{21} \lambda^1 \wedge \lambda^3
 - \tau^{12} \lambda^2 \wedge \lambda^4, \nonumber \nn\\
 &   \Updownarrow & \nn\\
 \pi_A  & = & ( \tau^{11}, \tau^{22}, 1, - {\rm det}(\tau), \tau^{21}, - \tau^{12}),
\end{eqnarray}
where the choice of basis in (\ref{eq:2form-basis-on-A}) is assumed in
the component description in the last line. 

Because $A$ is an Abelian surface (and not just a
complex torus), $\tau^{12} = \tau^{21}$. Hence the integral cycle 
$(C_{13}+ C_{24}) \longleftrightarrow (0,0,0,0,1,1)^T \in H_2(A; \Z)$ is 
part of $S_A$ for any Abelian surface, so that we generically have
$\rho_A = 1$.

The transcendental lattice of a generic Abelian surface $T_A$ is the orthogonal complement of 
$(C_{13}+C_{24})$ in $H_2(A; \Z)$, which is generated by the first four basis
vectors in (\ref{eq:2cycle-basis-on-A}) and $(C_{13}-C_{24})$.
Thus we find that $T_A = U \oplus U \oplus (-2)$. From this it follows that 
$T_X = U[2] \oplus U[2] \oplus (-4)$ for a Kummer surface $X = {\rm Km}(A)$ of 
a generic Abelian surface $A$.

\vspace{5mm}

When the Abelian surface has the special form, $A = E \times F$ (product
type) for some pair of elliptic curves $E$ and $F$, we can take 
$\tau^{12} = \tau^{21} = 0$. As each of the two elliptic curves $E$ and
$F$ allow for a projective embedding, $\rho_A=2$. This can be seen
explicitly by noting that both $C_{13}$ and $C_{24}$ are orthogonal to
the period vector $\pi_A$. Correspondingly, the transcendental lattices
of $A$ is $T_A = U \oplus U$, so that $\rho_X = 18$ and $T_X = U[2]
\oplus U[2]$ for the associated Kummer surface. 

\vspace{5mm}

Next we consider the case when the Abelian surface has the even more
special form of being the product of two mutually isogenous curves, 
$A=E\times E'$. Two elliptic curves $E$ and $E'$ are said to be {\it
isogenous} when their complex structure parameters $\tau$ and $\tau'$ 
in the upper half plane are related by a $GL(2; \Q)$ transformation, $\tau
= (a\tau'+b)/(c\tau'+d)$ where $a,b,c,d \in \Z$ do not necessarily
satisfy $ad - bc = 1$. In this case, we may exploit the modular group 
action on each of the two elliptic curves to choose $\tau^{11}=\tau$ and 
$\tau^{22}=\tau'$ such that 
\begin{equation}\label{condisogeny}
 \tau' = \frac{n_2}{n_1} \;  \tau \, , 
\end{equation}
with $n_1,n_2\in \Z$. In this case, also
$(n_2 C_{32} - n_1 C_{14}) \leftrightarrow (n_2,-n_1,0,0,0,0)^T$ 
is orthogonal to $\pi_A$ and hence algebraic. Now $\rho_A = 3$, and
hence $\rho_X = 19$. $T_A$ is generated by $n_2 C_{32}+ n_1 C_{14}$,
$C_{12}$ and $C_{43}$, and the intersection form is $U \oplus (2n_1n_2)$. Thus 
for the corresponding Kummer surface, $X = {\rm Km}(E \times E')$,  
$T_X = U[2] \oplus (4n_1n_2)$.

\vspace{5mm}

Finally, let us consider the case where $A = E\times E'$, with $E$ and
$E'$ isogenous, and furthermore suppose that they are mutually coprime 
integers $r,p,q$ such that $r \tau^2 + p \tau + q = 0$ holds for 
the complex structure parameter $\tau$ of $E$.\footnote{
Elliptic curves with this property are said to 
{\it have complex multiplication}. See \cite{Gukov:2002nw, Moore} and 
references therein for use of elliptic curves with complex
multiplication in string theory.}\raisebox{4pt}{,}\footnote{
One can also see from \eqref{condisogeny} that $E'$ also 
satisfies the same property when $E$ does, and vice versa.
} 
This implies that $\tau^2$ can be written as a $\Q$-coefficient linear
combination of $\tau$ and $1$. After using (\ref{condisogeny}), 
we may thus write
\begin{equation}
\pi_A = (\tau, n_2/n_1 \tau, 1, - n_2/n_1 \tau^2 , 0, 0) \, ,
\end{equation}
so that now 
$(pn_2 C_{32} + qn_2 C_{12} - rn_1 C_{43})  \leftrightarrow 
(p n_2,\; 0, \; q n_2, \; -rn_1, \; 0, \; 0)^T$ is an algebraic cycle. 
From this it follows that now $\rho_A=4$ and $\rho_X = 20$. 
The intersection form of $T_A$ can be obtained from this information, 
and that of $T_X$ for the corresponding Kummer surface 
is then again obtained by simply multiplying by 2.

Consider the example $\tau=\tau'=i$, for which $\tau^2 + 1 = 0$. 
In this case $C_{12}-C_{43}$ is an algebraic cycle of $A$, and the 
rank-2 transcendental lattice $T_A$ is generated by $C_{12}+C_{43}$ and 
$-(C_{32}+C_{14})$. $T_A = (2)\oplus (2)$, and $T_X = (4)\oplus (4)$.

As another example, consider $\tau=\tau'=\omega$, with 
$\omega=e^{2\pi i/3}$. Now we have $\tau^2 + \tau + 1 = 0$, 
$(C_{14}+C_{12}-C_{43})$ is an algebraic cycle of $A$, and the rank-2
lattice $T_A$ is generated by $C_{12}+C_{43}$ and
$-C_{32}-C_{14}+C_{43}$. 
Hence 
\begin{equation}
T_A = \left[ \begin{array}{cc} 2 & 1 \\ 1 & 2 \end{array}\right] \quad\mbox{and}\quad
T_X = \left[ \begin{array}{cc} 4 & 2 \\ 2 & 4 \end{array}\right] \, .
\end{equation}
%

%%%%%%%%%%%%%%%%%%%%%%%%%%%%%%%%%%%%%%%%%%%%%%%%%%%%%%%
\subsubsection{$T_X$, $S_X$ and their Symmetries for $\rho_X = 17$ 
Kummer Surfaces}
%%%%%%%%%%%%%%%%%%%%%%%%%%%%%%%%%%%%%%%%%%%%%%%%%%%%%%%%

In this section, we provide a brief summary of known results on the 
symmetries of the Neron--Severi and transcendental lattices of 
Kummer surfaces that we use in this article. 

Let us begin with the Kummer surface ${\rm Km}(A)$ of an Abelian surface 
$A$ with generic complex structure. As already stated above, 
the transcendental lattice $T_X$ is 
\begin{equation}
 T_X = U[2] \oplus U[2] \oplus (-4)
\end{equation}
generated by $\{ \overline{C}_{32}, \overline{C}_{14} \}$, 
$\{ \overline{C}_{12}, \overline{C}_{43} \}$ and 
$\overline{C}_{13}-\overline{C}_{24}$.
${\rm Isom}(T_X)^{({\rm Hodge})} = \{ \pm {\rm id.}\} \cong
\Z/2\Z\vev{-{\rm id.}}$.
It is much easier to compute the discriminant group $G_{T_X} \cong
G_{S_X}$ from $T_X$ than from $S_X$. As an Abelian group, 
\begin{equation}
\label{G_T_X_Kummer}
G_{T_X} \cong (\Z/2\Z)^4 \times (\Z/4\Z)   
\end{equation}
generated by $\overline{C}_{32}/2$, $\overline{C}_{14}/2$, 
$\overline{C}_{12}/2$, $\overline{C}_{43}/2$ and 
$(\overline{C}_{13}-\overline{C}_{24})/4$. Based on an explicit
computation, one can derive that 
\begin{equation}
 {\rm Isom}(q) \cong \Z/2\Z \vev{-_4} \times S_6 \, .
\end{equation}
Here, the $S_6$ acts on both the $(\Z/2\Z)^4$ and the $\Z/4\Z$ factor and $\vev{-_4}$ reverses the sign
of the $\Z/4\Z$. Since $p_T: (- {\rm id.}) \longmapsto -_4$, 
\begin{equation}
 p_T \left[{\rm Isom}(T_X)^{({\rm Hodge})}\right] \cong \Z/2\Z\vev{-_4}.
\label{eq:Tx-Hodge-image-JacobKumm}
\end{equation}

Let us now move on and provide a description of the Neron--Severi lattice 
$S_X$ for $X = {\rm Km}(A)$.
The involution acting on the Abelian surface has $16$ fixed points, which
give rise to $16$ singularities of type $A_1$. These fixed points are
located at the sixteen 2-torsion points\footnote{The 2-torsion points 
of an abelian variety $A$ are those points $p$ which give the unit element,
i.e. the origin, when added to themselves under the group law of $A$.} of $A$, 
\begin{equation}
  \mu_{\vec{s}} = 
    \frac{1}{2} \vec{s}\vec{\ell} = \sum_{i=1}^4 \frac{1}{2} s_i \ell_i
    \in A, 
\end{equation}
where $\vec{s} =(s_1,s_3,s_2,s_4) \in (\F_2)^4$. The exceptional cycles 
of $X={\rm Km}(A)$ obtained after the minimal resolution of the
$\C^2/(\Z/2\Z)$ singularities at these 2-torsion points are denoted by 
$G_{\vec{s}}$ for $\vec{s} \in (\F_2)^4$. They are algebraic and hence
elements of the Neron--Severi lattice $S_X$. 
There is one more independent divisor class in $S_X$, which corresponds 
to the algebraic cycle $(C_{13}+C_{24})$ in $A$ and is denoted by $H$.
$G_{\vec{s}}$'s and $H$ combined can be chosen as 
a set of $\Q$-coefficient generators of $S_X \otimes \Q$. The symmetric pairing 
on $S_X$ (and $S_X \otimes \Q$ also) is determined by 
$(G_{\vec{s}}, G_{\vec{s}'}) = -2\delta_{\vec{s},\vec{s}'}$, $(H, H)=4$, and 
$(H, G_{\vec{s}})=0$. The signature of $S_X$ is $(1,16)$.

Let $L$ denote (temporarily) the lattice generated by $H$ and $G_{\vec{s}}$; 
$L = {\rm Span}_\Z \{H, G_{\vec{s}}\}$. It must be an index $2^6$ sublattice 
of $S_X$, because ${\rm discr}(L)=2^{18}$, whereas 
${\rm discr}(S_X) = |{\rm discr}(T_X)| = 2^6$.
There are two equivalent ways to describe which elements of $G_L=L^*/L$ 
should be added to $L$ in order to obtain $S_X$, i.e. which 
rational linear combinations of the $G_{\vec{s}}$ and $H$ are integral
cycles of $X = {\rm Km}(A)$.

One way to do this is to say that the $(\Z/2\Z)^6 \subset G_L$ elements 
to be added are given by $\Z/2\Z$ generated by 
\begin{equation}
 T_{0000} = \frac{1}{2} \left(H 
      -G_{0000}-G_{1000}-G_{1100}-G_{0110}-G_{0111}-G_{0101} \right)
\label{eq:trope-T0-def}
\end{equation}
modulo $L$, and by $(\Z/2\Z)^5$ (\cite{Nik-Kummer} Cor.5) in the form of 
\begin{equation}
 [D^{(\bar{k})}] := \left[\frac{1}{2} \sum_{\vec{s} \in (\F_2)^4}
    (\vec{k} \cdot \vec{s} + k_0) G_{\vec{s}} \right] \in G_L.
\label{eq:Nik-Kummer-32elements}
\end{equation}
for any $\bar{k} := (\vec{k}, k_0) = (k_1,k_3,k_2,k_4,k_0) \in (\F_2)^5$. 

The other way to describe the elements of $L^*/L$ corresponding to 
$S_X/L$ deals with $T_{0000}$ in (\ref{eq:trope-T0-def}) in a way 
similar to how the sixteen $G_{\vec{s}}$'s are treated.
Let $I(0,0,0,0)$ denote the six elements $\vec{s} \in (\F_2)^4$ that
appear in the linear combination of (\ref{eq:trope-T0-def}).
Now we can rewrite (\ref{eq:trope-T0-def}) as follows:
\begin{equation}
 T_{0000} = \frac{1}{2} 
   \left(H - \sum_{\vec{s} \in I(0,0,0,0)} G_{\vec{s}} \right) .
\end{equation}
Similarly, we define fifteen other subsets of $(\F_2)^4$ by 
\begin{equation}
 I (\vec{s}) := \left\{ \vec{s} + \vec{t} \; | \vec{t} \in I(0,0,0,0) 
                \right\}, \qquad 
    {\rm for~} \vec{s} \in (\F_2)^4,
\end{equation}
and define fifteen other $\Q$-coefficient linear combinations of 
the generators of the lattice $L$ by:
\begin{equation}
  T_{\vec{s}} = \frac{1}{2} 
       \left( H - \sum_{\vec{t} \in I(\vec{s})} G_{\vec{t}} \right).
\end{equation}
The second description of the integral components of $S_X$ for 
$X = {\rm Km}(A)$ is to say that $S_X$ is generated by the sixteen 
$G_{\vec{s}}$'s and sixteen $T_{\vec{s}}$'s with $\Z$-coefficients; 
$G_{\vec{s}}$'s are called {\it nodes}, and $T_{\vec{s}}$'s {\it tropes}
in the literature.

The second description of $S_X$ is known to be equivalent to the first 
one using (\ref{eq:trope-T0-def}) 
and (\ref{eq:Nik-Kummer-32elements}). 
Symmetries can be made manifest in the second one, while practical
computation may be easier in the first one.

The discriminant group $G_{S_X} \cong S_X^*/S_X$ should be isomorphic 
to the one calculated from $T_X$. In the language of $S_X$, 
$G_{S_X} \cong (\Z/2\Z)^4 \times (\Z/4\Z)$, because 
$S_X \cong D_8 \oplus D_8 \oplus (4)$ \cite{Keum-auto-JacobKumm}.
A set of generators is explicitly given in \cite{Keum-auto-JacobKumm}:
\begin{eqnarray}
 B_1 = (G_{1100}+G_{0100}+G_{0110}+G_{1110})/2, &&
 B_2 = (G_{1000}+G_{0100}+G_{0111}+G_{1011})/2, \nonumber \\
 B_3 = (G_{1100}+G_{0100}+G_{0111}+G_{1111})/2, &&
 B_4 = (G_{1000}+G_{0100}+G_{0110}+G_{1010})/2, \nonumber \\
 B_0 = H/4 + (G_{0000}+G_{1000}+G_{1100}+G_{0100})/2. 
\label{eq:discr-grp-S-KummerGeneric}
\end{eqnarray}
In this basis, the discriminant bilinear form is
$b = \left[\begin{array}{cc} 0 & 1/2 \\ 1/2 & 0 \end{array} 
 \right]^{\oplus 2} \oplus (1/4)$, where diagonal (resp. off-diagonal) 
entries are evaluated mod $2 \Z$ (resp. $\Z$) \cite{Keum-auto-JacobKumm}.
The discriminant form is isomorphic to the one obtained from $T_X$, \eqref{G_T_X_Kummer}.

We will now describe a certain subgroup of the isometry group of this Neron--Severi 
lattice. 
First, there is a subgroup $(\Z/2\Z)^4 \subset {\rm Isom}(S_X)$ consisting of
translations. For any $\vec{s}_0 \in (\F_2)^4$, a translation is given by 
\begin{equation}
 \tau_{\vec{s}_0} : \left\{ \begin{array}{l}
      G_{\vec{s}} \longmapsto G_{\vec{s}+\vec{s}_0}, \\
      T_{\vec{s}} \longmapsto T_{\vec{s}+\vec{s}_0}, \\
      H \longmapsto H. 
      \end{array} \right.
\end{equation}
There is another subgroup $\Z/2\Z\vev{{\rm sw}} \subset {\rm Isom}(S_X)$, 
the generator of which, 
\begin{equation}
 {\rm sw}: \left\{ \begin{array}{l}
    G_{\vec{s}} \longmapsto T_{\vec{s}}, \\
    T_{\vec{s}} \longmapsto G_{\vec{s}}, \\
    H \longmapsto 3H - \sum_{\vec{s}\in (\F_2)^4} G_{\vec{s}}.
         \end{array} \right.
\end{equation}
is called switch. The switch and translations commute, and they form 
a $(\Z/2\Z)^5$ subgroup. 

%%%%%%%%%%%%%%%%%%%%%%%%%%%%%%%%%%%%%%%%%%%%%%%%%%%%%%%%%%%%%%%%%

\begin{figure}[tbp]
\begin{center}
     \scalebox{.5}{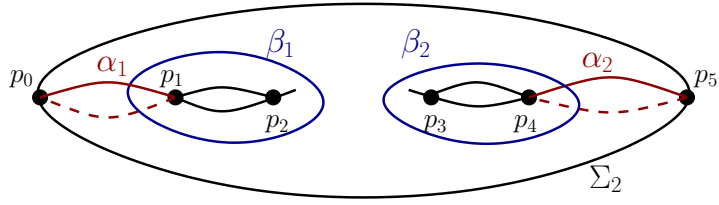} 
\caption{A curve of genus 2 with a symplectic choice of basis for 
$H_1(\Sigma_2;\Z)$. We have furthermore displayed the locations of 
the branch points (Weierstrass) $p_a$ which occur in the realization 
of the genus two curve as a hyperelliptic curve.
\label{fig:Sigma_2}} 
\end{center}
\end{figure}

Another subgroup $S_6 \subset {\rm Isom}(S_X)$ is understood better, 
when we introduce a different notation that involves more geometric 
aspects of Kummer surfaces. For this reason, we take a moment here 
to make a little digression (see \cite{Keum-auto-JacobKumm}). 
As we have already remarked, Abelian surfaces can be realized 
as Jacobians of a genus 2 curve. 
Choose a curve of genus 2 $\Sigma_{g=2}$ appropriately, so that 
$A = {\rm Jac}(\Sigma_{g=2})$. All genus 2 curves are known to be 
hyperelliptic, and we take its expression of the form 
$y^2 = \prod_{a=0}^5 (x-x_a)$.
The points $p_a = \{(x,y)=(x_a,0) \in \Sigma_2\}$ for $a=0,\cdots,5$ are 
called {\it Weierstrass points}. The
Abel--Jacobi map 
\begin{equation}
\mu:  \Sigma_2 \ni p \longmapsto (z^1, z^2) =
  \left( \int_{p_0}^p \omega^1,  \int_{p_0}^p \omega^2 \right) 
   \in \C^2 / {\rm Span}_\Z \{ \ell_1, \ell_2, \ell_3, \ell_4 \} \in A
\end{equation}
sends the six Weierstrass points $p_a$'s to six of the sixteen 
2-torsion points of $A$. Here, $\omega^1$ and $\omega^2$ are 
independent holomorphic 1-forms on $\Sigma_2$ normalized so that 
$(\omega^i,\alpha_j) = \delta^i_{\; j}$ for 1-cycles $\alpha_1$ and 
$\alpha_2$. Two other independent 1-cycles of $\Sigma_2$, $\beta_1$ and 
$\beta_2$ are chosen so that 
$(\alpha_i, \beta_j)= - (\beta_j,\alpha_i) = \delta_{ij}$, 
$(\alpha_i, \alpha_j)=0$ and $(\beta_i, \beta_j) = 0$, see
figure \ref{fig:Sigma_2}. 
%%%%%%%%%%%%%%%%%%%%%%%%%%%%%%%%%%%%%%%%%%%%%%%%%%%%%%%%%%%%%%%%%%%%%%%%%
%%%%%%%%%%%%%%%%%%%%%%%%%%%%%%%%%%%%%%%%%%%%%%%%%%%%%%%%%%%%%%%%%%%%%%%%%
It is easy to see from the figure that the six 2-torsion points 
are as follows:\footnote{They are contained in the image of the 
Abel--Jacobi map, $\mu(\Sigma_2) \subset A$. These six 2-torsion points 
precisely correspond to the exceptional cycles that appeared in the
linear combination of (\ref{eq:trope-T0-def}). $T_{0000}$ in $S_X$
corresponds to the curve $\mu(\Sigma_2) \subset A$. 
Similarly, $T_{\vec{s}}$'s for other $\vec{s} \in (\F_2)^4$ in $S_X$
correspond to curves in $A$ obtained by shifting $\mu(\Sigma_2)$ 
by $\mu_{\vec{s}}$.
}
\begin{eqnarray}
 p_0 \mapsto \int_{p_0}^{p_0} \omega = \mu_{0000} =: \mu_{00}, \quad &
 p_1 \mapsto \int_{p_0}^{p_1} \omega = \mu_{1000} =: \mu_{01}, & \quad 
 p_2 \mapsto \int_{p_0}^{p_2} \omega = \mu_{1100} =: \mu_{02}, \nonumber \\
 p_3 \mapsto \int_{p_0}^{p_3} \omega = \mu_{0110} =: \mu_{03}, \quad &
 p_4 \mapsto \int_{p_0}^{p_4} \omega = \mu_{0111} =: \mu_{04}, & \quad 
 p_5 \mapsto \int_{p_0}^{p_5} \omega = \mu_{0101} =: \mu_{05}. \nonumber 
\end{eqnarray}
%
%%%%%%%%%%%%%%%%%%%%%%%%%%%%%%%%%%%%%%%%%%%%%%%%%%%%%%%%%%%%%%
The remaining ten 2-torsion points correspond to period integrals 
$\int_{p_a}^{p_b} \omega$, so that they can 
be denoted also by $\mu_{ab}$ using the indices 
$a,b = 1,\cdots, 5$, $a \neq b$. Note that $\mu_{ab} = \mu_{ba}$ in $A$. 
Since one can also regard $\mu_a$'s for $a=1,\cdots,5$ as 
$\mu_{0a}=\mu_{a0}$, and $\mu_0$ as $\mu_{00}=\mu_{11} = \cdots =
\mu_{55}$, all the sixteen 2-torsion points can be labelled by using 
$(ab)$ with $a,b = 0,\cdots, 5$.
The sixteen nodes $X={\rm Km}(A)$ associated with 
$\mu_{\vec{s}}$ for $\vec{s} \in (\F_2)^4$---$G_{\vec{s}}$'s---can 
also be labelled by a pair of $a, b$; 
the following notation is introduced (as in \cite{Keum-auto-JacobKumm}):
\begin{align}
\label{NodesFromalphabeta}
&N_{00}   = G_{0000} \quad &N_{01}&   = G_{1000} \quad &N_{02}&   = G_{1100}    \quad&N_{03}&= G_{0110} \nn \\
&N_{04}   = G_{0111} \quad &N_{05}&   = G_{0101} \quad &N_{12}&= G_{0100} \quad&N_{13}&=G_{1110}\nn\\
&N_{14}= G_{1111} \quad &N_{15}&= G_{1101} \quad &N_{23}&= G_{1010} \quad&N_{24}&=G_{1011} \nonumber \\
&N_{25}= G_{1001} \quad &N_{34}&=G_{0001}  \quad &N_{35}&= G_{0011}
 \quad&N_{45}&=G_{0010}  % \nonumber 
\end{align}
The sixteen tropes $T_{\vec{s}}$ of $X = {\rm Km}(A)$ can also be
relabelled by $a,b=0, \cdots, 5$ by using the same dictionary between 
$\vec{s} \in (\F_2)^4$ and $(ab)$; 
$T_{00}=T_{0000}$, $T_{01}=T_{1000}$, $T_{12}=T_{0100}$ etc.
With this new notation, the $S_6$ subgroup of ${\rm Isom}(S_X)$ acts 
as\footnote{The $S_6$ subgroup of ${\rm Isom}(S_X)$ also acts as
permutation on the sixteen tropes, but not in a way that is as simple 
as on the the sixteen nodes. $S_6$ is generated by binary permutations 
of the form $(a,b)$ with $1 \leq a < b \leq 5$ and $(0,a)$ with $a \neq
0$, and under such permutations, 
\begin{eqnarray}
 (a,b): & &
 T_{0a} \longleftrightarrow T_{0b}, \qquad 
 T_{ac} \longleftrightarrow T_{bc}, \qquad 
 T_{00}, T_{0c}, T_{cd} \quad {\rm inv.}, \\
 (0,a): & & 
 T_{00} \longleftrightarrow T_{0a}, \qquad 
 T_{bc} \longleftrightarrow T_{de}, \qquad
 T_{0b}, T_{ab}, \quad {\rm inv.}, 
\end{eqnarray}
where $1\leq a,b,c,d,e \leq 5$ are all different. 
} 
\begin{equation}
 S_6 \ni \sigma: \left\{ 
   \begin{array}{l}
       N_{ab} \longmapsto N_{\sigma(a)\sigma(b)} \\
       H \longmapsto H,
   \end{array}
\right.
\end{equation}
where $S_6$ acts as permutation on the indices $a,b \in \{0,\cdots, 5\}$.

The $(\Z/2\Z)^5$ and $S_6$ subgroups of ${\rm Isom}(S_X)$ as a whole
form a subgroup $(\Z/2\Z)^5 \rtimes S_6$, or, equivalently 
$(\Z/2\Z)^5 : S_6$. This subgroup of ${\rm Isom}(S_X)$ is characterized 
as the one acting as permutation on the nodes ($G_{\vec{s}}$'s) and 
tropes ($T_{\vec{s}}$'s) while preserving intersection numbers among 
them \cite{Nik-analogue-Kumm}.
Reference \cite{Kondo-auto-JacobKumm} Lemma 4.5 further proved that this 
$(\Z/2\Z)^5 \rtimes S_6$ subgroup of ${\rm Isom}(S_X)$ does correspond 
to the ${\rm Aut}(D') \subset {\rm Isom}(S_X)^{({\rm Amp})}$ associated 
with the realization of $S_X$ for $X = {\rm Km}(A)$ using the embedding 
of $(T_0)_{\rm root} = A_3 \oplus A_1^{\oplus 6}$ into 
${\rm II}_{1,25}$.

The homomorphism 
$p_S: {\rm Isom}(S_X) \longrightarrow {\rm Isom}(G_{S_X},q)$ maps 
this subgroup ${\rm Aut}(D') \cong (\Z/2\Z)^5 \rtimes S_6$ to 
${\rm Isom}(G_{S_X},q)$. It is known that the kernel is the 
$(\Z/2\Z)^4$ subgroup corresponding to translations, 
$p_S({\rm sw}) = (-_4) \in {\rm Isom}(q)$, and the $S_6$ subgroup of 
${\rm Aut}(D')$ becomes the $S_6$ subgroup of ${\rm Isom}(q)$.
In particular, $p_S: {\rm Aut}(D') \longrightarrow {\rm Isom}(q)$ 
is surjective.\footnote{ 
The switch and $S_6$ do not commute in ${\rm Isom}(S_X)$, but images under $p_S$ 
do commute in ${\rm Isom}(q)$.}
Within the ${\rm Aut}(D')$ subgroup of ${\rm Isom}(S_X)^{({\rm Amp})}$, 
only the $(\Z/2\Z)^5$ subgroup---translations and switch---fall into 
the image of ${\rm Isom}(T_X)^{({\rm Hodge})}$, 
(\ref{eq:Tx-Hodge-image-JacobKumm}), so that they can be lifted to 
automorphisms of $X = {\rm Km}(A)$.

%%%%%%%%%%%%%%%%%%%%%%%%%%%%%%%%%%%%%%%%%%%%%%%%%%%%%%%
\subsubsection{The Neron--Severi Lattice of 
$X=X_{[2~2~2]} = {\rm Km}(E_\omega \times E_\omega)$}
\label{sssec:NS-for-KmEwEw}
%%%%%%%%%%%%%%%%%%%%%%%%%%%%%%%%%%%%%%%%%%%%%%%%%%%%%%%%

In this appendix, we give a minimal description 
of the Neron--Severi lattice of 
$X = {\rm Km}(E_\omega \times E_\omega) = X_{[2~2~2]}$ so notations 
to be used in the main text are explained. For a more systematic 
study, see \cite{KK-auto-prod}.

Let us begin with a more general K3 surface, $X = {\rm Km}(E \times F)$, 
where the Abelian surface is of product type $A = E \times F$, yet 
the elliptic curves 
$E = \left(\C /{\rm Span}_\Z \{ \ell_1, \ell_3 \} \right) = \C/(\Z+\tau
\Z)$ and 
$F = \left(\C / {\rm Span}_\Z \{ \ell_2, \ell_4 \} \right) = 
\C/(\Z+\tau'\Z)$ are not necessarily isogenous and their 
complex structure parameters $\tau, \tau'$ are generic. 

Let $\{p_1, \cdots, p_4\}$ be the four 2-torsion points in $E$, and 
$\{q_1, \cdots, q_4\}$ be those in $F$:
\begin{eqnarray}
 && p_1 = [0] \in E, \quad 
  p_2 = [\ell_1/2] \in E, \quad 
  p_3 = [\ell_3/2] \in E, \quad 
  p_4 = [\ell_1/2 + \ell_3/2] \in E, \\
 && q_1 = [0] \in F, \quad 
  q_2 = [\ell_2/2] \in F, \quad 
  q_3 = [\ell_4/2] \in F, \quad 
  q_4 = [\ell_2/2 + \ell_4/2] \in F,
\end{eqnarray}
$E_j$ for $j=1,\cdots,4$ are $(-2)$ curves of $X = {\rm Km}(E \times F)$ 
corresponding to $E \times q_j$, and $F_i$ for $j=1,\cdots, 4$ to 
$p_i \times F$. The exceptional curve at $(p_i,q_j)$ is now denoted by 
$G_{ij}$ instead of $G_{\vec{s}}$; $G_{0000} = G_{11}$, $G_{1000}=G_{21}$, 
$G_{0001}=G_{13}$ and $G_{0111}=G_{34}$, for example. 
Under the symmetric pairing of $S_X$, 
\begin{equation}
 (G_{ij}, G_{kl}) = -2 \delta_{ik} \delta_{jl}, \qquad 
  (G_{ij}, E_l) = \delta_{jl}, \qquad 
  (F_k, G_{ij}) = \delta _{ik}. 
\end{equation}
All of 
\begin{equation}
 2E_j + \sum_{i=1}^4 G_{ij}
\label{eq:lin-equiv-Es}
\end{equation}
are equivalent in $S_X$, as they correspond to the four singular fibres (all
of type $I_0^*$) of one of the two Kummer pencils 
(elliptic fibrations)\footnote{When $F_1$ is chosen as the zero section, 
$W_{\rm root} = \oplus_j D_4^{(j)} = \oplus_j {\rm Span}_\Z\{E_j,
G_{2j}, G_{3j}, G_{4j}\}$, and 
$MW = W_{\rm frame}/W_{\rm root} \cong \Z/2\Z\vev{[F_2-F_1]} \times 
\Z/2\Z \vev{[F_3-F_1]}$.} 
$\pi_X: {\rm Km}(E \times F) \longrightarrow F/(\Z/2\Z) \simeq \P^1$. 
For the other Kummer pencil, 
$\pi_X: {\rm Km}(E \times F) \longrightarrow E/(\Z/2\Z) \simeq \P^1$, 
all of the four singular fibres of type $I_0^*$, 
\begin{equation}
 2F_i + \sum_{j=1}^4 G_{ij}
\label{eq:lin-equiv-Fs}
\end{equation}
are similarly equivalent in $S_X$. 

The 16 + 4 + 4 = 24 $(-2)$ curves, 
$\{ G_{ij}, E_j, F_i \}$, can be chosen as a set of integer-coefficient 
generators of $S_X$ for $X = {\rm Km}(E \times F)$ with
non-isogenous elliptic curves $E$ and $F$ and generic complex structure parameters $\tau,
\tau'$, where the symmetry of $S_X$ is manifest. The set of 24 curves, 
however, is redundant as a $\Q$-coefficient basis of the rank-18 $S_X$.
We can choose $\{G_{ij}, E_1, F_1\}$ as a basis of $S_X$, though some of the 
symmetries of $S_X$ are not manifest in this basis. 

\vspace{.5cm}
Let us now turn to the case $X = {\rm Km}(E_\omega \times E_\omega)$.
This is a special case of $X = {\rm Km}(E \times F)$, 
in that $E$ and $F$ are isogenous and both have the complex structure 
$\omega=e^{2\pi i/3}$, so that the Picard number is enhanced to 20. 

All of the twenty-four $(-2)$-curves $\{G_{ij}, E_j, F_i\}$ remain in the $S_X$ lattice. 
As we consider the case where $E$ and $F$ are isogenous curves, however, 
the ``diagonal subset'' of $A = E \times F$ and its images under
the translations are also algebraic cycles of $A$, and their images 
in $S_X$ also remain algebraic.
Denoting the coordinates of % the two elliptic curves 
$A = E_\omega \times E_\omega$ by $(z^1,z^2) \in 
\C/(\Z + \omega \Z) \times \C/(\Z+\omega\Z)$, 
we embed a torus [resp. rational curve] 
(with complex coordinate $z$ [resp. $x/(-)$]) in $A$ [resp. $X = A/(-)$] 
by mapping $z$ to
\begin{align}
 D_1 :& \quad (z,z) &  D_2 :& \quad (\tfrac12+ z,z) \nn\\
 D_3 :& \quad (\tfrac12\omega+z,z) & D_4 :& \quad (\tfrac12(1+\omega)+z,z) \, .
\end{align}
Furthermore, due to the special complex structure $\tau = \omega$, we
may also embed rational curves (by exploiting complex multiplication) as
\begin{align}
C_1:&\quad (\omega z,-\omega^2 z) &  C_2: & \quad (\omega(z-\tfrac12),-\omega^2 z)  \nn \\
C_3:&\quad (\omega (z-\tfrac12\omega),-\omega^2 z) &  C_4: & \quad (\omega(z-\tfrac12(1+\omega)),-\omega^2 z)  \, .
\end{align}

The thirty-two $(-2)$ curves of $X = {\rm Km}(E_\omega \times E_\omega)$ 
are grouped into $\{G_{ij}\}$ and $\{ C_i, D_i, E_i, F_i\}$;  
there are no mutual intersection between the sixteen elements of 
$\{C_i,D_i,E_i,F_i\}$ and also between the sixteen elements
$\{G_{ij}\}$.  Just like each one of $E_i$'s or $F_i$'s finds four 
$(-2)$ curves (and precisely four) within the sixteen $\{ G_{ij}\}$'s 
to intersect, each one of $C_i$'s and $D_i$'s also have four of 
$\{ G_{ij}\}$'s to intersect:
\begin{align}
 C_1 \cdot G_{11}, G_{23}, G_{34}, G_{42} = 1\,, &\quad\quad D_1\cdot G_{11}, G_{22}, G_{33}, G_{44}=1, \nn \\
 C_2 \cdot G_{14}, G_{22}, G_{31}, G_{43} = 1\,, &\quad\quad D_2\cdot
 G_{12}, G_{21}, G_{34}, G_{43}=1,  \nn \\
 C_3 \cdot G_{12}, G_{24}, G_{33}, G_{41} = 1\,, &\quad\quad D_3\cdot G_{13}, G_{31}, G_{24}, G_{42}=1, \nn \\
 C_4 \cdot G_{13}, G_{21}, G_{32}, G_{44} = 1\,, &\quad\quad D_4\cdot
 G_{14}, G_{41}, G_{23}, G_{32}=1 \, .   
 \end{align}
All other intersection numbers between $(C_i, G_{kl})$ and $(D_i, G_{kl})$ vanish.

These thirty-two $(-2)$ curves are known to be a set of
integer-coefficient generators of the Neron--Severi lattice $S_X$. 
This set respects the symmetry of $S_X$, but is redundant as a
$\Q$-basis for the rank-20 lattice $S_X$. As a such a basis, we can take 
$\{G_{ij}, C_1, D_1, E_1, F_1\}$, for example. This is because 
there are linear equivalence relations such as 
\begin{equation}
 2D_1 + (G_{11}+G_{22}+G_{33}+G_{44}) \sim 
 2D_2 + (G_{21}+G_{12}+G_{34}+G_{43}) 
\end{equation}
for $C_i$'s and $D_i$'s, just like those for $E_i$'s and $F_i$'s in 
(\ref{eq:lin-equiv-Es}, \ref{eq:lin-equiv-Fs}).

\vspace{5mm}

Let us now describe the identification between $S_X$ for 
$X = {\rm Km}(E_\omega \times E_\omega)$ and the orthogonal complement 
of a primitive embedding $\phi_{T_0}: T_0 = D_4 \oplus A_2
\hookrightarrow {\rm II}_{1,25}$ defined 
by (\ref{eq:D4-A2-embed-II-1-25-1}--\ref{eq:D4-A2-embed-II-1-25-q}, 
\ref{eq:Ku}, \ref{eq:L24-Pi-relation}). {\it Type 0) Leech root} is 
a Leech root $\lambda \in \Pi$ that is orthogonal to any one of 
$\{ \lambda_{1,2,3,4}, \lambda_p, \lambda_q\}$, 
which are the $\phi_{T_0}$-images of the simple roots of $D_4 \oplus A_2$.
Under the choice of the 
embedding (\ref{eq:D4-A2-embed-II-1-25-1}--\ref{eq:D4-A2-embed-II-1-25-q}), 
all the type 0) Leech roots are of the form $\lambda = (2\nu_K,1,1)$
for a codeword $K \in {\cal C}_{24}(8)$, and the codeword $K$ satisfies 
the following conditions \cite{KK-auto-prod}:
\begin{equation}
 \{0, \infty \} \subset K \quad {\rm but} \quad 2 \nin K, \quad 
{\rm and~either} 
  \left\{
\begin{array}{ll}
   1 \in K & {\rm and~}|K \cap K_U|=4, \; {\rm or} \\
   1 \nin K & {\rm and~} |K \cap K_U|=2.
\end{array}
  \right.
\label{eq:type0-cond--D4A2-embed}
\end{equation}

Under the choice of $K_U$ in (\ref{eq:Ku}), there are sixteen codewords 
$K \in {\cal C}_{24}(8)$ of the first type, $1 \in K$ and $|K \cap K_U|=4$.
They are constructed from the sextet decompositions in
(\ref{eq:five-sextets}) in the form of $K_{\alpha a}$, $K_{\beta a}$, 
$K_{\gamma a}$ and $K_{\delta a}$ with 
$a \in \{ \heartsuit, \clubsuit, \spadesuit, \diamondsuit \}$.
We identify the sixteen type 0) Leech roots with the $\{ G_{ij}
\}$'s. To be more explicit, the correspondence between the sixteen 
$\{G_{ij}\}$'s and the sixteen codewords are as follows:
\begin{eqnarray}
&& G_{24} \leftrightarrow K_{\alpha \clubsuit}, \qquad 
 G_{43} \leftrightarrow K_{\alpha \spadesuit}, \qquad 
 G_{11} \leftrightarrow K_{\alpha \heartsuit}, \qquad 
 G_{32} \leftrightarrow K_{\alpha \diamondsuit},\nn \\ 
&& G_{13} \leftrightarrow K_{\beta \clubsuit}, \qquad 
 G_{41} \leftrightarrow K_{\beta \spadesuit}, \qquad
 G_{34} \leftrightarrow K_{\beta \heartsuit}, \qquad 
 G_{22} \leftrightarrow K_{\beta \diamondsuit},\nn \\
&& G_{42} \leftrightarrow K_{\gamma \clubsuit}, \qquad
 G_{14} \leftrightarrow K_{\gamma \spadesuit}, \qquad 
 G_{33} \leftrightarrow K_{\gamma \heartsuit}, \qquad 
 G_{21} \leftrightarrow K_{\gamma \diamondsuit}, \nn\\
&& G_{44} \leftrightarrow K_{\delta \clubsuit}, \qquad 
 G_{31} \leftrightarrow K_{\delta \spadesuit}, \qquad 
 G_{12} \leftrightarrow K_{\delta \heartsuit}, \qquad 
 G_{23} \leftrightarrow K_{\delta \diamondsuit}.
 \label{gijvsK}
\end{eqnarray}
There are also sixteen codewords satisfying the second type of
condition---$1 \nin K$ and 
$|K \cap K_U| = 2$---in (\ref{eq:type0-cond--D4A2-embed}). Under the 
choice of $K_U$ in (\ref{eq:Ku}), they are\\

\begin{minipage}{2.5cm}
\begin{align}
K_p \nn 
\end{align}
\end{minipage}
\begin{minipage}{4cm}
\begin{align}
\begin{array}{|cc|cc|cc|}
\hline
 * &  * & {\color{white} *} & {\color{white} *}  & {\color{white} *}  & {\color{white} *}  \\
  * &  * &  & & &\\
  * &  * &  & & &\\
  * &  * &  & & &\\
\hline
\end{array} \nn
\end{align}
\end{minipage}\hspace{1.2cm}
\begin{minipage}{2.5cm}
\begin{align}
K_{q1} \nn 
\end{align}
\end{minipage}
\begin{minipage}{4cm}
\begin{align}
\begin{array}{|cc|cc|cc|}
\hline
 * &   *& {\color{white} *} & {\color{white} *}  & {\color{white} *}  & {\color{white} *}  \\
   * &   * &  & & &\\
&&   * &   * &  & \\
&&   * &   * &  & \\
\hline
\end{array} \nn
\end{align}
\end{minipage}
%%%%%%%%%%%%%%%%%%%%%%%%%%%%%%%%%%%%%%%%%%%%%%%%%%%%%%%%%%%%%%%%

%%%%%%%%%%%%%%%%%%%%%%%%%%%%%%%%%%%%%%%%%%%%%%%%%%%%%%%%%%%%%%%%
\begin{minipage}{2.5cm}
\begin{align}
K_{q\omega} \nn 
\end{align}
\end{minipage}
\begin{minipage}{4cm}
\begin{align}
\begin{array}{|cc|cc|cc|}
\hline
 * &  * & {\color{white} *} & {\color{white} *}  & {\color{white} *}  & {\color{white} *}  \\
 &  &  *  &  * & &\\
 * & *&   &  & &\\
    &   & * & *& &\\
\hline
\end{array} \nn
\end{align}
\end{minipage}\hspace{1.2cm}
\begin{minipage}{2.5cm}
\begin{align}
K_{q\bar{\omega}} \nn 
\end{align}
\end{minipage}
\begin{minipage}{4cm}
\begin{align}
\begin{array}{|cc|cc|cc|}
\hline
 * &  * & {\color{white} *} & {\color{white} *}  & {\color{white} *}  & {\color{white} *}  \\
 &  &  *  &  * & &\\
  & & *  &  *& &\\
   * &  * &  & & &\\
\hline
\end{array} \nn
\end{align}
\end{minipage}
%%%%%%%%%%%%%%%%%%%%%%%%%%%%%%%%%%%%%%%%%%%%%%%%%%

%%%%%%%%%%%%%%%%%%%%%%%%%%%%%%%%%%%%%%%%%%%%%%%%%%%
\begin{minipage}{2.5cm}
\begin{align}
K_{r\bar{\omega}} \nn 
\end{align}
\end{minipage}
\begin{minipage}{4cm}
\begin{align}
\begin{array}{|cc|cc|cc|}
\hline
 * &  * & {\color{white} *} & {\color{white} *}  & {\color{white} *}  & {\color{white} *}  \\
  * &   &   *  & &  *  &\\
   &   *   &  *   & &  *  &\\
   &  &  & & &\\
\hline
\end{array} \nn
\end{align}
\end{minipage}\hspace{1.2cm}
\begin{minipage}{2.5cm}
\begin{align}
K_{\bar{r}\bar{\omega}} \nn 
\end{align}
\end{minipage}
\begin{minipage}{4cm}
\begin{align}
\begin{array}{|cc|cc|cc|}
\hline
 * &   * & {\color{white} *} & {\color{white} *}  & {\color{white} *}  & {\color{white} *}  \\
  &   * &   &  * &  * &\\
 *  &  &  & * & * &\\
   &  &  & &  &\\
\hline
\end{array} \nn
\end{align}
\end{minipage}

\begin{minipage}{2.5cm}
\begin{align}
K_{r \omega} \nn 
\end{align}
\end{minipage}
\begin{minipage}{4cm}
\begin{align}
\begin{array}{|cc|cc|cc|}
\hline
 * &  * & {\color{white} *} & {\color{white} *}  & {\color{white} *}  & {\color{white} *}  \\
  &  * & *  & & * &\\
  &  &  & & &\\
  * &  &  * & &  * &\\
\hline
\end{array} \nn
\end{align}
\end{minipage}\hspace{1.2cm}
\begin{minipage}{2.5cm}
\begin{align}
K_{\bar{r}\omega} \nn 
\end{align}
\end{minipage}
\begin{minipage}{4cm}
\begin{align}
\begin{array}{|cc|cc|cc|}
\hline
 *&   * & {\color{white} *} & {\color{white} *}  & {\color{white} *}  &  {\color{white}*}  \\
* &   &     & *&  *  &\\
    &    &     & &   &\\
    & * &   & * &* &\\
\hline
\end{array} \nn
\end{align}
\end{minipage}

\begin{minipage}{2.5cm}
\begin{align}
K_{r 1} \nn 
\end{align}
\end{minipage}
\begin{minipage}{4cm}
\begin{align}
\begin{array}{|cc|cc|cc|}
\hline
 * &  * & {\color{white} *} & {\color{white} *}  & {\color{white} *}  & {\color{white} *}  \\
   &   &     & &    &\\
   * &   &  *   & &  *  &\\
   &  * & *  & &  *&\\
\hline
\end{array} \nn
\end{align}
\end{minipage}\hspace{1.2cm}
\begin{minipage}{2.5cm}
\begin{align}
K_{\bar{r}1} \nn 
\end{align}
\end{minipage}
\begin{minipage}{4cm}
\begin{align}
\begin{array}{|cc|cc|cc|}
\hline
  * &   * & {\color{white} *} & {\color{white} *}  & {\color{white} *}  & {\color{white} *}  \\
  & & & & & \\
  &  *  &   & * &  * &\\
  * &  &  &  *&  *&\\
\hline
\end{array} \nn
\end{align}
\end{minipage}
%%%%%%%%%%%%%%%%%%%%%%%%%%%%%%%%%%%%%%%%%%%%%%%%%%%%%%%%%%

%%%%%%%%%%%%%%%%%%%%%%%%%%%%%%%%%%%%%%%%%%%%%%%%%%%%%%%%
\begin{minipage}{2.5cm}
\begin{align}
K_{\bar{s}1} \nn 
\end{align}
\end{minipage}
\begin{minipage}{4cm}
\begin{align}
\begin{array}{|cc|cc|cc|}
\hline
 * &   * & {\color{white} *} & {\color{white} *}  & {\color{white} *}  & *  \\
  &    &   &   &  * &\\
   &  *&  & * &  &\\
   & * &  *& &  &\\
\hline
\end{array} \nn
\end{align}
\end{minipage}\hspace{1.2cm}
\begin{minipage}{2.5cm}
\begin{align}
K_{s1}\nn 
\end{align}
\end{minipage}
\begin{minipage}{4cm}
\begin{align}
\begin{array}{|cc|cc|cc|}
\hline
 *&   * & {\color{white} *} & {\color{white} *}  & {\color{white} *}  &  *  \\
 &   &     & &  *  &\\
 *   &    &  *   & &   &\\
  *  &  &   & * & &\\
\hline
\end{array} \nn
\end{align}
\end{minipage}

\begin{minipage}{2.5cm}
\begin{align}
K_{s \omega} \nn 
\end{align}
\end{minipage}
\begin{minipage}{4cm}
\begin{align}
\begin{array}{|cc|cc|cc|}
\hline
  * &   * & {\color{white} *} & {\color{white} *}  & {\color{white} *}  &  *\\
   * &   &    & *   &    &\\
   &      &     & &   *  &\\
 *   &  &  * & & &\\
\hline
\end{array} \nn
\end{align}
\end{minipage}\hspace{1.2cm}
\begin{minipage}{2.5cm}
\begin{align}
K_{\bar{s}\omega} \nn 
\end{align}
\end{minipage}
\begin{minipage}{4cm}
\begin{align}
\begin{array}{|cc|cc|cc|}
\hline
 * &   * & {\color{white} *} & {\color{white} *}  & {\color{white} *}  &  *  \\
   &   * &   *   & &    &\\
   &   &     & &   *  &\\
   &   * &   & * &  &\\
\hline
\end{array} \nn
\end{align}
\end{minipage}

\begin{minipage}{2.5cm}
\begin{align}
K_{s\bar{\omega}} \nn 
\end{align}
\end{minipage}
\begin{minipage}{4cm}
\begin{align}
\begin{array}{|cc|cc|cc|}
\hline
 * &   * & {\color{white} *} & {\color{white} *}  & {\color{white} *}  &  *  \\
*   &    &   *   & &    &\\
*   &   &     & *&     &\\
   &    &   &  &  *&\\
\hline
\end{array} \nn
\end{align}
\end{minipage}\hspace{1.2cm}
\begin{minipage}{2.5cm}
\begin{align}
K_{\bar{s}\bar{\omega}} \nn 
\end{align}
\end{minipage}
\begin{minipage}{4cm}
\begin{align}
\begin{array}{|cc|cc|cc|}
\hline
  * &   * & {\color{white} *} & {\color{white} *}  & {\color{white} *}  &  * \\
   & *   &     &  *&    &\\
    &    *&   *   & &    &\\
   &   &   & &   *&\\
\hline
\end{array} \nn
\end{align}
\end{minipage}

\vspace{.5cm}

We set up the identification between the sixteen $(-2)$-curves $\{C_i, D_i, E_i, F_i\}$ in $S_X$ and 
these sixteen codewords as follows:
\begin{eqnarray}
&& E_1 \leftrightarrow K_{p}, \qquad 
 F_1 \leftrightarrow K_{q\bar{\omega}}, \qquad 
 C_1 \leftrightarrow K_{q1}, \qquad 
 D_1 \leftrightarrow K_{q\omega}, \nn\\ 
&& E_{2} \leftrightarrow K_{r\bar{\omega}}, \qquad 
 F_{2} \leftrightarrow K_{s1}, \qquad
 C_{2} \leftrightarrow K_{s\omega}, \qquad 
 D_{2} \leftrightarrow K_{\bar{r}\omega},\nn \\
&& E_{3} \leftrightarrow K_{r\omega}, \qquad
 F_{3} \leftrightarrow K_{\bar{r}\bar{\omega}}, \qquad 
 C_{3} \leftrightarrow K_{\bar{r}1}, \qquad 
 D_{3} \leftrightarrow K_{\bar{s}1}, \nn\\
&& E_{4} \leftrightarrow K_{r1}, \qquad 
 F_{4} \leftrightarrow K_{\bar{s}\omega}, \qquad 
 C_{4} \leftrightarrow K_{\bar{s}\bar{\omega}}, \qquad 
 D_{4} \leftrightarrow K_{s\bar{\omega}}.
 \label{efcdvsK}
\end{eqnarray}
\noindent
As a check, we can recover the intersection form by noting that two Leech roots associated to an
octad have intersection number one if they share two entries in the MOG and intersection number zero
if they share four entries in the MOG. 

As explained in section \ref{ssec:KmEoxEo-multiplicity}, we can embed 
the Neron--Severi lattice into ${\rm II}_{1,25}$ and obtain the 
root lattice $D_4\oplus A_2$ as its orthogonal complement.
Restricting the closure of the fundamental chamber $\overline{C_\Pi}$ of 
the reflection group of ${\rm II}_{1,25}$ to the positive cone of $S_X$, 
we find $D'$ which is bounded by a finite number of faces. 
The group of automorphisms of of $D'$ has the following 
structure:
\begin{equation}
 {\rm Aut}(D') \cong (\Z/2\Z)^4 . A_4 . (S_3\times \Z/2\Z) \, .
\end{equation}
The kernel of $p_S: {\rm Aut}(D') \longrightarrow {\rm Isom}(q) \cong
(S_3 \times \Z/2\Z)$ is 
\begin{equation}
 {\rm Aut}(D')_0 \cong (\Z/2\Z)^4 . A_4, 
\end{equation}
and $p_S$ is surjective \cite{KK-auto-prod}.

In the case of $X = {\rm Km}(E_\omega \times E_\omega)$, 
tasks I--IV and the verification of assumptions (as-1)--(as-7) were carried
out in \cite{KK-auto-prod}, specifically in 
Lemma 3.1 (task I), Lemma 3.2 (task II type 0)), 
Lemma 3.6 (task II type non-0), Lemma 3.4 for (as-3), 
\S 4.3 (task III and (as-4)), Lemma 2.1 for (as-6), 
Lemma 3.4 for (as-7), and Lemma 3.5 (task IV). 

Since we would like to use the action of ${\rm Aut}(D')$ on $S_X$
and compute the orbit of elliptic divisors under the action of this
group, let us extract more details of the action of the group 
from \cite{KK-auto-prod}.

As before, we may consider translations on $A = E_\omega \times
E_\omega$ by any of the two-torsion points $\mu_{ij}$ which generate 
the normal subgroup $(\Z/2\Z)^4$ of ${\rm Aut}(D')$. 
These will not mix the $C,D,E,F,G$ but leave each class of cycles
separate. We use the following names for the translations:
\begin{align}
 t_1 \sim \mu_{21}&=(\tfrac12,0)\nn\\
 t_2 \sim \mu_{31}&=(\tfrac{\omega}{2},0)\nn\\
 t_3 \sim \mu_{12}&=(0,\tfrac12)\nn\\
 t_4 \sim \mu_{13}&=(0,\tfrac{\omega}{2})
\end{align}
They act in the following way on the generators of the Neron--Severi lattice:
\begin{align}
t_1: &&&& \quad t_2: && \nn\\ 
&E_i \leftrightarrow E_i & & C_1 \leftrightarrow C_4 & &E_i \leftrightarrow E_i & & C_1 \leftrightarrow C_2\nn \\
&&&                          C_2 \leftrightarrow C_3 & &&& 			      C_3 \leftrightarrow C_4 \nn \\ 
&F_1 \leftrightarrow F_2 & & D_1 \leftrightarrow D_2 & &F_1 \leftrightarrow F_3 & & D_1 \leftrightarrow D_3 \nn \\ 
&F_3 \leftrightarrow F_4 & & D_3 \leftrightarrow D_4 & &F_2 \leftrightarrow F_4 & & D_2 \leftrightarrow D_4 \nn \\
&G_{1j}\leftrightarrow G_{2j} & & G_{3j} \leftrightarrow G_{4j} & & G_{1j}\leftrightarrow G_{3j} & & G_{2j} \leftrightarrow G_{4j}\nn \\ 
&&&&&&\nn\\
t_3: &&&& t_4: && \nn\\ 
&F_i \leftrightarrow F_i & & C_1 \leftrightarrow C_3 & &F_i \leftrightarrow F_i & & C_1 \leftrightarrow C_4\nn \\
&&&                          C_2 \leftrightarrow C_4 & &&& 			      C_2 \leftrightarrow C_3 \nn \\ 
&E_1 \leftrightarrow E_2 & & D_1 \leftrightarrow D_2 & &E_1 \leftrightarrow E_3 & & D_1 \leftrightarrow D_3 \nn \\ 
&E_3 \leftrightarrow E_4 & & D_3 \leftrightarrow D_4 & &E_2 \leftrightarrow E_4 & & D_2 \leftrightarrow D_4 \nn \\
&G_{i1}\leftrightarrow G_{i2} & & G_{i3} \leftrightarrow G_{i4} & & G_{i1}\leftrightarrow G_{i3} & & G_{i2} \leftrightarrow G_{i4} 
\end{align}
\noindent
Using the embedding of the Neron--Severi lattice into Leech roots, 
we may realize these maps as transformations in 
$M_{24} \subset Co_0 \cong {\rm Isom}(\Lambda_{24})$. 
% automorphisms of the Leech lattice. 
Written in MOG form they are given by
\begin{align}
t_1: &  &  t_2:&\nn \\
& \begin{array}{|cc|cc|cc|}
\hline
  &   &  & a &  & b  \\
c  &    &   e  & f & e   & h\\
d  & c   &     & b &     & a \\
   &  d  & g  & f  &  g & h\\
\hline
\end{array} 
&&
\begin{array}{|cc|cc|cc|}
\hline
  &   &  & g &  & e  \\
  & a & d& f & d& h \\
b &   &c & f & c & h \\
a &b  &  & e &  & g\\
\hline
\end{array} \nn \\
t_3: &  &  t_4:&\nn \\
& \begin{array}{|cc|cc|cc|}
\hline
  &   &  & g &  &   \\
e  & d  & b &  & c   & h\\
c  & e   &a &  &d & h \\
 a  & b  & f  &   &  f & g\\
\hline
\end{array} 
&&
\begin{array}{|cc|cc|cc|}
\hline
  &   &  & f &  & d  \\
e  & & b& a & c& f \\
c & h  & & a &  & g \\
 &b  & h & d & e & g\\
\hline
\end{array}
\end{align}
Here, any two entries which share a letter are exchanged under the $\Z/2\Z$, all others are left fixed.

The action by the alternating group $A_4$ is generated\footnote{
When we see $A_4$ as a subgroup of $S_4
\vev{\sigma_{12}, \sigma_{23}, \sigma_{34}}$, it is 
\begin{equation}
 A_4 = \langle a := (\sigma_{12} \cdot \sigma_{23}), \; b := (\sigma_{23}
  \cdot \sigma_{34}) \; | \; a^3 = b^3 = 1, \; (ab)^2 = 1 \rangle .
\end{equation}
Each of the generators $g_1$ and $g_2$ also satisfy $(g_1)^3 = {\rm
id.}$ and $(g_2^{-1})^3 = {\rm id.}$ on $A/(-{\rm id.})$, and
furthermore, $(g_2^{-1} \cdot g_1): (z^1, z^2) \longmapsto 
(\omega z^1 +z^2, \omega(z^1-z^2))$ is an involution, 
$(g_2^{-1} \cdot g_1)^2 = {\rm id}.$
Thus, $g_1$ and $g_2^{-1}$ satisfy the same set of relations 
as the generators $a$ and $b$ of the group $A_4$. } 
by the automorphisms on $A/(-id.)$:
\begin{align}
 g_1:\quad &(z^1,z^2) \mapsto (z^2,-z^1+z^2) \nn \\
 g_2:\quad&(z^1,z^2) \mapsto (\omega(z^2-z^1),-\omega^2 z^2)
\end{align}
Given the explicit expression for the various algebraic curves, we can
work out the action of $g_1$: $(i=2,3,4)$
\begin{align}           
\sta{G_{i1}}{G_{1i}}{G_{ii}}&  & \sta{G_{23}}{G_{34}}{G_{42}}& & \sta{G_{24}}{G_{43}}{G_{32}}    \nn \\
&&&& \nn \\
\sta{C_{1}}{C_{1}}{C_{1}}&  & \sta{C_{2}}{C_{4}}{C_{3}}&  & \sta{E_{i}}{F_{i}}{D_{i}}    
\end{align}
Similarly, the action of $g_2$ is 
\begin{align}           
\sta{G_{11}}{G_{11}}{G_{11}}&  & \sta{G_{22}}{G_{14}}{G_{23}}&  & \sta{G_{33}}{G_{12}}{G_{34}} \nn \\
&&&&\nn\\
\sta{G_{44}}{G_{13}}{G_{42}}&  & \sta{G_{21}}{G_{31}}{G_{41}}&  & \sta{G_{32}}{G_{24}}{G_{43}} \nn \\  
&&&&\nn\\
\sta{E_{1}}{E_{1}}{E_{1}}&  & \sta{E_{2}}{E_{4}}{E_{3}}&  & \sta{F_{1}}{C_{1}}{D_{1}}    \nn \\
&&&&\nn\\
\sta{F_{2}}{C_{2}}{D_{4}}& &  \sta{F_{3}}{C_{3}}{D_{2}} & &  \sta{F_{4}}{C_{4}}{D_{3}}
\end{align}
% these generate $G_0$
%
This $A_4$ subgroup of ${\rm Aut}(D')_0$ can also be realize 
as transformations in $M_{24}$. In the MOG representation, 
\begin{align}
g_1: &  &  g_2:&\nn \\
& \begin{array}{|cc|cc|cc|}
\hline
  &   &  &  &  & d_1  \\
b_1  & a_2  & e_2   & c_2 &  & f_1\\
e_3  & c_3   &b_2   & a_3 & d_3    & f_2 \\
 c_1 & e_1  & a_1  & b_3  & d_2 & f_3 \\
\hline
\end{array} 
&&
\begin{array}{|cc|cc|cc|}
\hline
  &   &  &  &  &  \\
a_1  & b_1 & c_1 & d_1 & e_1 & f_1 \\
a_2 & b_2  &c_2 & d_2 & e_2 & f_2 \\
a_3 &b_3  & c_3 & d_3 & e_3 & f_3\\
\hline
\end{array} \nn
\end{align}
Here, the maps $g_1$ and $g_2$ act by sending $\xi_1\rightarrow \xi_2 \rightarrow \xi_3 \rightarrow \xi_1$ for
any letter.
\\
\\
The remaining maps in ${\rm Aut}(D')$ are 
$S_3 \times \Z/2\Z \cong {\rm Aut}(D')/{\rm Aut}(D')_0$. 
Let $S_3 = \langle s_1, s_2 \; | (s_1)^3 = 1, \; (s_2)^2 = 1, \; 
s_2 \cdot s_1 \cdot s_2 = (s_1)^2 \rangle$.
The order-3 element $s_1$ can be represented by an automorphism 
\begin{equation}
s_1 :\quad (z^1,z^2) \mapsto (\omega z^1, \omega z^2) \, ,
\end{equation}
which acts on $\{ C_i,D_i,E_i,F_i,G_{ij}\}$ such that it leaves the 
index $1$ fixed and cyclically permutes all other indices:
\begin{equation}
 \sta{2}{3}{4} \, .
\end{equation}
The other generator $s_2$ is represented by an isometry of $S_X$, 
\begin{align}
s_2: & \nn \\
     & E_2 \leftrightarrow E_3 &  & F_3 \leftrightarrow F_4  & &  G_{12} \leftrightarrow G_{13} & & G_{22} \leftrightarrow G_{23} \nn\\
     & G_{32} \leftrightarrow G_{43} &  & G_{33} \leftrightarrow G_{42}  & &  G_{31} \leftrightarrow G_{41} & & G_{34} \leftrightarrow G_{44} \nn\\
     & C_1 \leftrightarrow D_1 &  & C_2 \leftrightarrow D_4  & &  C_3 \leftrightarrow D_{3} & & C_{4} \leftrightarrow D_{2} \, ,
\end{align}
which leaves $E_1,E_4,F_1,F_2,G_{11},G_{14},G_{21},G_{24}$ fixed.
\\
\noindent
Let $\sigma$ be the generator of the remaining $\Z/2\Z$ factor; 
it acts on $S_X$ as follows:
\begin{align}
\sigma: & \nn \\
&E_1 \leftrightarrow G_{24} &  & E_2 \leftrightarrow G_{43} & & E_3 \leftrightarrow G_{32} & & E_4 \leftrightarrow G_{11} \nn \\
&F_1 \leftrightarrow G_{44} &  & F_2 \leftrightarrow G_{31} & & F_3 \leftrightarrow G_{23} & & F_4 \leftrightarrow G_{12} \nn \\
&D_1 \leftrightarrow G_{14} &  & D_2 \leftrightarrow G_{42} & & D_3 \leftrightarrow G_{21} & & D_4 \leftrightarrow G_{33} \nn \\ 
&C_1 \leftrightarrow G_{34} &  & C_2 \leftrightarrow G_{22} & & C_3 \leftrightarrow G_{41} & & C_4 \leftrightarrow G_{13}.  
\end{align}

Since the image of ${\rm Isom}(T_X)^{({\rm Hodge})} \cong \Z/6\Z\vev{\theta_6}$
under $p_T$ is $\Z/2\Z \vev{\sigma} \times \Z/3\Z \vev{s_1} \subset  
{\rm Isom}(q)$, the subgroup ${\rm Aut}(D')^{({\rm Hodge})}$ is given 
by \cite{KK-auto-prod} 
\begin{equation}
  {\rm Aut}(D')^{({\rm Hodge})} \cong (\Z/2\Z)^4 . A_4. 
     (\Z/3\Z\vev{s_1} \times \Z/2\Z \vev{\sigma}).
\end{equation}
%

%%%%%%%%%%%%%%%%%%%%%%%%%%%%%%%%%%%%%%%%%%%%%%%%%%%%%%%%%%%%%%%%

\end{document}

%% file: dnchain.pdf_t
\begin{picture}(0,0)%
\includegraphics{dnchain.pdf}%
\end{picture}%
\setlength{\unitlength}{4144sp}%
\begingroup\makeatletter\ifx\SetFigFont\undefined%
\gdef\SetFigFont#1#2#3#4#5{%
  \reset@font\fontsize{#1}{#2pt}%
  \fontfamily{#3}\fontseries{#4}\fontshape{#5}%
  \selectfont}%
\fi\endgroup%
\begin{picture}(5403,3150)(598,-4396)
\put(5761,-2401){\makebox(0,0)[lb]{\smash{{\SetFigFont{20}{24.0}{\rmdefault}{\mddefault}{\updefault}{\color[rgb]{0,0,0}$Y$}%
}}}}
\put(3961,-3391){\makebox(0,0)[lb]{\smash{{\SetFigFont{20}{24.0}{\rmdefault}{\mddefault}{\updefault}{\color[rgb]{0,0,0}$X$}%
}}}}
\put(3061,-3391){\makebox(0,0)[lb]{\smash{{\SetFigFont{20}{24.0}{\rmdefault}{\mddefault}{\updefault}{\color[rgb]{0,0,0}$\vec{v}_5$}%
}}}}
\put(2161,-3391){\makebox(0,0)[lb]{\smash{{\SetFigFont{20}{24.0}{\rmdefault}{\mddefault}{\updefault}{\color[rgb]{0,0,0}$\vec{v}_6$}%
}}}}
\put(1261,-3391){\makebox(0,0)[lb]{\smash{{\SetFigFont{20}{24.0}{\rmdefault}{\mddefault}{\updefault}{\color[rgb]{0,0,0}$\vec{v}_7$}%
}}}}
\put(3061,-1501){\makebox(0,0)[lb]{\smash{{\SetFigFont{20}{24.0}{\rmdefault}{\mddefault}{\updefault}{\color[rgb]{0,0,0}$\vec{v}_6 '$}%
}}}}
\put(2161,-1501){\makebox(0,0)[lb]{\smash{{\SetFigFont{20}{24.0}{\rmdefault}{\mddefault}{\updefault}{\color[rgb]{0,0,0}$\vec{v}_7 '$}%
}}}}
\put(4861,-4291){\makebox(0,0)[lb]{\smash{{\SetFigFont{20}{24.0}{\rmdefault}{\mddefault}{\updefault}{\color[rgb]{0,0,0}$\vec{v}_4 '$}%
}}}}
\put(5131,-3301){\makebox(0,0)[lb]{\smash{{\SetFigFont{20}{24.0}{\rmdefault}{\mddefault}{\updefault}{\color[rgb]{0,0,0}$Z$}%
}}}}
\put(4861,-1501){\makebox(0,0)[lb]{\smash{{\SetFigFont{20}{24.0}{\rmdefault}{\mddefault}{\updefault}{\color[rgb]{0,0,0}$U$}%
}}}}
\put(3961,-1501){\makebox(0,0)[lb]{\smash{{\SetFigFont{20}{24.0}{\rmdefault}{\mddefault}{\updefault}{\color[rgb]{0,0,0}$\vec{v}_5 '$}%
}}}}
\end{picture}%

%% file: enchain.pdf_t
\begin{picture}(0,0)%
\includegraphics{enchain.pdf}%
\end{picture}%
\setlength{\unitlength}{4144sp}%
\begingroup\makeatletter\ifx\SetFigFont\undefined%
\gdef\SetFigFont#1#2#3#4#5{%
  \reset@font\fontsize{#1}{#2pt}%
  \fontfamily{#3}\fontseries{#4}\fontshape{#5}%
  \selectfont}%
\fi\endgroup%
\begin{picture}(7500,2646)(301,-3442)
\put(3961,-2401){\makebox(0,0)[lb]{\smash{{\SetFigFont{20}{24.0}{\rmdefault}{\mddefault}{\updefault}{\color[rgb]{0,0,0}$X$}%
}}}}
\put(5761,-2401){\makebox(0,0)[lb]{\smash{{\SetFigFont{20}{24.0}{\rmdefault}{\mddefault}{\updefault}{\color[rgb]{0,0,0}$Y$}%
}}}}
\put(4861,-3346){\makebox(0,0)[lb]{\smash{{\SetFigFont{20}{24.0}{\rmdefault}{\mddefault}{\updefault}{\color[rgb]{0,0,0}$Z$}%
}}}}
\put(5221,-1951){\makebox(0,0)[lb]{\smash{{\SetFigFont{20}{24.0}{\rmdefault}{\mddefault}{\updefault}{\color[rgb]{0,0,0}$U$}%
}}}}
\put(5221,-1051){\makebox(0,0)[lb]{\smash{{\SetFigFont{20}{24.0}{\rmdefault}{\mddefault}{\updefault}{\color[rgb]{0,0,0}$\vec{v}_6 '$ }%
}}}}
\put(6661,-2401){\makebox(0,0)[lb]{\smash{{\SetFigFont{20}{24.0}{\rmdefault}{\mddefault}{\updefault}{\color[rgb]{0,0,0}$\vec{v}_6$}%
}}}}
\put(7561,-2401){\makebox(0,0)[lb]{\smash{{\SetFigFont{20}{24.0}{\rmdefault}{\mddefault}{\updefault}{\color[rgb]{0,0,0}$\vec{v}_7 '$}%
}}}}
\put(3061,-2401){\makebox(0,0)[lb]{\smash{{\SetFigFont{20}{24.0}{\rmdefault}{\mddefault}{\updefault}{\color[rgb]{0,0,0}$\vec{v}_5$}%
}}}}
\put(2161,-2401){\makebox(0,0)[lb]{\smash{{\SetFigFont{20}{24.0}{\rmdefault}{\mddefault}{\updefault}{\color[rgb]{0,0,0}$\vec{v}_7$}%
}}}}
\put(1261,-2401){\makebox(0,0)[lb]{\smash{{\SetFigFont{20}{24.0}{\rmdefault}{\mddefault}{\updefault}{\color[rgb]{0,0,0}$\vec{v}_8$}%
}}}}
\put(361,-2401){\makebox(0,0)[lb]{\smash{{\SetFigFont{20}{24.0}{\rmdefault}{\mddefault}{\updefault}{\color[rgb]{0,0,0}$\vec{v}_8 '$}%
}}}}
\end{picture}%

%% file: dn7e8.pdf_t
\begin{picture}(0,0)%
\includegraphics{dn7e8.pdf}%
\end{picture}%
\setlength{\unitlength}{4144sp}%
\begingroup\makeatletter\ifx\SetFigFont\undefined%
\gdef\SetFigFont#1#2#3#4#5{%
  \reset@font\fontsize{#1}{#2pt}%
  \fontfamily{#3}\fontseries{#4}\fontshape{#5}%
  \selectfont}%
\fi\endgroup%
\begin{picture}(8103,2856)(-1202,-4306)
\put(3961,-2401){\makebox(0,0)[lb]{\smash{{\SetFigFont{20}{24.0}{\rmdefault}{\mddefault}{\updefault}{\color[rgb]{0,0,0}$X$}%
}}}}
\put(5761,-2401){\makebox(0,0)[lb]{\smash{{\SetFigFont{20}{24.0}{\rmdefault}{\mddefault}{\updefault}{\color[rgb]{0,0,0}$Y$}%
}}}}
\put(5131,-1681){\makebox(0,0)[lb]{\smash{{\SetFigFont{20}{24.0}{\rmdefault}{\mddefault}{\updefault}{\color[rgb]{0,0,0}$U$}%
}}}}
\put(4861,-3301){\makebox(0,0)[lb]{\smash{{\SetFigFont{20}{24.0}{\rmdefault}{\mddefault}{\updefault}{\color[rgb]{0,0,0}$Z$}%
}}}}
\put(3061,-2401){\makebox(0,0)[lb]{\smash{{\SetFigFont{20}{24.0}{\rmdefault}{\mddefault}{\updefault}{\color[rgb]{0,0,0}$S'$}%
}}}}
\put(6661,-2401){\makebox(0,0)[lb]{\smash{{\SetFigFont{20}{24.0}{\rmdefault}{\mddefault}{\updefault}{\color[rgb]{0,0,0}$Y_R$}%
}}}}
\put(2161,-2401){\makebox(0,0)[lb]{\smash{{\SetFigFont{20}{24.0}{\rmdefault}{\mddefault}{\updefault}{\color[rgb]{0,0,0}$P$}%
}}}}
\put(2161,-4201){\makebox(0,0)[lb]{\smash{{\SetFigFont{20}{24.0}{\rmdefault}{\mddefault}{\updefault}{\color[rgb]{0,0,0}$\vec{v}_7 '$}%
}}}}
\put(1261,-4201){\makebox(0,0)[lb]{\smash{{\SetFigFont{20}{24.0}{\rmdefault}{\mddefault}{\updefault}{\color[rgb]{0,0,0}$\vec{v}_8 '$}%
}}}}
\put(1261,-2401){\makebox(0,0)[lb]{\smash{{\SetFigFont{20}{24.0}{\rmdefault}{\mddefault}{\updefault}{\color[rgb]{0,0,0}$Q_2$}%
}}}}
\put(361,-4201){\makebox(0,0)[lb]{\smash{{\SetFigFont{20}{24.0}{\rmdefault}{\mddefault}{\updefault}{\color[rgb]{0,0,0}$\vec{v}_9 '$}%
}}}}
\put(361,-2401){\makebox(0,0)[lb]{\smash{{\SetFigFont{20}{24.0}{\rmdefault}{\mddefault}{\updefault}{\color[rgb]{0,0,0}$\vec{v}_8$}%
}}}}
\put(-539,-2401){\makebox(0,0)[lb]{\smash{{\SetFigFont{20}{24.0}{\rmdefault}{\mddefault}{\updefault}{\color[rgb]{0,0,0}$\vec{v}_9$}%
}}}}
\end{picture}%

%% file: CoxEw.pdf_t
\begin{picture}(0,0)%
\includegraphics{CoxEw.pdf}%
\end{picture}%
\setlength{\unitlength}{4144sp}%
\begingroup\makeatletter\ifx\SetFigFont\undefined%
\gdef\SetFigFont#1#2#3#4#5{%
  \reset@font\fontsize{#1}{#2pt}%
  \fontfamily{#3}\fontseries{#4}\fontshape{#5}%
  \selectfont}%
\fi\endgroup%
\begin{picture}(7230,7363)(1966,-7811)
\put(3331,-2761){\makebox(0,0)[lb]{\smash{{\SetFigFont{20}{24.0}{\rmdefault}{\mddefault}{\updefault}{\color[rgb]{0,0,0}$F_2$}%
}}}}
\put(5581,-6001){\makebox(0,0)[lb]{\smash{{\SetFigFont{20}{24.0}{\rmdefault}{\mddefault}{\updefault}{\color[rgb]{0,0,0}$G_{33}$}%
}}}}
\put(3421,-6361){\makebox(0,0)[lb]{\smash{{\SetFigFont{20}{24.0}{\rmdefault}{\mddefault}{\updefault}{\color[rgb]{0,0,0}$C_2$}%
}}}}
\put(2791,-5641){\makebox(0,0)[lb]{\smash{{\SetFigFont{20}{24.0}{\rmdefault}{\mddefault}{\updefault}{\color[rgb]{0,0,0}$G_{31}$}%
}}}}
\put(4231,-7081){\makebox(0,0)[lb]{\smash{{\SetFigFont{20}{24.0}{\rmdefault}{\mddefault}{\updefault}{\color[rgb]{0,0,0}$G_{43}$}%
}}}}
\put(4951,-7711){\makebox(0,0)[lb]{\smash{{\SetFigFont{20}{24.0}{\rmdefault}{\mddefault}{\updefault}{\color[rgb]{0,0,0}$E_3$}%
}}}}
\put(6121,-7171){\makebox(0,0)[lb]{\smash{{\SetFigFont{20}{24.0}{\rmdefault}{\mddefault}{\updefault}{\color[rgb]{0,0,0}$G_{13}$}%
}}}}
\put(7021,-6451){\makebox(0,0)[lb]{\smash{{\SetFigFont{20}{24.0}{\rmdefault}{\mddefault}{\updefault}{\color[rgb]{0,0,0}$F_1$}%
}}}}
\put(7921,-5821){\makebox(0,0)[lb]{\smash{{\SetFigFont{20}{24.0}{\rmdefault}{\mddefault}{\updefault}{\color[rgb]{0,0,0}$G_{12}$}%
}}}}
\put(9181,-4471){\makebox(0,0)[lb]{\smash{{\SetFigFont{20}{24.0}{\rmdefault}{\mddefault}{\updefault}{\color[rgb]{0,0,0}$E_2$}%
}}}}
\put(5356,-691){\makebox(0,0)[lb]{\smash{{\SetFigFont{20}{24.0}{\rmdefault}{\mddefault}{\updefault}{\color[rgb]{0,0,0}$E_4$}%
}}}}
\put(5581,-3796){\makebox(0,0)[lb]{\smash{{\SetFigFont{20}{24.0}{\rmdefault}{\mddefault}{\updefault}{\color[rgb]{0,0,0}$D_1$}%
}}}}
\put(3961,-1861){\makebox(0,0)[lb]{\smash{{\SetFigFont{20}{24.0}{\rmdefault}{\mddefault}{\updefault}{\color[rgb]{0,0,0}$G_{24}$}%
}}}}
\put(6436,-1771){\makebox(0,0)[lb]{\smash{{\SetFigFont{20}{24.0}{\rmdefault}{\mddefault}{\updefault}{\color[rgb]{0,0,0}$G_{34}$}%
}}}}
\put(7336,-2581){\makebox(0,0)[lb]{\smash{{\SetFigFont{20}{24.0}{\rmdefault}{\mddefault}{\updefault}{\color[rgb]{0,0,0}$C_1$}%
}}}}
\put(8236,-3571){\makebox(0,0)[lb]{\smash{{\SetFigFont{20}{24.0}{\rmdefault}{\mddefault}{\updefault}{\color[rgb]{0,0,0}$G_{42}$}%
}}}}
\put(4141,-5101){\makebox(0,0)[lb]{\smash{{\SetFigFont{20}{24.0}{\rmdefault}{\mddefault}{\updefault}{\color[rgb]{0,0,0}$G_{41}$}%
}}}}
\put(5761,-5101){\makebox(0,0)[lb]{\smash{{\SetFigFont{20}{24.0}{\rmdefault}{\mddefault}{\updefault}{\color[rgb]{0,0,0}$D_4$}%
}}}}
\put(7381,-5101){\makebox(0,0)[lb]{\smash{{\SetFigFont{20}{24.0}{\rmdefault}{\mddefault}{\updefault}{\color[rgb]{0,0,0}$G_{32}$}%
}}}}
\put(1981,-4561){\makebox(0,0)[lb]{\smash{{\SetFigFont{20}{24.0}{\rmdefault}{\mddefault}{\updefault}{\color[rgb]{0,0,0}$E_1$}%
}}}}
\put(5626,-2311){\makebox(0,0)[lb]{\smash{{\SetFigFont{20}{24.0}{\rmdefault}{\mddefault}{\updefault}{\color[rgb]{0,0,0}$G_{44}$}%
}}}}
\put(2476,-3706){\makebox(0,0)[lb]{\smash{{\SetFigFont{20}{24.0}{\rmdefault}{\mddefault}{\updefault}{\color[rgb]{0,0,0}$G_{21}$}%
}}}}
\end{picture}%

%% file: e6dynkin.pdf_t
\begin{picture}(0,0)%
\includegraphics{e6dynkin.pdf}%
\end{picture}%
\setlength{\unitlength}{4144sp}%
\begingroup\makeatletter\ifx\SetFigFont\undefined%
\gdef\SetFigFont#1#2#3#4#5{%
  \reset@font\fontsize{#1}{#2pt}%
  \fontfamily{#3}\fontseries{#4}\fontshape{#5}%
  \selectfont}%
\fi\endgroup%
\begin{picture}(3900,1816)(2101,-3410)
\put(2161,-3301){\makebox(0,0)[lb]{\smash{{\SetFigFont{20}{24.0}{\rmdefault}{\mddefault}{\updefault}{\color[rgb]{0,0,0}$\alpha_{1}$}%
}}}}
\put(3061,-3301){\makebox(0,0)[lb]{\smash{{\SetFigFont{20}{24.0}{\rmdefault}{\mddefault}{\updefault}{\color[rgb]{0,0,0}$\alpha_{2}$}%
}}}}
\put(3961,-3301){\makebox(0,0)[lb]{\smash{{\SetFigFont{20}{24.0}{\rmdefault}{\mddefault}{\updefault}{\color[rgb]{0,0,0}$\alpha_{3}$}%
}}}}
\put(4861,-3301){\makebox(0,0)[lb]{\smash{{\SetFigFont{20}{24.0}{\rmdefault}{\mddefault}{\updefault}{\color[rgb]{0,0,0}$\alpha_{4}$}%
}}}}
\put(5761,-3301){\makebox(0,0)[lb]{\smash{{\SetFigFont{20}{24.0}{\rmdefault}{\mddefault}{\updefault}{\color[rgb]{0,0,0}$\alpha_{5}$}%
}}}}
\put(4321,-1861){\makebox(0,0)[lb]{\smash{{\SetFigFont{20}{24.0}{\rmdefault}{\mddefault}{\updefault}{\color[rgb]{0,0,0}$\alpha_{6}$}%
}}}}
\end{picture}%

%% file: X3cox.pdf_t
\begin{picture}(0,0)%
\includegraphics{X3cox.pdf}%
\end{picture}%
\setlength{\unitlength}{4144sp}%
\begingroup\makeatletter\ifx\SetFigFont\undefined%
\gdef\SetFigFont#1#2#3#4#5{%
  \reset@font\fontsize{#1}{#2pt}%
  \fontfamily{#3}\fontseries{#4}\fontshape{#5}%
  \selectfont}%
\fi\endgroup%
\begin{picture}(5564,5985)(141,-4891)
\put(4989,-3931){\makebox(0,0)[lb]{\smash{{\SetFigFont{20}{24.0}{\rmdefault}{\mddefault}{\updefault}{\color[rgb]{0,0,0}$\lambda_{17}$}%
}}}}
\put(4131,-4426){\makebox(0,0)[lb]{\smash{{\SetFigFont{20}{24.0}{\rmdefault}{\mddefault}{\updefault}{\color[rgb]{0,0,0}$\lambda_{18}$}%
}}}}
\put(780,-61){\makebox(0,0)[lb]{\smash{{\SetFigFont{20}{24.0}{\rmdefault}{\mddefault}{\updefault}{\color[rgb]{0,0,0}$\lambda_2$}%
}}}}
\put(1559,389){\makebox(0,0)[lb]{\smash{{\SetFigFont{20}{24.0}{\rmdefault}{\mddefault}{\updefault}{\color[rgb]{0,0,0}$\lambda_1$}%
}}}}
\put(2339,839){\makebox(0,0)[lb]{\smash{{\SetFigFont{20}{24.0}{\rmdefault}{\mddefault}{\updefault}{\color[rgb]{0,0,0}$\lambda_{22}$}%
}}}}
\put(1793,-1726){\makebox(0,0)[lb]{\smash{{\SetFigFont{20}{24.0}{\rmdefault}{\mddefault}{\updefault}{\color[rgb]{0,0,0}$\lambda_{25}$}%
}}}}
\put(1793,-2266){\makebox(0,0)[lb]{\smash{{\SetFigFont{20}{24.0}{\rmdefault}{\mddefault}{\updefault}{\color[rgb]{0,0,0}$\lambda_7$}%
}}}}
\put(4131,-1726){\makebox(0,0)[lb]{\smash{{\SetFigFont{20}{24.0}{\rmdefault}{\mddefault}{\updefault}{\color[rgb]{0,0,0}$\lambda_{24}$}%
}}}}
\put(4131,-2356){\makebox(0,0)[lb]{\smash{{\SetFigFont{20}{24.0}{\rmdefault}{\mddefault}{\updefault}{\color[rgb]{0,0,0}$\lambda_{16}$}%
}}}}
\put(3352,-826){\makebox(0,0)[lb]{\smash{{\SetFigFont{20}{24.0}{\rmdefault}{\mddefault}{\updefault}{\color[rgb]{0,0,0}$\lambda_{19}$}%
}}}}
\put(157,-1411){\makebox(0,0)[lb]{\smash{{\SetFigFont{20}{24.0}{\rmdefault}{\mddefault}{\updefault}{\color[rgb]{0,0,0}$\lambda_4$}%
}}}}
\put(157,-2311){\makebox(0,0)[lb]{\smash{{\SetFigFont{20}{24.0}{\rmdefault}{\mddefault}{\updefault}{\color[rgb]{0,0,0}$\lambda_5$}%
}}}}
\put(156,-511){\makebox(0,0)[lb]{\smash{{\SetFigFont{20}{24.0}{\rmdefault}{\mddefault}{\updefault}{\color[rgb]{0,0,0}$\lambda_3$}%
}}}}
\put(157,-3211){\makebox(0,0)[lb]{\smash{{\SetFigFont{20}{24.0}{\rmdefault}{\mddefault}{\updefault}{\color[rgb]{0,0,0}$\lambda_6$}%
}}}}
\put(2573,-3076){\makebox(0,0)[lb]{\smash{{\SetFigFont{20}{24.0}{\rmdefault}{\mddefault}{\updefault}{\color[rgb]{0,0,0}$\lambda_{20}$}%
}}}}
\put(5690,-2356){\makebox(0,0)[lb]{\smash{{\SetFigFont{20}{24.0}{\rmdefault}{\mddefault}{\updefault}{\color[rgb]{0,0,0}$\lambda_{14}$}%
}}}}
\put(5690,-3256){\makebox(0,0)[lb]{\smash{{\SetFigFont{20}{24.0}{\rmdefault}{\mddefault}{\updefault}{\color[rgb]{0,0,0}$\lambda_{15}$}%
}}}}
\put(5690,-1456){\makebox(0,0)[lb]{\smash{{\SetFigFont{20}{24.0}{\rmdefault}{\mddefault}{\updefault}{\color[rgb]{0,0,0}$\lambda_{13}$}%
}}}}
\put(5690,-556){\makebox(0,0)[lb]{\smash{{\SetFigFont{20}{24.0}{\rmdefault}{\mddefault}{\updefault}{\color[rgb]{0,0,0}$\lambda_{12}$}%
}}}}
\put(4131,614){\makebox(0,0)[lb]{\smash{{\SetFigFont{20}{24.0}{\rmdefault}{\mddefault}{\updefault}{\color[rgb]{0,0,0}$\lambda_{10}$}%
}}}}
\put(4911,164){\makebox(0,0)[lb]{\smash{{\SetFigFont{20}{24.0}{\rmdefault}{\mddefault}{\updefault}{\color[rgb]{0,0,0}$\lambda_{11}$}%
}}}}
\put(2494,-4786){\makebox(0,0)[lb]{\smash{{\SetFigFont{20}{24.0}{\rmdefault}{\mddefault}{\updefault}{\color[rgb]{0,0,0}$\lambda_{23}$}%
}}}}
\put(1715,-4336){\makebox(0,0)[lb]{\smash{{\SetFigFont{20}{24.0}{\rmdefault}{\mddefault}{\updefault}{\color[rgb]{0,0,0}$\lambda_9$}%
}}}}
\put(936,-3841){\makebox(0,0)[lb]{\smash{{\SetFigFont{20}{24.0}{\rmdefault}{\mddefault}{\updefault}{\color[rgb]{0,0,0}$\lambda_8$}%
}}}}
\end{picture}%

%% file: e8dynkin.pdf_t
\begin{picture}(0,0)%
\includegraphics{e8dynkin.pdf}%
\end{picture}%
\setlength{\unitlength}{4144sp}%
\begingroup\makeatletter\ifx\SetFigFont\undefined%
\gdef\SetFigFont#1#2#3#4#5{%
  \reset@font\fontsize{#1}{#2pt}%
  \fontfamily{#3}\fontseries{#4}\fontshape{#5}%
  \selectfont}%
\fi\endgroup%
\begin{picture}(6600,1816)(-599,-3410)
\put(-539,-3301){\makebox(0,0)[lb]{\smash{{\SetFigFont{20}{24.0}{\rmdefault}{\mddefault}{\updefault}{\color[rgb]{0,0,0}$\alpha_{-\theta}$}%
}}}}
\put(361,-3301){\makebox(0,0)[lb]{\smash{{\SetFigFont{20}{24.0}{\rmdefault}{\mddefault}{\updefault}{\color[rgb]{0,0,0}$\alpha_{1}$}%
}}}}
\put(1261,-3301){\makebox(0,0)[lb]{\smash{{\SetFigFont{20}{24.0}{\rmdefault}{\mddefault}{\updefault}{\color[rgb]{0,0,0}$\alpha_{2}$}%
}}}}
\put(2161,-3301){\makebox(0,0)[lb]{\smash{{\SetFigFont{20}{24.0}{\rmdefault}{\mddefault}{\updefault}{\color[rgb]{0,0,0}$\alpha_{3}$}%
}}}}
\put(3061,-3301){\makebox(0,0)[lb]{\smash{{\SetFigFont{20}{24.0}{\rmdefault}{\mddefault}{\updefault}{\color[rgb]{0,0,0}$\alpha_{4}$}%
}}}}
\put(3961,-3301){\makebox(0,0)[lb]{\smash{{\SetFigFont{20}{24.0}{\rmdefault}{\mddefault}{\updefault}{\color[rgb]{0,0,0}$\alpha_{5}$}%
}}}}
\put(4861,-3301){\makebox(0,0)[lb]{\smash{{\SetFigFont{20}{24.0}{\rmdefault}{\mddefault}{\updefault}{\color[rgb]{0,0,0}$\alpha_{6}$}%
}}}}
\put(5761,-3301){\makebox(0,0)[lb]{\smash{{\SetFigFont{20}{24.0}{\rmdefault}{\mddefault}{\updefault}{\color[rgb]{0,0,0}$\alpha_{7}$}%
}}}}
\put(4321,-1861){\makebox(0,0)[lb]{\smash{{\SetFigFont{20}{24.0}{\rmdefault}{\mddefault}{\updefault}{\color[rgb]{0,0,0}$\alpha_{8}$}%
}}}}
\end{picture}%

%% file: Sigma_2.pdf_t
\begin{picture}(0,0)%
\includegraphics{Sigma_2.pdf}%
\end{picture}%
\setlength{\unitlength}{4144sp}%
\begingroup\makeatletter\ifx\SetFigFont\undefined%
\gdef\SetFigFont#1#2#3#4#5{%
  \reset@font\fontsize{#1}{#2pt}%
  \fontfamily{#3}\fontseries{#4}\fontshape{#5}%
  \selectfont}%
\fi\endgroup%
\begin{picture}(8220,2384)(2506,-4532)
\put(7201,-2761){\makebox(0,0)[lb]{\smash{{\SetFigFont{25}{30.0}{\rmdefault}{\mddefault}{\updefault}{\color[rgb]{0,0,.56}$\beta_2$}%
}}}}
\put(8551,-3661){\makebox(0,0)[lb]{\smash{{\SetFigFont{20}{24.0}{\rmdefault}{\mddefault}{\updefault}{\color[rgb]{0,0,0}$p_4$}%
}}}}
\put(9451,-4381){\makebox(0,0)[lb]{\smash{{\SetFigFont{25}{30.0}{\rmdefault}{\mddefault}{\updefault}{\color[rgb]{0,0,0}$\Sigma_2$}%
}}}}
\put(4321,-3121){\makebox(0,0)[lb]{\smash{{\SetFigFont{20}{24.0}{\rmdefault}{\mddefault}{\updefault}{\color[rgb]{0,0,0}$p_1$}%
}}}}
\put(2521,-3121){\makebox(0,0)[lb]{\smash{{\SetFigFont{20}{24.0}{\rmdefault}{\mddefault}{\updefault}{\color[rgb]{0,0,0}$p_0$}%
}}}}
\put(9361,-2941){\makebox(0,0)[lb]{\smash{{\SetFigFont{25}{30.0}{\rmdefault}{\mddefault}{\updefault}{\color[rgb]{.56,0,0}$\alpha_2$}%
}}}}
\put(10711,-3121){\makebox(0,0)[lb]{\smash{{\SetFigFont{20}{24.0}{\rmdefault}{\mddefault}{\updefault}{\color[rgb]{0,0,0}$p_5$}%
}}}}
\put(3556,-2986){\makebox(0,0)[lb]{\smash{{\SetFigFont{25}{30.0}{\rmdefault}{\mddefault}{\updefault}{\color[rgb]{.56,0,0}$\alpha_1$}%
}}}}
\put(5581,-2761){\makebox(0,0)[lb]{\smash{{\SetFigFont{25}{30.0}{\rmdefault}{\mddefault}{\updefault}{\color[rgb]{0,0,.56}$\beta_1$}%
}}}}
\put(5581,-3661){\makebox(0,0)[lb]{\smash{{\SetFigFont{20}{24.0}{\rmdefault}{\mddefault}{\updefault}{\color[rgb]{0,0,0}$p_2$}%
}}}}
\put(7471,-3661){\makebox(0,0)[lb]{\smash{{\SetFigFont{20}{24.0}{\rmdefault}{\mddefault}{\updefault}{\color[rgb]{0,0,0}$p_3$}%
}}}}
\end{picture}%